\documentclass[12pt,reqno]{amsart}
\usepackage[T1]{fontenc}
\usepackage{lmodern}
\usepackage{amsmath,amssymb,amsthm,array}
\usepackage[a4paper,top=28mm,bottom=28mm,left=32mm,right=32mm,headheight=14pt]{geometry}
\usepackage{xcolor}
\definecolor{inkblue}{RGB}{23,66,90}
\usepackage[unicode,colorlinks,allcolors=inkblue,bookmarksnumbered=true,bookmarksdepth=2]{hyperref}
\hypersetup{pdftitle={Modularity of Heegner Cycles and Arithmetic Inner Products on Kuga\textendash Sato Varieties},pdfauthor={Tuoping Du}}
\makeatletter
\renewcommand*\l@section{\@tocline{1}{2pt}{1pc}{}{\bfseries}}
\renewcommand*\l@subsection{\@tocline{2}{0pt}{2.4pc}{5pc}{}}
\makeatother
\numberwithin{equation}{section}
\newtheorem{theorem}{Theorem}[section]
\newtheorem{proposition}[theorem]{Proposition}
\newtheorem{corollary}[theorem]{Corollary}
\newtheorem{lemma}[theorem]{Lemma}
\newtheorem{conjecture}[theorem]{Conjecture}
\newtheorem*{auxlemma}{Lemma}
\theoremstyle{definition}
\newtheorem{remark}[theorem]{Remark}
\newtheorem*{introremark}{Remark}

\newcommand{\Hom}{\operatorname{Hom}}
\newcommand{\Pic}{\operatorname{Pic}}
\newcommand{\CH}{\operatorname{CH}}

\begin{document}

\title[Heegner cycles and arithmetic heights]
{Modularity of Heegner Cycles and Arithmetic Inner Products on Kuga--Sato Varieties}
\author[Tuoping Du]{Tuoping Du}
\address{School of Mathematics and Physics, North China Electric Power University, Beijing 102206, P.R. China}
\email{dtp1982@163.com}
\date{24 November 2026}
\begin{abstract}
For \(\kappa\ge2\), the principal part of a weakly holomorphic modular
form of weight \(\tfrac32-\kappa\) determines a linear combination of
codimension-\(\kappa\) Heegner cycles on a compact Kuga--Sato variety.
We use isogeny correspondences to prove that every such combination
vanishes in the Chow group with complex coefficients. These rational
equivalences imply that the generating series of Heegner cycles is a
Chow-valued cusp form of weight \(\kappa+\tfrac12\) at every 
level and for every compatible fundamental twist.

We prove an arithmetic inner product formula for theta lifts on the
modular curve \(X_0(N)\) and on higher-dimensional Kuga--Sato varieties.
It expresses their height pairings in terms of central derivatives of
quadratic twists of newform \(L\)-functions, with an explicit product
of local operators on oldspaces. On each primitive newform line
selected by the twist, the height of one projected Heegner cycle at a
fundamental index with nonzero Fourier coefficient determines all
pairings of projected cycles on that line. This specialization
recovers the Gross--Zagier formula in weight two and the
genus-character case of the Gross--Zagier--Zhang formula in higher weight.
\end{abstract}

\maketitle
{\small\linespread{1.00}\selectfont\tableofcontents}

\section{Introduction}\label{rd-intro}

The Gross--Kohnen--Zagier theorem identifies Heegner divisor classes
on a modular curve as the Fourier coefficients of a Jacobi form
\cite{GKZ}. For higher Heegner cycles the corresponding question
concerns rational equivalence on a Kuga--Sato variety. Does the polar
part of every weakly holomorphic modular form give a sum of Heegner
cycles rationally equivalent to zero? In the untwisted case Du and Peng
reduced Chow-valued modularity to this question
\cite[Conjecture~0.4 and Proposition~0.5]{DuPeng}.
Theorem~\ref{intro-chow-vanishing} gives the vanishing for all
compatible fundamental twists.
The principal-part pairing then gives the Chow-valued cusp form of
Theorem~\ref{med-target}.

The second question, from Kudla's program, asks for the arithmetic
inner product of the resulting theta lifts. Theorems~\ref{med-inner}
and~\ref{ks-arithmetic-inner-product} express this pairing in terms
of a central \(L\)-derivative. On a primitive newform line the
projected Heegner cycles are all multiples of one class. Evaluating
the pairing at a fundamental Fourier index and applying Waldspurger's
formula recovers the Gross--Zagier height in the divisor case and the
Gross--Zagier--Zhang height for a genus character in higher codimension,
under the hypotheses of \eqref{med-gz} and
Corollary~\ref{p19-higher-gzz}, respectively.

\subsection{Heegner cycles}
A CM point \(y=(E,C)\) on \(X_0(N)\) consists of an elliptic curve
with complex multiplication and a cyclic subgroup of order \(N\).
For \(\kappa\ge2\) put \(a=\kappa-1\). On a fine-level cover the fibre
of \(\mathcal Y_\kappa(N)\) at \(y\) is \(E^{2a}\). An oriented CM
vector \(x\) belongs to the rational space of trace-zero
\(2\times2\) matrices with \(Q(x)=N\det x\). It gives a trace-zero endomorphism \(u\) of \(E\) after
clearing denominators. Its graph has a mixed part
\[
 D_u=[\Gamma_u]-[E\times\{0\}]-(\deg u)[\{0\}\times E]
\]
with intersection zero against either axis. Apply the Scholl projector
to \(D_u^{\boxtimes a}\) and normalize its fibre self-intersection
to have degree \((-1)^a\). The projector takes the part odd in every
factor and alternating under permutation. Pushing the normalized
class forward from the fibre gives the \emph{Heegner cycle}
\(S_\kappa(x)\) of codimension \(\kappa\). Reversing \(x\) multiplies
it by \((-1)^a\).

Fix a fundamental discriminant \(\Delta\), allowing \(\Delta=1\), and
\(r\bmod 2N\) with \(r^2\equiv\Delta\bmod 4N\). For \(n>0\) and
\(\mu\in L^\sharp/L\simeq\mathbb Z/2N\mathbb Z\), the coefficient
\(Z_\kappa(n,\mu)=Z^{(N)}_{\Delta,r,\kappa}(n,\mu)\) is the sum of
the \(S_\kappa(x)\) over \(\Gamma_0(N)\)-orbits with
\(Q(x)=N\det x=|\Delta|n\) and \(x\bmod L=r\mu\). Each orbit has
weight \(2n^{a/2}\chi_\Delta(x)/|\Gamma_x|\), with genus character
\(\chi_\Delta\) and stabilizer \(\Gamma_x\).
Section~\ref{rd-cycles} gives the lattice \(L\) and the descent from
fine level. For \(\kappa\ge2\) the coefficients belong to
\(A_N=\CH^\kappa_{\rm hom}(\mathcal Y_\kappa(N)_{\overline{\mathbb Q}})_{\mathbb C}\),
the group of homologically trivial cycles modulo rational equivalence
with complex coefficients. For \(\kappa=1\) we subtract the degree
of the orbit divisor times the cusp at infinity and take its Jacobian
class. All tensor products are algebraic.

\subsection{Vanishing and Chow-valued modularity}
Assume \(\kappa\ge2\), and put \(k=\kappa+\tfrac12\) and
\(w=2-k=\tfrac32-\kappa\). Write \(e_\mu\) for the basis vector
indexed by \(\mu\) in \(\mathbb C[L^\sharp/L]\). The output Weil
representation \(\widetilde\rho_L\) is \(\rho_L\) for \(\Delta>0\)
and \(\bar\rho_L\) for \(\Delta<0\). The input has the dual type
\(\overline{\widetilde\rho_L}\). For a harmonic Maass form \(F\) of
weight \(w\) with cuspidal shadow \(\xi_w F\), write \(c_F^+(n,\mu)\)
for the Fourier coefficients of its holomorphic part. Its principal
part, here the sum of its strictly negative powers of \(q=e^{2\pi i\tau}\),
selects the finite combination
\[
 Z_\kappa(F)=\sum_{\mu,\,n>0}c_F^+(-n,\mu)\,Z_\kappa(n,\mu).
\]
Write \(M^!\) for the weakly holomorphic forms, which are
holomorphic on the upper half-plane and meromorphic at the cusps.

\begin{theorem}[Vanishing for weakly holomorphic forms]
\label{intro-chow-vanishing}
For every \(f\in M^!_{w,\overline{\widetilde\rho_L}}\),
\begin{equation}\label{intro-weak}
 Z_\kappa(f)=0
 \qquad\text{in }\CH^\kappa
       (\mathcal Y_\kappa(N)_{\overline{\mathbb Q}})_{\mathbb C}.
\end{equation}
\end{theorem}

\begin{introremark}[Geometry of the vanishing theorem]\leavevmode
\begin{enumerate}
\item \emph{Untwisted cycles.} Restricting the isogeny graphs to
the cyclic isogeny locus gives the CM sums
(Theorem~\ref{graph-to-moments}). On a family of Kummer surfaces
these graphs satisfy a recurrence modulo classes supported over
finitely many pairs of base points. Restriction kills those
classes, giving the untwisted relations
(Corollary~\ref{all-untwisted-modularity}). In codimension two,
the action of elliptic homomorphisms on \(H^1\) detects the mixed
divisors on \(E^2\) (Proposition~\ref{restricted-divisor-graph-errors}).
In higher codimension Proposition~\ref{fixed-evaluation} places the
evaluation class on the line spanned by the isogeny graph in Chow.
The graph has nonzero cohomology class, so a numerical pairing then
determines its coefficient.

\item \emph{Twists and compactification.} Transfer to auxiliary
untwisted levels reduces the obstruction to the Hecke eigenvalues
\(\pm(p+1)p^{\kappa-1}\). At a prime common to \(N\) and \(\Delta\),
trace to level \(N/p\) makes \(p\) a good prime.
Proposition~\ref{marked-cm-spectrum} excludes both eigenvalues.
Restriction is injective on the Scholl summand
(Lemma~\ref{scholl-open-chow}), giving the compact result in
Theorem~\ref{rd-m-vanishing}.

\item \emph{Abel--Jacobi images.} On parabolic cohomology the
same graphs act by the Eichler--Shimura operators. This proves
vanishing of the Hodge and \(\ell\)-adic Abel--Jacobi images
(Theorem~\ref{general-aj-modularity}). Here injectivity concerns
restriction to the open target. After the Tate twist, the boundary
has neither classes of Hodge type \((0,0)\) nor \(\ell\)-adic
Galois invariants.
\end{enumerate}
\end{introremark}

The kernel of the shadow operator \(\xi_w\) consists of weakly
holomorphic forms \cite{BFGeometricTheta}. Two harmonic lifts with
the same shadow therefore give the same Chow class by
Theorem~\ref{intro-chow-vanishing}. Write
\(S=S_{k,\widetilde\rho_L}^{\Pi_{\Delta,N}}\) for the cusp forms
of the indicated weight and type fixed by the residue symmetries
at primes dividing both \(N\) and \(\Delta\)
(Section~\ref{rd-cycles}).

Write \(\Phi_{\Delta,r,\kappa}^{(N)}\) for the generating series of
\(Z_{\Delta,r,\kappa}^{(N)}(n,\mu)\). We abbreviate it to \(\Phi_\kappa\)
when \(N,\Delta,r\) are fixed and retain the superscript \((N)\) when
comparing levels.

\begin{theorem}[Chow-valued modularity]\label{med-target}
The Heegner series \(\Phi_\kappa\) is a cusp form of weight \(k\)
with values in \(A_N\). For any Petersson-orthogonal basis \((g_j)\)
of \(S\) and harmonic lifts with \(\xi_w F_j=g_j/\|g_j\|_{\rm Pet}^2\),
\[
 \Phi_\kappa(\tau)=\sum_{\mu,\,n>0}q^n e_\mu\otimes Z_\kappa(n,\mu)
             =\sum_j g_j(\tau)\otimes Z_\kappa(F_j)
             \qquad\text{in }S\otimes A_N .
\]
\end{theorem}
\begin{introremark}[How Theorem~\ref{intro-chow-vanishing} gives modularity]
An invariant principal part \(P\) defines a Fourier functional
\[
 \beta(P)(g)=\sum_{\mu,\,n>0}c_P(-n,\mu)c_g(n,\mu).
\]
These functionals span \(S^\vee\), since the Fourier coefficients
determine \(g\). The principal-part pairing identifies \(\ker\beta\)
with the polar parts of the weakly holomorphic forms
(Proposition~\ref{med-exact}). By Theorem~\ref{intro-chow-vanishing},
the cycle assignment vanishes on this kernel and hence defines a map
\(S^\vee\to A_N\). Its tensor in \(S\otimes A_N\) has coefficient
\(Z_\kappa(n,\mu)\) at \((n,\mu)\). The basis dual to the \(g_j\)
gives the expansion above.
\end{introremark}

The Petersson pairing in the modular variable now defines the theta
lift \(\theta_{\Delta,r,\kappa}(f)=(\Phi_\kappa,f)_{\rm Pet}\in A_N\)
for \(f\in S\). Theorem~\ref{ks-arithmetic-inner-product} computes
the inner product of its admissible arithmetic lifts.

Finitely many cycles determine the whole series.
Choose Fourier evaluations \(\ell_{a_1},\ldots,\ell_{a_s}\) forming a basis of
\(S^\vee\), all with exponent at most \(k/12\). Writing
\(\ell_i=\sum_j u_j(i)\ell_{a_j}\) gives
\[
                  Z_\kappa(i)=\sum_{j=1}^{s}u_j(i)Z_\kappa(a_j).
\]
The \(u_j(i)\) are rational. For \(i\) outside this basis the
corresponding difference of polar monomials is the principal part
of a weakly holomorphic form. These forms constitute a rational
basis of all invariant weakly holomorphic inputs
(Theorem~\ref{med-complete-basis}). Further relations among the
selected cycles form the kernel of \(S^\vee\to A_N\). On the dual
of a retained Hecke component with nonzero central derivative,
finitely many height pairings in one CM field separate vectors,
so the map is injective there \eqref{eq:height-detects-weak}.
For a harmonic input in the
invariant sector with shadow on that component, \(Z_\kappa(F)=0\)
is therefore equivalent to weak holomorphy.

The supporting fibres retain information lost in the total
space. Their Betti classes detect every nonzero invariant principal
part (Theorem~\ref{fibre-detection}). For a weakly holomorphic
input the resulting tuple on the fibres is a boundary in Bloch's
localization sequence (Proposition~\ref{med-boundary}).
At integral pole bound \(B\ge k/12\) these boundaries span a space
of dimension \(B d_\Pi-\dim S\), where \(d_\Pi\) counts the residue
components after the central and shared-prime symmetries
(Corollary~\ref{bounded-fibre-boundaries}).

\subsection{Relation to earlier work}
Gross and Zagier computed the height of a Heegner point and described
Deligne's higher-weight vectors above CM points \cite[Ch.~V, \S4]{GrossZagier}.
Gross, Kohnen and Zagier related the variation of the divisor class with
its discriminant to the coefficients of a Jacobi form \cite{GKZ}.
Borcherds produced divisor relations from automorphic products with
prescribed weakly holomorphic principal parts \cite{BorcherdsGKZ}.
Bruinier and Ono incorporated genus twists and harmonic Maass forms
\cite{BO}.

For higher cycles, Zhang proved a Gross--Zagier formula \cite{Zhang}.
Qiu's arbitrary-level formula retains the local toric periods, and his
unramified Hecke decomposition concerns CM cycles modulo the radical of
the height pairing \cite[\S2.3.1, Theorem~2.3.4 and Corollary~2.4.7]{QiuCM}.
Xue computed mixed heights for distinct CM fields and deduced a
higher-weight Gross--Kohnen--Zagier theorem under a positivity
hypothesis \cite{XueGKZ}. Nekov\'a\v{r} used \(\ell\)-adic
Abel--Jacobi images of individual cycles in an Euler system
\cite{NekovarKuga}. Lilienfeldt expressed their complex Abel--Jacobi
images as cusp-form integrals and proved infinite rank in the
Griffiths group for a fixed imaginary quadratic field, extending
Schoen's result for elliptic modular threefolds
\cite{SchoenCM,LilienfeldtGriffiths}. Here the coefficients are
genus-weighted orbit sums. Corollary~\ref{fixed-index-infinite-rank}
holds the Fourier index fixed and varies the fundamental twist.
The resulting cycles lie in distinct quadratic Galois character
spaces in ordinary Chow. Their independence follows from this
character law and the nonvanishing of their projected heights.

Sreekantan constructed rational equivalences in codimension two
as boundaries in higher Chow groups \cite{SreekantanHeegner}.
Zemel proved higher-codimension modularity modulo the strongly
principal relations from a theta lift \cite[Theorem~4.6]{ZemelGKZ}.
Du and Peng proved untwisted modularity modulo the relations obtained
from residue divisors of algebraic meromorphic modular forms of weight
\(2\kappa\). They reduced modularity in the ordinary Chow group to
the vanishing of the cycles of weakly holomorphic forms
\cite[Theorem~0.3, Conjecture~0.4 and Proposition~0.5]{DuPeng}.
The target in Theorem~\ref{intro-chow-vanishing} is the ordinary
Chow group of the compact variety. Its relations therefore pass to
every height pairing defined on Chow.

\subsection{The arithmetic inner product on a modular curve}
The modular-curve formula and the higher-weight CM height calculation
originate in the author's \href{https://arxiv.org/abs/2412.00688}{2024 preprint}.
Their proofs are included in Sections~\ref{rd-h-lift}--\ref{rd-h-height}
with the height and Petersson normalizations used below.
Let \(G\) be a normalized weight-two newform of level \(N\) and
trivial character. Let \(\Delta\) be a fundamental discriminant
prime to \(N\), allowing \(\Delta=1\), and fix any root
\(r^2\equiv\Delta\bmod 4N\). Suppose \(G\otimes\chi_\Delta\)
has root number \(-1\). Let \(g\) be a weight-\(3/2\) Shimura
newform corresponding to \(G\), with real Fourier coefficients and type
\(\rho_L\) for \(\Delta>0\) and \(\bar\rho_L\) for \(\Delta<0\).
The Petersson projection of the degree-zero Heegner divisor series
\(A_{\Delta,r}\) defines the Jacobian class
\(\theta_{\Delta,r}(g)=(A_{\Delta,r},g)_{\rm Pet}\).
We use the absolute N\'eron--Tate height and the unaveraged
Petersson norm.

\begin{theorem}[Arithmetic inner product on \(X_0(N)\)]\label{med-inner}
Let \(\widehat\theta_{\Delta,r}(g)\) be any admissible extension of
\(\theta_{\Delta,r}(g)\) to a regular proper arithmetic model,
with arithmetic degree divided by the degree of the field of definition.
Then
\begin{equation}\label{rd-h-nt}
 -\langle\widehat\theta_{\Delta,r}(g),\widehat\theta_{\Delta,r}(g)\rangle_{\rm GS}
 =\|\theta_{\Delta,r}(g)\|_{\rm NT}^{2}
 =\frac{2\sqrt{N|\Delta|}}{\pi}\,\|g\|_{\rm Pet}^{2} L'(G,\chi_\Delta,1).
\end{equation}
\end{theorem}

\begin{introremark}[Primes dividing the level and the twist]
Let \(\Delta\ne1\) be fundamental with \(r^2\equiv\Delta\bmod 4N\).
Write \(\mathcal R=\{p:p\mid(\Delta,N)\}\), \(t=\#\mathcal R>0\), and
\(\sigma=\operatorname{sgn}(\Delta)\). Suppose the exact-level
newform \(G\) has root number \(-\sigma\) and \(a_p(G)=-1\) at
every \(p\in\mathcal R\). For its real Shimura newform \(g\) of the type
specified above, \(G\otimes\chi_\Delta\) has root number \(-1\) and
\begin{equation*}\tag{\ref*{rd-h-nt}a}\label{eq:curve-shared-height}
 \|\theta_{\Delta,r}(g)\|_{\rm NT}^2
 =2^t\,\frac{2\sqrt{N|\Delta|}}{\pi}\,
       \|g\|_{\rm Pet}^2 L'(G,\chi_\Delta,1).
\end{equation*}
If instead \(a_p(G)=+1\) at any common prime, with
\(\varepsilon(G)=-\sigma\) as above, the projected lift vanishes.
Section~\ref{rd-h-lift} gives both calculations, including at \(2\).
\end{introremark}

Kudla, Rapoport and Yang related arithmetic theta lifts on compact
Shimura curves to central derivatives \cite[Chapter~9]{KRYBook}.
On \(X_0(N)\), Du and Yang computed the pairings of arithmetic
theta series with the metrized Hodge bundle \cite{DY1,DY2}.

At squarefree level the Petersson projection of the full arithmetic
theta series gives the admissible lift on each cuspidal Shimura
eigenspace, for both signs of \(\Delta\)
(Theorems~\ref{p19a-lift} and~\ref{squarefree-signed-theta}).
For a negative twist the first horizontal divisor has opposite
nonzero intersections with the two components at each split bad
place. The component series is unary (Proposition~\ref{unary-split-components})
and hence orthogonal to the cuspidal Shimura eigenspace. For an
arbitrary weight-\(3/2\) cusp form its contribution is
\eqref{vertical-arithmetic-pairing}. In the untwisted case the
lift is nonzero precisely when \(L'(G,1)\ne0\).

Let \(y^G_{\Delta,r}(a)\) denote the projected divisor at index
\(a\), and \(b_g(a)\) the coefficient of \(g\). Put \(t=0\) in
Theorem~\ref{med-inner} and \(t=\#\mathcal R\) under the sign
conditions of \eqref{eq:curve-shared-height}. Whenever \(b_g(a)\ne0\),
\begin{equation}\label{eq:curve-single-derivative}
 L'(G,\chi_\Delta,1)=
 \frac{\pi\|g\|_{\rm Pet}^2}{2^{t+1}\sqrt{N|\Delta|}}\,
 \frac{\|y^G_{\Delta,r}(a)\|_{\rm NT}^{2}}{|b_g(a)|^2}.
\end{equation}
Indeed, \(y^G_{\Delta,r}(v)=b_g(v)y^G_{\Delta,r}(a)/b_g(a)\)
for every \(v\) by \eqref{med-proj}. One may take
\(a=(|d|/(4N),\mu)\) with \(d\) fundamental, prime to \(2N\Delta\),
of sign opposite to \(\Delta\), and \(d\equiv\mu^2\bmod 4N\).
In the twisted case we fix the height normalization using Gross--Zagier
\eqref{med-gz} and Waldspurger \eqref{med-wald}. Dividing the two
identities at this index cancels \(L(G,\chi_d,1)\) and gives
\eqref{rd-h-nt} by projection. Conversely, substituting the
Waldspurger value of \(\lvert b_g(a)\rvert^2\) into the projected
inner product recovers \eqref{med-gz} wherever both coefficient
formulas apply. The mixed-height identity is \eqref{p19-divisor-mixed}.

Over a totally real multiquadratic field \(K/\mathbb Q\) of degree \(d\), the lifts
for distinct characters are orthogonal. If all twists satisfy
Theorem~\ref{med-inner}, the determinant of their normalized heights
is \((2\sqrt N/\pi)^d L^{(d)}(G/K,1)/d!\).
When primes of \(N\) ramify in \(K\), Theorem~\ref{thm:multiquadratic-height}
gives the additional power of two in terms of their inertia groups
under the local sign conditions above.

\subsection{From a mixed height to the Gross--Zagier--Zhang formula}
Bruinier, Ehlen and Yang expressed heights of higher Heegner cycles
in terms of CM values of Green functions. For odd \(\kappa>1\) their
formula gives the untwisted mixed height below and identifies its
coprime coefficients with those of a cusp form
\cite[Theorem~6.5, Corollary~6.6 and Theorem~6.12]{BEY}.
For the remainder of the introduction let \(\kappa>1\) and let
\(H=(-1)^\kappa\langle\,\cdot\,,\,\cdot\,\rangle_{\rm GS}\) be the
height of admissible extensions, with absolute arithmetic degree
\cite[\S\S2.2.4--2.2.5]{QiuCM}. Choose \(d_0\) fundamental
(possibly \(1\), with trivial character) such that \(D_0=\Delta d_0\)
is negative and fundamental. A root \(\mu_0^2\equiv d_0\bmod4N\)
fixes the auxiliary CM datum \(U\) of index \(a_U=(m_0,\mu_0)\),
where \(m_0=|d_0|/(4N)\). We also write \(d_U=d_0\), \(D_U=D_0\).
The points of \(Z(U)\) have CM by \(\mathbb Q(\sqrt{D_0})\).
In the CM integral, \(U=x_0^\perp\) for a primitive vector of norm
\(|\Delta|m_0\) in component \(r\mu_0\) with \(\chi_\Delta(x_0)\ne0\).
The integral uses both orientations. The cycle \(Z_\kappa(U)=Z_\kappa(a_U)\)
retains the residue label. Write \(\ell_U(g)=b_g(U)=c_g(a_U)\) and
\(S_U(F)=\operatorname{Sh}_{m_0,\mu_0}(\xi_w F)\) as in
\eqref{merged-shimura-coeff}, and put
\[
 C_{\kappa,N,x}=
 \frac{2\sqrt{Nx}\,\Gamma(\kappa-\tfrac12)}
      {(4\pi)^{\kappa-1}\pi^{3/2}}\qquad(x>0).
\]

For the CM integral assume \((\Delta,2N)=1\). Take \(F\) with real
Fourier data, cuspidal shadow and polar support disjoint from the
primitive CM line, as in Propositions~\ref{p19g-cm-derivative}
and~\ref{p19g-shimura}. Write \(j=\lfloor(\kappa-1)/2\rfloor\).
The CM value \(\Phi^j_{\Delta,r}(Z(U),F)\) uses the genus weights.
The boundary \(B_U(F)\) in \eqref{p19g-boundary} is normalized per
unit of the unsigned mass \(\deg Z(U)=4h(D_0)/|\mathcal O_{D_0}^{\times}|\).

\begin{theorem}[CM values and mixed heights]\label{merged-height-theorem}
For either sign of \(\Delta\) and either parity of \(\kappa\),
\[
 \frac{(-1)^\kappa m_0^{(\kappa-1)/2}}{2}
 \bigl(\Phi^j_{\Delta,r}(Z(U),F)-\deg Z(U)\,B_U(F)\bigr)
 = C_{\kappa,N,|\Delta|}\,L'(S_U(F),\chi_\Delta,\kappa).
\]
If \((D_0,2N)=1\), adding the boundary and finite intersections
gives the full height in \eqref{height-residual-finite}. If also
\(\Delta=r=1\) and \((4Nn,D_0)=1\) at each polar index \(n\), then
\begin{equation}\label{med-untwistedheight}
 H(Z_\kappa(F),Z_\kappa(U))=C_{\kappa,N,1}\,
 L'(S_U(F),\kappa).
\end{equation}
\end{theorem}

\begin{introremark}[The two parities]
Substituting \eqref{p19g-central} in \eqref{p19g-stokes} gives the
CM integral identity. Du and Peng state the untwisted height identity
\eqref{med-untwistedheight} for both parities \cite[Theorem~0.6]{DuPeng}.
\begin{enumerate}
\item \emph{Odd \(\kappa\).}
The CM integral uses the Siegel kernel. For the untwisted height,
the finite height cancels the boundary term in the archimedean height
\cite[Theorem~6.5 and Corollary~6.6]{BEY}.

\item \emph{Even \(\kappa\).}
Reversing a CM vector fixes its base point and changes its sign in
the oriented Millson kernel. Keeping both vectors oriented in the
Green pairing \eqref{arch-signed-double-sum} gives the required
cancellation. The linear factor in the unary theta series produces
the same Legendre polynomial as the intersection of the two fibre
classes, with opposite coefficients in the archimedean and finite heights.
\end{enumerate}
\end{introremark}

\begin{introremark}[Twisted mixed heights]
Let \(\Delta\ne1\) be prime to \(2N\). Proposition~\ref{p19-height-comparison}
gives, for either sign of \(\Delta\) and either parity of \(\kappa\),
\[
\begin{aligned}
 H(Z_\kappa(F),Z_\kappa(U))
   &= C_{\kappa,N,|\Delta|}\,L'(S_U(F),\chi_\Delta,\kappa)+R_U(F),\\
 R_U(F)&=\frac{(-1)^\kappa m_0^{(\kappa-1)/2}\deg Z(U)}{2}\,B_U(F)
                 +\sum_p H_p(Z_\kappa(F),Z_\kappa(U)).
\end{aligned}
\]
Here the CM line is primitive and the polar support is disjoint
from \(U\). For fundamental \(D_0\) prime to \(2N\), the oriented
Green pairing gives the archimedean comparison. For other auxiliary
data the proposition assumes this comparison.
\begin{enumerate}
\item \emph{The boundary.} The finite sum \eqref{p19g-boundary}
computes \(B_U(F)\) using both lattice residues in the theta kernel.
A change of genus sign cancels a pair of summands when it preserves
their other factors.
\item \emph{Finite places.} At a good prime split in
\(\mathbb Q(\sqrt{D_0})\), distinct CM fields give disjoint horizontal
closures. At a nonsplit prime, a sign-reversing pairing cancels the
local sum if it preserves the other genus weight, the residue
conditions and the full intersection weights. Any vertical terms
on the arithmetic model also enter \(H_p\)
(Proposition~\ref{finite-genus-cancellation}).
\end{enumerate}
\end{introremark}

We now allow any \(N\ge1\), any fundamental \(\Delta\equiv r^2\bmod 4N\)
and complex Fourier coefficients for \(F\). Let \(t\) be the number
of primes dividing both \(N\) and \(\Delta\). Write \(\Pi_{\Delta,N}\)
for the average over the corresponding residue involutions
(Section~\ref{rd-cycles}). For \((N,\Delta)=1\) this is the identity.

Call a fundamental CM datum \(U\) \emph{admissible for \(F\)} if
\(d_0\) is prime to \(2N\Delta\) and to every \(4Nn\) with a nonzero
polar coefficient \(c_F^+(-n,\mu)\), and if its CM field
\(\mathbb Q(\sqrt{D_0})\) is distinct from each polar CM field
\(\mathbb Q(\sqrt{-4N|\Delta|n})\). Here the polar coefficients are taken before passage to Chow.
Let \(S_{\rm prim}^{(1)}\) be the sum in \(S\) of the good Hecke
eigenspaces spanned by a single primitive Jacobi newform whose
Shimura correspondent has exact level \(N\). Only the operators
away from \(2N\Delta\) enter this definition.

\begin{theorem}[Auxiliary cycles and the mixed height]\label{intro-newform-mixed}
For either sign of \(\Delta\), every fundamental CM datum \(U\)
with \((d_U,2N\Delta)=1\) satisfies, for all harmonic inputs \(F\),
\begin{equation}\label{intro-newform-height}
 H(Z_\kappa(F),Z_\kappa(U))
 =2^t C_{\kappa,N,|\Delta|}\,
 \overline{L'\bigl(S_U(\Pi_{\Delta,N}F),\chi_\Delta,\kappa\bigr)}.
\end{equation}
For any finite-dimensional space of inputs, infinitely many such \(U\)
are admissible for every \(F\) in that space.
\end{theorem}

The global identity allows \(U\) to be fixed before \(F\) is chosen,
even for a self-pairing. Admissibility for \(F\) is needed to compute
the local terms with its given polar representatives. Vanishing at
a particular finite prime requires the disjointness or cancellation
conditions of Section~\ref{rd-h-height}, or a separate local calculation.
For an oldspace
the proof uses a height matrix of cycles with the same CM field
(Proposition~\ref{coprime-toric-block}). Its fundamental column is
determined by the untwisted height. For negative twists the skew
Jacobi basis and its Petersson matrix are computed in
Section~\ref{sec:jacobi-petersson}.

Let \(g\) be a real generator of a retained primitive Hecke line of
exact index \(N\), so that \(\Pi_{\Delta,N}g=g\). Its normalized
Shimura correspondent \(G\) has weight \(2\kappa\) and trivial
character. Choose a real lift with \(\xi_wF_g=g/\|g\|_{\rm Pet}^2\).
Write
\[
 X_g=Z_\kappa(F_g),\qquad b_0=c_g(m_0,\mu_0),\qquad
 A_g=C_{\kappa,N,|\Delta|}\frac{L'(G,\chi_\Delta,\kappa)}
                                  {\|g\|_{\rm Pet}^2}.
\]
The geometric Hecke action gives the self-adjoint projection
\(Z_\kappa(U)^G=e_GZ_\kappa(U)=b_0X_g\)
(Section~\ref{sec:newform-heights}, see also \cite[\S5]{ABS}).
Since \(S_U(F_g)=b_0G/\|g\|_{\rm Pet}^2\), the mixed height is
\(H(X_g,Z_\kappa(U))=2^t b_0A_g\). If an ambient newform has
\(\Pi_{\Delta,N}g=0\), its projected lift is weakly holomorphic.
Then \(X_g=0\) and both sides of \eqref{intro-newform-height}
vanish by Theorem~\ref{intro-chow-vanishing} and
Corollary~\ref{med-antisector}.

Let \(\chi\) be the genus character of \(K=\mathbb Q(\sqrt{D_0})\)
attached to \(D_0=\Delta d_0\).

\begin{corollary}[Gross--Zagier--Zhang for the genus character \(\chi\)]
\label{p19-higher-gzz}
Let \(U\) be a fundamental CM datum with \((d_0,2N\Delta)=1\).
Suppose \(G\) has exact level \(N\) and \(g\) is primitive and retained.
Then
\begin{equation}\label{eq:intro-gzz}
 H(Z_\kappa(U)^G,Z_\kappa(U)^G)=
 \frac{2^t(2\kappa-2)!\sqrt{|D_0|}\,m_0^{\kappa-1}}
 {2^{4\kappa-2}\pi^{2\kappa}\|G\|_{\rm Pet}^2}\,
 L'_K(G,\chi,\kappa).
\end{equation}
\end{corollary}

The height formulas are related as follows.
\[
 \begin{gathered}
  \text{Theorem~\ref{merged-height-theorem}}\\[-2pt]
  \scriptstyle\text{CM integral}
 \end{gathered}
 \xrightarrow[\text{toric heights and Gram matrices}]{\text{untwisted comparison}}
 \begin{gathered}
  \text{Theorem~\ref{intro-newform-mixed}}\\[-2pt]
  \scriptstyle\text{mixed height}
 \end{gathered}
 \xrightarrow[\text{Waldspurger}]{\text{Hecke projection}}
 \begin{gathered}
  \text{Corollary~\ref{p19-higher-gzz}}\\[-2pt]
  \scriptstyle\text{self-height}
 \end{gathered}
\]
In the untwisted coprime range the finite intersection cancels
the Stokes boundary. To pass to arbitrary compatible twists, the
first arrow uses this reference height together with the primitive
toric height and the oldspace Petersson matrix
(Theorem~\ref{positive-full-mixed}). For the second arrow,
Hecke projection gives \(Z_\kappa(U)^G=b_0X_g\), and Waldspurger's
formula evaluates \(b_0^2\)
(see also \cite[Proposition~0.9]{DuPeng} for the untwisted case).

\begin{introremark}[Choosing \(U\) and recovering the heights]\leavevmode
\begin{enumerate}
\item \emph{Existence.}
Choose an odd prime \(q\) avoiding \(2N\Delta\) and the polar data
of \(F\), with \(d_0=-\operatorname{sgn}(\Delta)q\equiv1\bmod 4N\).
Then \(\mu_0=1\) gives an admissible index. For a primitive \(g\),
\(\gg X/\log X\) admissible discriminants \(|d_0|<X\) have
\(b_g(U)\ne0\) \cite[Corollary~3]{OS}
(Lemma~\ref{signed-fundamental-norm}). At squarefree level every
nonzero cusp form has \(\gg X^\eta\) such discriminants for each
\(0<\eta<1\), with the coefficient bound of
Lemma~\ref{fundamental-separation}. In a retained oldspace of
dimension \(s\) and nonzero derivative, fundamental height tests
have a common kernel of dimension \(s-1\), independent of the
compatible twist with nonzero derivative. Imprimitive indices
supply the other directions (Corollary~\ref{twisted-oldspace-rank}).

\item \emph{One height on a line.}
On a retained Hecke line it suffices that
\(\operatorname{Sh}_{m_0,\mu_0}(g)=b_0G\),
\(L(G,\chi_\Delta,\kappa)=0\), and \eqref{eq:gzz-waldspurger}
hold to deduce \eqref{eq:intro-gzz} from \eqref{intro-newform-height}.
The deduction includes \(D_0=-3,-4\) whenever the congruence and
coprimality conditions hold. The Fourier coefficient already carries
the CM stabilizer weight.
Conversely, when \(b_g(U)\ne0\) the projected self-height divided
by \(b_g(U)\) gives the mixed height. For
\(\mathbb Cg\subset S_{\rm prim}^{(1)}\) this determines the
derivative and, for all Fourier indices \(v,w\), gives
\[
 H(e_GZ_\kappa(v),e_GZ_\kappa(w))
 =b_g(v)b_g(w)\,\frac{H(X_g,Z_\kappa(U))}{b_g(U)} .
\]
This includes coincident CM fields and zero derivatives. If
\(S=\mathbb Cg\) it gives the full pairing. Across several primitive
lines the scaled derivatives are generalized eigenvalues of the
height matrix (Theorem~\ref{primitive-height-pencil} and
Corollary~\ref{primitive-height-count}).

\item \emph{Zero coefficients.}
On such a line \(b_g(U)=0\) gives a zero projected cycle. When \(N\) is squarefree,
\((\Delta,2N)=1\) and
\(\operatorname{sgn}(\Delta)(-1)^{\kappa-1}=-1\), reflection also gives
indices at which the whole cycle and the Shimura lift vanish. Here
\(D_0\) may meet \(2N\) (Section~\ref{sec:primitive-pairings}).
\end{enumerate}
\end{introremark}

\subsection{The arithmetic inner product on a Kuga--Sato variety}
On the regular fine models of Appendix~\ref{sec:ks-integral}, each
Heegner coefficient has a unique admissible arithmetic representative
in the cyclic Scholl summand of character \(\chi_\Delta\). These lifts
give an arithmetic Chow-valued cusp form \(\widehat\Phi_\kappa\).
Write \(\widehat\theta_{\Delta,r,\kappa}(f)=(\widehat\Phi_\kappa,f)_{\rm Pet}\).
We take absolute arithmetic degree and divide by the degree of the
fine cover. Both theta lifts are conjugate-linear. Since the Hermitian
intersection and Petersson products are linear in the first variable,
the input \(h\) occurs first on the right of \eqref{eq:ks-two-lifts}.

\begin{theorem}[Arithmetic inner product on a Kuga--Sato variety]
\label{ks-arithmetic-inner-product}
Let \(\kappa\ge2\), and let \(g\) span a retained primitive line
of exact level \(N\) corresponding to the newform \(G\). Then
\[
 (-1)^\kappa
 \langle\widehat\theta_{\Delta,r,\kappa}(g),
        \widehat\theta_{\Delta,r,\kappa}(g)\rangle_{\rm GS}
 =2^{\omega((N,\Delta))} C_{\kappa,N,|\Delta|}\,
        L'(G,\chi_\Delta,\kappa)\,\|g\|_{\rm Pet}^2 .
\]
More generally, if \(f,h\) lie in the same retained eigenspace
of primitive level \(M\mid N\), there is an explicit positive
Petersson self-adjoint operator \(E_{G,N}\) such that
\begin{equation}\label{eq:ks-two-lifts}
 (-1)^\kappa
 \langle\widehat\theta_{\Delta,r,\kappa}(f),
        \widehat\theta_{\Delta,r,\kappa}(h)\rangle_{\rm GS}
 =2^{\omega((M,\Delta))} C_{\kappa,N,|\Delta|}\,
       L'(G,\chi_\Delta,\kappa)\,(E_{G,N}h,f)_N .
\end{equation}
Distinct good Hecke eigenspaces are orthogonal. The operator is
independent of the compatible twist retaining the eigenspace.
It factors over the primes dividing \(N/M\) as in
\eqref{eq:ks-local-tensor}. For squarefree \(N/M\) it equals
\(2^{\omega(N/M)}\) times the identity.
\end{theorem}

On a primitive line Theorem~\ref{intro-newform-mixed} determines the
inner product from one nonzero fundamental coefficient. Conversely,
the inner product determines every projected self-height on that line.
At a fundamental index, Waldspurger's formula then gives
Corollary~\ref{p19-higher-gzz}. The normalization in
Section~\ref{sec:mixed-inner-product} uses the toric height formulas of
Zhang and Qiu.

On an oldspace all fundamental indices test the same column of the
height matrix. Imprimitive indices in one CM field determine the
relative entries of the full matrix and hence the operator \(E_{G,N}\)
in \eqref{eq:ks-old-operator} (Proposition~\ref{coprime-toric-block}).
Once these entries are known, a single projected self-height at a
fundamental index with \(b_G(U)\ne0\) fixes every pairing on the block,
even when the central derivative vanishes
\eqref{eq:calibrated-projected-heights}.

The intersection number is independent of the admissible
representatives. To lift the series itself one must also lift its
rational equivalences. At good primes the two Frobenius branches
of the isogeny correspondence exclude vertical classes. At level
primes the quadratic character permits descent along the Igusa
cover to a good fibre at lower level. Proposition~\ref{ks-admissible-restriction}
uses this to identify the zero-curvature arithmetic classes.
Corollary~\ref{higher-multiquadratic-determinant} gives the height
determinant over a totally real multiquadratic field.

Related intersection formulas have been proved for special divisors
and arithmetic theta lifts on unitary Shimura varieties
\cite{BHY,LL1,LL2}. Over function fields, Yun and Zhang relate
intersections of shtuka cycles to higher derivatives \cite{YunZhang}.

\subsection{Vanishing of individual coefficients}
For the trivial twist every coefficient vanishes for odd \(\kappa\ge3\) with
\(\kappa+2N\le9\), and even \(\kappa\ge2\) with \(N=1\) or
\(\kappa+2N\le12\) (Corollary~\ref{bounded-graph-rows}). At level
one there is a uniform criterion: if
\(S_{2\kappa}(\mathrm{SL}_2(\mathbb Z))=0\), every signed class from
a smooth CM fibre has zero pushforward, for all compatible twists.
Indeed the identity acts through this cuspidal Hecke algebra by
Proposition~\ref{marked-cm-spectrum}. For \(2\le\kappa\le5\) the
same vanishing follows directly from the Chow groups of pointed
elliptic curves (Theorem~\ref{pointed-vanishing}).

Symmetries of the residues give a different cancellation. If a
shared-prime involution sends
\(f\) to \(-f\), the residue terms in \(Z_\kappa(f)\) cancel in pairs
(Corollary~\ref{med-antisector}). If \(N\mid r\) and
\(\operatorname{sgn}(\Delta)(-1)^{\kappa-1}=-1\), the terms for
\(x\) and \(-x\) cancel directly in \eqref{rd-coeff}.

\subsection{Organization of the paper}
The proof of modularity occupies Sections~\ref{rd-cycles}--\ref{rd-m-twist}.
Section~\ref{rd-cycles} gives the principal-part criterion.
Section~\ref{rd-m-untwisted} restricts full-level isogeny graphs
to cyclic level, reducing the untwisted case to \(N+1\) or \(N-1\)
scalar series according to parity. The Kummer surface calculation
in Section~\ref{rd-g-section} yields a recurrence at the four branch
fibres. Section~\ref{med-u-five} solves it using the characters of a
two-torsion frame. Section~\ref{rd-m-twist} treats the twists by
Hecke transfer and proves modularity of the Abel--Jacobi images.

Section~\ref{rd-h-lift} proves the arithmetic inner product formula
on the modular curve and computes its multiquadratic height determinant.
For higher cycles, Section~\ref{p19-cm-section} computes the CM value
and Stokes boundary of the regularized theta lift.
Section~\ref{rd-h-height} compares the boundary with the finite
intersections and computes the primitive height by toric periods.
Section~\ref{merged-conclusion} gives the squarefree height series,
the oldspace height matrices at arbitrary cyclic level for compatible
twists, and the Kuga--Sato arithmetic inner product formula.
Appendix~\ref{med-all-signs} constructs rational bases of weakly
holomorphic forms of both Weil types. Appendix~\ref{app-fibres}
computes the fibre boundaries and identifies the admissible
arithmetic representatives on the compact Scholl summand.

\section{CM cycles and the principal-part pairing}\label{rd-cycles}

\subsection{The geometric coefficients}
On the space of rational trace-zero matrices put \(Q(x)=N\det x\).
Take the lattice and discriminant module
\begin{align*}
 L&=\left\{\begin{pmatrix} b&-a/N\\ c&-b\end{pmatrix}:a,b,c\in\mathbb Z\right\},\\
 A_L&=L^\sharp/L\simeq\mathbb Z/2N\mathbb Z,\qquad
 Q(\mu)=-\frac{\mu^2}{4N}\pmod{\mathbb Z}.
\end{align*}
The dual is taken with the bilinear form
\((x,y)=Q(x+y)-Q(x)-Q(y)\). The component \(\mu\) is represented by
the diagonal matrix with entries \(\mu/(2N),-\mu/(2N)\).
The space has signature \((1,2)\). We use the Weil representation
for which \(T\) acts on \(e_\mu\) by \(e^{2\pi iQ(\mu)}\).
Set \(\widetilde\rho_L=\rho_L\) for \(\Delta>0\) and
\(\widetilde\rho_L=\bar\rho_L\) for \(\Delta<0\).
The input representation is \(\rho_{\rm in}=\overline{\widetilde\rho_L}\).

To define the cycle on a fibre, work on a smooth fine-level cover of
the modular curve and its smooth Kuga--Sato compactification as in
\cite{SchollMotives,Zhang}. At cyclic level we take the invariant
part of the compact Scholl summand on a normal full-level cover
under the group preserving the cyclic subgroup. Pullback and
degree-normalized trace identify these summands on a common
refinement, as in Appendix~\ref{app-fibres}. For the CM coefficients take geometric
Chow groups over \(\overline{\mathbb Q}\) with complex coefficients.
Each finite collection of their defining cycles descends to a common
number field. The height calculations below use such a field and
divide arithmetic degrees by its degree over \(\mathbb Q\).
For a CM elliptic curve \(E\) let \(u\)
be the trace-zero endomorphism of the chosen orientation, after
clearing denominators. The mixed part of its graph is the divisor on
\(E^2\) obtained by removing the two axis classes. For multiplication
by \(\sqrt D\), \(D<0\), our convention is
\[
       [\Gamma_{\sqrt D}]-[E\times0]+D[0\times E].
\]
This divisor has intersection zero with both axes. As a pullback
of the Poincar\'e bundle, with the sign of our graph convention, its
class in \(\Pic(E^2)_{\mathbb Q}\) is additive in \(u\).

Put \(a=\kappa-1\). On \(E^{2a}\) take the external product of
\(a\) mixed divisors and project to the part odd in each coordinate
and alternating under permutations. On cohomology the signs from
the degree-one factors give the symmetric power. With the chosen
orientation, normalize the fibre self-intersection to have degree
\((-1)^a\) \cite[\S6.3]{BEY}. Pushing forward from the fibre,
a divisor in the total space, gives \(S_\kappa(x)\) of codimension
\(a+1=\kappa\), with \(S_\kappa(-x)=(-1)^a S_\kappa(x)\).

Here \(S_\kappa(x)\) represents a single oriented lift from the
normal fine cover. For a torsion-free normal subgroup
\(\Gamma'\subset\Gamma=\Gamma_0(N)\), the orbit of \(x\) splits
into \([\Gamma:\Gamma']/|\Gamma_x|\) fine orbits. Multiplying their
pushforward by \(2/[\Gamma:\Gamma']\) gives the orbit weight
\(2/|\Gamma_x|\).

Write \(L_\nu[t]=\{x\in\nu+L:Q(x)=t\}\), let \(\Gamma_x\) be the
stabilizer in \(\Gamma_0(N)\), and let \(\chi_\Delta\) on vectors
denote the genus character with the usual value zero on inadmissible
vectors \cite{BO}. For \(n>0\) congruent to
\(\operatorname{sgn}(\Delta)Q(\mu)\bmod\mathbb Z\), define
\begin{equation}\label{rd-coeff}
 Z^{(N)}_{\Delta,r,\kappa}(n,\mu)=n^{a/2}
 \sum_{x\in\Gamma_0(N)\backslash L_{r\mu}[n|\Delta|]}
          \frac{2\chi_\Delta(x)}{|\Gamma_x|}\,S_\kappa(x).
\end{equation}
The normalization in \eqref{rd-coeff} has a common scalar for all
Fourier indices of a given twist. Indeed, \(u_x=2Nx\) is an integral endomorphism
of the CM elliptic curve: its matrix preserves the period lattice and
has determinant \(4N|\Delta|n\). Let \(D_{u_x}\) be its mixed graph
divisor, \(e\) the signed projector on the fibre and \(i_x\) the fibre
inclusion. With the orientation fixed above, put
\begin{equation}\label{eq:common-chow-scalar}
 c_{a,N,\Delta}=\frac{\sqrt{\binom{2a}{a}}}{2^a(4N|\Delta|)^{a/2}},
 \qquad
 n^{a/2} S_\kappa(x)
 = c_{a,N,\Delta}\,i_{x*}e\bigl(D_{u_x}^{\boxtimes a}\bigr).
\end{equation}
To compute the scalar, choose a symplectic eigenbasis of \(H^1(E)\)
for \(u_x\), with eigenvalues \(\lambda,-\lambda\), so
\(\lambda^2=-\deg u_x\). The mixed graph tensor is, up to a common
sign, \(\lambda(\alpha\otimes\beta+\beta\otimes\alpha)\).
The signed projector symmetrizes the \(2a\) degree-one factors in
cohomology. Every word with \(a\) copies of each basis vector then
has the same coefficient, of magnitude
\(|\lambda|^a2^a/\binom{2a}{a}\). Only complementary words pair.
The Koszul sign in the product and the \(a\) reversed symplectic
pairings cancel, leaving self-intersection
\((-1)^a4^a(\deg u_x)^a/\binom{2a}{a}\). This proves the formula.

After division by this common scalar, every orbit coefficient belongs
to the rational Chow group over \(\overline{\mathbb Q}\). In particular,
if the principal coefficients of \(f\) lie in a subfield \(E_0\subset
\mathbb C\), its normalized relation in Theorem~\ref{intro-chow-vanishing}
already holds with coefficients in \(E_0\): extension of scalars from
\(E_0\) to \(\mathbb C\) is injective. The genus character
satisfies \(\chi_\Delta(-x)=\operatorname{sgn}(\Delta)\chi_\Delta(x)\):
changing the sign of a represented integer gives precisely the
factor \(\chi_\Delta(-1)\). Hence the component symmetry is
\[
 Z_\kappa(n,-\mu)=\operatorname{sgn}(\Delta)(-1)^a Z_\kappa(n,\mu).
\]
This agrees with the metaplectic central character at weight \(k\).
A self-opposite component with negative sign has zero coefficient.
For \(\kappa=1\), formula \eqref{rd-coeff} gives a divisor on the
base. To pass to the Jacobian one subtracts its degree times
the cusp \(P_\infty\).

There is a second symmetry when \(\Delta\) and \(N\) share a prime.
Put \(J=\{p:p\mid(\Delta,N)\}\). For an odd \(p\in J\),
compatibility gives \(p\mid r\). If \(p^2\mid N\) it would then give
\(p^2\mid\Delta\), contrary to fundamentality. If \(2\in J\),
\(\Delta=4d\) with \(d\equiv2,3\bmod4\). Thus \(\Delta\) is \(8\)
or \(12\bmod16\), whereas the squares modulo \(16\) are \(0,1,4,9\).
Hence \(4\nmid N\), and \(p\Vert N\) for every \(p\in J\).
The congruences
\[
  \sigma_p=-1\pmod{2p},\qquad \sigma_p=1\pmod{2N/p}
\]
define a unit modulo \(2N\). They imply
\(\sigma_p^2=1\bmod 4N\), so multiplication by \(\sigma_p\)
preserves the discriminant form. They also give
\(r\sigma_p\mu=r\mu\bmod 2N\), since \(p\mid r\).
Consequently the two sums \eqref{rd-coeff} have identical summands.
Let
\[
             \Pi_{\Delta,N}=\prod_{p\in J}\frac{1+\sigma_p}{2}.
\]
This real projector commutes with the Weil action and with the shadow
operator. Both \(\Phi_\kappa\) and the test space may be restricted to
its image. The product is the identity if \(J\) is empty.

\begin{corollary}[Relations at shared primes in ordinary Chow]\label{med-antisector}
For every input \(F\) one has \(Z_\kappa(F)=Z_\kappa(\Pi_{\Delta,N}F)\).
In particular, any weakly holomorphic input on which at least one
shared-prime involution acts by \(-1\) gives a rationally trivial
cycle in the full compact Chow group, for every pole order.
\end{corollary}
\begin{proof}
The coefficient sum of \(F-\sigma_pF\) vanishes term by term after
reindexing the residue by \(\sigma_p\), by the equality of vector
cosets in \eqref{rd-coeff}. Apply the commuting averages in succession.
\end{proof}

\subsection{The principal part and the Petersson pairing}
For the remainder of this section assume \(\kappa\ge2\), so \(w<0\).
For \(\kappa=1\) the kernel of taking the strictly negative part on
\(M^!_{1/2,\rho_{\rm in}}\) is \(M_{1/2,\rho_{\rm in}}\); the
divisor argument in Section~\ref{rd-h-lift} uses the Jacobian-valued
series directly.
Use the Petersson inner product linear in its first entry. For
\(\tau=u+iv\), the conjugate-linear operator
\[
 \xi_w F=2iv^w\,\overline{\frac{\partial F}{\partial\bar\tau}}
\]
takes a harmonic weak Maass form with cuspidal shadow to \(S_{k,\widetilde
\rho_L}\). Write \(P_F\) for the strictly negative terms of the holomorphic part
\(F^+\). The constant term does not pair with cusp forms. A weakly
holomorphic form is holomorphic on the upper half-plane and meromorphic
at the cusps. By the exact sequence of Bruinier and Funke, \(\xi_w\)
is onto the cusp space and its kernel consists of the weakly holomorphic
forms \cite[Theorem~3.7]{BFGeometricTheta}. We work in the invariant
sector from now on, since \(Z_\kappa\) kills its complement.

Let \(\mathcal P\) be the vector space of finite polar polynomials
of weight \(w\) and type \(\rho_{\rm in}\) with the indicated
symmetries. Denote by \(S^\vee\) the complex-linear dual. The coefficient
functionals define
\begin{equation}\label{rd-principal}
 \beta:\mathcal P\longrightarrow S^\vee,\qquad
 \beta(P)(g)=\sum_{n>0,\mu}c_P(-n,\mu)c_g(n,\mu).
\end{equation}
This coefficient pairing is complex bilinear. In its Petersson
description below, conjugation in the second entry cancels the
conjugate-linearity of \(\xi_w\).

\begin{proposition}[Which principal parts are weakly holomorphic?]\label{med-exact}
In the invariant sector there is an exact sequence of complex vector spaces
\[
 0\longrightarrow M^!_{w,\rho_{\rm in}}{}^\Pi
 \xrightarrow{\,f\mapsto P_f\,}
 \mathcal P \xrightarrow{\ \beta\ } S^\vee\longrightarrow0.
\]
\end{proposition}
\begin{proof}
Given \(P\in\mathcal P\), choose a harmonic lift \(F\) with polar
part \(P\) and average it under \(\Pi\). The principal-part pairing gives
\[
                 \beta(P)(g)=\langle g,\xi_w F\rangle_{\rm Pet}
\]
\cite[Propositions~3.5 and~3.11]{BFGeometricTheta}. Thus \(\beta(P)=0\)
exactly when \(F\) is weakly holomorphic. A weakly holomorphic form
with zero polar part is holomorphic at every cusp. On a finite cover
trivializing the type, a suitable positive power of each component is a
section of a negative power of the Hodge bundle. The Hodge bundle has positive
degree by the divisor of the discriminant, so each component vanishes.
Finally the coefficient functionals span \(S^\vee\), since their common
kernel consists of cusp forms with zero Fourier expansion.
\end{proof}

\subsection{The obstruction in the ordinary Chow group}
\begin{theorem}[The obstruction and a bounded choice of coefficients]\label{rd-main}
Let \(W_N\subset A_N\) be the span of the coefficients. For
\(f\in M^!_{w,\overline{\widetilde\rho_L}}\) put
\[
 Z_\kappa(f)=\sum_{\mu,n>0}c_f(-n,\mu)Z_\kappa(n,\mu),
 \qquad D_N=\{Z_\kappa(f): f\in M^!_{w,\overline{\widetilde\rho_L}}\}.
\]
Then the image of \(\Phi_\kappa\) in \(W_N/D_N\) is a cusp form of weight
\(k\). A linear map from \(W_N\) to any complex vector space takes
\(\Phi_\kappa\) to a cusp form precisely when it kills \(D_N\).
If \(s=\dim S\), one may choose \(s\) coefficient positions with
\(0<n\le k/12\) which span this quotient and determine all its Fourier
coefficients. Subtracting the corresponding interpolation formulas in
\(W_N\) gives a generating set of \(D_N\). In particular,
Chow-valued modularity is equivalent to \(D_N=0\).
\end{theorem}
\begin{proof}
Extend \eqref{rd-coeff} linearly to \(Z_\kappa:\mathcal P\to W_N\)
by summing each polar coefficient against its cycle. Thus
\(Z_\kappa(F)=Z_\kappa(P_F)\) for a harmonic lift \(F\).
To see that \(Z_\kappa\)
is onto, project \(q^{-n}e_\mu\) to the allowed central symmetry and
the \(\Pi\)-invariant sector. Denote the resulting polynomial by \(p_i\)
where \(i=(n,\mu)\). The averaging idempotents are self-transpose for
the component pairing. The symmetries already proved therefore give
\[
       Z_\kappa(p_i)=Z_\kappa(n,\mu),\qquad \beta(p_i)(g)=c_g(n,\mu).
\]
For a self-opposite component on which the sign is negative both sides
are zero. Choose one signed representative of each remaining orbit,
ordered by exponent and then by residue. These \(p_i\)'s form a basis
of \(\mathcal P\).

For a linear map \(b:W_N\to V\), the Fourier values
\(bZ_\kappa(p_i)\) are those of a tensor in \(S\otimes V=\Hom(S^\vee,V)\)
precisely when \(bZ_\kappa\) factors through \(\beta\).
Since \(\beta\) is onto, this is equivalent to
\(bZ_\kappa(\ker\beta)=0\). Proposition~\ref{med-exact} and the
invariance of the cycles identify \(Z_\kappa(\ker\beta)=D_N\).
For the quotient map \(b:W_N\to W_N/D_N\), denote the resulting map
by \(T:S^\vee\to W_N/D_N\). It is onto because the coefficients span
\(W_N\), giving \(\dim(W_N/D_N)\le\dim S\).

By Lemma~\ref{med-eta} below, choose \(p_{a_1},\ldots,p_{a_s}\) of
exponent at most \(k/12\) whose images form a basis of \(S^\vee\).
For each remaining \(i\), subtract from \(p_i\) the unique linear
combination of the \(p_{a_j}\) with the same image under \(\beta\).
These differences form a basis of \(\ker\beta\). Indeed, subtracting
them with the prescribed coefficients from any relation leaves a
relation among the basis vectors \(\beta(p_{a_j})\). Applying
\(Z_\kappa\) gives the generators of \(D_N\), and passing to the
quotient gives the interpolation formulas.
\end{proof}

\subsection{Choosing the Fourier coefficients}\label{med-basis}

The relations in Theorem~\ref{rd-main} admit a basis indexed by Fourier
positions. The following bound makes the choice independent of the
pole order of the input.

\begin{lemma}[Low coefficients separate cusp forms]\label{med-eta}
The functionals \(g\mapsto c_g(n,\mu)\), taken over all components
\(\mu\) and all Fourier exponents \(0<n\le k/12\) of the original
vector-valued series, span \(S^\vee\).
\end{lemma}
\begin{proof}
If all these coefficients vanish, \(g/\eta^{2k}\) tends to zero
at infinity since \(\eta^{2k}\) has order \(k/12\). The quotient
has weight zero and finite-image type. On a finite cover
trivializing the type, each component is holomorphic and tends to
zero at every cusp by the transformation law. It is therefore zero
by compactness. Duality gives the assertion.
\end{proof}

Let \(p_i\) run through the signed orbit representatives in the
order fixed above. Put \(\ell_i=\beta(p_i)\) and \(z_i=Z_\kappa(p_i)\).
Retain \(\ell_i\) whenever it is independent of its predecessors.
By Lemma~\ref{med-eta}, the resulting basis
\(\ell_{a_1},\ldots,\ell_{a_s}\) is complete by exponent \(k/12\),
where \(s=\dim S\). For each \(i\), write
\[
               \ell_i=\sum_{j=1}^s u_j(i)\ell_{a_j}.
\]
In any basis \(g_1,\ldots,g_s\) of \(S\), the coordinates \(u(i)\)
solve \(B u(i)=(c_{g_t}(i))_t\), with \(B_{tj}=c_{g_t}(a_j)\).
Changing basis multiplies both sides by the same invertible matrix.

\begin{theorem}[The complete family of cycle tests]\label{med-complete-basis}
For each index \(i\) outside the selected set there is a unique
weakly holomorphic form \(f_i\) in the invariant sector with
\begin{equation}\label{med-pivot}
          P_{f_i}=p_i-\sum_{j=1}^s u_j(i)p_{a_j}.
\end{equation}
The coefficients \(u_j(i)\) are rational. The forms \(f_i\) have
rational Fourier coefficients and form a basis of
\(M^!_{w,\rho_{\rm in}}{}^\Pi\). Their images in the compact Chow
group are
\begin{equation}\label{med-remainders}
          R_i=z_i-\sum_{j=1}^s u_j(i)z_{a_j}.
\end{equation}
One has \(D_N=\sum_i\mathbb C R_i\). In particular, the modularity assertion of Theorem~\ref{med-target} is equivalent to the
simultaneous vanishing of \eqref{med-remainders}.
\end{theorem}
\begin{proof}
Let \(P\in\ker\beta\). For each nonselected \(i\), subtract its
coefficient in \(P\) times the right side of \eqref{med-pivot}.
Only finitely many subtractions occur. What remains is supported
at the selected positions and pairs to zero. Independence of the
selected Fourier functionals makes this remainder zero. Conversely,
each row in \eqref{med-pivot} pairs to zero and has coefficient one
at its own nonselected position and zero at every other such position.
These rows form a basis of \(\ker\beta\). Proposition~\ref{med-exact}
identifies their unique weakly holomorphic lifts, and applying
\(Z_\kappa\) gives \eqref{med-remainders}.

McGraw's theorem applied to the lattice or its negative gives a
rational Fourier basis for the chosen Weil type \cite{McGraw}.
The constant-term conditions for cuspidality and the projector
\(\Pi\) are rational. We may thus compute \(u(i)\) in a rational
basis of \(S\), which makes each \(P_{f_i}\) rational.

In fact a weakly holomorphic form \(h\) of weight \(w<0\) is defined
over the field \(E_0\subset\mathbb C\) of its polar coefficients.
For an integral pole bound \(t\), multiplication by
\(\Delta_{\rm mod}^t\) identifies the bounded space with holomorphic
forms of weight \(w+12t\). Since \(\Delta_{\rm mod}=q+O(q^2)\)
has rational coefficients, McGraw's theorem supplies a rational
basis of the bounded space. The principal-part matrix in this basis
is injective because a holomorphic form of negative weight is zero.
Solving for \(h\) over \(E_0\) proves the claim. Take \(h=f_i\).
\end{proof}

The basis respects every pole bound because \(u_j(i)=0\) for
\(a_j\) later than \(i\). If \(s=0\), each orbit monomial gives a
test. Proposition~\ref{filtered-weak-basis} gives the least pole
bound for generation under multiplication by the modular invariant.

\begin{corollary}[Interpolation in the Chow group]\label{med-candidate}
There is a unique \(\Psi\in S\otimes W_N\) whose coefficients at the
selected positions agree with the original cycles. At position \(i\)
its coefficient is \(\sum_j u_j(i)z_{a_j}\). The difference
\(\Phi_\kappa-\Psi\) has coefficients \(R_i\), and their span is independent
of the selected basis of functionals.
\end{corollary}
\begin{proof}
Send \(\ell_{a_j}\) to \(z_{a_j}\) and extend linearly on \(S^\vee\).
Tensor duality gives \(\Psi\). At an arbitrary position, evaluate this
map on the unique expansion of \(\ell_i\). The difference spans
\(Z_\kappa(\ker\beta)\), which was defined before the choice.
\end{proof}

\medskip
\noindent\hypertarget{integral-reconstruction}{\textit{Integral reconstruction.}}
Fix \(N,\kappa\) and one of the two Weil types. Write \(M=4N\),
\(\psi(m)=[{\rm SL}_2(\mathbb Z):\Gamma_0(m)]\), and set
\[
 B_{\mathbb Z}=\frac{k}{12}\max_{p\mid M}\psi(p^{v_p(M)}).
\]
Every Fourier functional is an integral combination of those of
exponent at most \(B_{\mathbb Z}\). Consequently the same integers
reconstruct every normalized coefficient \(Z_\kappa(i)/c_{a,N,\Delta}\)
in rational Chow from coefficients in this range, simultaneously
for all \(\Delta,r\) of the fixed type. Over \(\mathbb Z[1/(2N)]\)
the smaller range \(n\le k/12\) suffices.

Indeed, take the full cusp space of the chosen type before applying
\(\Pi_{\Delta,N}\). Its forms with integral Fourier coefficients
constitute a full lattice \(L_{\mathbb Z}\) by McGraw's theorem
and Lemma~\ref{med-eta}. Fix a prime \(p\) and put \(M_p=p^{v_p(M)}\).
On \(\Gamma_0(M_p)\) the Weil matrices and their inverses are integral
at a prime above \(p\) of a common coefficient field. In the explicit
matrix formula the only denominator is a square root of
\(|D|/|D[c]|\), where \(c\) is the lower left entry, \(D=L^\sharp/L\),
and \(D[c]\) is its \(c\)-torsion. This index is prime to \(p\)
because \(D[c]\) contains the \(p\)-primary subgroup
\cite[Theorem~6.4 and Lemma~5.6]{StrombergWeil}.

Suppose \(g\in L_{\mathbb Z}\) has nonzero reduction modulo \(p\).
Take the finite union of the orbits of its components under the
weight-\(k\) slash action of the inverse image of \(\Gamma_0(M_p)\).
The eighth powers of the resulting functions are permuted in
integral weight \(8k\). If \(j\) of these powers have nonzero
reduction, their product is the reduction of the \(j\)-th elementary
symmetric function of the whole family. This scalar form has
weight \(8kj\) on \(\Gamma_0(M_p)\), with integral exponents
and coefficients integral at the chosen prime. Were every component
of \(g\) zero through exponent \(k\psi(M_p)/12\) modulo \(p\),
each nonzero factor in that product would have order strictly
greater than \(8k\psi(M_p)/12\). This contradicts Sturm's
congruence bound \cite{SturmCongruence}. Thus the indicated
functionals separate \(L_{\mathbb Z}/pL_{\mathbb Z}\) at every
prime, including those dividing \(2N\). By duality they generate
\(L_{\mathbb Z}^{\vee}\) integrally through \(B_{\mathbb Z}\), and
away from \(2N\) through \(k/12\). Restricting these identities
to each invariant sector and applying Theorem~\ref{med-target}
gives the asserted Chow relations with the normalization
\eqref{eq:common-chow-scalar}.

\subsection{The cycle on a newform line}\label{rd-components}
Choose an orthogonal basis \(g_1,\ldots,g_s\) of \(S\), and choose
harmonic lifts \(F_j\) with
\(\xi_w F_j=g_j/\langle g_j,g_j\rangle_{\rm Pet}\). Write
\(\overline z\) for the image of \(z\) in \(W_N/D_N\).

\begin{corollary}[Normalized shadow expansion]\label{rd-expansion}
One has
\begin{equation}\label{rd-basis}
    \sum_{n,\mu}q^n e_\mu\otimes\overline{Z_\kappa(n,\mu)}
       =\sum_{j=1}^s g_j\otimes\overline{Z_\kappa(F_j)} .
\end{equation}
The cycle on the right is independent of the choice of \(F_j\) in
this quotient. The same assertion holds in ordinary Chow if
\(D_N=0\).
\end{corollary}
\begin{proof}
The difference between two lifts is weakly holomorphic, so its cycle
maps to zero. The functionals \(\beta(P_{F_j})\) are the basis of
\(S^\vee\) dual to the \(g_i\)'s by the Petersson formula. Applying
the map \(T\) from the proof of Theorem~\ref{rd-main} to these
functionals gives \eqref{rd-basis}.
\end{proof}

A good Hecke eigenspace contributes the sum of its basis terms
to \eqref{rd-basis}. One term gives the \(G\)-component when
the Shimura eigenspace is a line. In higher codimension
Section~\ref{sec:newform-heights} realizes the projector on the
Chow coefficients by cyclic isogenies of degree prime to \(2N\Delta\).

\subsection{The boundary on the CM fibres}\label{med-localization}
Work on a smooth compactification \(X\) at neat fine level, with
characteristic-zero coefficients and finite-group projector \(e\).
Let \(T\) be a finite stable set of smooth CM base points, let
\(i_T:S_T=\coprod_{x\in T}X_x\hookrightarrow X\) be its fibre union,
and put \(U_T=X\setminus S_T\). For a principal part \(P\) supported
there, write \(b_T(P)\in e\CH^{\kappa-1}(S_T)\) for the weighted
tuple in \eqref{rd-coeff} before pushing forward. Its image
\(i_{T*}b_T(P)\) is the cycle on \(X\).

\begin{proposition}[Support and rational equivalence]\label{med-boundary}
The cycle \(i_{T*}b_T(P)\) vanishes in \(e\CH^\kappa(X)\) if and
only if \(b_T(P)=\partial_Tu\) for some \(u\in e\CH^\kappa(U_T,1)\).
This equivalence also holds after enlarging \(T\) and extending
the tuple by zero.
\end{proposition}
\begin{proof}
Bloch's localization sequence \cite{BlochLoc} for the closed immersion \(i_T\) has
the segment
\[
 \CH^\kappa(U_T,1)\xrightarrow{\partial_T}
 \CH^{\kappa-1}(S_T)\xrightarrow{i_{T*}}\CH^\kappa(X).
\]
The fibres have codimension one, which accounts for the shift.
Taking the \(e\)-summand preserves exactness by equivariance.
Enlarging \(T\) leaves the pushforward of the extended tuple
unchanged, so the same sequence applies.
\end{proof}

A representative of \(u\) is a finite sum of rational functions
\(f_j\) on codimension-\((\kappa-1)\) subvarieties \(V_j\) of
\(U_T\). Their divisors sum to zero on \(U_T\). On taking closures
in \(X\), the remaining divisor represents \(b_T(P)\) on \(S_T\).
Its entries may be nonzero on individual fibres.
Corollary~\ref{bounded-fibre-boundaries} counts the independent
boundary tuples at every integral pole bound \(B\ge k/12\).
For a finite flat change of level of degree \(d\), pushforward
after pullback multiplies by \(d\). Vanishing on the smooth
cover therefore descends to the invariant Scholl summand.

\subsection{The field of a twisted coefficient}
Use the cyclic-level model over \(\mathbb Q\), or the complete set of
determinant components on a common fine cover before descent. Galois
acts on the geometric Chow group and fixes the coefficient scalars.

\begin{proposition}[The Galois character of a Heegner cycle]
\label{coef-galois-character}
For each compatible \((\Delta,r)\), every Fourier index and
\(\sigma\in G_{\mathbb Q}\),
\[
 \sigma Z^{(N)}_{\Delta,r,\kappa}(n,\mu)
       =\chi_\Delta(\sigma)Z^{(N)}_{\Delta,r,\kappa}(n,\mu).
\]
After division by the common scalar of Section~\ref{rd-cycles},
the coefficient descends to the rational Chow group over
\(\mathbb Q(\sqrt\Delta)\). For fixed \(N,\kappa\), the spans for
distinct fundamental twists form a direct sum in geometric Chow.
Any Galois-invariant bilinear height pairs these spaces orthogonally.
\end{proposition}
\begin{proof}
Let \(u_x=2Nx\) on the elliptic curve with cyclic subgroup \(C\).
If \(b_0\) represents \(r\mu\bmod 2N\), the lattice and residue
condition says that \((b_0I+u_x)/2\) is integral and kills \(C\).
Indeed, in the period basis \((1,z_x)\) the matrix of the endomorphism
is \(-Ju_xJ^{-1}\), with \(J=\operatorname{diag}(-1,1)\), and \(C\)
is generated by \(1/N\). This gives \(b\equiv b_0\bmod 2N\).
Galois transport preserves this condition and the norm.

Choose bases of the full Tate modules of the curve and its conjugate
adapted to \(C\), with the same compatible system of Weil pairings.
Transport by \(\sigma\) has matrix \(k\in K_0(N)\) with determinant
\(\chi_{\rm cyc}(\sigma)\). In the convention for \(x\), the conjugating
matrix is \(JkJ^{-1}\) and has the same determinant. The local
genus-character law proved in Section~\ref{p19-translated} multiplies
the weight by \(\prod_{p<\infty}(\Delta,\det k_p)_p=\chi_\Delta(\sigma)\).
The law includes imprimitive vectors and units at \(2\) and at primes
common to \(N,\Delta\). Transport takes the graph of \(u_x\) to that
of its conjugate, preserves the axis corrections, and commutes with
the projector. The stabilizer multiplicity is preserved as well.
Summing proves the character formula, which in codimension one is
\cite[Lemma~5.1]{BO}.

After dividing by the common scalar, take a finite Galois field of
definition containing \(\mathbb Q(\sqrt\Delta)\). The normalized
trace to that quadratic field pulls back to the original geometric
class, since its conjugates over the field agree. For any finite
family of twists, the character idempotents of their multiquadratic
compositum separate the coefficient spaces. Finally, for distinct
characters choose \(\sigma\) on which their product is \(-1\).
Galois invariance equates the mixed height with its negative.
This proof uses no modularity or nondegeneracy of heights.
\end{proof}

Theorem~\ref{med-target} thus has a rational form. Put
\(K=\mathbb Q(\sqrt\Delta)\), and let \(S_{\mathbb Q}\) denote the forms
in \(S\) with rational Fourier coefficients. With the scalar in
\eqref{eq:common-chow-scalar} one has
\[
 c_{a,N,\Delta}^{-1}\Phi_\kappa
 \ \in\ S_{\mathbb Q}\otimes_{\mathbb Q}
 \bigl(\CH^\kappa_{\rm hom}(\mathcal Y_\kappa(N)_K)_{\mathbb Q}
       \bigr)^{\chi_\Delta}.
\]
Indeed, restriction from \(K\) to \(\overline{\mathbb Q}\) on rational
Chow is injective. A geometric rational equivalence is defined over
a finite extension, and its trace recovers the relation over \(K\).
The superscript refers to the action of \(\operatorname{Gal}(K/\mathbb Q)\)
on the descended classes.
The normalized coefficients belong to the indicated character space.
Their rational span is finite-dimensional, since its complexification
is the coefficient span of the modular form. In a rational basis of
this span each coordinate series has rational coefficients and belongs
to \(S\) by Theorem~\ref{med-target}.

\section{Isogeny graphs and the untwisted coefficients}
\label{rd-m-untwisted}

Fix \(N\geq1\), \(a=\kappa-1\geq1\), and initially \(\Delta=r=1\).
A theta basis reduces the principal-part test at every pole order
to finitely many scalar series. We identify their coefficients
by restricting full-level isogeny graphs to the cyclic level.

\subsection{The scalar series}
Write \(Z_\kappa(m,\mu)=Z^{(N)}_{1,1,\kappa}(m,\mu)\) and put
\(\epsilon\in\{0,1\}\) with \(\epsilon\equiv a\bmod 2\). The coefficients
have reflection symmetry \(Z_\kappa(m,-\mu)=(-1)^a Z_\kappa(m,\mu)\): reversing the
CM orientation changes the mixed divisor by a sign and the cycle uses
\(a\) such factors. Equivalently the criterion below starts with any
sequence having this symmetry and the exponents allowed by the Weil
representation. Define
\[
 I=\begin{cases}\{0,\ldots,N\},&\epsilon=0,\\
                  \{1,\ldots,N-1\},&\epsilon=1,\end{cases}
 \qquad h=|I|,\qquad s_i=\epsilon+2i,\quad b_i=a+s_i.
\]
Thus every \(b_i\) is positive and even. Define the homogeneous trace
polynomial by coefficient extraction in a new variable \(x\):
\[
   P_{a,s}(T,D)=[x^s](1-Tx+Dx^2)^{-a-1}.
\]
For \(0\leq i<h\) form a power series \(F_i=\sum_{n\geq1} F_i(n)q^n\)
with coefficients in \(A_N\) by the explicit formula
\begin{equation}\label{med-u-moment}
 F_i(n)=\sum_{\substack{T\in\mathbb Z\\T^2<4Nn}}
 P_{a,s_i}(T,Nn)\,
 Z_\kappa\bigl(n-T^2/(4N),T\bmod 2N\bigr).
\end{equation}
Every sum is finite. The strict inequality is important: it is
exactly the condition that the selected Fourier coefficient of a
weakly holomorphic input is polar.

\begin{theorem}[Scalar criterion]\label{med-u-scalar}
The untwisted cycle series belongs to
\(S_{a+3/2,\rho_L}\otimes A_N\) if and only if for each \(0\leq i<h\)
\[
             F_i\in S_{b_i+2}(\mathrm{SL}_2(\mathbb Z))\otimes A_N .
\]
In particular, the number of scalar series depends only on the level
and the parity of the weight. For \(h=0\) both the allowed input and
the cycle series vanish by reflection symmetry.
\end{theorem}

\subsection{A basis for weakly holomorphic inputs}
The input has weight \(w=1/2-a\) and Weil type
\(\rho_+\) of \(\mu\mapsto\mu^2/(4N)\). For \(R e_\mu=e_{-\mu}\),
the Weil action of the central element \((-I,i)\) is \(-iR\).
Its weight factor is \(i^{2w}=i(-1)^a\), so \(Rf=(-1)^af\).
At this parity put
\(v=1/2+\epsilon\) and take the vector
\[
  \theta_\mu(\tau)=\sum_{T\equiv\mu\ (2N)}
         T^\epsilon q^{T^2/(4N)},\qquad 0^0=1.
\]
Translation and Poisson summation on the residue classes give the
two generators of the Weil action. Thus \(\theta\) has weight \(v\)
and type \(\rho_+\). Its components in \(I\) determine it by reflection.
In odd parity the two self-opposite components, \(0\) and \(N\), vanish.

Let \(D=q\,d/dq\), \(E_2=1-24\sum_{n\geq1}\sigma_1(n)q^n\), and
\(\mathcal D_t=D-tE_2/12\), the Serre derivative on forms of weight
\(t\). Define \(U_0=\theta\) and \(U_{j+1}=\mathcal D_{v+2j}U_j\).
These vectors have weights \(v+2j\). The following determinant also
controls poles of coordinates at every point of the upper half-plane.

\begin{lemma}[The theta Wronskian and the scalar coefficients]\label{med-u-wronskian}
For \(h>0\) the matrix \(M=(U_j(u))_{u\in I,\,0\leq j<h}\) satisfies
\[
              \det M=c\,\eta^{\,2W},
       \qquad W=h(v+h-1),\qquad c\in\mathbb C^\times.
\]
Every weakly holomorphic input \(f\) has a unique expression
\(f=\sum_{j=0}^{h-1}H_jU_j\) with scalar full-level forms
\(H_j\in M^!_{w-v-2j}=M^!_{-b_j}\).
\end{lemma}
\begin{proof}
The leading term of \(\theta_u\) is \(c_uq^{\lambda_u}\), where
\(\lambda_u=u^2/(4N)\) and \(c_u\ne0\). At the even endpoint \(N\)
the two minimal representatives have the same sign. The Serre
recursion gives the coefficient of \(q^{\lambda_u}\) in \(U_j(u)\)
as \(c_u\prod_{k=0}^{j-1}(\lambda_u-(v+2k)/12)\).
The polynomials in \(\lambda_u\) are monic of successive degrees,
so their determinant is Vandermonde and is nonzero. Summing
\(u^2/(4N)\) for \(0\le u\le N\) in even parity gives
\((N+1)(2N+1)/24\). For \(1\le u<N\) in odd parity it gives
\((N-1)(2N-1)/24\). Both expressions equal \(W/12\).
Thus \(\det M/\eta^{2W}\) has
weight zero, finite-order character and a finite nonzero limit at
the cusp. A power of the quotient is holomorphic on \(X(1)\), hence
constant. Cramer's rule now gives the unique \(H_j\), of weight
\(w-v-2j\), with poles only at the cusp.
\end{proof}

For the cycle pairing a second basis is better suited. Given
\(u\in M^!_{-b_i}\), put
\(\mathcal B_i(u)=(-1)^i[u,\theta]_i\), where the bracket has the
convention of \cite[(3.6)]{BEY}. Explicitly,
\[
 \mathcal B_i(u)=\sum_{\ell=0}^i(-1)^\ell
   \binom{-b_i+i-1}{i-\ell}\binom{v+i-1}{\ell}
                 D^\ell u\,D^{i-\ell}\theta .
\]
Here the binomial coefficients are generalized. The Rankin--Cohen
transformation law \cite[(3.6)]{BEY} gives weight \(-b_i+v+2i=w\)
and type \(\rho_+\). It applies to weakly holomorphic inputs because
differentiation preserves finite pole order at the cusp.
Since \(D^j\theta=U_j\) plus lower columns, the highest column in
\(\mathcal B_i(u)\) is \(U_i\), with coefficient
\(\binom{-b_i+i-1}{i}u=(-1)^i\binom{b_i}{i}u\).
This scalar multiplier is nonzero for \(b_i\ge i\).
Applying Lemma~\ref{med-u-wronskian} and eliminating the highest column
successively writes every input uniquely as a sum of brackets,
with no bound on its pole order.

\subsection{Pairing with a bracket}
Fix \(i\) and abbreviate \(b=b_i,\ s=s_i\). In the bracket the
coefficient \(c_u(-n)\) and the theta summand indexed by \(T\)
meet at exponent \(-n+T^2/(4N)\). The sign from \(D^\ell u\)
cancels the sign in the bracket. Hence before imposing the strict
trace bound the multiplier of this term is
\[
 T^\epsilon\sum_{\ell=0}^i
  \binom{-b+i-1}{i-\ell}\binom{\epsilon+i-1/2}{\ell}
                 n^\ell (T^2/(4N))^{i-\ell}.
\]
It equals \(c_{a,i,N}P_{a,s}(T,Nn)\), where
\[
       c_{a,i,N}=\frac{\binom{-b+i-1}{i}}{(4N)^i\binom ba}\ne0 .
\]
Here is a coefficient check for the identity in all degrees.
The coefficient of \(T^{s-2\ell}D^\ell\) in \(P_{a,s}(T,D)\),
divided by its coefficient \(\binom{b}{a}\) at \(T^s\), is
\[
       (-1)^\ell\frac{(b-\ell)!\,s!}{b!\,\ell!\,(s-2\ell)!}.
\]
In the bracket polynomial the same ratio follows from
\(\binom{-b+i-1}{i-\ell}=(-1)^{i-\ell}\binom{b-\ell}{i-\ell}\)
and
\(\binom{\epsilon+i-1/2}{\ell}
 =s!(i-\ell)!/(4^\ell\ell!(s-2\ell)!i!)\).
Summing the terms of negative exponent and retaining their full
component multiplicities therefore gives the identity in \(A_N\)
\begin{equation}\label{med-u-constant-term}
             Z_\kappa(\mathcal B_i(u))=c_{a,i,N}\operatorname{CT}(uF_i).
\end{equation}
It is a Fourier coefficient identity and requires no convergence
assumption on \(F_i\).

\begin{proof}[Proof of Theorem~\ref{med-u-scalar}]
By the principal-part criterion and the bracket basis,
Chow-valued cuspidality is equivalent to
\(\operatorname{CT}(uF_i)=0\) for every \(i\) and \(u\in M^!_{-b_i}\).
A cusp form \(F_i\) satisfies this condition by the residue theorem
for \(uF_i\,d\tau\) on \(X(1)\). Its sole possible pole is the cusp,
with residue \(\operatorname{CT}(uF_i)/(2\pi i)\).
Conversely, the condition says that \(q^{-n}\mapsto F_i(n)\)
kills the principal parts of all such \(u\). Scalar principal-part
duality in weights \(-b_i,b_i+2\) thus gives a map
\(S_{b_i+2}^\vee\to A_N\), whose tensor in \(S_{b_i+2}\otimes A_N\)
has precisely the prescribed coefficients. For any other compatible
root of the trivial twist, the isometry \(\mu\mapsto r\mu\)
permutes the components, since \(r^2\equiv1\bmod 4N\).
\end{proof}

\subsection{Explicit coefficient identities}
The scalar criterion can be tested entirely in the Chow group.
For even \(k=b_i+2\geq4\) put
\[
  d_k=\left\lfloor k/12\right\rfloor-\mathbf1_{\,k\equiv2\ (12)} .
\]
Choose the normalized echelon basis \(g_1,\ldots,g_{d_k}\) of
\(S_k(\mathrm{SL}_2(\mathbb Z))\), characterized by
\([q^t]g_j=\delta_{jt}\) for \(1\leq t\leq d_k\).
For \(1\le j\le d_k\), write \(k-12j=4u+6v\) with \(u,v\ge0\).
Then \(\Delta_{\rm mod}^jE_4^uE_6^v\), with
\(\Delta_{\rm mod}=\eta^{24}\), has initial exponent \(j\).
These \(d_k\) forms are a basis by the dimension formula, and triangular
elimination gives the stated normalization. Thus \(\Phi_\kappa\) is cuspidal
if and only if, for every \(i\) and \(n>d_k\),
\begin{equation}\label{med-u-residual}
               F_i(n)=\sum_{j=1}^{d_k}([q^n]g_j)\,F_i(j),
                 \qquad k=b_i+2 .
\end{equation}
The right side gives the coefficients of the unique candidate matching
the first \(d_k\) entries. For \(d_k=0\) the condition reads
\(F_i(n)=0\) for all \(n\ge1\).

\subsection{Restricting a graph to a cyclic isogeny}
Let \(U_b\) be the \(b\)-fold power of the universal elliptic curve
over \(\mathcal M_{1,1}\), and let \(e_b\) take the part odd in each
coordinate and alternating under permutations. On \(U_b\times U_b\)
write \(H_b(n)\) for the \((e_b\times e_b)\)-projection of the
isogeny graph: for every subgroup \(K\subset E\) of order \(n\), take
the graph of \(E^b\to(E/K)^b\), and push forward from the stack of
such subgroups. Both projections of its support are finite over the
corresponding elliptic powers. All Chow groups below have rational
coefficients, extended to \(\mathbb C\) for the normalized CM classes.
On a fine cover these signed graphs have unique compact Scholl
lifts by Lemma~\ref{scholl-open-chow} applied to the product.

\begin{theorem}[Restriction of the full-level graphs]\label{graph-to-moments}
Fix \(s=s_i\), \(b=a+s\). There is a linear operation on ordinary Chow
classes on \(U_b\times U_b\), independent of \(n\), whose value on
\(H_b(n)\) is a fixed nonzero multiple of the restriction of
\(F_i(n)\) to the smooth-base locus. Consequently, for a given \(b\),
the cusp-form coefficient identities of weight \(b+2\) for \(H_b(n)\)
imply \eqref{med-u-residual} for every cyclic level and every row
with \(a+s_i=b\). The same implication holds modulo classes supported
over finitely many pairs of base points.
\end{theorem}
\begin{proof}
\smallskip\noindent\emph{Restrict the graph to a cyclic isogeny.}
On the base of pairs \((E,L)\) with \(L\) cyclic of order \(N\),
write \(q_L:E\to E/L\) for the quotient. Pull the first elliptic
power back by \((E,L,z)\mapsto(E/L,q_L^b z)\). This map is finite
flat. On the base it is the forgetful map composed with the Fricke
involution, and on the fibres it is an isogeny. Refined pullback
along the graph of the forgetful map to \(\mathcal M_{1,1}\) now
identifies the target base elliptic curve with \(E\). On the resulting
relative \(E^b\times E^b\), intersect each of the last \(s\) pairs
with the mixed Poincar\'e divisor. Forget the \(2s\) coordinates by
proper pushforward and apply \(e_{2a}\) to the remaining ones.
All these operations are independent of \(n\) and can be computed
on fine charts. They kill a class over a base point pair because,
where it meets the base graph, its excess normal line comes from a point.

\smallskip\noindent\emph{Identify the CM incidences.}
An incidence is an endomorphism \(\gamma:E\to E\) of degree \(Nn\)
killing \(L\), since every such endomorphism factors uniquely
through \(q_L\), even for a noncyclic kernel. Write
\(E=\mathbb C/(\mathbb Z+\tau\mathbb Z)\), \(L=\langle1/N\rangle\).
An integral matrix \(M\) fixing \(\tau\) acts by multiplication by
\(\lambda=M_{21}\tau+M_{22}\). It kills \(L\) precisely when its
bottom row is divisible by \(N\). For trace \(T\) we obtain the
mutually inverse assignments
\[
 M=\tfrac T2 I+Nx,\qquad
 x=\begin{pmatrix}B/(2N)&-A/N\\ C&-B/(2N)\end{pmatrix}.
\]
They give \(B\equiv T\bmod 2N\) and \(Q(x)=n-T^2/(4N)\), and
commute with change of cyclic marking. If \(M\) is nonscalar, its
fractional linear transformation has derivative
\(\bar\lambda/\lambda\ne1\) at \(\tau\). The fixed point thus has
multiplicity one on a fine chart. The fibre product records the
isomorphism from the quotient curve to its target. Thus descent
weights each oriented orbit by the reciprocal of its stabilizer
order, also at \(j=0,1728\), with the normalization of
\eqref{rd-coeff}. A scalar branch
instead meets along the whole base, with excess line pulled back
from \(T_{\mathcal M_{1,1}}\). Its first Chern class vanishes
rationally on the smooth modular stack because the modular
discriminant trivializes the twelfth power of the Hodge line.
Only \(T^2<4Nn\) contributes, as in \eqref{med-u-moment}.

\smallskip\noindent\emph{Compute the fibre contribution.}
Let \(\lambda_1,\lambda_2\) be the eigenvalues of a nonscalar
\(\gamma\) in a symplectic eigenbasis of \(H^1(E)\). Put
\[
 Q_\gamma(X,Y)=\lambda_2 X_1Y_2-\lambda_1 X_2Y_1,\qquad
 \mathcal D=\partial_{X_1}\partial_{Y_2}
             -\partial_{X_2}\partial_{Y_1}.
\]
In the two symmetric powers the signed graph is \(Q_\gamma^b\).
Intersecting a mixed divisor and integrating over its coordinate pair
contracts with the alternating form. The \(s\) contractions therefore
act by \(\mathcal D^s\) up to a nonzero combinatorial factor, and
the final symmetrization sets \(Y=X\). With
\(J(v)=(1-\lambda_1v)(1-\lambda_2v)\) we have
\[
 \left.(e^{v\mathcal D} e^{Q_\gamma(X,Y)})\right|_{Y=X}
       =J(v)^{-1}\exp\bigl(Q_\gamma(X,X)/J(v)\bigr).
\]
Indeed, the two differentiations act on disjoint pairs. On either
pair, expansion gives \(e^{t\partial_x\partial_y}e^{cxy}
=(1-ct)^{-1}\exp(cxy/(1-ct))\). Taking degree \(2a\) in \(X\)
and degree \(s\) in \(v\) yields
\[
 \left.\mathcal D^s Q_\gamma^b\right|_{Y=X}
  =\frac{b!\,s!}{a!} P_{a,s}(T,Nn)\,Q_\gamma(X,X)^a .
\]
On the diagonal only the mixed tensor of
\(Nx=\gamma-(T/2)\operatorname{id}\) remains. For odd \(T\) this
is a quasi-endomorphism, to which the mixed Poincar\'e divisor extends
by additivity and \(\mathbb Q\)-linearity. On the abelian fibre, the
graph divisors and the mixed Poincar\'e divisors are symmetrically
distinguished. Their signed products and the homomorphism pushforwards
in the contraction remain so. O'Sullivan's theorem therefore lifts
the cohomological identity to rational Chow \cite{OSullivanCycles}.

Since \(D_{2Nx}=2D_{Nx}\), the normalization
\eqref{eq:common-chow-scalar} at \(\Delta=1\) identifies the
pushforward of the remaining graph tensor with a fixed nonzero
multiple of \(Q(x)^{a/2}S_\kappa(x)\). Summing the
incidences with the factor \(P_{a,s}(T,Nn)\) gives the asserted
multiple of \(F_i(n)\) on the smooth base. Lemma~\ref{scholl-open-chow}
carries the resulting identities to compact Chow.
\end{proof}

\begin{corollary}[A uniform range of Chow relations]\label{bounded-graph-rows}
At every cyclic level, \(F_i(n)=0\) in ordinary compact Chow for
every \(n\ge1\) and every row with \(b_i\le8\). In particular,
for the trivial twist and every compatible root, the entire series
\(\Phi_\kappa\) vanishes in ordinary compact Chow whenever \(h=0\) or
\[
                  a+\epsilon+2(h-1)\le8.
\]
At any other level the weakly holomorphic inputs generated by these
rows still give rationally trivial cycles for arbitrary pole order.
\end{corollary}
\begin{proof}
Remove the zero coordinates and the pair diagonals in the second
factor of \(U_b\times U_b\). This leaves
\(U_b\times\mathcal M_{1,b+1}\). For \(b+1\le10\) the rational
Chow group of the pointed stack is generated by its fundamental
class, and its exterior products generate the Chow group of a
product with any such stack \cite[Theorem~5.1 and \S5.2]{CanningLarson}.
Every restricted class therefore lifts to a class
\(\alpha\boxtimes[U_b]\), which is killed by inversion in one
coordinate of the second factor after odd projection. By localization,
the difference from the original class is supported on the omitted
loci. Inversion at a zero coordinate or transposition on a pair
diagonal fixes that locus pointwise and has negative character.
Hence \(1\times e_b\) kills the difference as well. In particular
\(H_b(n)=0\) for every positive even \(b\le8\), for all \(n\).
Theorem~\ref{graph-to-moments} now gives the first assertion.
Equation~\eqref{med-u-constant-term} treats all polar parts in these
rows. To recover the coefficients themselves, put
\(C_j=\langle D^j\theta,\Phi_\kappa\rangle\), pairing over all residues.
Expand \(P_{a,s_i}(T,Nn)\) in the monomials \(T^{s_i-2l}(Nn)^l\).
In the convolution defining \(F_i\), \(D^l\) on the whole series
gives the factor \(n^l\), while \(D^{i-l}\theta\) provides
\(T^{s_i-2l}/(4N)^{i-l}\). Thus \(F_i\) is a linear combination of
\(D^l C_{i-l}\) for \(0\le l\le i\), with coefficient
\((4N)^i\binom{b_i}{a}\ne0\) at \(C_i\). If all rows vanish,
induction gives \(C_0=\cdots=C_{h-1}=0\). The columns \(D^j\theta\)
have the same determinant as the Serre derivative columns in
Lemma~\ref{med-u-wronskian}, by unitriangular change of basis.
Pairing on the reflection subspace is nondegenerate, including the
self-opposite residues in even parity. Inverting over the field of
formal Laurent series in \(q^{1/(4N)}\) gives \(\Phi_\kappa=0\)
coefficientwise (equivalently, apply each linear functional first).
For \(h=0\) reflection already gives this conclusion. The isometry
\(\mu\mapsto r\mu\) transfers the conclusion to every compatible root.
\end{proof}

\begin{proposition}[The diagonal and the cycles of weakly holomorphic inputs]
\label{diagonal-hecke-defect}
Fix \(a=\kappa-1\ge1\) and a cyclic level \(N\), with the trivial
twist and root one. For positive even \(b\), let \(\mathbb T_b\)
be the abstract full-level Hecke algebra in weight \(b+2\), and put
\[
 I_b^{\rm cusp}=\ker\!\left(\mathbb T_b\longrightarrow
        \operatorname{End}_{\mathbb C} S_{b+2}(\mathrm{SL}_2(\mathbb Z))\right).
\]
The geometric Hecke operators on the target act on the span of the
graphs modulo \(J_b\), with \(T_nH_b(1)=H_b(n)\). Write \(R_i\)
for the operation of Theorem~\ref{graph-to-moments}, divided by its
nonzero normalization constant so that \(R_i H_{b_i}(n)=F_i(n)|_{\rm sm}\).
Let \(D_{1,1}(N)\) be the span of the cycles of all untwisted
weakly holomorphic inputs. On the smooth base one has the equality
of subspaces
\[
            D_{1,1}(N)|_{\rm sm}
       =\sum_{0\le i<h} R_i\big(I_{b_i}^{\rm cusp} H_{b_i}(1)\big).
\]
In particular, untwisted modularity in compact Chow is equivalent
to the vanishing of each summand.
\end{proposition}
\begin{proof}
Take the basis \(T_n\) of the abstract algebra with
\(T_pT_n=T_{pn}+p^{b+1}T_{n/p}\) and \(T_1=1\). The target
Hecke operators commute: flags of kernels of coprime orders
are identified by reordering the isogenies. The kernel count in
Lemma~\ref{full-graph-hecke} gives the displayed recursion on graphs.
Applying it to \(H_b(1)\) therefore gives \(T_nH_b(1)=H_b(n)\)
by induction on \(n\).

For \(T=\sum_n u_nT_n\) and a cusp form \(g\), the first
coefficient of \(Tg\) is \(\sum_nu_n[q^n]g\). If this vanishes
for every \(g\), then for every \(j\)
\[
         [q^j]Tg=[q]T_jTg=[q]T(T_jg)=0 .
\]
It follows that \(T\in I_b^{\rm cusp}\), and the converse is
immediate. By scalar principal-part duality these \(T\) are
exactly those whose \(u_n\) form the negative Fourier coefficients
of some \(u\in M^!_{-b}\). For \(b=b_i\),
\eqref{med-u-constant-term} identifies
\(R_i(TH_b(1))\) with \(c_{a,i,N}^{-1}Z_\kappa(\mathcal B_i(u))|_{\rm sm}\).
All inputs are sums of such brackets and every bracket is an input.
This proves the equality. Lemma~\ref{scholl-open-chow} and
Theorem~\ref{rd-main} give the last assertion.
\end{proof}

\section{Evaluation on the Kummer surface}\label{rd-g-section}

Along an isogeny curve in the product of two modular curves, the
isogeny graph gives an additional divisor on the product of the
elliptic curves. We use the elliptic pencil on its Kummer surface to
study relations among exterior powers of this divisor. On a fixed
surface we identify the signed evaluation of stable maps with a
multiple of the projected graph. In the family we determine the
support and multiplicity of that contribution.

The fixed-surface argument in Proposition~\ref{fixed-evaluation}
requires a Chow decomposition with all marked points and node flags
present. Lemma~\ref{vertex-derived-cut} gives the disjoint-position
comparison used for that decomposition. When at most one smooth
fibre can occur in the vertical image, the decomposition also follows from the
support of evaluation.
Proposition~\ref{family-supported-diagonal} compares the refined
pullback in \eqref{rechecked-family-diagonal} with the reduced
classes of the individual vertices before their flag evaluations
are matched. Section~\ref{med-u-five} uses these comparisons to compute the
unmarked graph relations. Section~\ref{rd-m-twist} treats the marked
pencils and transfers the relations to the twisted cycle series.

\subsection{Distinct positions in the pencil}
Let \(\Bbbk\) be an algebraically closed field of characteristic zero.
Let \(E_1,E_2\) be the elliptic curves over the geometric generic point of
an isogeny correspondence. In what follows one can take any isogenous pair
with \(\operatorname{End}(E_1)_{\mathbb Q}=\mathbb Q\). Choose a generator
\(\eta\) of \(\operatorname{Hom}(E_1,E_2)\); its sign amounts to a choice
of orientation. Put \(A=E_1\times E_2\), \(T=A^b\), and write its
points as \(((z_1,w_1),\ldots,(z_b,w_b))\). All Chow groups in this
section have rational coefficients. For \(i=1,2\) let \(G_i\) be the
group of signed permutations on the \(b\) coordinates belonging to
\(E_i\). Define the character
\[
 \chi_i(\epsilon_1,\ldots,\epsilon_b;\sigma)
       =\epsilon_1\cdots\epsilon_b\,\operatorname{sgn}(\sigma).
\]
Write \(e_i\) for the corresponding character idempotent and
\(e=e_1e_2\) for their product acting on \(\CH^*(T)\). Thus \(e_i\)
first takes the odd part in each coordinate, then alternates the
permutations. For the quotient map \(\lambda:E_1\to P=E_1/\{\pm1\}\)
put
\[
 U=\big\{\lambda(z_i)\ne\lambda(z_j)\ (i\ne j),\
                          z_i\notin E_1[2]\ \text{for all }i\big\}\subset T .
\]
The open set is stable under \(G_1\times G_2\), so restriction
commutes with both character projectors.

\begin{lemma}[Restriction to distinct pencil positions]\label{med-g-open}
Restriction induces an isomorphism
\[
                       e\CH^q(T)\ \xrightarrow{\ \sim\ }\ e\CH^q(U)
                       \qquad(q\geq 0).
\]
\end{lemma}
\begin{proof}
Localization gives surjectivity on the full Chow groups, and
restriction commutes with \(e\). To prove injectivity, write \(D=T\setminus U\).
The kernel before projection is generated by pushforwards of integral
cycles in \(D\). Each lies in one of the loci
\[
 z_i=z_j,\qquad z_i=-z_j,\qquad z_i=t\quad (t\in E_1[2]).
\]
Transposing \(z_i,z_j\) fixes the first locus pointwise. The same
transposition followed by inversion of both coordinates fixes the
second. Inversion of \(z_i\) fixes the third. Each of these operations
has character \(-1\) in \(G_1\). If \(C\) is an integral cycle on
the corresponding locus and \(\gamma\) the indicated operation, then
\(e_1[C]=e_1\gamma_*[C]=-e_1[C]\). Thus every generator of the kernel
vanishes after projection with rational coefficients, including cycles
in the intersection of several components of \(D\).
\end{proof}

Let \(q_A=\tfrac14(1-[-1,1]^*)(1-[1,-1]^*)\) on \(A\), where the
signs denote the two coordinate inversions, and put
\(D_\eta=q_A[\Gamma_\eta]\) for the mixed part of the isogeny graph.
We next determine the signed part of the exterior products that
will occur in the evaluation.

\begin{lemma}[Products of single-slot classes]\label{med-g-products}
For the isogenous pair fixed above, with
\(\operatorname{End}(E_1)_{\mathbb Q}=\mathbb Q\), the images under
\(e\) of exterior products in \(\CH^b(T)\) span
the single line \(\mathbb Q\, e(D_\eta^{\boxtimes b})\). The same
statement holds for their restrictions to \(U\).
\end{lemma}
\begin{proof}
Take homogeneous classes in each slot, with codimensions adding to \(b\).
The projector on the two inversions in one slot kills \(\CH^0(A)\),
since the fundamental class is fixed. If one slot has codimension
larger than one, some other slot has codimension zero. Hence only
products of divisors survive.

Given a line bundle on \(E_1\times E_2\), tensor it with pullbacks
from the factors to trivialize its restrictions to the two zero axes.
The resulting normalized bundle defines a homomorphism
\(E_1\to\operatorname{Pic}^0(E_2)\) by restriction to
\(\{z\}\times E_2\) as \(z\) varies. Pullback of the normalized
Poincar\'e bundle gives the inverse construction by the see-saw
principle. The projector \(q_A\) kills the axis terms, and inversion in either coordinate negates the
homomorphism. The principal polarization of \(E_2\) therefore gives
\[
 q_A\CH^1(A)=\operatorname{Hom}(E_1,E_2)_{\mathbb Q}
             =\mathbb Q D_\eta .
\]
For the last equality note that, after composing with a rational inverse
to the isogeny, every homomorphism becomes a scalar endomorphism of
\(E_1\). Applying the two alternating projectors proves containment
in the stated line. Its generator is nonzero: in cohomology the graph
gives an isomorphism on \(H^1\), and its projected power induces the
isomorphism on the \(b\)-th symmetric powers of \(H^1\). Here the
sign from permuting degree-one cohomology factors cancels the sign
in the alternating projector. One may check this after passage to
complex cohomology by descent to a finitely generated field, or in
characteristic-zero \(\ell\)-adic cohomology. Lemma~\ref{med-g-open}
gives the assertion on \(U\).
\end{proof}

The two lemmas reduce the test for membership in the graph line to
the open locus \(U\). It suffices to express the signed restriction
of the middle-dimensional cycle there as a sum of restricted exterior
products in Chow. All collisions between two retained first
coordinates disappear on this open. Other collisions arise when the
stable map carries additional marks.

\subsection{The squarefree double ramification relation}
Let \(S\) be the minimal resolution of \(A/\{\pm1\}\), where the sign
acts simultaneously on both coordinates. Projection from \(A\) to
\(E_1\) induces an elliptic pencil \(p:S\to P\simeq\mathbb P^1\)
with fibre class \(F\). The four special pencil positions are the
images of \(E_1[2]\). To transfer a cycle on \(S\) to \(A\), use the
correspondence given by the blowup of \(A\) at its two-torsion points,
with its finite map to \(S\) and its blowdown to \(A\). Use this
correspondence separately at each retained mark. Equality of pencil
positions then means \(\lambda(z_i)=\lambda(z_j)\) on its support.

For the family construction the two curves at first vary independently
on fine modular curves with labelled two-torsion. Write \(b=2d\ge2\)
and use the genus \(d-1\) double ramification relation in degree \(d\).
Give its \(2b\) marked points weights
\((x_1,\ldots,x_b,-x_1,\ldots,-x_b)\). Extract the coefficient
of \(x_1\cdots x_b\). This gives a Chow relation on stable curves:
the relation vanishes for \(d> d-1\), and its coefficients are
polynomial in the weights
\cite[Theorem~1.2]{CladerJandaDR}\cite[Theorem~4.1]{SpelierDR}. Pull back along
stabilization on the stack of maps to the Kummer surfaces in the
chosen total class, and cap with the labelled ordinary family class of
Proposition~\ref{family-supported-diagonal}.
Retaining the evaluations at the positive marks and inserting \(F\)
at the negative marks gives a projected relation. We describe its
boundary terms in terms of isogeny graphs.

Each vertex represents a connected map. Each edge has two flags,
one on either branch of a node, whose evaluations in \(S\) must agree.
Take a total class of degree two over the pencil with zero mixed
part on the abelian surface. In the contributing terms over the
generic point of an isogeny curve, two nonconstant vertices have
opposite nonzero mixed parts. Both have positive fibre degree because vertical classes
extend, so each has degree one. The family comparison eliminates
the extending contributions and restricts the two vertex classes
to the reduced virtual classes on the fixed surface. We retain
both flags at every original node for the diagonal pullback. The
resulting class on \(A^b\) has codimension \(b\), with one further
codimension in the family from the isogeny curve.

For a fixed vertex let \(B\) be its effective class on \(S\), so
\(B\cdot F=1\). Write \(I\) for the set of all marks and flags on the
vertex, and \(K=\overline M_{g,I}(S,B)\) for its space of connected
stable maps. The evaluation class, allowing nonnegative cotangent
exponents \(a_i\), is
\[
             V_{B,I}=\operatorname{ev}_*
                  \left(\prod_{i\in I}\psi_i^{a_i}\cap[K]^{\rm red}\right)
                  \in\CH_*(S^I).
\]
Here \(\psi_i\) is the first Chern class of the marked cotangent line.
Divisor insertions at legs can be made on \(S^I\) by the projection
formula. Proposition~\ref{fixed-evaluation} also allows powers of
the hyperplane class of the image linear system, inserted before
pushforward. Its exponent is zero in the graph relations of
Sections~\ref{med-u-five} and~\ref{rd-m-twist}. For \(\varnothing\ne J\subset I\), set
\[
 S_P^J=\{(v_j)_{j\in J}\in S^J:
                 p(v_j)=p(v_{j'})\ \text{for every }j,j'\in J\}.
\]
Let \(\iota_J:S_P^J\hookrightarrow S^J\) be the inclusion.

\begin{lemma}[The term with a single nonconstant vertex]\label{family-single-positive}
On the framed product of the two smooth modular curves, every term of
the squarefree relation with a single nonconstant vertex has zero
signed Chow evaluation. This holds for all \(b=2d\).
\end{lemma}
\begin{proof}
If both marks weighted \(x_i,-x_i\) are on the same vertex of the
original graph, their sum drops out of every edge weighting and
their leg factors depend only on \(x_i^2\). That graph contributes
zero to the squarefree coefficient. Thus all \(b\) pairs have their
ends on different original vertices. Passing to prestable curves
pulls back vertex decorations separately; correcting a stabilized
cotangent class introduces a rational-tail edge in place of one
cotangent factor \cite[Propositions~3.13--3.14]{BaeSchmitt}.
Such refinements cannot reunite separated marks, and the total
codimension \(d\) bounds the number of specified edges by \(d\).

On a graph with just one nonconstant vertex, at least one mark from
each pair is therefore on a constant vertex. Two retained marks at
that vertex have the same first pencil position, killed after
Kummer transfer by signed transposition. A retained and an inserted
mark there give a torsion first coordinate: represent \(F\) by a
branch fibre. Inversion kills this term. For two inserted marks the
product is \(F^2=0\) on the family; the framed quotient pencil is
a projective line bundle with three disjoint named sections, hence
trivialized by them, and disjoint constant fibres represent \(F\).
Thus each of the \(b\) pairs requires a different constant vertex.
Root a spanning tree at the nonconstant vertex. Each constant
vertex has a different edge toward the root. Hence the
graph would have at least \(b=2d\) edges, whereas it has at most
\(d\) and \(d\ge1\).
\end{proof}

\paragraph{The fixed-surface equality.}
For the vertex classes and insertions appearing above, the needed
factorization is the inclusion
\begin{equation}\label{med-g-factor}
 V_{B,I}\ \in\
 \sum_{\pi\ {\rm partition\ of}\ I}
       \operatorname{span}_{\mathbb Q}
       \left\{\ \mathop{\boxtimes}_{J\in\pi}(\iota_J)_*\alpha_J:
                               \alpha_J\in\CH_*(S_P^J)\ \right\}.
\end{equation}
Within each block \(J\), the evaluations may vary together on the
fibre product \(S_P^J\). Across blocks the class must be a sum of
exterior products. The next proposition shows why this distinction
matters for the isogeny graph. Proposition~\ref{fixed-evaluation}
then establishes the decomposition in Chow.

\begin{proposition}[From vertex evaluation to an isogeny graph]\label{rd-g-graph}
Form a codimension-\(b\) evaluation on \(A^b\) from the two
degree-one vertices on the fixed surface, the constant vertices, the
node matchings and the slotwise insertions described above. If the
nonconstant vertices satisfy \eqref{med-g-factor} with all flags
included, its image under \(e\) is a rational multiple of
\(e(D_\eta^{\boxtimes b})\).
\end{proposition}
\begin{proof}
Choose a summand of \eqref{med-g-factor} at each nonconstant vertex.
Represent each block by a vertex of a new graph. Add one vertex for
each constant component, whose flags all evaluate at the same point.
For every original node join its two flags, allowing loops.

Refined pullback along the surface diagonals commutes with proper
pushforward from the block supports \(S_P^J\). The resulting class,
including any excess Chern classes, lies on their fibre product.
The two evaluations at a node have the same pencil position.
Following a path gives the same position at any two marks in one
connected component. Different components have disjoint slots, so matching the
nodes, inserting the leg classes and pushing forward the evaluations
all preserve the exterior product across components. Kummer transfer
also preserves it.

On \(U\), a component with two retained outputs contributes zero
because their pencil positions coincide. A component with no output
pushes forward to a scalar, possibly zero. Each remaining component
contributes a class on one copy of \(A\). Thus on \(U\) we obtain a sum
of products of single-slot classes. Lemmas
\ref{med-g-open} and \ref{med-g-products} identify the signed class
on all of \(A^b\). When the graph and the input cycle are defined
over the function field of the oriented sheet, the rational
equivalence over its algebraic closure descends to a finite extension.
Taking its trace and dividing by the degree proves the assertion
over the function field.
\end{proof}

\subsection{Cutting a degree-one map along its section}
The line bundle of \(B\) restricts to degree
one on the generic elliptic fibre. It has a unique effective
representative there. Its closure is a section \(M\) of the pencil
by properness over the smooth base curve. Every divisor in \(|B|\)
therefore contains this section with multiplicity one. Two vertical
divisors in the same linear equivalence class differ by the divisor
of a rational function from \(P\): a function with vertical divisor
restricts to an invertible regular function on the complete generic
fibre and hence belongs to \(\Bbbk(P)^*\). At each reducible fibre remove
from the vertical part as many full scheme fibres as its multiplicities
permit. The residual effective divisor \(V_t\) at each such \(t\) is
independent of the representative, and for some \(h\ge0\) one obtains
\[
                      |B|=M+\sum_t V_t+p^*|\mathcal O_P(h)|.
\]
This equality determines the possible image divisors. Several stable
maps can have the same image divisor.
If a scheme fibre is \(\sum_j m_j C_j\) and the vertical divisor
has coefficients \(a_j\) there, the number of full fibres removed
is \(\min_j\lfloor a_j/m_j\rfloor\). Its entire nonzero remainder
\(V_t\) is one fixed label at \(t\), including all its components.
If \(k=\min_j\lfloor a_j/m_j\rfloor>0\) and \(V_t\ne0\), the
moving divisor has \(k\) labels at the same point as the fixed label.
All these labels stay in one group on a disjoint-position chart.
Their coincidence is retained in the diagonal term below.

There are finitely many pencil positions where a map can differ from
its section. The image-divisor morphism records every vertical
component of positive degree with its multiplicity
\cite[\S4]{KoolThomas}. The branch morphism of the stabilized
degree-one map to \(P\) has degree \(2g\). A contracted tail of
arithmetic genus \(r\) attached at one point contributes \(2r\)
there by the component and node formula \cite{FPBranch}. We also
record every marked point and each fixed residual position \(t\).
Away from these positions, a contracted tail could only be an
unmarked rational tree with a single attachment. A leaf farthest
from that attachment would violate stability. The map therefore
agrees with the section there. On a geometric fibre we can cut the map
into pieces with disjoint sets of recorded positions and recover it by
gluing along the common section.

Away from the recorded positions, cutting duplicates the horizontal
section. In the ordinary deformation theory, normalization matches
the two tangent values at every node. In the reduced theory each
nonconstant piece also has a semiregularity quotient
\(H^2(S,\mathcal O_S)\). To compare the two reduced pieces, we
keep the determinant of the image complex during gluing.

\begin{lemma}[Cutting the fixed-surface map at disjoint positions]
\label{vertex-derived-cut}
Suppose the moving image divisor, branch divisor, markings and fixed
vertical remainders are split into two labelled groups with disjoint
supports. For the resulting distribution of classes, genera and
markings, the reduced derived map space on this chart is the product
of the two reduced map spaces restricted to disjoint positions.
The equivalence respects cotangent lines and refinements of the split.
\end{lemma}
\begin{proof}
\smallskip\noindent\emph{The determinant on the overlap.}
Write \(M_0\) for the horizontal section, \(V\) for the sum of the
fixed vertical remainders, and \(\mathcal D(f)=\det Rf_*\mathcal O_C\).
Over a derived affine \(T\), the reduced determinant construction
\cite[Proposition~4.6 and Definition~4.7]{STV} specifies a stable map
\(f:C\to S_T\), a degree-\(h\) line \(L\) on \(P_T\), and an equivalence
\[
 \mathcal D(f)\simeq\mathcal O_{S_T}(M_0+V)\otimes p_T^*L .
\]
The reduced stack is the fibre over zero of the derived coordinate
in the splitting of \(\mathbf R\operatorname{Pic}(S)\)
of \cite[Proposition~4.6]{STV}. Its ordinary Picard component
retains the scalar automorphisms of a line bundle. Since \(P\simeq
\mathbb P^1\), pullback from \(\mathbf R\operatorname{Pic}^{h}(P)
\simeq B\mathbb G_m\) followed by tensoring with the fixed line
identifies this component with the stack of lines used above.
The pointed determinant fibre is quasi-smooth for every
nonzero class \cite[Theorem~4.11]{STV}. The ordinary mapping complex
relative to the pointed prestable stack is \(R\pi_{C*}f^*T_S\), with
the same semiregularity quotient in the reduced theory
\cite[Theorem~2.4 and Lemma~2.5]{KoolThomas}.

On the classical truncation let \(D_1,D_2\) be the two groups of
records. Addition of effective divisors of fixed degrees on \(P\) is
\'{e}tale where their supports are disjoint: the idempotent separating
the two divisors lifts uniquely across nilpotent thickenings. This
allows arbitrary collisions within either group. Pull back this
classical chart and take its unique derived \'{e}tale lift. Only the
chart is needed to define \(U_1=P_T\setminus |D_2|\) and
\(U_2=P_T\setminus |D_1|\).

Over \(U_{12}=U_1\cap U_2\) the map \(h=pf\) is an isomorphism on
geometric fibres. On truncations it is finite by properness and
quasi-finiteness, and flat by the fibrewise criterion for schemes flat
over \(t_0T\) (\href{https://stacks.math.columbia.edu/tag/039E}{Stacks Project,
Tag 039E}). Thus \(t_0h\) is finite \'etale of degree one. Derived
flatness over \(T\) identifies the higher homotopy sheaves on source
and target with the corresponding pullbacks from \(T\). Hence \(h\)
itself is an isomorphism there and \(f\) is a section.

On any such open \(U\), we claim that the space of sections with their
determinant data is contractible. The equation of \(M_0\) generates
\(Rp_*\mathcal O_S(M_0)=\mathcal O_P\). To see this, K3 Riemann--Roch
and \(h^0(S,F)=2\) give \(H^1(S,\mathcal O_S(-F))=0\).
The fibre sequence then gives \(H^0(F_t,\mathcal O_{F_t})=\Bbbk\),
even for a nonreduced fibre. Its only global functions are constants,
and \(M_0\) meets it in a point. Since \(\omega_{F_t}\simeq\mathcal O_{F_t}\),
duality gives \(H^1(F_t,\mathcal O_{F_t}(M_0))=0\). The space of
sections of this degree-one line has dimension one, so cohomology
and base change give the assertion. The pencil \(p\) is flat.
Derived base change therefore gives the same identity for
\(Rp_*\mathcal O(M_0)\) after pullback to \(P_T\) and restriction to \(U\).

The graph of a section \(q\) in the smooth locus is a regular immersion
of codimension one. For its ideal \(I_q\), the triangle
\(I_q\to\mathcal O\to\mathcal O_{\Gamma_q}\) identifies the graph
determinant with \(I_q^\vee\) and gives the defining section. The derived
zero locus of this section recovers the virtual Cartier divisor
\cite[Proposition~3.2.6]{KhanRydhCartier}. By the projection formula and
derived base change, the space of sections of \(\mathcal O(M_0)\otimes p^*L\)
over \(S_U\) is the space of sections of \(L\) over \(U\).
A section \(t\) vanishing at a geometric point of \(U\) would put the
entire fibre in its divisor. An invertible \(t\), conversely, recovers
the graph with its determinant equivalence as the derived zero locus
of \(s_{M_0}\,p^*t\), where \(s_{M_0}\) cuts out \(M_0\). Thus the
graph and determinant data form the stack of lines with a trivializing
section. This stack is contractible over every base, compatibly with
restriction. In particular, there is no residual motion of the
section with fixed relative determinant. Infinitesimally, the
sequence for \(M_0\subset S\) gives the boundary isomorphism
\(N_{M_0/S}\simeq R^1p_*\mathcal O_S\) on \(M_0\simeq P\),
since \(Rp_*\mathcal O_S(M_0)=\mathcal O_P\). This is the
differential of \(q\mapsto\mathcal O(q-M_0)\) in the relative Picard.

\smallskip\noindent\emph{Cutting the pointed map.}
Retain \(f\) on \(U_1\) and glue it to \(M_0\) on \(U_2\) using the
trivialization on \(U_{12}\). The omitted fixed vertical divisor has its
canonical trivialization on \(U_1\). The same overlap identification
glues \(L|_{U_1}\) to the trivial line on \(U_2\), and open base
change for \(Rf_*\mathcal O_C\) glues the determinant equivalences.
Interchanging the indices gives the second map. Properness over \(P_T\)
is local on \(P_T\). Flatness over \(T\) and the nodal condition are
local on the source. The glued domain is proper over \(T\) with
connected fibres.

To check that the glued domain belongs to the pointed prestable
stack of \cite[\S2.2]{STV}, work on a bounded prestable chart of
the derived moduli of flat proper curves
\cite[Theorem~3.32 and Example~3.36]{PridhamModuli}. On the nodal
open take the \(n\)-fold fibre product of the smooth locus of the
universal curve over the chart, and remove its pairwise diagonals.
This gives the chart for \(n\) distinct ordered smooth markings.
The ordinary prestable stack maps there by its universal curve and
agrees on truncations. At a pointed curve \((C,\Sigma)\) the derived
tangent is \(R\operatorname{Hom}(L_C,\mathcal O_C(-\Sigma))[1]\).
It has cohomology only in degrees \(-1,0\): the curve has dimension
one and its local \(\mathcal Ext^1\) is supported at the nodes.
The universal curve identifies these two groups with the
infinitesimal automorphisms and deformations on the smooth ordinary
stack. The tangent complexes are perfect, so Nakayama gives zero
relative cotangent complex. Square-zero lifting
\cite[Proposition~8.1.2]{LurieDerived} and nilcompleteness give
equivalence on derived affines. Stability of the map is fibrewise;
the domain may contain unstable components. For the corresponding
analytic comparison see \cite[Proposition~4.5 and Corollary~4.9]{PortaYu}.

\smallskip\noindent\emph{Gluing and refinement.}
Conversely, glue the two maps across the open where both equal \(M_0\).
Their auxiliary lines tensor together, and the determinant is
\[
 \mathcal D(f)\simeq
 \mathcal D(f_1)\otimes\mathcal D(f_2)\otimes\mathcal O_{S_T}(-M_0).
\]
On each \(U_i\), the other map contributes the section determinant.
The graph trivialization identifies these cancellations on \(U_{12}\).
Only this overlap uses the defining section. On each occupied open the
determinant equivalence is part of the original data.

Descent for curves, maps, lines and determinant equivalences expresses
the glued data as a homotopy fibre product over the space of bare
data on \(U_{12}\). For two lines with frames \(\alpha,\beta\), the
transition carrying one frame to the other is \(\beta\alpha^{-1}\).
These transitions compose on a triple overlap. Over a derived test
base the argument uses the contractibility of the full space of bare
data proved above. The homotopy fibre product is therefore the
product of the two spaces of occupied data, with their automorphisms.
For a finite disjoint split, the intersections of distinct
occupied opens are bare, and refining the cover gives the same descent
limit. Open restriction also preserves the relative dualizing line near
each marked section, hence its cotangent line, including at node flags.
\end{proof}

The order in which we forget positions matters. A vertical component
can attach at the first retained mark while a second mark moves freely.
Forgetting the attachment then leaves two distinct retained positions,
so restriction to the open of Lemma~\ref{med-g-open} does not remove
this collision. We apply inclusion and exclusion while the image and
branch positions still carry their labels.
For a labelled set \(A'\) and classes \(F_B\) attached to its subsets,
the relevant expression is
\[
 C_{A'}=\sum_{\pi\ {\rm partition\ of}\ A'}\!
       (-1)^{|\pi|-1}(|\pi|-1)!\ \mathop{\boxtimes}_{B\in\pi} F_B .
\]
M\"obius inversion expresses \(F_{A'}\) as a sum of products of the
\(C_B\). We show that \(C_B\) comes from the diagonal on which the
positions in \(B\) coincide.

\begin{lemma}[Compatible sheaves and collisions]\label{med-g-mobius}
Let \(Y=\prod_{i\in I}Y_i\), where the \(Y_i\) are varieties with
position maps to \(P\). For each partition \(\pi\) of \(I\), give a
coherent sheaf \(H_\pi\) on \(Y\). Write \(\pi\le\sigma\) for
refinement. Suppose \(H_\pi\) and \(H_\sigma\) are identified on
the open where any two different \(\pi\)-blocks in a common
\(\sigma\)-block have disjoint position sets. Require these
identifications to compose for chains of refinements on the common
opens. With \(1\) the one-block partition and \(\mu\) the partition
M\"obius function, the class
\[
                       \sum_{\pi}\mu(\pi,1)[H_\pi]\quad\text{in }G_0(Y)
\]
comes from the closed locus where all positions agree.
\end{lemma}
\begin{proof}
Induct on the closed union \(Z\) of the supports. Fix an irreducible
component \(Z_0\) and let \(\eta\) be the partition by equality of
positions at its generic point. Take an ambient open \(j:V\hookrightarrow Y\)
meeting \(Z_0\) densely, missing the other components, on which indices
from different \(\eta\)-blocks have distinct positions. On \(V\) the
sheaf \(H_\pi\) is identified with \(H_r\) for \(r=\pi\wedge\eta\).
In \(j_*(H_r|_V)\), take the images of all the corresponding \(H_\pi\);
denote their intersection by \(J_r\). The open immersion is
quasi-compact. The images are coherent quotients of \(H_\pi\), and
their finite intersection is coherent by noetherianity. Each \(J_r\)
is supported on \(Z_0\).

For each \(\pi\), let \(K_\pi\) be the kernel of the indicated
map from \(H_\pi\) and \(Q_\pi\) its image modulo \(J_r\). Both new
systems are supported on the proper closed subset \(Z\setminus V\).
They inherit the required identifications. To see this for a merge
\(\pi\le\sigma\), first suppose \(\pi\wedge\eta=\sigma\wedge\eta=r\).
The two maps into \(j_*(H_r|_V)\) agree on the merge open by
composition of the identifications on its intersection with \(V\)
and adjunction for \(j\). Passing to the common intersection \(J_r\)
preserves this compatibility. If the two meets differ, the merge
requires two positions equal throughout \(Z_0\) to be disjoint.
Its open misses \(Z_0\). All the images there are zero, so the
original identification identifies the kernels, and both quotients
are zero. This proves the induction hypotheses for the new systems.

If \(\eta=1\), the common sheaves \(J_r\) are already supported on the
small diagonal. Otherwise for \(s\le\eta\) one has
\[
   \sum_{\pi:\,\pi\wedge\eta\ge s}\mu(\pi,1)
                    =\sum_{\pi\ge s}\mu(\pi,1)=0
\]
since \(s<1\). Inversion on the interval of refinements of \(\eta\)
shows that the coefficient of each fixed meet \(r\) is zero.
All terms with that meet contain the same class \([J_r]\), so they cancel. The identity
\([H_\pi]=[K_\pi]+[J_r]+[Q_\pi]\) and noetherian induction now leave
only sheaves supported on the small diagonal. Filter each such
sheaf by powers of the ideal of that diagonal. The filtration is
finite by coherence, and each quotient is a sheaf on the diagonal,
as required in \(G_0(Y)\).
\end{proof}

\begin{proposition}[Fixed-surface evaluation at degree one]\label{fixed-evaluation}
For every effective \(B\) of fibre degree one on the fixed Kummer
pencil, the Chow decomposition \eqref{med-g-factor} holds in every
genus, with arbitrary powers of the marked cotangent classes and of
the image hyperplane. Consequently the evaluation of
Proposition~\ref{rd-g-graph} lies on the signed graph line.
\end{proposition}
\begin{proof}
The target is a point if \(I=\varnothing\). Assume \(I\ne\varnothing\)
and write \(|B|=M_0+V+p^*|\mathcal O_P(h)|\).
If the map covers at most one smooth fibre, its image already gives
the required decomposition. In general we keep track of the fibre
positions before forgetting any marks. Lemma~\ref{vertex-derived-cut}
gives a product over disjoint positions. Lemma~\ref{med-g-mobius}
then expresses the remaining terms by classes on the diagonals where
positions coincide. We carry this out for coherent sheaves and
recover the reduced evaluation cycle by Riemann--Roch.

\smallskip\noindent\emph{Support when at most one smooth fibre occurs.}
Each smooth pencil fibre covered by the map contributes at least one
to \(h\) and requires a separate component of geometric genus at least
one by Riemann--Hurwitz. Hence the map covers at most \(\min(g,h)\)
smooth fibres.

Suppose \(\min(g,h)\le1\). Let \(R\subset S\) be the union of \(M_0\)
and the four reduced singular fibres. Its components are rational.
Over a branch point the central component is \(E_2/\{\pm1\}\), and
resolution adds exceptional rational curves. The connected image is
the union of the section and the nonconstant vertical images. All
marks outside \(R\) therefore lie on the same smooth fibre, so
evaluation is supported on
\[
 \bigcup_{J\subseteq I} S_P^J\times R^{I\setminus J},
 \qquad S_P^\varnothing=\operatorname{Spec}\Bbbk.
\]
Every integral cycle in this union lifts with rational coefficients
to \(S_P^J\times(\mathbb P^1)^{I\setminus J}\) for some \(J\), by
normalizing the chosen components of \(R\) and using finite pushforward.
The projective bundle formula separates the line factors. Pushing down
gives the block \(J\), if nonempty, and singleton blocks elsewhere,
as required by \eqref{med-g-factor}.

\smallskip\noindent\emph{Recording and separating positions.}
For the remaining range \(\min(g,h)>1\), keep all the position labels.
For a label set \(J\), let \(h_J\) count its moving image labels and
\(2g_J\) its branch labels. A set with an odd branch count contributes
zero. The corresponding map has genus \(g_J\) and class
\[
       M_0+\sum_{t\,\text{labelled in }J}V_t+h_JF .
\]
The record map \(a_J\) sends its moving image and branch divisors to
\(\operatorname{Sym}^{h_J}P\times\operatorname{Sym}^{2g_J}P\)
and its marked points to the corresponding surface factors. Each
fixed remainder has a point for its position space. Finite flat
pullback orders each moving divisor. The empty label set gives \(M_0\).

These records are additive under cutting, including over nonreduced
bases. The determinant and canonical section of the image complex
restrict to each occupied open \cite[\S4.1, (46)]{KoolThomas}.
The equation of the fixed part \(M_0+V\) identifies
\(H^0(P,\mathcal O_P(h))\) with \(H^0(S,\mathcal O_S(B))\), so its
removal commutes with restriction. To obtain the branch divisor,
forget the marks and stabilize \(pf\). The horizontal component
has degree one, and each vertical tail meets it at just one position.
Cutting away other positions therefore preserves the special points
of the tail. The stabilized pieces glue to \(\bar u\). Its branch complex
\(R\bar u_*[\bar u^*\omega_P\to\omega_{\bar C/T}]\) commutes with
restriction and base change and is canonically acyclic on the overlap
\cite[\S\S2.3,3.2]{FPBranch}. Its divisor is the sum of the branch
divisors of the pieces. Each marked evaluation stays on its own piece.

Write \(o:Y\to Y'\) for the ordering of moving image and branch
positions. Addition of disjoint partial divisors is \'{e}tale, and
ordering within each part is finite flat even at collisions. On a
disjoint-position chart, Lemma~\ref{vertex-derived-cut} gives the
K\"unneth decomposition of the homotopy sheaves of the reduced derived
map space \(K_B^{\rm red,der}\). Flat base change and K\"unneth for
the proper record maps identify
\[
 R^i(a_B)_*\bigl(\mathcal H^{-j}(\mathcal O_{K_B^{\rm red,der}})
                          \otimes\mathcal V_B\bigr)
\]
with the direct sum of exterior products for \(i_1+i_2=i\) and
\(j_1+j_2=j\), including the summands of the split insertion
\(\mathcal V_B\). Finite flat pullback preserves this decomposition
at collisions within a block. For the stack pushforward, the finite
stabilizers in characteristic zero are linearly reductive. Pushforward
to the coarse space is exact, and on quotient charts invariants
for a product of stabilizers commute with exterior products
\cite[Theorem~3.2 and Corollary~3.3]{AOVTame}. The ordinary
K\"unneth formula over the ground field then applies on affine opens
of the record targets. Functoriality of cutting and of the
K\"unneth maps makes it compatible with refinements. Each attachment
contributes twice its genus, which accounts for the empty support of
blocks with an odd number of branch labels on this open.

\smallskip\noindent\emph{A faithful bundle compatible with cuts.}
For the insertion \(\mathcal V_B\), take a basepoint-free divisor
\(H\) contracting \(M_0\) and ample off it. A sufficiently divisible
multiple of \(2A+(A\cdot M_0)M_0\) for \(A\) ample has these
properties. Choose a section of \(H\) nonvanishing on \(M_0\).
For \(u=p f\) and marking divisor \(\Sigma\), set
\[
 L'=\omega_{\mathcal C/K}(\Sigma)\otimes u^*\omega_P^{-1}
                         \otimes f^*\mathcal O_S(3H).
\]
The differential of \(u\), the marks and the cube of the chosen
section give a section \(s\) of \(L'\) invertible off the recorded
positions. The degree on the horizontal component is the number of
special points. On a nonconstant vertical component the last factor
offsets any negative dualizing degree, and on a constant component
stability gives positive degree. Hence \(L'\) is ample on each
nonbare domain. The bare section has trivial stabilizer.

The finitely many label subsets give bounded families. Choose a
common \(k\) so that \((L')^k\) is very ample on each nonbare fibre
and its higher direct images under \(u\) vanish on every geometric
fibre. Over a parameter scheme \(T\), consider the cone of the
section \(\mathcal O_{P_T}\to Ru_*(L')^k\). The graph of \(u\) in
\(C_T\times_T P_T\) is a regular immersion since \(P_T/T\) is smooth.
Projection to \(P_T\) is a base change of the proper flat nodal
curve \(C_T/T\), so \(u\) is proper and perfect and the cone is
perfect. Flatness of \(C_T\) and \(P_T\) over \(T\) gives
Tor-independent base change to the geometric fibres.

On a geometric fibre the section is nonzero on the horizontal component,
so the map of sheaves is injective. Its cokernel is supported at the
recorded positions. The fibre criterion for relatively perfect complexes
makes the cone a coherent sheaf flat over \(T\). Its support is finite
over \(T\), being proper and quasi-finite. Pushforward thus gives a vector
bundle \(W\) with fibre \(H^0((L')^k)/\langle s^k\rangle\), since
\(H^1(P,\mathcal O_P)=0\).

A stabilizer element acts on \(H^0((L')^k)\) by an operator \(T\)
fixing \(s^k\). If it acts trivially on the quotient \(W\), then
\((T-1)H^0((L')^k)\subseteq\langle s^k\rangle\subseteq\ker(T-1)\).
For its order \(m\) we therefore have \(1=T^m=1+m(T-1)\).
Characteristic zero gives \(T=1\), and very ampleness forces the
element to fix the domain. Thus \(W\) is faithful. Over disjoint position sets its finite
cokernel splits canonically, giving \(W=W_1\oplus W_2\).

\smallskip\noindent\emph{From disjoint positions to diagonal classes.}
Choose a common multiple \(m\) of the stabilizer orders and a
primitive \(m\)-th root of unity \(\zeta\). In the Grothendieck
ring of bundles, set
\[
 m^{-\operatorname{rk}W}\prod_{j=1}^{m-1}\lambda_{-\zeta^j}(W),
       \qquad \lambda_t(W)=\sum_i t^i[\textstyle\bigwedge^i W],
\]
where the displayed class is denoted by \(D(W)\).
On a line with character \(x\) it is
\((1+x+\cdots+x^{m-1})/m\). It is multiplicative under direct sums
and has rank one. At a nonidentity stabilizer element faithfulness
gives an eigenvalue \(x\ne1\) with \(x^m=1\), so its character
there is zero.

Write \(Y'_J\) for the target of \(a_J\), \(o_J:Y_J\to Y'_J\)
for the ordering map, and \(I_J\) for the marks in \(J\).
With the same \(k,m,l\) for every block, set
\[
 F_J=o_J^*(a_J)_*\!\left(
 [\mathcal O_J^{\rm vir}]\,D(W_J)^l
 \prod_{i\in I_J}(1-[L_i^{-1}])^{a_i}\right)
 \quad\text{in }G_0(Y_J)_{\mathbb C}.
\]
Here \([\mathcal O_J^{\rm vir}]\) is the alternating class of the
homotopy sheaves on the reduced map space for \(J\). A block with
an odd number of branch labels has \(F_J=0\).
Expand each factor. Fix the exterior degree in each
copy of each \(\lambda_{-\zeta^j}(W)\), the power of \(L_i^{-1}\)
at each mark, and the total homotopy and direct-image degrees.
For each partition, distribute these degrees among its blocks,
keeping each marked factor on its block. Take the direct sum of the
resulting exterior products of coherent sheaves. The splitting of
\(W\) and the preceding K\"unneth isomorphisms identify these sheaves
on the required opens, compatibly with further splits. Thus
Lemma~\ref{med-g-mobius} applies before restoring the fixed signs
and scalar weights. The normalization also factors, since the rank
\(k(2g+|\Sigma|+3H\cdot f_*[\mathcal C])-g\) from Riemann--Roch
adds under cutting.

Partition inversion gives \(F_{A'}=\sum_\pi\boxtimes_{J\in\pi}C_J\)
because \(\sum_{\pi\ge\rho}\mu(\rho,\pi)\) vanishes except at
the one-block partition, where it is one.
For the inclusion \(i_J:Z_J\to Y_J\) of the locus
where the positions of \(J\) coincide, Lemma~\ref{med-g-mobius} gives
\(C_J=(i_J)_*\gamma_J\) for some \(\gamma_J\in G_0(Z_J)\).
Proper covariance gives \(\tau_{Y_J}(C_J)=(i_J)_*\tau_{Z_J}(\gamma_J)\).
Multiplicativity on exterior products then preserves all block supports.

\smallskip\noindent\emph{Recovering the reduced evaluation class.}
Set \(v=g+|I|-\sum_i a_i\). The alternating sum of the homotopy
sheaves of the quasi-smooth reduced map stack is its virtual
structure sheaf. On the identity sector its Riemann--Roch class is
the Todd class times the reduced virtual cycle
\cite[Remark~6.3, Theorem~6.12 and Remark~6.13]{KhanVGRR}.
The projective-frame presentation of stable maps has a linearized
ample line and embeds equivariantly as a closed Hilbert locus in
an invariant open of projective space. Its resolution property gives
global two-term resolutions of the perfect obstruction theories, as
required for the quotient-stack formula \cite[Theorem~1.1]{RaviSreedhar}.
For the proper record maps to \(Y'\) we use the inertia
Riemann--Roch transformation for quotient stacks
\cite[Theorem~1.2]{KrishnaSreedhar}. On the sector of an automorphism
\(\gamma\), multiplication by a bundle class acts by the Chern character
of its \(\gamma\)-weighted eigenspaces, by the projection formula
and \cite[Corollary~1.4 and \S4.1, (4.2)]{RaviSreedhar}.

Let \(q:K\to\mathfrak M^{\rm pre}\) forget the map and retain its
pointed prestable domain. The ordinary absolute virtual tangent
class is \([R\pi_{C*}f^*T_S]+q^*[T_{\mathfrak M^{\rm pre}}]\).
For the determinant fibre, the tangent triangle of the homotopy
pullback adds one copy of \(H^2(S,\mathcal O_S)\) in \(K^0(K)\)
\cite[Proposition~4.6 and Theorem~4.11]{STV}.
For the fixed-surface theory the dual of the semiregularity cone
adds the same line \cite[(15) and Theorem~2.4]{KoolThomas}.
Thus both absolute theories on the same classical stack \(K\) have class
\[
 [T^{\rm vir}]=[R\pi_{C*}f^*T_S]+q^*[T_{\mathfrak M^{\rm pre}}]
                  +[H^2(S,\mathcal O_S)\otimes\mathcal O_K].
\]
The frame presentation gives global two-term resolutions on \(K\).
Apply Siebert's formula to \(K\) over the ground field. It takes
the component of \(c(T^{\rm vir})^{-1}\cap c_F(K)\) in the virtual
dimension. Both the expression and its dimension agree, giving
the fixed-surface class \([K]^{\rm red}\) for the determinant fibre
\cite[Theorem~4.6]{SiebertVirtual}, \cite[Theorem~9.2.2]{LuWebbVirtual}.

On the identity sector, \(D(W)\) has rank one and the first term of
\(1-[L_i^{-1}]\) is \(\psi_i\) in codimension one. After these
insertions the identity-sector class has dimension at most \(v\).
Its dimension-\(v\) term is
\(a_*(\prod_i\psi_i^{a_i}\cap[K]^{\rm red})\).
On every nonidentity sector the traced Chern character of \(D(W)\)
has zero constant term. The moving Euler factor in the inertia
formula is invertible in degree zero, since every moving eigenvalue
differs from one. Its inverse and the Todd factor have only
nonnegative codimensions. Multiplication by \(D(W)^l\) therefore
lowers dimension by at least \(l\). Choose \(l\) greater than the
dimensions of the inertia components of all the map stacks occurring
in the records.
Their rational Chow groups agree with those of their
coarse spaces and vanish in negative dimensions
\cite[\S10.1]{KrishnaSreedhar}, so only the identity sector contributes.

The ordering map \(o\) is finite flat between smooth varieties, hence
locally a complete intersection. Its relative Todd class has constant
term one, so the dimension-\(v\) term on \(Y\) is the flat pullback of
the cycle above. The image hyperplane pulls back to a sum of point
divisors on the ordered image factors. Capping with its \(e\)-th
power gives a sum of block products in dimension \(v-e\).
The finite flat map has rank \(h!(2g)!\), including over collisions.
By push--pull, forgetting the order and the positions gives this
rank times the desired evaluation. Dividing proves the decomposition
with complex coefficients.

The quotient of \(\CH_*(S^I)_{\mathbb Q}\) by the right-hand side
of \eqref{med-g-factor} injects into its complexification, so
the decomposition holds with rational coefficients.

For an algebraically closed ground field of characteristic zero,
choose a finitely generated field of definition for the pencil,
the insertions and the reduced obstruction theory, and embed it
in \(\mathbb C\). Over its algebraic closure, spread the finitely
many block cycles and rational equivalences from \(\mathbb C\)
over a smooth integral parameter scheme. Specialization at a closed
point preserves proper pushforward and exterior products. The
evaluation specializes to the original class, since its theory
and insertions are pulled back from the field of definition.
Extending to the given ground field proves \eqref{med-g-factor}.
Proposition~\ref{rd-g-graph} gives the signed graph assertion.
\end{proof}

Once Proposition~\ref{fixed-evaluation} places the evaluation on the
signed graph line, pairing with a class that is nonzero on its
generator determines the scalar. Here \(B\cdot F=1\) makes \(B\) primitive, so
\cite[Corollary~3]{BuellesK3Classes} applies to this numerical comparison.
Extending the generic family identities from their supporting isogeny
curves can leave classes over finitely many base point pairs. On a
fine cover, frame the cyclic subgroup on the source. Its quotient
defines a graph in the product with the smooth target curve, so each
branch is a Cartier divisor even where branches meet over a CM pair.
A point outside the chosen divisor has empty pullback. For a point on
it, the self-intersection formula caps a point-fibre class with
the first Chern class of a line pulled back from that point, which is
zero. Thus restriction to a cyclic level kills all these error terms.

\subsection{The two vertex classes in the family}
\begin{lemma}[The period equation on infinitesimal families]\label{picard-artin-support}
Let \(S\to B\) be a smooth projective family of complex K3 surfaces.
On a marked local period chart at \(b\), write \(s_\gamma=(\gamma,\omega)\)
for the equation of a fixed integral class \(\gamma\), using a local
generator \(\omega\) of the Hodge line. The classical stack of stable
maps of class \(\gamma\) factors through \(s_\gamma=0\), with its
scheme structure. This remains true with arbitrary marked flags.

In the framed product Kummer family, let \(X=(s_L=0)\) be a smooth
total-class sheet in the three-dimensional period base. At a generic
isogeny point let \(z\) be the extra mixed direction. Then
\(ds_L,ds_z\) are independent. If two vertex labels sum to \(L\)
and have opposite nonzero mixed parts, their two period sheets meet
with intersection multiplicity one on a transverse surface slice.
\end{lemma}
\begin{proof}
Test the first assertion on a local Artin \(\mathbb C\)-scheme \(T\)
at \(b\). For a family \(f:C\to S_T\), the line bundle
\(\det Rf_*\mathcal O_C\) commutes with base change. Indeed the
graph of \(f\) is a regular immersion into \(C\times_T S_T\), and
projection to \(S_T\) is the base change of the proper flat family
of nodal curves. This gives a proper perfect pushforward and
Tor-independent base change. On the special fibre its first Chern
class is \(\gamma\): at each prime divisor the determinant records
the generic length of the image complex, hence the multiplicity of
the pushed-forward curve. Contracted components contribute only in
codimension two. This is the determinant construction for stable
maps in \cite{STV}.

The analytic exponential sequence on \(S_T\) is exact also on a
nilpotent thickening, by the finite logarithm on nilpotent units.
The analytic spaces \(S_T\) and \(S_b\) have the same underlying
topological space. By cohomology and base change for the smooth K3
family, \(H^1(S_T,\mathcal O_{S_T})=0\) and
\(H^2(S_T,\mathcal O_{S_T})\) is free of rank one over
\(\Gamma(T,\mathcal O_T)\). Exactness of the exponential sequence
and GAGA for the proper scheme \(S_T\) therefore identify
\[
 \operatorname{Pic}(S_T)\simeq
 \ker\!\left(H^2(S_b,\mathbb Z)\longrightarrow
                         H^2(S_T,\mathcal O_{S_T})\right).
\]
Under the relative de Rham comparison and Serre duality, pairing
the image of \(\gamma\) with \(\omega\) gives \(s_\gamma|_T\).
Thus a line bundle of class \(\gamma\) exists on \(S_T\) exactly
when this period vanishes, and its isomorphism class is then unique.
In particular, the determinant of a stable map supplies the required
line bundle on every Artin test of the map stack. To see the
scheme structure, take the local ring \((A,\mathfrak m)\) at any
complex point of a Noetherian analytic atlas. The pulled-back period vanishes in
\(A/\mathfrak m^j\) for every \(j\), hence lies in
\(\bigcap_{j\ge1}\mathfrak m^j=0\) by the Krull intersection theorem.
The factorization descends from the atlas. On an algebraic Picard
chart, faithful flatness of completion detects the same vanishing
of the defining ideal. Retaining marked flags leaves the determinant
of the map unchanged.

On \(X\), a branch of the isogeny locus is given in elliptic period
coordinates by
\[
       (c\tau_1+e)\tau_2-(a\tau_1+b')=0,\qquad ae-b'c>0 .
\]
This is the mixed period equation, up to a unit and a nonzero
constant after choosing the scale of \(z\). Its derivative in
\(\tau_2\) is \(c\tau_1+e\ne0\). Thus \(ds_z\) restricts
nontrivially to \(X\). Modulo the polarization lattice, the two
vertex labels have coordinates \((u_1,t)\), \((1-u_1,-t)\) in
the \(L,z\) directions. The determinant is \(-t\ne0\). Their
periods are consequently regular parameters on the transverse
slice, giving the stated length one.
\end{proof}

\begin{lemma}[Restriction to the generic points of the support]
\label{generic-labelled-support}
Let \(X\) be a smooth surface of finite type over \(\Bbbk\) and let
\(Y\to X\) be of finite type. Let
\(C=\bigcup_i C_i\) be a finite union of distinct closed integral curves in
\(X\), with generic points \(\eta_i\). Modulo classes pushed forward
from fibres over closed points of \(C\), restriction induces an
isomorphism
\[
       \CH_l(Y_C)\ \big/\!
        \sum_{x\in C\ {\rm closed}}\operatorname{im}\CH_l(Y_x)
          \ \simeq\ \bigoplus_i\CH_{l-1}(Y_{\eta_i}).
\]
The dimension on the right is relative to \(k(\eta_i)\).
\end{lemma}
\begin{proof}
Remove the intersection points and the singular points of the
curves. The remaining base is a disjoint union of smooth curves.
A generic integral cycle extends by closure, increasing its dimension
by one. For the kernel, write a generically trivial restriction as
a finite sum of principal divisors on each generic fibre. Take the
closures of their supports and the same rational functions on those
closures. After subtracting the resulting divisors, the cycle is
supported on finitely many closed fibres. Conversely, every class
pushed forward from such a fibre restricts to zero at every \(\eta_i\).
Localization also accounts for the intersection and singular points
removed at the outset.
The same proof applies with rational coefficients on a fine cover.
\end{proof}

In particular, a comparison made on the labelled support can be
pushed into \(Y\) modulo point-fibre classes. A relation already
pushed into \(Y\) need not determine separate relations on the
curves \(C_i\).

\begin{proposition}[The family class on its supporting curves]\label{family-supported-diagonal}
Let \(X\) be the product of the two smooth modular curves with
two-torsion framed, using a fine cover. In the first pencil of the
product Kummer family fix an extending effective total class \(L\)
of first degree two and zero mixed part. Use the component of the ordinary three-dimensional family theory
labelled by \(L\), with flat pullback along the framed cover of its
period sheet as described below.
In a decorated boundary term
with \(r\ge2\) nonconstant vertices, the part on which all their
classes extend vanishes in rational Chow on the family. Modulo
classes over finite base point pairs, every remaining contribution
comes from an isogeny curve. At its labelled generic point only
\(r=2\) contributes, with multiplicity one from the base intersection.
The two factors are the fixed-surface reduced classes with every
original flag retained, before matching node evaluations.
The same support and multiplicity assertions hold on the marked
cyclic Kummer family of Section~\ref{rd-m-twist} for the total
classes used there. In each case, the pushforwards from the
labelled curves sum to the decorated graph term modulo point
fibres.
\end{proposition}
\begin{proof}
\smallskip\noindent\emph{The period base.}
Fix a relatively ample integral class \(H_0\). The two axes and the
sixteen exceptional curves span the generic rational Picard lattice
\(\Lambda_{\mathbb Q}\), and the framing fixes its saturation.
Only finitely many extending classes \(D\) can occur at a
nonconstant vertex. Indeed, if \(\sum m_j C_j\) represents \(D\),
then \(\sum m_j\le H_0D\le H_0L\). Adjunction and nonnegative
intersections of distinct components give \(D^2\ge-2(H_0L)^2\).
It follows that the projection of \(D\) onto the negative definite
space \(H_0^\perp\cap\Lambda_{\mathbb R}\) lies in a bounded set.
Perturb \(H_0\) in the relatively ample cone and clear denominators
to obtain an integral ample class \(H\) outside the lines spanned by
\(L\) and these \(D\). Choose a rational functional \(u\) with
\(u(H)=0\), \(u(L)=1\), and \(u(D)\ne0\) for all such \(D\).

In the fine moduli space of \(H\)-polarized K3 surfaces, mark sixteen
further independent classes of \(\ker u\cap\Lambda\) in the
rigidified Picard space. Together with \(H\) they span this rank-seventeen
lattice rationally. Along the framed Kummer family the sixteen period
equations have independent differentials by local Torelli: their
span together with \(H\) is nondegenerate of signature \((1,16)\).
The labelled locus \(B\) is therefore smooth of dimension three
near \(X\), with K3 family \(\pi:S\to B\) \cite{MPPeriods}.
Pass to a common fine
level on \(B\) and \(X\). Push--pull for this finite \'{e}tale map
descends the resulting rational Chow identities. Write
\(\mathcal P=\operatorname{Pic}_{S/B}\) for the rigidified relative
Picard space and \(\mathcal L=R^2\pi_*\mathcal O_S\). The sheet
\(P_L\) specifying the total class has equation \(s_L=0\) by
Lemma~\ref{picard-artin-support}. The map \(X\to P_L\) is
\'{e}tale because the two derivatives of the mixed period vector
\((1,\tau_1,\tau_2,\tau_1\tau_2)\) are independent modulo the
period line. The determinant map of the ordinary component
\(K_{B,L}\) factors through \(P_L\). We therefore use on
\(K_{X,L}=K_{B,L}\times_{P_L}X\) the flat pullback of the ordinary
virtual class. Any operational relation on stable curves still
caps this class to zero.

\smallskip\noindent\emph{Gluing with the vertex classes recorded.}
Order the flags. Pulling the decorations to the bounded locus of
the smooth stack of prestable curves met by the graph gives the
prestable gluing of \cite[Proposition~3.13, formula~(25)]{BaeSchmitt}.
The proper birational map from this gluing space to the pullback
of the stable gluing space accounts for the chains contracted in
stabilization. Label constant components by zero and
positive-degree components by one, with \(1+1=1\).
For stabilized cotangent classes the rational-tail expansion
retains the original flags \cite[Proposition~3.14]{BaeSchmitt}.

For each of the \(k\) vertices, let \(V_i\) be the ordinary virtual
class on its stack \(K_i\) over a separate copy of \(B\).
The virtual tangent complex relative to the prestable curve and the
base is \(R\pi_{i*}f_i^*T_{S/B}\). Before proper pushforward the
ordered gluing class, with its decorations understood, is
\[
 \Delta_{\rm ev}^!\,(\Delta_B^{(k)})^!
                          \bigl(\boxtimes_{i=1}^k V_i\bigr).
\]
First \((\Delta_B^{(k)})^!\) identifies the base parameters. Then
\(\Delta_{\rm ev}^!\) matches the evaluations at both flags of each
edge, including a loop. Normalizing
the source at the nodes gives the triangle
\[
 R\pi_{C*}f^*T_{S/B}\ \longrightarrow\
 \bigoplus_i R\pi_{i*}f_i^*T_{S/B}\ \longrightarrow\
 \bigoplus_e {\rm ev}_e^*T_{S/B}\ \longrightarrow\
 (R\pi_{C*}f^*T_{S/B})[1].
\]
The middle arrow takes the difference at the two flags.
Compatibility of the ordinary obstruction theories gives the
gluing formula \cite[\S3]{BehrendCuttingEdges}. It preserves the
total determinant label \(L\). On the node-matching stack \(T\), the
normalization sequence has node term \(s_*\mathcal O_T\), where
\(s:T\to S_T\) is the corresponding section of the relative surface.
The section is a regular immersion of codimension two. Its Koszul
resolution has canonically trivial determinant, since the determinants
of the two conormal terms cancel. The same resolution remains exact
after any base change of the node-matching parameter stack \(T\).
Hence the node term contributes no determinant, and the vertex labels
add under gluing.

The tuples of determinant labels with this sum are parametrized by
\[
  \mathcal Q=(\mathcal P^{\times_B k})\times_{\mathcal P} X ,
\]
where the maps to \(\mathcal P\) are addition and the section for
\(L\). Smooth geometrically integral K3 fibres with
\(H^1(\mathcal O)=0\) make \(\mathcal Q\to X\) separated and
unramified. Each extending tuple gives a section. Indeed, the
first-axis direction on the exceptional curve over the origin gives
a rational point of the generic Kummer surface, so its invariant
generic classes descend to line bundles. Closure extends a generic
divisor on the smooth total space. Vertical prime divisors are
pulled back from the base, so the extension defines a relative
Picard section. Each section is open and closed by the diagonal
criterion.

There are finitely many vertex degree distributions. For each one,
the product of the vertex stacks with common base and total label
\(L\) is proper over \(X\), since the total determinant condition
is open and closed in the proper fixed-degree stack over \(X\).
Its image in \(\mathcal Q\) is finite over \(X\), being
quasi-finite, separated and universally closed. Away from the
extending sections an additional Picard class requires an isogeny
between the two elliptic factors. By normalization and generic
separability, after removing finitely many points of \(X\) this
part of the reduced support consists of disjoint smooth curves,
each finite \'{e}tale over its image. Their preimages split the
matching stack there into open and closed summands. On any such
summand \(W\), reduction induces an isomorphism
\(i_*:\CH(W_{\rm red})\to\CH(W)\). Transport the class by
\(i_*^{-1}\) to retain its virtual multiplicities. Proper
covariance of refined pullback then applies at the nodes.

\smallskip\noindent\emph{The excess on the common period locus.}
For the \(r\) nonconstant vertices let \(P_i\to B\) be their smooth
labelled sheets in the relative Picard space. In \'{e}tale charts
they are algebraic closed immersions whose completed ideals are the
period equations. Put \(P_i=B\) at the constant vertices. The
factorization \(K_i\to P_i\) holds scheme-theoretically by
Lemma~\ref{picard-artin-support}. Let \(Z\) be a smooth component
where the sheets meet cleanly in codimension \(d\) in \(B\), and
write \(j_i:Z\to P_i\). Below \(j_i^!V_i\) is the refined pullback
along \(j_i\) for the map \(K_i\to P_i\); throughout this comparison
the vertex class \(V_i\) is kept on \(K_i\). Set
\[
  Q_Z=\operatorname{coker}\!\left(N_{Z/B}\longrightarrow
                      \bigoplus_{i=1}^r N_{P_i/B}|_Z\right).
\]
Excess intersection on the matching stack over \(Z\) gives
\begin{equation}\label{eq:family-clean-excess}
 (\Delta_B^{(k)})^!(\boxtimes_i V_i)
  =c_{r-d}(Q_Z)\cap
     (\Delta_Z^{(k)})^!\bigl(\boxtimes_i j_i^!V_i\bigr).
\end{equation}
Indeed the normal map for \(Z\to\prod_iP_i\) is
\[
 (\textstyle\bigoplus_i T_{P_i}|_Z)/T_Z
       \ \hookrightarrow\ (T_B|_Z)^{\oplus k}/T_B|_Z .
\]
Its cokernel is \(Q_Z\). Factoring \(Z\to\prod_iP_i\) through \(Z^k\)
gives \eqref{eq:family-clean-excess}. These normal sequences glue
because the sheets lie in the rigidified Picard space and their
intersection is the pullback of the base diagonal.

For an extending tuple, each \(D_i-u(D_i)L\) belongs to the
rational polarization lattice. After clearing denominators, \(D_i\)
and \(L\) thus define the same sheet, including on Artin tests.
The normal lines admit a common algebraic description. On a sheet
\(P_D\), cup product of the Kodaira--Spencer class with the de Rham
class \(c_1(D)\) induces \(N_{P_D/B}\to\mathcal L|_{P_D}\).
This is the first-order obstruction to extending \(D\), so the
smooth period equation identifies the two lines. The map is additive
in \(D\) and descends from charts carrying a representative line bundle.
Thus on \(Z=X\), where \(d=1\), the normal map has entries \(u(D_i)\), giving
\(Q_X\simeq\mathcal L^{\oplus(r-1)}|_X\). The two elliptic
discriminants trivialize the twelfth power of the Kummer Hodge line
on \(X\). Thus \(c_1(\mathcal L|_X)=0\) in rational Chow, and for
\(r\ge2\) the excess factor vanishes everywhere on this component.

On a remaining curve the period rows at the generic point are
\(a_i s_L+c_i s_z\). They sum to \(s_L\), with \(\sum_i a_i=1\)
and \(\sum_i c_i=0\). At least one \(c_j\ne0\), and an extending
label has nonzero row by the choice of \(u\). The sum of the rows is
\(s_L\). Subtract \(a_j s_L\) from row \(j\) and divide by the
nonzero scalar \(c_j\) to recover \(s_z\). Thus the rows generate
the ideal \((s_L,s_z)\), including its scheme structure. By
Lemma~\ref{picard-artin-support} its zero locus is smooth of codimension
two here. For \(r>2\) the excess bundle has positive rank, so its top Chern class
vanishes over the function field. Only closed points can contribute.
For \(r=2\) the
relations \(a_1+a_2=1\) and \(c_1+c_2=0\) give
\(a_1c_2-a_2c_1=-c_1\ne0\). Thus the two period equations cut the
curve with multiplicity one at its generic point.

\smallskip\noindent\emph{The two classes on an isogeny curve.}
Extend a separating parameter of the labelled curve to a smooth
map \(B\to T\) near its generic point. The base diagonal factors
through \(B^{\times_T k}\). By functoriality of refined pullback,
flat restriction to \(k(T)\) gives the diagonal of the smooth surface
\(B'\) on the restricted exterior product. After a finite separable
field extension the two period curves \(D_1,D_2\subset B'\) meet
transversely at a rational point \(y\). With \(j_i:y\to D_i\),
the excess in \eqref{eq:family-clean-excess} is zero, and hence
\begin{equation}\label{rechecked-family-diagonal}
 \Delta_{B'}^!(V_1\boxtimes V_2)
                   = j_1^!V_1\boxtimes j_2^!V_2 .
\end{equation}
Here \(T_yD_1\oplus T_yD_2\to T_yB'\) is the difference isomorphism.
Transversality makes the square of the two curves over \(B'\)
Tor-independent. Compatibility of refined Gysin maps with this
base change identifies pullback along \(D_2\hookrightarrow B'\)
on classes of the first vertex stack with \(j_1^!\), for any map
of that stack to \(D_1\). By the scheme-theoretic period factorization
the resulting classical stack \(K\) is the fixed-surface map stack
with its original flags.

Let \(\mathfrak M^{\rm pre}\) denote the smooth pointed prestable stack.
Along the regular immersion
\(D_2\times\mathfrak M^{\rm pre}\hookrightarrow
B'\times\mathfrak M^{\rm pre}\), proper base change for the flat
universal curve pulls back the mapping complex. Virtual base change
thus identifies \(j_1^!V_1\) with the ordinary family class over
\(D_2\) \cite[Proposition~7.2]{BehrendFantechiCone}. Relative to
\(\mathfrak M^{\rm pre}\) its virtual tangent class is
\[
 [R\pi_{C*}f^*T_{S_y}]+[T_yD_2\otimes\mathcal O_K].
\]
The period derivative, given by cup product of the Kodaira--Spencer
class with the first vertex class, has kernel \(T_yD_1\). It identifies
\(T_yD_2\) with the line \(H^2(S_y,\mathcal O_{S_y})\) added by the
reduced theory \cite[Theorem~2.4 and Lemma~2.5]{KoolThomas}.
Passing to the absolute family theory adds \(q^*[T_{\mathfrak M^{\rm pre}}]\),
where \(q\) forgets the map. The resulting class in \(K^0(K)\)
is exactly the absolute reduced tangent class computed in
Proposition~\ref{fixed-evaluation}. The two global two-term theories
on \(K\) have the same virtual dimension and the same Fulton class
\(c_F(K)\). Siebert's formula over the ground field therefore gives
\(j_1^!V_1=[K]^{\rm red}\) in Chow. To use the complex
comparison over the labelled generic point, descend the family,
the slice and both theories to a finitely generated field of
characteristic zero and embed it in \(\mathbb C\). The resulting
rational equivalence spreads over a finite-type model; specialization
at a closed point and normalized trace descend the equality to the
field of definition. Flat extension gives it over the labelled field.
Interchanging the two curves gives the second reduced class. Constant
vertices keep the ordinary obstruction bundle
\(\mathbb E^\vee\otimes T_{S_y}\), with \(\mathbb E\) the Hodge bundle.
Matching the original flags and applying the decorations now gives
the graph contribution, with its automorphism and labelled-cover factors.

\smallskip\noindent\emph{Pushforward from the supporting curves.}
In the unmarked first pencil every vertical component has extending
class. Nefness of the fibre shows that the two mixed degrees are
positive and sum to two, hence both equal one. For fixed \(L\) the
second-pencil degree bounds the isogeny degrees. On a labelled
support curve \(D\), write \(\alpha_D\) for the evaluation after
Kummer transfer, proper pushforward and the signed projector \(e\).
In the comparison below, \([\Gamma_{D,\eta_D}]\) denotes the class
of the isogeny graph with the same projector on the generic fibre.
Over the geometric generic point of \(D\),
Proposition~\ref{fixed-evaluation} puts the signed evaluation on the
graph line in Chow. Its generator has nonzero realization by
Lemma~\ref{med-g-products}. Pairing with a class nonzero on that
generator determines \(c_D\), using the numerical comparison for
primitive classes and matching insertions \cite[Corollary~3]{BuellesK3Classes}.
Compute this ratio of intersection degrees after embedding a finitely
generated field of definition in \(\mathbb C\). It is unchanged by
field extension. Both classes are defined over \(k(D)\), and the
rational equivalence descends to a finite extension. Normalized trace
gives \((\alpha_D)_{\eta_D}=c_D[\Gamma_{D,\eta_D}]\) in Chow over
\(k(D)\). By Lemma~\ref{generic-labelled-support} the difference of
their closures on the normalization is supported on finitely many
closed fibres. Its proper pushforward belongs to \(J_b\).

For the marked cyclic family the isogeny
\(E\times E'\to (E\times E')/\operatorname{graph}(\phi)\)
identifies the rational \(H^1\)-variations. The two axes and sixteen
exceptional classes give the same extending rank. The additional
period has the same simple zero up to nonzero scale. For the
prescribed total class \(T\), use the bound \(H_0T\) in place of
\(H_0L\) when choosing the polarization. Scaling a vertex class
by a nonzero integer scales its period by that integer and preserves
its zero scheme in characteristic zero. The excess calculation and
the identification of the two fixed-surface classes thus carry over.
The degree-one conclusion above uses nefness in the unmarked first
pencil.
\end{proof}

\section{A recurrence for the isogeny graphs}
\label{med-u-five}

The four branch fibres of the Kummer pencil give one recurrence.
We identify its solutions with cusp forms to apply Theorem~\ref{med-u-scalar}.

For even \(b=2d\geq2\), put \(K_b=b+2\) and use the full-level isogeny graphs
\(H_b(e)\) of Theorem~\ref{graph-to-moments}.
Work modulo classes over finitely many base point pairs, which
cyclic restriction kills by the self-intersection formula. For a functional
\(\ell\) on the span of the graph classes put
\(f_\ell=\sum_{e\geq1}\ell(H_b(e))q^e\) and write \(V_b\) for
their space of formal series. Define level-two weight-two forms
\[
 A=(\vartheta_3^4+\vartheta_4^4)/2,\qquad B=\vartheta_2^4/8,
 \qquad R_j=\mathbb C[A,B]_j .
\]
The subscript denotes homogeneous degree, with \(R_j=0\) for \(j<0\).

\subsection{Hecke operators and the space of solutions}
\begin{lemma}\label{full-graph-hecke}
The space \(V_b\) is stable under every full-level Hecke operator
\(T_p\) at a prime \(p\) in weight \(b+2\).
\end{lemma}
\begin{proof}
Postcomposition by the degree-\(p\) correspondence defines a finite
pull-push \(\Psi_p\) preserving support over finite sets of base point
pairs. A terminal kernel of order \(pn\) has one index-\(p\) subgroup
if its \(p\)-primary part is cyclic, and \(p+1\) otherwise. In the latter
case the kernel contains \(E[p]\), so the isogeny factors through \([p]\)
and an isogeny of degree \(n/p\). On each odd elliptic factor both
\([p]^*\) and \([p]_*\) act by \(p\) \cite[Theorem~3.1]{DeningerMurre}.
The \(p\) extra factorizations thus give
\[
 \Psi_p H_b(n)=H_b(pn)+p^{b+1}H_b(n/p),
\]
with the second term zero if \(p\nmid n\). Applying \(\ell\) gives
\(T_p f_\ell=f_{\ell\circ\Psi_p}\).
\end{proof}

\begin{proposition}[A criterion for full-level cusp forms]
\label{med-u-five-prop}
Let \(W_b\) be the space of series \(y=\sum_{t\ge2}y_tq^{t/2}\)
satisfying
\begin{equation}\label{rd-m-rec}
 P_b(t-2)y_t+
 \sum_{1\le j\le(t-2)/2}C_b(j,t-2-j)y_{t-2j}=0
 \qquad(t\ge2),
\end{equation}
where \(\deg P_b=d=b/2\). Suppose \(R_4 V_b\subset W_b\), and,
when \(d\ge5\), that
\(\dim\bigl(W_b/(W_b\cap\Delta_{\rm mod}R_{d-1})\bigr)\le4\).
Then \(V_b\subset S_{b+2}(\mathrm{SL}_2(\mathbb Z))\).
\end{proposition}
\begin{proof}
Two solutions agreeing at every index \(t\) with \(P_b(t-2)=0\)
agree everywhere. At a first differing index the earlier terms in
\eqref{rd-m-rec} cancel, forcing \(P_b(t-2)=0\). Restriction to
the root indices is thus injective, so \(\dim W_b\le d\).
For \(0\ne f\in V_b\), the five products \(A^{4-j}B^jf\) have
distinct initial exponents and hence require \(d\ge5\).
For \(d\ge5\), the map
\[
 R_4\longrightarrow W_b/(W_b\cap\Delta_{\rm mod}R_{d-1}),\qquad P\longmapsto Pf
\]
has a nonzero kernel since \(\dim R_4=5\). A nonzero multiplier in
the kernel expresses \(f\) as a modular ratio of weight \(b+2\)
on \(\Gamma(2)\). Its integral exponents give invariance under
translation. Multiplying numerator and denominator by the translate
of the denominator under \(\tau\mapsto\tau+1\) expresses the ratio
at level \(\Gamma_0(2)=\Gamma_1(2)\). Multiplication
by \(A^4\) embeds \(V_b\) into \(W_b\), so \(V_b\) is finite-dimensional.
The coefficient formula at \(3\) is the same at levels one and two.
Lemma~\ref{full-graph-hecke} therefore makes \(V_b\) stable under
\(T_3\), and Schaeffer's theorem gives holomorphy at level two
\cite[Theorem~1.1(a)]{SchaefferStability}.

The formal full-level translate \(T_2^{(1)}f\) belongs to \(V_b\).
The level-two operator \(U_2^{(2)}\) also preserves modularity there.
Comparing coefficients gives
\(T_2^{(1)}f-U_2^{(2)}f=2^{b+1}f(2\tau)\). Invariance of
\(f(2\tau)\) under the lower unipotent of entry two gives invariance
of \(f\) under the lower unipotent of entry one. Together with
translation this gives full level, and the expansion vanishes at the cusp.
\end{proof}

\subsection{Counting the four branch fibres}
Let \(J_b\) denote the classes supported over finite base point pairs.
Write \(Y_{j,t}=[q^{t/2}]A^{4-j}B^j\sum_{n\ge1}H_b(n)q^n\).

\begin{proposition}[The four branch fibres]\label{family-unmarked-rows}
For the unmarked product Kummer pencil there are a polynomial
\(P_b\) of degree \(d=b/2\) and coefficients
\(C_b(v,\nu)\), independent of the framing, such that for
\(0\le j\le4\) and \(t\ge2\) one has in
\(\CH^{b+1}(U_b\times U_b)/J_b\)
\begin{equation}\label{graph-family-system}
 P_b(t-2)Y_{j,t}+
 \sum_{1\le v\le(t-2)/2} C_b(v,t-2-v)Y_{j,t-2v}=0 .
\end{equation}
On the two-torsion cover, fix the selected branches and write \(D\)
for the tuple of prescribed target labels at the other branches.
One may multiply the contribution of each tuple by a scalar
\(\lambda(D)\) before summing. The resulting coefficients satisfy
\eqref{graph-family-system} with the same \(P_b\) and \(C_b\),
provided the same function \(\lambda\) is used at every index.
\end{proposition}
\begin{proof}
\emph{The two sections.}
Label the two-torsion over the product of fine modular curves.
On the first Kummer pencil write \(O,F\) for the zero section and
fibre and \(A_{u,v}\) for the exceptional curve at \((u,v)\). Select
\(j\) first labels \(J\) and prescribe a target label \(D(u)\) for
each \(u\notin J\). The extending integral class for \(n\ge1\) is
\[
 L=2O+nF+\sum_{u\notin J}(A_{u,0}-A_{u,D(u)})
                         +\sum_{u\in J}A_{u,0}.
\]
It has fibre degree two, exceptional intersections two at the
prescribed labels and zero elsewhere, and square \(4n-8+2j\).
For \(2n+j\ge3\), K3 Riemann--Roch gives effectivity. Apply the
squarefree double ramification relation of Section~\ref{rd-g-section}.
Modulo \(J_b\), Lemma~\ref{family-single-positive} and
Proposition~\ref{family-supported-diagonal} give two reduced vertices
over generic isogeny points. At each node the two evaluations are
matched along the diagonal of the K3 surface.
Each vertex has fibre degree one. Pulling its section back to
\(E_1\) gives a map intertwining inversion. The two sections
\(M_1,M_2\) thus have affine maps \(\eta+t_1,-\eta+t_2\)
with \(e=\deg\eta>0\) and two-torsion translations \(t_i\).
A vertical divisor changes the four outer intersection parities
at a branch either all together or not at all, since the central
component meets each outer component once. The intersections of
\(L\) are even, so the two sections meet the same outer component
at each branch. Write \(t_1=t_2=t_0\) and \(G(u)=\eta(u)+t_0\).

\smallskip\noindent\emph{The vertex classes and their coefficients.}
At each branch the central component \(C_u\) and the three outer
components other than \(A_{u,G(u)}\) span a negative \(D_4\)
orthogonal to both sections. Let \(R\) be the orthogonal sum of the
four lattices, with \(\|v\|^2=-v^2\). Since
\(F=2C_u+\sum_vA_{u,v}\), the vertical lattice is \(\mathbb ZF\oplus R\).
Under the double cover of the Kummer surface, the pullbacks of
\(M_1,M_2\) are their strict graph transforms, each of multiplicity
one. Before blowing up, their mixed parts cancel and the translation
is torsion, so their sum in rational Picard is
\(2[E_1\times0]+2e[0\times E_2]\). Each transform subtracts the exceptional
divisor over each \((u,G(u))\) once. As \(A_{u,v}\) pulls back to
twice the exceptional divisor, descent gives
\[
 M_1+M_2=2O+eF+
                \sum_u\bigl(A_{u,0}-A_{u,G(u)}\bigr).
\]
Consequently the contribution at \(u\) to \(L-M_1-M_2-(n-e)F\) is
\[
  A_{u,G(u)}-\mathbf1_{u\notin J}A_{u,D(u)} .
\]
If \(c\) counts the mismatches at unselected branches, their total
over all branches is \((c+j)F+w\) with \(w\in R\) and
\(\|w\|^2=4c+2j\). All integral splittings are thus of the form
\[
 B_1=M_1+xF+v,\quad
 B_2=M_2+(n-e+c+j-x)F+w-v,\qquad v\in R .
\]
Set \(h_i=1+B_i^2/2\). Completing the square gives
\[
 n+j/2=e+h_1+h_2+\|v-w/2\|^2,\qquad
 B_1B_2=2n+j-2-h_1-h_2 .
\]
Conversely \(x=h_1+\|v\|^2/2\) recovers both classes. For
\(h_i\ge0\) their positive fibre degree gives effectivity by
Riemann--Roch. For \(h_i<0\) primitivity gives zero reduced virtual
class \cite[Proposition~3]{PandharipandeYinK3}.

To compute the scalar on the Chow graph line of
Proposition~\ref{fixed-evaluation}, normalize a real decomposable
mixed cohomology class \(v_0\) by \(v_0B_1=1\). Then \(v_0B_2=-1\),
\(v_0^2=Fv_0=F^2=0\) and \(FB_i=1\). The pairings of \(F,v_0\)
with \(B_1,B_2\) are invertible. The Gram matrix on their four-space
is therefore nondegenerate and determined by the two squares and
\(\nu=B_1B_2\). An isometry of such four-spaces extends to K3
cohomology. Use this isometry at both flags of every node. The
primitive descendant formula \cite[Corollary~3]{BuellesK3Classes}
then identifies the numerical vertex evaluations for the two sets of
data. For a positive-degree rational vertex with fewer than three
flags and markings, insert the primary divisor \(F\). Its divisor
equation has leading factor \(FB=1\), and each other term lowers a
positive cotangent exponent and cups the insertion with \(F\).
Induction on the number of missing marks and then on total cotangent
degree extends the comparison to these vertices. At a constant
vertex the isometry preserves the pairing and \(c_2(S)\).
Summing at fixed \(h_1+h_2=v'\) now defines \(C_b(v',\nu)\)
independently of the framing.

\smallskip\noindent\emph{The leading coefficient.}
For \(P_b(\nu)=C_b(0,\nu)\), both vertices have square \(-2\).
Compute each primitive root tensor on a smooth rational \((-2)\)-curve
using the preceding isometry. Its evaluation class is a linear
combination of products of point and fundamental classes in cohomology.
Thus \(\nu\) occurs only
when the two roots meet. A direct edge supplies at most one factor
of \(\nu\). A connected cluster of constant vertices evaluates at one
point and uses at least two edges for a single pairing of divisor
classes. A loop contributes the fixed root square. With at most \(d\)
edges, a term in \(\nu^d\) must therefore have exactly \(d\) parallel
edges, no other vertices and no descendants. The genus formula forces
both vertices to have genus zero. A stabilization correction would
introduce a bridge, which the parallel-edge graph lacks for \(d>1\).
For \(d=1\) the one-vertex term has zero squarefree coefficient.
The unmarked degree-one map to each smooth root
is rigid with reduced multiplicity one, and distinct markings sweep a
dense open of its power. Each evaluation slot thus carries the
fundamental class.

For the parallel edges put \(f_r(u)=u(r-u)\) on representatives
\(0,\ldots,r-1\) of \(\mathbb Z/r\). Up to Pixton's common factor
the flow sum for total leg weight \(x\) on the first root is
\(f_r^{*d}(x)=r^{1-d}\sum_{u_1+\cdots+u_d=x}\prod_i f_r(u_i)\).
The centred second difference of \(f_r\) is \(2r\mathbf1_0-2\).
For \(d\ge2\), convolution gives \(2f_r^{*(d-1)}\) plus a constant
as the second difference of \(f_r^{*d}\). Hence the flow sum on
\(0<x<r\) is a polynomial of degree \(2d\). Its leading coefficient satisfies
\((2d)(2d-1)c_d=2c_{d-1}\) and \(c_1=-1\), giving
\(c_d=-2^d/(2d)!\) for every \(r\) sufficiently large.
At the \(2d+1\) fixed values \(x=1,\dots,2d+1\), the normalized flow
sums are polynomial in \(r\) for all large \(r\)
\cite[Proposition~2.4]{SpelierDR}. Interpolation makes the coefficients
of the polynomial in \(x\) polynomial in \(r\) as well.
Reflection of all edge weights gives the symmetry \(x\mapsto r-x\).
Its constant term in \(r\) is therefore an even polynomial in \(x\),
valid for signed leg weights. If both marks \(+x_i,-x_i\) lie on the
same root, \(x_i\) drops out and the squarefree coefficient is zero.
Otherwise \(x=\sum_i\sigma_i x_i\), with \(\sigma_i=\pm1\) recording
which mark lies on the first root. Substitution gives the squarefree
coefficient \(-2^d\prod_i\sigma_i\). Exchanging the two marks also
reverses their mixed graph sign. Thus all \(2^{2d}\) surviving
assignments contribute with the same sign, even after division by
the edge symmetry factor \(d!\).
The disjoint strict graph sections of slopes \(1,-1\) on the resolved
surface cannot be joined in this root contribution. Hence \(P_b(0)=0\),
which gives the initial case \((n,j)=(1,0)\).

\smallskip\noindent\emph{Summing the labels.}
At an unselected branch the four \(w_u/2\) give \(D_4^\vee/D_4\).
Indeed, the three mismatches pair integrally with the simple roots and
represent distinct norm-one cosets. At a selected branch the
shift is half a root. Their theta series in \(q^{\|\cdot\|^2}\)
are respectively \(A\) and \(B\), since
\(D_4^\vee=\mathbb Z^4\cup(\mathbb Z+\tfrac12)^4\) and reflection
in the first coordinate gives for \((1,1,0,0)/2\) half the series
of \((\mathbb Z+\tfrac12)^2\times\mathbb Z^2\). Theta duplication
gives the stated constants.

Sum the norm equation over \(D\) and put \(t=2n+j\). The four
translations and interchange of the roots have multiplicities
independent of \(e\), and a slope divisible by \(l\) scales its
signed graph by the factor \(l^b\) occurring in \(H_b(e)\).
This gives \eqref{graph-family-system}. The minimum selected
norm is \(1/2\), so the equations for \(t<2+j\) or opposite
parity are empty. All the comparisons preceded the sum on the
labelled family. Each prescribed tuple \(D\) specifies its own total class \(L\)
and hence its own double ramification relation. Multiplying these
relations by \(\lambda(D)\) and summing, with the same weights at
each index, gives the asserted identity.
\end{proof}

\subsection{Comparison on the two-torsion cover}
\begin{proposition}\label{framed-comparison-solutions}
Under the framed identities of Proposition~\ref{family-unmarked-rows},
let \(W_b\) be the solution space of \eqref{rd-m-rec}.
For every \(d=b/2\ge1\), its complete space of solutions and
its leading polynomial are
\[
 W_b=\Delta_{\rm mod} R_{d-1},\qquad
 P_b(X)=c_b\prod_{i=0}^{d-1}(X-i),\quad c_b\ne0 .
\]
\end{proposition}
\begin{proof}
\smallskip\noindent\emph{Two labelled branches.}
Let \(M\) act on two-torsion and fix a nonzero target covector \(v\).
At an unselected branch \(u\), weighting \(D(u)\) by \((-1)^{v(D(u))}\)
gives \((-1)^{v(Mu+t)}C\). Indeed the three nonzero cosets have equal
theta series by triality and their signs sum to \(-1\), leaving the
difference \(C=\vartheta_4(2\tau)^4\) of the even- and odd-sum integer
coset series. For two distinct unselected branches \(u_1,u_2\), put
\(w=u_1+u_2\). Weighting both cancels the translation \(t\) and gives
\(C^2 A^{2-j}B^j\) when \(j\) of the remaining branches are selected,
with sign \((-1)^{v(Mw)}=2\mathbf1_{M\langle w\rangle\subset\ker v}-1\).
Retain the
source and target lines and take the induced correspondence on
parabolic cohomology. Base point pairs act by zero in degree one,
as is seen by pairing with compactly supported lifts. For a target
cusp form \(f\) of trace zero to level one, summing over all target
lines gives zero. This removes both the unrestricted term and the
kernels containing the source line. Every other kernel has a unique
compatible target line. Up to a common nonzero factor the operator
therefore sums the isogenies with kernel disjoint from the source
line. At degree \(e\) the upper triangular representatives have
diagonal \(a,d'>0\), \(ad'=e\), \(a\) odd, and upper entry
\(0\le r<d'\). The base differential and the \(b\) elliptic
differentials supply the factor \(a^{b+1}/d'\). In the resulting
sum of translates \(f((a\tau+r)/d')\), extraction at \(q^1\)
leaves only \(a=1\) and gives \([q^e]f\).

On level-one pullbacks the same extraction for unsigned graphs
has unrestricted \(a\). The identity \(C^2=A^2-16B^2\) expresses
the three multipliers above as linear combinations of the five
unsigned multipliers. Pullback splits the trace after division by
three, so for all \(f\in S_{2d+2}(\Gamma_0(2))\) we have
\(C^2 A^{2-j}B^j f\in W_b\), \(0\le j\le2\).

On \(X(2)\) translation fixes \(A\) and negates \(B\).
The zero of \(B\) is the cusp at infinity, and the two factors of
\(C^2\) vanish at the other cusps. In the even-weight ring
\(\mathbb C[A,B]\) this gives
\[
 S_{2d+2}(\Gamma_0(2))=B^2 C^2\mathbb C[A,B^2]_{d-3},
 \qquad \Delta_{\rm mod}\in\mathbb C^\times B^2 C^4.
\]
For \(d\ge3\) it follows that \(\Delta_{\rm mod}A^{d-1-i}B^i\in W_b\) for
\(0\le i\le2\lfloor(d-3)/2\rfloor+2\).
When \(d\ge3\) these monomials give all of
\(\Delta_{\rm mod}R_{d-1}\) for odd \(d\). For even \(d\ge4\),
only the monomial \(\Delta_{\rm mod}B^{d-1}\) is missing.

\smallskip\noindent\emph{The sign of a two-torsion frame.}
Keep the source frame fixed. At three unselected branches choose
distinct nonzero target covectors \(v_i\). Translate their source
labels to \(0,e_1,e_2\) and pair them respectively with
\(e_1^\vee+e_2^\vee,e_1^\vee,e_2^\vee\). The covectors sum to
zero, so the translation \(t\) drops out. The product of the three
character weights contributes \(C^3(-1)^{\operatorname{tr}M}\).
Let the remaining branch be selected or unselected and unweighted,
giving a further factor \(B\) or \(A\).

Take a full-level form \(g\) of weight \(2d-4\) and put
\(f=\eta^{12}g\), where \(\eta\) is the Dedekind eta function.
The transformation laws of \(\eta^{12}\) on \(S,T\) give the sign
character of \(\mathrm{GL}_2(\mathbb F_2)\simeq S_3\).
Thus \(f\) is a cusp form of weight \(b+2\) on \(\Gamma(2)\) and
\(f=\sum_{m\ {\rm odd}}c(m)q^{m/2}\). Summing over the six
target frames weights an isogeny by
\[
 \sum_{D\in\mathrm{GL}_2(\mathbb F_2)}\operatorname{sgn}(D)
                  (-1)^{\operatorname{tr}(DM)}.
\]
For singular \(M\) an odd permutation fixes its image pointwise,
so the sum vanishes. For invertible \(M\) it equals
\(-4\operatorname{sgn}(M)\), by summing at the identity and
changing variables. Only odd isogeny degrees \(e\) remain.

Write \(ad'=e\) and \(0\le r<d'\) as in the two-branch trace.
In the two frames the isogeny has diagonal entries \(a,d'\) and
off-diagonal entry \(-r\), so at odd \(e\) its sign is \((-1)^r\).
For odd \(m\) the sum
\(\sum_{r=0}^{d'-1}(-1)^r\exp(\pi i m r/d')\) is \(d'\) if
\(d'\mid m\) and zero otherwise. The coefficient at \(q^{1/2}\)
then requires \(a=1\), \(m=d'=e\), and the factor \(1/d'\)
cancels. Thus the series in the isogeny degree is \(f(2\tau)\)
up to one nonzero constant. Both \(C^3A f(2\tau)\) and
\(C^3B f(2\tau)\) belong to \(W_b\).

The product formula \(C=\vartheta_4(2\tau)^4
=\eta(\tau)^8/\eta(2\tau)^4\) gives
\(C^3\eta(2\tau)^{12}=\Delta_{\rm mod}\). For \(d=2\), take
\(g=1\) to get both \(\Delta_{\rm mod}A\) and \(\Delta_{\rm mod}B\).
For even \(d\ge4\), take \(g=E_4^{(d-2)/2}\).
At \(\tau_0=(1+i)/2\) the weight-two form \(A\) vanishes: the
matrix \(\left(\begin{smallmatrix}1&-1\\2&-1\end{smallmatrix}\right)\)
fixes \(\tau_0\) and its weight-two factor there is \(-1\).
But \(B(\tau_0)\ne0\) and \(E_4(2\tau_0)=E_4(i)>0\).
Thus \(\Delta_{\rm mod}B\,g(2\tau)\) is independent of the
\(d-1\) monomials already obtained, all divisible by \(A\).
Since \(g(2\tau)\in R_{d-2}\), this supplies the missing direction.
The bound \(\dim W_b\le d\) proves equality for \(d\ge2\).

\smallskip\noindent\emph{The genus-zero row.}
For \(d=1\) only one edge joins the two nonconstant vertices in the
squarefree relation. Each vertex has two primary divisor insertions,
and dimension forces its node insertion to be a divisor as well.
The divisor equation gives the product of the three divisor
pairings with its primitive class times the reduced genus-zero
invariant \(N_h\) for square \(2h-2\). Contraction at the node
gives \(B_1B_2=\nu\); the mark-distribution factor is the same for
all \(h_1,h_2\). Thus
\[
 C_2(v,\nu)=c_2\nu\,[q^v]G(q)^2,\qquad
 G(q)=\sum_{h\ge0}N_hq^h=q/\Delta_{\rm mod}.
\]
Here the last identity is the primitive Yau--Zaslow formula
\cite[\S1.4, (7)]{PandharipandeYinK3}. With \(D=q\,d/dq\),
the recurrence is \(G^2(2D-2)y+D(G^2)y=0\), or
\((D-1)(Gy)=0\). Hence \(W_2=\mathbb C\Delta_{\rm mod}\).

Finally \(\Delta_{\rm mod}A^{d-1-i}B^i\) starts with
\(2^i q^{(2+i)/2}\). At \(t=2+i\) the recurrence gives
\(P_b(i)=0\) for \(0\le i<d\). Proposition~\ref{family-unmarked-rows}
gives its degree and nonzero leading coefficient.
\end{proof}

\begin{corollary}[The untwisted Chow series]\label{all-untwisted-modularity}
For every cyclic level \(N\), every \(\kappa\ge2\) and every
\(r\) with \(r^2\equiv1\bmod 4N\), one has \(D_{1,r}(N)=0\) in the ordinary compact
Chow group. The untwisted generating series belongs to \(S\otimes A_N\).
\end{corollary}
\begin{proof}
Let \(f\) be the series of a functional on the graphs modulo
\(J_b\). By Propositions~\ref{family-unmarked-rows}
and~\ref{framed-comparison-solutions}, both \(A^4f\) and \(B^4f\)
belong to \(\Delta_{\rm mod}R_{d-1}\). The two weight-two forms
\(A,B\) have no common zero on \(X(2)\): \(B\) vanishes only at
infinity, where \(A=1\). Dividing on these two opens gives a cusp
form \(f\) of weight \(b+2\). Its integral exponents give level
\(\Gamma_0(2)\), and the \(T_2-U_2\) argument in
Proposition~\ref{med-u-five-prop} gives full level.
Functionals separate points, so every linear relation among the Fourier
coefficient functionals on that cusp space also holds among the
\(H_b(n)\) modulo \(J_b\). The restriction in
Theorem~\ref{graph-to-moments} kills \(J_b\) and carries these relations
to the \(F_i(n)\). Theorem~\ref{med-u-scalar} and
Lemma~\ref{scholl-open-chow} give \(D_N=0\) in compact Chow.
Multiplication by any compatible root permutes the residues isometrically.
\end{proof}

The graph calculation also identifies all the relations among the graphs
before restriction. Let \(\mathcal G_b\) be the span of \(H_b(n)\)
modulo \(J_b\). The preceding proof gives a surjection
\(S_{b+2}(\mathrm{SL}_2(\mathbb Z))^\vee\to\mathcal G_b\) carrying
the \(n\)-th Fourier evaluation to \(H_b(n)\). It is an isomorphism.
Indeed, parabolic realization carries the signed graph to \(T_n\)
and kills point-pair errors (Section~\ref{rd-m-twist}). If a finite
combination \(\sum_n c_n H_b(n)\) is zero modulo \(J_b\),
\(\sum_n c_n T_n\) acts by zero on the holomorphic cusp space.
Taking its first coefficient gives \(\sum_n c_n[q^n]g=0\) for
every cusp form \(g\). This proves injectivity, and hence
\(V_b=S_{b+2}(\mathrm{SL}_2(\mathbb Z))\). Moreover, composition
of graphs satisfies the prime recurrence of
Lemma~\ref{full-graph-hecke} and multiplies coprime degrees.
Thus \(H_b(n)\mapsto T_n\) identifies \(\mathcal G_b\), as an
algebra with identity \(H_b(1)\), with the cuspidal Hecke algebra.
In particular, the simultaneous eigenprojectors are represented
modulo point pairs by linear combinations of graphs of degree
\(n\le\lfloor(b+2)/12\rfloor\), since these Fourier evaluations
span the dual by the full-level valence formula.

\begin{proposition}[A weaker graph relation for codimension two]
\label{restricted-divisor-graph-errors}
Work on a fine cover of the smooth product of modular curves and
take \(a=1\), so the output is a codimension-two cycle on the
Kuga--Sato threefold. In degree \(b\), enlarge \(J_b\) by the
following classes. Take a normalized isogeny curve \(C\) in the
base product, allowing a finite cover to record the isogeny. On
the product of the two elliptic \(b\)-fold powers over \(C\), take
a codimension-\(b\) class in the \(e_b\times e_b\) summand whose
fibre Betti class is zero on a dense open subset of \(C\). Include
its proper pushforward to the product over the two modular curves.
Denote the resulting complex subspace by \(J_b^{\,\mathrm{fib}}\).

Every cyclic restriction \(R_i\) of Theorem~\ref{graph-to-moments}
with \(a=1\) kills \(J_b^{\,\mathrm{fib}}\) in ordinary Chow.
This subspace is stable under the target Hecke operators. Thus for
\(\kappa=2\), the coefficient relations of weight \(b+2\) need hold
only modulo \(J_b^{\,\mathrm{fib}}\), at every cyclic level. The
framed recurrence and scalar argument of Section~\ref{med-u-five}
also suffice with this larger error space on each framed cover.
\end{proposition}
\begin{proof}
The first base change in \(R_i\) is finite flat. After this base
change a smooth curve again represents the branch over \(C\) and maps finitely
to the source curve. One may normalize each component. Modifications
at finitely many base points give classes already in \(J_b\).
On a smooth isogeny branch the Betti class along the fibres is a
section of the relative cohomology local system. Indeed it is the
restriction of the cycle class on the total space. Its vanishing on
a dense open implies its vanishing on every fibre of that branch.

The next operation is refined intersection with the graph of the
forgetful map on the bases. If a branch lies in this graph, the
excess factor is \(c_1\) of the pullback of the graph's normal line.
This line comes from \(T_{\mathcal M_{1,1}}\) and has zero rational
first Chern class on the smooth base. On every other branch the graph
cuts out a finite Cartier divisor on the smooth normalization.
Functoriality of refined pullback expresses its contribution as
the sum of the fibre restrictions with their intersection
multiplicities. Each has zero Betti class. The remaining \(b-1\)
Poincar\'e contractions and pushforward therefore give a divisor of
zero Betti class on \(E^2\). It lies in the part odd in both
coordinates. Removing the restrictions of a line bundle to the two
axes identifies this part of \(\Pic(E^2)_{\mathbb Q}\) with
\(\operatorname{Hom}(E,E^\vee)_{\mathbb Q}\) by the Poincar\'e
bundle. A homomorphism of elliptic curves acts faithfully on \(H^1\).
The contracted divisor is consequently zero in the Chow group of
the fibre. Summing the pushforwards proves the first assertion,
also after finite descent.

Postcomposition by a degree-\(p\) isogeny pulls the presentation
back to the finite cover of the isogeny curve recording the
\(p\)-subgroup, then pushes it forward along isogenies. After
normalization the generic fibre class remains zero by Betti
functoriality, and point-pair classes stay in \(J_b\). This proves
Hecke stability. For the framed comparison by parabolic cohomology,
the sign projector on the \(2b\) elliptic coordinates leaves relative
cohomology only in degree \(2b\). Hence on \(C\) the projected
Leray spectral sequence identifies the projected \(H^{2b}\) of
the total space with \(H^0(C,R^{2b})\) in that summand. The added
classes consequently contribute zero to this comparison as well.
Lemma~\ref{full-graph-hecke} and
Propositions~\ref{med-u-five-prop} and~\ref{framed-comparison-solutions}
therefore apply to the graph series in this quotient.
Theorem~\ref{graph-to-moments} gives the claimed cyclic relations.
\end{proof}

In codimension two one can combine the labelled family class of
Proposition~\ref{family-supported-diagonal} with the numerical vertex
pairings from Proposition~\ref{family-unmarked-rows}. Modulo
\(J_b^{\,\mathrm{fib}}\) the fixed-surface evaluation can be
computed directly from its Betti class.
Indeed, along a normalized isogeny curve the chosen isogeny
identifies the two rational \(H^1\) local systems. The signed
fibre tensor then gives an endomorphism of
\(\operatorname{Sym}^b H^1\) commuting with finite-index modular
monodromy. This group is Zariski dense in \(\operatorname{SL}_2\);
the symmetric power is irreducible. Its invariant endomorphisms
are scalars, and the projected graph gives a nonzero generator.
The numerical pairing against that generator therefore determines
the fibre class. The difference has a supported Chow presentation
with zero fibre Betti class, as required. For \(\kappa=2\)
Proposition~\ref{restricted-divisor-graph-errors} annihilates
each such class after cyclic restriction at every level.

For any \(a=\kappa-1\ge1\) one can identify what these same
comparisons give in the smooth total space. For a smooth CM point
\(x\) write \(i_x\) for its fibre inclusion. Put
\[
 \mathcal R_{a,N}=
 \sum_{x\ {\rm CM}}(i_x)_*\!
      \left(e_{2a}\CH^a(E_x^{2a})_{\rm hom}\right),
\]
using complex coefficients and descent from a fine cover.
The same supported decomposition and numerical comparisons imply
that every untwisted weakly holomorphic input satisfies
\[
               Z_\kappa(f)|_{Y_0(N)}\in\mathcal R_{a,N}.
\]
Indeed the framed recurrences hold modulo \(J_b^{\,\mathrm{fib}}\)
in every degree \(b\), by the monodromy argument. The proof of
Proposition~\ref{restricted-divisor-graph-errors} through the
refined intersection applies verbatim for \(b-a\) contractions.
A branch that meets the cyclic restriction properly gives a
nonscalar endomorphism, so its intersection points are CM.
Each remaining fibre class still has zero Betti class after the
contractions and projection. Thus \(R_i(J_b^{\,\mathrm{fib}})
\subset\mathcal R_{a,N}\) for every row. Proposition~\ref{diagonal-hecke-defect}
gives the assertion for all principal parts. When \(a=1\),
\(\mathcal R_{1,N}=0\) by the divisor argument above. In higher
codimension the graph-line comparison gives a way to remove this
subspace from the modularity argument.

\section{Twisted cycles and the kernel of transfer}\label{rd-m-twist}

Transfer expresses each twisted relation, after applying a product of
Hecke operators, in terms of the untwisted relations proved in
Corollary~\ref{all-untwisted-modularity}. It remains to show that the
product has trivial kernel on the relevant cycles in Chow.
Qiu's unramified Hecke decomposition applies to CM cycles modulo the
radical of the height pairing, away from a set containing \(2\), the
level primes and the primes ramified in the CM field, with at least
three elements \cite[Theorem~1.2.1]{QiuCM}. Here we need injectivity
before passing to the height quotient.
Write \(d=\kappa-1\geq1\), \(m=2d\), \(K_m=m+2\), and \(t_p=p^d\).
Let \(D_{\Delta,r}(N)\) denote the image of the map taking a weakly
holomorphic principal part to its cycle in compact Chow. Keep the root
\(r\) fixed. The operator in the transfer is
\begin{equation}\label{md-t-product}
 L_{\Delta,N}=
 \prod_{\substack{p\mid\Delta\\p\nmid N}}
      \bigl(T_p^2-(p+1)^2t_p^2\bigr)
 \prod_{p\mid(\Delta,N)}(B_p-2pt_p).
\end{equation}
Here \(T_p\) is the geometric Hecke operator at a good prime. At a
shared prime \(B_p=U_p+U_p^*\) is the sum of the two oriented
correspondences. We will prove the transfer inclusion
\begin{equation}\label{md-t-inclusion}
 L_{\Delta,N}D_{\Delta,r}(N)\ \subset\
       \sum_i C_i D_{1,1}(M_i),
\end{equation}
where \(C_i\) are correspondences from finitely many auxiliary cyclic
levels \(M_i\) to \(N\), independently of the principal part.

\subsection{Transfer through finitely many cyclic levels}

Write \(x=[A,B,C]\), \(D(x)=B^2-4AC\), and \(c=|\Delta|\).
The action \(x|g=(\det g)^{-1}x\circ g\) preserves \(D\).
In these coordinates the exponent in \eqref{rd-coeff} is
\(q_*(x)=-D(x)/(4Nc)\). Let \(v_\mu\) be the finite test consisting
of the coset \(r\mu+L\), the congruence
\(q_*(x)\equiv\operatorname{sgn}(\Delta)Q(\mu)\pmod{\mathbb Z}\),
and the factor \(\chi_\Delta(x)\). These tests are finite linear
combinations of lattice cosets. Their real matrix of coefficients is the twist
intertwiner of \cite[Definition~3.1 and Proposition~3.2]{AETraces}.
We track its complex conjugate below when passing to the finite
Weil action with inverse additive character.

At \(p\) let \(S_p\) be the locally constant compactly supported
functions on rational binary forms and \(G_p=\operatorname{PGL}_2(\mathbb Q_p)\).
The ordinary cyclic coset of depth \(t\ge0\) has
\[
 v_p(A)\ge t,\qquad v_p(C)\ge0,\qquad
 B\equiv b\pmod{p^{t+\delta}},\qquad \delta=v_p(2).
\]
Let \(E_p\) be the span of these cosets after changes of frame by \(G_p\),
scalar dilations, and multiplication by additive characters of \(D\).
For \(H=\{\operatorname{diag}(u,1):u\in\mathbb Z_p^\times\}\) we claim
that \(E_p=G_p\cdot S_p^H\), where the right side denotes the linear span.
The cosets, their dilations and their discriminant multipliers are
\(H\)-invariant. Conversely, after dilation an invariant test depends
on integral triples modulo \(p^k\). If the valuations \(i,j\) of
\(A,C\) are both less than \(k\), its unit orbits are distinguished by
these valuations, \(B\bmod p^k\), and
\(p^{-i-j}AC\bmod p^{\min(k-i,k-j)}\). Indeed the two congruences on
the scaling unit are compatible exactly when these products agree.
If an outer residue is zero only the truncated valuations and \(B\)
are needed. For lower bounds \(l,m\) on the outer valuations and
large middle precision \(n\), put
\[
 s=n-\delta-l-m,\qquad j=m+s,\qquad t=l+m+2s .
\]
Applying the cyclic coset to
\((p^{s+j}A,p^s B,p^{s-j}C)\) gives the desired bounds and precision.
Differences impose exact valuations. Refine the middle precision so
that \(B^2\) is fixed modulo \(p^{i+j+\min(k-i,k-j)+v_p(4)}\);
finite Fourier inversion in \(D=B^2-4AC\) then imposes the product
condition, also at two. This proves the claim. Compact averaging of
a lift now gives \((S_p/E_p)^H=0\).

\begin{lemma}\label{md-t-unit}
Let \(H=\{\operatorname{diag}(u,1):u\in\mathbb Z_p^\times\}\).
A smooth complex representation \(V\) of
\(\operatorname{PGL}_2(\mathbb Q_p)\) with \(V^H=0\) factors through
\(\mathbb Q_p^\times/\mathbb Q_p^{\times2}\) by the determinant.
\end{lemma}
\begin{proof}
Let \(n(x)\) be upper unipotent. If its action is nontrivial,
a compact subgroup \(n(p^{-R}\mathbb Z_p)\) has a nonzero eigenvector
\(v\) for a nontrivial character. Indeed each vector generates a
finite-dimensional representation of this compact abelian subgroup,
and this representation is a sum of characters. Choose \(c\) so that
\(1+p^c\mathbb Z_p\) fixes \(v\) under \(H\). Using an additive
character of conductor \(\mathbb Z_p\), the frequency of \(v\) has
valuation \(a<R\). Decompose \(v\) into character components for
\(n(p^{-R'}\mathbb Z_p)\), where \(R'\ge R\) and \(R'-a\ge c\).
Every nonzero component has frequency of the same valuation \(a\),
since its restriction to the smaller subgroup is prescribed. Choose
one such component \(w\). Its character stabilizer in \(H\) consists
of the units congruent to one modulo \(p^{R'-a}\). This subgroup
fixes \(v\) and the character projection producing \(w\). The sum of
the translates of \(w\) over the cosets of the stabilizer is
\(H\)-invariant. It is nonzero since the translates are in distinct
character spaces. This contradicts \(V^H=0\). The upper unipotent
therefore acts trivially. The lower unipotent is conjugate to it.
Together they generate the image of \(\operatorname{SL}_2(\mathbb Q_p)\),
the kernel of the determinant square-class map.
\end{proof}

\begin{proposition}[Transfer of weakly holomorphic relations]\label{md-transfer-proof}
In compact Chow the inclusion \eqref{md-t-inclusion} holds for every
compatible fundamental discriminant and root. The levels and
correspondences are independent of the input principal part. Each
coefficient of the twisted series is fixed by normalized edge
reversal \(W_p\) at every \(p\mid(\Delta,N)\).
\end{proposition}
\begin{proof}
\emph{Trace on a cycle coefficient.}
For a finite \(\Gamma\)-invariant test \(v\), let \(F_q(v)\) sum the
positive-norm orbits with exponent \(q(x)\), weight
\(2v(x)/|\Gamma_x|\), and fibre class \(Y_d(x)\), the projected
product of \(d\) mixed graph divisors. Since \(n=-D(x)/(4Nc)\)
and the normalized fibre class divides this product by a fixed
multiple of \(|D(x)|^{d/2}\), \(F_{q_*}(v_\mu)\) agrees with
\eqref{rd-coeff} up to one nonzero scalar depending on level and twist.
An isogeny \(\phi\) of degree \(e\) transports \(u\) to
\(u'=\phi u\phi^{-1}\). Pullback of the mixed divisor of \(u'\)
corresponds to \(\widehat\phi u'\phi=eu\), so the projection formula
on \(E^2\) gives pushforward factor \(e\) on one divisor and \(e^d\)
on the product.

For a positive-determinant rational frame \(g\) from \(\Gamma'\) to
\(\Gamma\), put \(J=\Gamma\cap g\Gamma'g^{-1}\). The trace sums the
changed test over \(\Gamma/J\). The \(\Gamma_x\)-orbit of \(\gamma J\)
has size \([\Gamma_x:\Gamma_x\cap\gamma J\gamma^{-1}]\), so its weight
is \(2/|\Gamma_x\cap\gamma J\gamma^{-1}|\), the intermediate-level
weight. Isogenies preserve cusps, so the same trace holds in compact Chow.

\smallskip\noindent\emph{Local Hecke factors.}
The full test \(v_\mu\) is fixed by the product of the groups
\[
 K_\ell=\{k\in\operatorname{SL}_2(\mathbb Z_\ell):
                         k_{21}\in N\mathbb Z_\ell\}.
\]
Indeed its stabilizer is open and contains \(\Gamma_0(N)\), which
is dense by surjectivity of reduction of \(\mathrm{SL}_2(\mathbb Z)\)
modulo every multiple of \(N\).
For \(p\mid c,\ p\nmid N\), choose the integral branches of \(T_p\)
with determinant \(p\) and lower left entry divisible by \(N\).
Each of the \((p+1)^2\) double branches has a matrix \(g\) of
determinant \(p^2\). Its projective representative \(g/p\) lies in
\(K_\ell\) for every \(\ell\ne p\). At \(p\) its determinant class
is trivial, so it acts as the identity modulo \(E_p\) by
Lemma~\ref{md-t-unit}. Including the graph factor \(t_p^2\) gives
\((T_p^2-(p+1)^2t_p^2)v_\mu\in E_p\otimes\bigotimes'_{\ell\ne p}S_\ell\).

For \(p\mid(c,N)\) one has \(p\Vert N\). Denote by \(P\) pullback
along the map forgetting the \(p\)-line, by \(P^*\) its trace, and
by \(W_p\) normalized edge reversal. The \(p+1\) branches of \(PP^*\)
change that line by integral determinant-one matrices with trivial
transport away from \(p\), hence \(PP^*\equiv(p+1)I\bmod E_p\).
Changing to one of the other \(p\) lines and reversing, in either
order, gives
\[
 B_p=t_p\{W_p(PP^*-I)+(PP^*-I)W_p\}
       \equiv 2pt_pW_p\pmod{E_p}.
\]
The exact Atkin--Lehner involution preserves the genus character,
including its zero weights \cite[\S4]{BO}, and takes the middle
component to \(r\sigma_p\mu=r\mu\) (Corollary~\ref{med-antisector}).
This identity on rational vectors holds on the finite adelic test
by density and local constancy. Thus \(W_p\) fixes the full \(v_\mu\)
before summation, including at \(p=2\). Since \(W_p\) normalizes
\(E_p\), the shared-prime factor may be applied as \(B_p-2pt_pW_p\).

The factors at distinct primes commute, since quotients by
coprime subgroups give the same isogeny and bijective choices of
flags in either order. Apply each exceptional factor first in turn.
Every subsequent branch preserves \(E_p\) by \(G_p\)-stability.
At \(\ell\nmid c\) the component fixes \(x\bmod cL_\ell=L_\ell\),
and the genus character depends on \(x\bmod\Delta L\)
\cite[\S3]{AETraces}. Diagonal units preserve the middle residue,
outer valuation bounds and discriminant, so the local test already
lies in the unit-invariant subspace. Almost all are standard lattice
indicators. Choosing complements to the \(E_p\)'s at the finitely
many exceptional primes shows that the intersection of the
single-place conditions is their tensor product.

\smallskip\noindent\emph{Realizing the local tests at cyclic levels.}
Take \(S\) to contain \(2\), the primes of \(Nc\), and every place
with a nonstandard summand. For one product of local boxes choose
frame lifts \(g_p\in\mathrm{GL}_2(\mathbb Q_p)\). Set
\(a=\prod_{p\in S}p^{v_p(\det g_p)}\), \(A=\operatorname{diag}(a,1)\)
and \(h_p=\operatorname{diag}(\det g_p/a,1)\). Then \(h_p\) fixes
the box and \(g_ph_p^{-1}A^{-1}\) has determinant one. Strong
approximation gives \(s\in\operatorname{SL}_2(\mathbb Z[S^{-1}])\)
close enough to these matrices that \(sA\) acts like \(g_p\) on
each box. This is equality of test functions: the stabilizer of a
compactly supported locally constant function is open. Outside \(S\)
the frame preserves the standard lattice.
One positive rational dilation likewise realizes all prescribed
valuations, absorbing units into the middle residues. If the box
at \(p\) has depth \(s_p\), Chinese remaindering gives a single
residue modulo \(2M\), where \(M=\prod_p p^{s_p}\).

Use the finite additive character \(\psi_f(x)=e^{-2\pi i x}\)
on \(\mathbb Q\). Compact support bounds \(v_p(D)\) below, so the
discriminant multipliers require only finite precision at \(p\).
Approximation in \(\mathbb Z[S^{-1}]\) realizes them as
\(\psi_f(bD)\) for one rational \(b\). We have therefore expressed
the transferred test as a finite sum of rational frame changes,
dilations and additive norm multipliers of ordinary component tests.
Averaging each term over the output group gives, by the trace
calculation, a correspondence \(C_g\) from cyclic level \(M\) on
geometric Chow groups over \(\overline{\mathbb Q}\), up to a fixed
scalar. The orthogonal frame action commutes with the
Weil action, so \(C_g\) applies to every component and every finite
Weil translate with the same source invariance group. Modulo the
sum of the \(C_gD_{1,1}(M)\), Theorem~\ref{rd-main} makes these
ordinary component series after \(C_g\) components of cusp forms.

\smallskip\noindent\emph{The Weil transformation law.}
Write \(k=d+3/2\). For an upper triangular rational metaplectic
element \(u\) with positive diagonal and compatible real lift,
reindexing gives
\[
       F_q(\omega_q(u)v)=F_q(v)|_k u^{-1}.
\]
Indeed \(n(b)\) gives the multiplier \(e^{-2\pi i bq(x)}\) on a
rational vector. For \(\operatorname{diag}(a,a^{-1})\), changing
\(x\) to \(ax\) gives \(a^{-d}\) from the graph and
\(|a|_f^{3/2}=a^{-3/2}\) from the measure. The finite Weil-index
ratios have product one since \(a>0\) preserves the real signature.
For \(q_M=-D/(4M)\), the ordinary Weil law gives the same covariance
for integral metaplectic elements in the Chow quotient. Every
rational element factors as \(u\gamma\) by choosing the primitive bottom
row of \(\gamma\in\operatorname{SL}_2(\mathbb Z)\) in the given
direction with sign making \(u\) positive on the diagonal. A
central sign is absorbed by the integral lift. For a translated
component \(\omega_q(a)v\), apply this factorization to \(\gamma a\).
The right slash action gives
\(F_q(v)|_k(\gamma a)^{-1}=(F_q(v)|_k a^{-1})|_k\gamma^{-1}\),
so covariance holds on the rational Weil span.

Set \(\lambda=M/(Nc)>0\), so \(q_*=\lambda q_M\) and
\(F_{q_*}(v)(\tau)=F_{q_M}(v)(\lambda\tau)\). Conjugation by the
rational frame \(\operatorname{diag}(\lambda,1)\) takes \(n(b)\)
to \(n(\lambda b)\) and \(S\) to
\(\operatorname{diag}(\lambda,\lambda^{-1})S\). On the finite side,
scaling the Fourier pairing by \(\lambda\) scales the self-dual
measure by \(|\lambda|_f^{3/2}\), and the product formula identifies
the Weil indices because the real signatures agree. Thus covariance
holds also for \(q_*\). Rational translates and positive rational
rescalings of cusp forms vanish at every rational cusp.

Finally let \(v'_\mu\) denote the transferred full test. The
intertwiner of \cite[Proposition~3.2]{AETraces} has real
coefficients and source type \(\widetilde\rho_L=\rho_L\) for
\(\Delta>0\), \(\widetilde\rho_L=\overline{\rho_L}\) for \(\Delta<0\).
Our inverse additive character conjugates its transformation law.
Since transfer uses orthogonal frame changes, for every integral
lift \(\gamma\) this gives
\[
 \omega_{q_*}(\gamma)v'_\mu
       =\sum_\nu\overline{\widetilde\rho_L(\gamma)}_{\nu\mu}\,v'_\nu.
\]
Covariance with \(\gamma^{-1}\) and unitarity give
\[
 F_{q_*}(v'_\mu)|_k\gamma
       =\sum_\nu\widetilde\rho_L(\gamma)_{\mu\nu}\,F_{q_*}(v'_\nu).
\]
The transferred column is thus a cusp form of the required type
in the output quotient. The principal-part criterion gives
\eqref{md-t-inclusion}, with levels and correspondences determined
by the finite tests alone.
\end{proof}

\subsection{Excluding the exceptional eigenvalues}

To recover the cycle from its transfer, we use a spectral criterion
for modules of arbitrary dimension. It suffices that a power of the
error ideal annihilate the module.

\begin{proposition}[Exclusion of the exceptional eigenvalues]\label{md-t-spectrum}
Let \(\mathcal A\) be a unital complex algebra acting on \(V\) and
\(J\) a two-sided ideal with \(J^s V=0\) for some \(s\).
Fix an integer \(d\ge1\) and a cyclic level \(h\ge1\). For \((n,h)=1\) assume there are
operators \(T_n\in\mathcal A\) with \(T_1=1\) and
\[
 T_p T_n=T_{pn}+p^{\,K_m-1}T_{n/p}\qquad(p\nmid h),
\]
where a nonintegral subscript gives zero. Suppose for one \(M>1\),
divisible by \(h\), the scalar series
\begin{equation}\label{md-t-series}
 F_\ell(\tau)=\sum_{\substack{n\ge1\\(n,h)=1}}
               \ell(T_n\bmod J)q^n,\qquad \ell\in(\mathcal A/J)^\vee,
\end{equation}
are modular ratios of weight \(K_m=2d+2\) for \(\Gamma_1(M)\) and span a
finite-dimensional space. Put \(M_h=\prod_{p\mid h}p^{v_p(M)}\).
Then all these series lie in \(M_{K_m}(\Gamma_0(M_h))\). For every
\(p\nmid h\), let \(P_p\) be the monic polynomial whose simple roots
are the distinct classical \(T_p\)-eigenvalues on that space. Then
\(P_p(T_p)\in J\) and \(P_p(T_p)^sV=0\). In particular, both
\(T_p\pm(p+1)p^d\) have polynomial inverses on \(V\).
Every eigenvalue on \(V\) lies in
\([-2p^{d+1/2},2p^{d+1/2}]\) or has the form
\(\zeta p^{2d+1}+\zeta^{-1}\) for a root of unity \(\zeta\).
\end{proposition}
\begin{proof}
Let \(U_p\) extract the coefficients at multiples of \(p\), put
\(V_pf(\tau)=f(p\tau)\), and write \(T_p^0=U_p+p^{K_m-1}V_p\).
For \(p\nmid h\) the recursion gives \(T_p^0 F_\ell=F_{\ell'}\),
where \(\ell'(a)=\ell(T_pa)\). Both sides have zero coefficient when
\((n,h)>1\). Choose a prime \(q\equiv1\bmod M\). Its diamond operator
is trivial, so the finite-dimensional span of the \(F_\ell\) is stable
under the classical \(T_q=T_q^0\). Schaeffer's theorem makes this
span holomorphic at level \(\Gamma_1(M)\), including at the cusps
\cite[Theorem~1.1(a)]{SchaefferStability}.

The same recursion determines the level. At weight \(K_m>2\) a basis
at level \(\Gamma_1(L)\) consists of \(g(t\tau)\) with
\(t\,\operatorname{cond}(g)\mid L\). Here \(g\) is a primitive cusp
eigenform or an Eisenstein series attached to two primitive characters,
whose conductors have product \(\operatorname{cond}(g)\)
\cite[Theorems~5.8--5.9 and \S9.2]{SteinModularForms}.
In the common basis at level \(pL\), \(V_p\) replaces \(t\) by \(pt\).
If \(p\nmid L\), the spans of the original basis and its dilation
have zero intersection. If \(p\mid L\) and \(V_pf\) still has level \(L\),
each term of \(f\) satisfies \(t\,\operatorname{cond}(g)\mid L/p\).

Apply these observations to \(f=F_\ell\). For \(p\nmid M\), the
forms \(T_p^0f\) and \(T_pf\) both have level \(M\), so their difference
\(p^{K_m-1}V_p(f-\langle p\rangle f)\) belongs to both subspaces.
Thus \(\langle p\rangle f=f\). These primes generate the diamond
operators, giving level \(\Gamma_0(M)\). Now suppose the entire
family has level \(\Gamma_0(L)\). If \(p\mid L\) and \(p\nmid h\),
\(U_pf\) has level \(L\), as does \(T_p^0f=F_{\ell'}\). Their
difference gives \(V_pf\) at level \(L\). Thus every \(F_\ell\)
has level \(L/p\). Iteration gives level \(M_h\).

For \(p\nmid h\), \(T_p^0\) is now the classical \(T_p\), diagonal
in the basis above. Thus \(P_p(T_p^0)F_\ell=0\). Its first Fourier
coefficient is \(\ell(P_p(T_p)\bmod J)\) by the recursion and
\(T_1=1\). The functionals on \(\mathcal A/J\) separate points,
so \(P_p(T_p)\in J\) and \(P_p(T_p)^sV=0\).
The cuspidal eigenvalues are real and satisfy Deligne's bound
\cite{DeligneWeil}. The two Eisenstein characters have product one,
giving eigenvalues \(\alpha=\zeta p^{2d+1}+\zeta^{-1}\) with
\(\lvert\zeta\rvert=1\). For every such eigenvalue we have
\[
 2p^{d+1/2}<(p+1)p^d<p^{2d+1}-1\le\lvert\alpha\rvert .
\]
The first gap is \(p^d(\sqrt p-1)^2\). For the second, set \(u=p^d\ge2\).
Then \(pu^2-1-(p+1)u=u\bigl(p(u-1)-1\bigr)-1>0\). The last
inequality follows from the reverse triangle inequality.
B\'ezout's identity for \(P_p^s\)
and \(X\pm(p+1)p^d\) gives the polynomial inverses.
\end{proof}

\subsection{The marked graph relations on CM fibres}
\begin{lemma}[Two actions of a supported correspondence]
\label{labelled-two-actions}
Let \(Y\) be a smooth fine modular curve and \(\pi:\mathcal E\to Y\)
its universal elliptic curve. Suppose a correspondence on the signed
\(m\)-fold elliptic power is presented by classes on the fibre products
over normalized oriented isogeny curves, together with terms over base
point pairs. If every relative class has zero signed fibre realization,
its action on \(H^1_{\mathrm{par}}(Y,\operatorname{Sym}^m R^1\pi_*\mathbb C)\)
is zero. The product of any \(m+1\) such correspondences acts by zero
on the span \(V_h\) of signed middle-codimension classes pushed
forward from smooth CM fibres after descent to cyclic level \(h\).
\end{lemma}
\begin{proof}
The two endpoint maps of each normalized isogeny curve are finite
\'{e}tale on the interior. The action factors through pullback by
the first endpoint, action of the relative correspondence, and
trace by the second endpoint. After applying the projector the
relative direct image of the elliptic power is concentrated in
degree \(m\), with its only cohomology sheaf the indicated local
system. A degree-zero map between two such complexes is determined
by the induced map of local systems. That map is zero on a dense
open set by assumption and therefore zero on each connected
component. Terms over base point pairs act by zero on parabolic
degree one: testing against compactly supported degree-one lifts
gives the restriction of the base classes to a point. This proves
the first assertion and also identifies the parabolic action of
a scalar fibre graph with its usual pullback and trace.

For the second assertion a base point pair acts on a source fibre
by either the empty intersection or the self-intersection of a
point in a smooth curve; its excess line on the fibre is trivial.
Compose the other terms over the finite curves of composable
isogenies and specialize at the source CM point. The restricted
classes have zero realization by smooth proper base change. After
transport along the isogenies one is composing endomorphisms of
the same signed CM motive. Over \(\mathbb C\) coefficients it is
the direct sum of \(m+1\) invertible CM lines. Indeed the two
CM eigenprojectors on \(h^1\) have rank one, and for a rank-one
summand \(U\) the unit splits off \(U^\vee\otimes U\) by the trace
pairing. The complementary projector has zero realization and
vanishes by Kimura nilpotence \cite[Proposition~7.5]{KimuraFinite}.
Each of the \(m+1\) summands in the signed power is a tensor product
of these two lines. On an invertible line an endomorphism is a scalar
and its realization detects it. Any matrix entry of a product of
\(m+1\) errors is a sum of paths through the \(m+1\) lines. A path
has \(m+2\) vertices and revisits a line. Its intervening nonempty
subproduct has zero realization, hence vanishes on that line.
Summing the paths, taking the finite trace, and using compact
restriction proves the assertion.
\end{proof}

\begin{corollary}[Nilpotence on middle CM fibre classes]
\label{cm-middle-error-bound}
Let \(m=2d\ge2\). A product of \(m-1\) correspondences satisfying
the hypotheses of Lemma~\ref{labelled-two-actions} acts by zero on
\(V_h\). In particular, for \(m=2\) each such correspondence
annihilates \(V_h\).
\end{corollary}
\begin{proof}
Write \(h^1(E)=M_+\oplus M_-\) for a smooth CM fibre \(E\), with
complex coefficients, and let \(\mathbb L\) denote the Lefschetz
motive. The cup product pairs \(M_+\) perfectly with \(M_-\), so
\(M_+\otimes M_-\simeq\mathbb L\). After the common twist
\(\mathbb L^{-d}\), the signed power has one line
\[
   N_j=M_+^{\otimes(d+j)}\otimes M_-^{\otimes(d-j)}
          \otimes\mathbb L^{-d},\qquad -d\le j\le d .
\]
Its signed fibre classes of codimension \(d\) are the direct sum
of \(\operatorname{Hom}(\mathbf1,N_j)\). For \(j=\pm1\) the line
is \(M_\pm^{\otimes2}\otimes\mathbb L^{-1}\). Its group of cycles
is therefore a summand of the divisors on \(E^2\) odd under inversion
in each factor. By the Poincar\'e bundle this divisor group identifies
with \(\operatorname{Hom}(E,E^\vee)_{\mathbb C}\), on which Betti
realization is injective. Before the twist, \(M_+^{\otimes2}\)
and \(M_-^{\otimes2}\) have types \((2,0)\) and \((0,2)\), in some
order. Neither contains a divisor class. Hence
\(\operatorname{Hom}(\mathbf1,N_{1})=
\operatorname{Hom}(\mathbf1,N_{-1})=0\).

Along a chain of isogenies, transport the CM embeddings from the
first fibre. Each relative correspondence has codimension \(m\)
on the product of the two \(m\)-dimensional fibres, so it preserves
codimension \(d\) and the twist above. A path of matrix entries
passing through \(N_1\) or \(N_{-1}\) kills the input cycle.
Only \(m-1\) lines remain. A product of \(m-1\) entries visits
\(m\) lines, so a line repeats. Transport identifies the intervening
composite with an endomorphism of that invertible line. Its
realization is zero and hence it is zero. Applying this to each
isogeny chain and tracing gives the asserted bound on \(V_h\).
\end{proof}

The proof of the good-prime spectral bound uses the following fact
about theta series to recover the contribution of each discriminant label.

\goodbreak\hypertarget{frame-interpolation}{}
\begin{auxlemma}[Interpolation on the frame lattice]
Let \(W\) be a positive definite integral lattice of rank sixteen,
\(\mathcal B\subset W_{\mathbb R}\) a finite set, and \(m\ge2\) even.
Write \(\mathcal O=(W^\vee/W)/\{\pm1\}\) and, for \(o\in\mathcal O\), put
\[
 \Theta_{o,P}(\tau)=\sum_{\alpha\in o}\sum_{a\in W+\alpha}
                         P(a)q^{\|a\|^2/2}.
\]
In any nonempty open subset of \(W_{\mathbb R}\) one can choose a
rational vector \(a_H\) and a totally isotropic \(7\)-plane
\(U_0\subset(a_H^\perp)_{\mathbb C}\) with the following property.
There are homogeneous products \(P_j\) of linear forms
\((u,\,\cdot\,)\), \(u\in U_0\), vanishing on \(\mathcal B\), such that
each \(P_j\) has positive even degree and \((\Theta_{o,P_j})_{j,o}\)
has a nonzero minor of order \(|\mathcal O|\).
\end{auxlemma}
\begin{proof}
\emph{Choosing an isotropic plane.}
For each \(o\in\mathcal O\), choose a nonzero \(a_o\in W^\vee\) with label
in \(o\), off the finitely many lines spanned by the nonzero vectors
of \(\mathcal B\). This is possible in every coset of the full
lattice \(W\). Take \(R=\max_o\|a_o\|^2/2\) and let \(\mathcal S\)
be the set of sign-orbits of the \(x\in W^\vee\) with
\(\|x\|^2/2\le R\), including the zero orbit.

Choose \(a_H\) in the given open set, rational and nonzero, outside
the lines spanned by \(x\pm y\) for distinct sign-orbits in \(\mathcal S\),
and outside the planes spanned by \(a_o\) and each nonzero
\(v\in\mathcal B\). Orthogonal projection to \(a_H^\perp\)
therefore preserves these distinctions and these independent pairs.
Choose an isotropic \(7\)-plane \(U_0\subset(a_H^\perp)_{\mathbb C}\)
on which each independent real pair \((a_o,v)\) still has rank two
and each \(x\pm y\) for distinct sign-orbits has nonzero
restriction. These are finitely many open
conditions on the irreducible orthogonal Grassmannian. Each rank
condition is attained by putting the real pair in the span of
\(e_1,e_3\) and taking the plane spanned by
\(e_{2j-1}+i e_{2j}\), \(1\le j\le7\), in an orthonormal frame.
For the other conditions the projected \(x\pm y\) is nonzero
and isotropic vectors span the complexified space. Thus each open
is nonempty and their intersection is nonempty.

\smallskip\noindent\emph{Polynomials separating the labels.}
For a general \(u\in U_0\), the values \(\lambda(x)^2=(u,x)^2\)
are distinct on \(\mathcal S\). For each exceptional line choose a
generator \(v\) and \(u_v\in U_0\) with \((u_v,v)=0\) and
\((u_v,a_o)\ne0\) for every \(o\). The independence just proved
excludes only finitely many proper hyperplanes in the kernel of
\((\,\cdot\,,v)\). Put \(B_*(x)=\prod_v(u_v,x)^2\), taking
\(B_*(x)=\lambda(x)^2\) if there are no exceptional lines. The products
\[
                P_j(x)=B_*(x)\lambda(x)^{2j},
                   \qquad 0\le j<|\mathcal S|,
\]
have positive even degree. They vanish on \(\mathcal B\)
and are harmonic because the Laplacian contracts pairs of vectors
in the totally isotropic plane \(U_0\).

For each \(o\), interpolate a polynomial \(L_o(t)\) of degree
less than \(|\mathcal S|\) taking the value \(1/B_*(a_o)\) at
\(\lambda(a_o)^2\) and zero at the other interpolation points.
For \(Q_o=B_*L_o(\lambda^2)\), the row
\((\Theta_{o',Q_o})_{o'\in\mathcal O}\) through exponent \(R\)
has a single term, \(2q^{\|a_o\|^2/2}\) in column \(o\).
The factor two counts \(a_o\) and \(-a_o\), including when their
discriminant label is self-opposite. The determinant of these
rows has leading term \(2^{|\mathcal O|}q^{\sum_o\|a_o\|^2/2}\).
Since each interpolated row is a constant linear combination of
the homogeneous rows, some square minor of the homogeneous rows
has nonzero determinant. This proves the assertion.
\end{proof}

\begin{proposition}[The good-prime spectrum on CM fibres]\label{marked-cm-spectrum}
Let \(m=2d\ge2\) and let \(V_h\) be the complex span in compact
Chow of the signed middle-codimension classes pushed forward from
smooth CM fibres of the elliptic \(m\)-fold power at cyclic level \(h\).
The prime-to-\(h\) Hecke action on \(V_h\) admits a filtration of
length at most \(m-1\). Its graded pieces are direct sums of joint
eigenspaces occurring in \(M_{m+2}(\Gamma_0(H_h))\), where one
may take \(H_h=h^2\operatorname{rad}(h)\) for odd \(h\) and
\(H_h=4h^2\operatorname{rad}(h)\) for even \(h\).
At level one the action on \(V_1\) factors through the full-level
cuspidal Hecke algebra. In particular, for every \(p\nmid h\) the
operators \(T_p\pm(p+1)p^d\) are invertible on \(V_h\).
\end{proposition}
\begin{proof}
We apply the spectral criterion of Proposition~\ref{md-t-spectrum}.
We first specify the level of the correspondences. On a normal fine
cover of the open cyclic modular stack let \(\widetilde{\mathcal A}\)
be the algebra of signed correspondences presented on normalized
isogeny curves or over finitely many base point pairs. The endpoints
carry independent fine markings, so each marked primitive graph is
included. Let \(\widetilde J\) consist of the correspondences admitting,
modulo point pairs, a presentation with zero signed fibre realization
on every supporting curve.

Write \(U\) for the smooth elliptic power over the open fine modular
curve. Consider two composable correspondences on \(U\). Write
\(\alpha\) for the class of the first in \(\CH(U\times U)\), and
present the second by a class \(\gamma\in\CH_*(T)\) with proper
endpoint maps \(p_1,p_2:T\to U\). Their composite, in this order, is
\[
 (1\times p_2)_*
 \bigl((1\times p_1)^*\alpha\cap([U]\times\gamma)\bigr).
\]
Here smoothness identifies \(\alpha\) with an operational class.
Thus the formula depends only on the Chow class of the first
correspondence. To see the same for the second, take a proper
presentation of the first and apply the formula to the transposes
in reverse order. Transposing the result gives the original
composite: both constructions pull back along the diagonal of the
middle copy of \(U\) and push forward to the two outer copies.
Since the endpoint maps on the base isogeny curves are finite,
composition on either side preserves support over finitely many
point pairs. On the remaining curves, fibre realization respects
composition. Thus \(\widetilde J\) is a two-sided ideal.

Let \(P\) denote the graph of the fine cover of the elliptic powers,
of degree \(\delta\). Then \(PP^t=\delta\,1\) downstairs, and
\(e=P^tP/\delta\) averages the deck transformations upstairs.
Cyclic-level correspondences act on the invariant summand by
\[
  I(A)=P^tAP/\delta,\qquad I(AB)=I(A)I(B),\qquad I(1)=e.
\]
Use \(\mathcal A=e\widetilde{\mathcal A}e\) and \(J=e\widetilde Je\)
in Proposition~\ref{md-t-spectrum}. Pullback and normalized trace
identify \(V_h\) with the invariant part of the fine-level CM space.
Since \(J^{m-1}\subset e\widetilde J^{\,m-1}e\), the fibrewise
argument of Corollary~\ref{cm-middle-error-bound} gives \(J^{m-1}V_h=0\).
In the linear graph relations at fine level, a functional \(\ell\)
on \(\mathcal A/J\) evaluates a graph \(G_D\) as
\(\widetilde\ell(G_D)=\ell(eG_De\bmod J)\). This makes sense for
each branch even if it does not preserve the invariant summand.
We construct a fixed finite-dimensional space for \eqref{md-t-series},
first for \(h\ge2\).
\smallskip\noindent\emph{1. The lattice of the marked Kummer surface.}
On the product of two fine modular curves write
\((E,L),(E',L')\) for the two marked curves. Take the finite
\'etale cover of isomorphisms \(\phi:L\to L'\) and two-torsion
frames of \(A_\phi=(E\times E')/\operatorname{graph}(\phi)\).
On a connected component \(X\) of this cover, put
\(S=\operatorname{Km}(A_\phi)\). The two abelian axes and the
sixteen labelled exceptional divisors are defined over \(X\).
They span the Picard lattice \(\Lambda\) of the geometric generic
fibre rationally. On that fibre choose a Jacobian elliptic pencil,
which exists since the Picard rank is \(18\)
\cite[Proposition~2.3]{BertinFibrations}. Its fibre \(F\) and
zero section \(O\) give a unimodular plane \(U=\langle F,O+F\rangle\).
The lattice splits as \(U\oplus W\), with positive norm
\(\|a\|^2=-a^2\) on the rank-sixteen frame \(W\).

Both classes \(F\) and \(O\) are Galois invariant. Indeed, a positive multiple
of each is an integral combination of the labelled divisors and
\(\Lambda\) is torsion free. The tangent to the first axis gives a
rational point of the exceptional line above the origin, so Picard
descent gives line bundles on the generic surface. Close their
divisors on the smooth total space. The difference between a suitable
multiple of either closure and the corresponding labelled combination
is vertical, hence pulled back from \(X\) by smoothness and geometric
integrality of the fibres. Thus the labelled combinations represent
the desired classes on every fibre. The pencil itself is needed only
on the geometric generic fibre.

Over a generic marked isogeny curve \(D\) of primitive degree \(i\),
choose an orientation \(\eta:E\to E'\). The extra rational Picard
direction \(z\perp\Lambda\) is fixed by prescribing mixed entry
\(\eta\) on \(E\times E'\). The pure mixed divisor there has square
\(-2i\). The Kummer quotient and the abelian isogeny have degrees
\(2,h\), respectively, giving \(z^2=-i/h\). On any framed lift
\(D'\) the integral divisor classes of mixed part \(sz\) have
\[
 s\in d_D\mathbb Z,\qquad
 d_D=\#((\eta L+L')/L'),\qquad
 a\in W+\alpha_{D'}(s),\quad \alpha_{D'}(s)\in W^\vee/W,
\]
where \(a\) is their frame component. To see the condition on \(s\),
choose a generator \(P\) of \(L\) and put \(Q=\phi(P)\). The
symmetric matrix of a line bundle on \(E\times E'\) with mixed
entries \(s\eta,s\widehat\eta\) kills \((P,Q)\) precisely when its
integer diagonal entries \(c_1,c_2\) satisfy
\[
 c_1P=-s\widehat\eta Q,\qquad c_2Q=-s\eta P .
\]
The second equation is soluble exactly when \(s\eta L\subset L'\).
By adjointness for the Weil pairings, this is equivalent to
\(s\widehat\eta L'\subset L\), since each cyclic subgroup of order
\(h\) is its own annihilator. Thus both diagonal entries can be
chosen, and the possible \(s\) are precisely the multiples of \(d_D\).

Translation by \((P,Q)\) fixes the line bundle. Over the geometric
generic point a lift of this generator can be rescaled to have \(h\)-th
power one, so the bundle descends to \(A_\phi\). Divisibility of
\(\operatorname{Pic}^0\) gives a symmetric representative. Rescale its
inversion isomorphism to square to the identity and descend on the free
locus. The closure of a representing divisor on the smooth Kummer
resolution is integral. Its pullback and blowdown recover the original
class since the omitted two-torsion has codimension two. Both descents
also apply for even \(h\). The integral
kernel of the mixed component is the saturated generic Hodge
lattice \(\Lambda\). Pairing with \(W\) puts the frame component in
\(W^\vee\), and changing the lift translates it by an element of
\(W\). The splitting by \(U\) realizes every such translate.

\smallskip\noindent\emph{2. Bounding the exceptional classes.}
We need one polynomial insertion that removes all contributions
from vertices of nonpositive \(F\)-degree.
Choose an extending relatively ample class \(H_0\). If an effective
nonextending vertex \(B=n\delta\), with \(\delta\) primitive, has
nonzero reduced class, then \(\delta^2\ge-2\)
\cite[Proposition~3]{PandharipandeYinK3}.
If \(FB\le0\), the positive cone and \(H_0\delta>0\)
force \(\delta^2=-2\). Its extending projection has square
\(-2+s_0^2i/h\le0\), where \(s_0z\) is its mixed part, by the same
cone argument. Thus \(i\le2h\). For \(F\delta=0\) its frame vector
has norm squared \(2-s_0^2i/h\). Otherwise put \(c=H_0F\),
\(j=-F\delta>0\), and write \(\delta=-jH_0/c+xF+y\) with
\(y\perp H_0,F\). Then
\[
     2=j^2H_0^2/c^2+2j(H_0\delta)/c+\|y\|^2 .
\]
This bounds \(j,x,y\) in finitely many Picard lattices. Hence the
primitive frame directions form a finite set \(\mathcal B\).
Include their negatives. If \(FB=1,\ B^2\ge-2\) and \(B\) is
noneffective, Riemann--Roch and the positive cone imply that
\(-B\) is an effective root of negative \(F\)-degree, also covered
by this bound.

Choose a rational ample class \(H\) and a totally isotropic complex
\(7\)-plane \(U_0\) perpendicular in \(W_{\mathbb C}\) to its frame
projection. Take \(P(a)=\prod_{j=1}^l(u_j,a)\) in the positive frame
pairing, with \(u_j\in U_0\), \(l>0\) even, and \(P|_{\mathcal B}=0\).
Set \(b=m+l=2d'\). Use the squarefree Pixton coefficient of
Proposition~\ref{family-unmarked-rows} in genus \(d'-1\) for
\(T=2O+tF,\ t\ge2\). Insert \(F\) at all negative marks,
retain \(m\) positive evaluations, and insert the \(u_j\) at the
others. The family reduced dimension \(d'-1+2b\), less the relation
degree \(d'\) and \(2b-m\) divisor conditions, is \(m-1\)
relative to the base surface.

\smallskip\noindent\emph{3. The contribution of a pair of vertices.}
First consider a single nonconstant vertex.
Choose extending divisor representatives by the two
coordinate curves and exceptional divisors. Any two have total
Chow product their fibre pairing times one fixed relative point
class: upstairs the normal lines of the coordinate elliptic
subgroups have torsion first Chern class, their transverse
intersection is torsion, and torsion translations act trivially by
invertibility of multiplication on the rational abelian motive.
The exceptional divisors are disjoint, have zero product with the
total transforms of the axes, and their squares are minus the
pullback of their centres since the first Chern class of the
normal bundle of each centre is zero. Descent by finite trace
gives the assertion on \(S\). Thus two fixed marks at a constant
vertex contribute zero by isotropy. A retained mark there with a
fixed mark has a torsion coordinate \(t\) after transfer. The
involution \(x\mapsto 2t-x\) fixes it, and torsion translation
acts trivially in rational Chow. Hence its odd projection vanishes.
Two retained marks there are killed by the signed
transposition. The count in Lemma~\ref{family-single-positive}
puts at least \(b\) marks on at most \(d'\) constant vertices.

Apply Proposition~\ref{family-supported-diagonal} on the marked family
before inserting the divisors. Its auxiliary period polarization may
be chosen separately for each total class \(T\). The divisors
\(F,u_j\) are inserted after pullback to the framed base, so the
same frame polynomials work for every \(t\). For fixed \(t\), \(T\) is
effective with \(T^2=4t-8\ge0\). The two abelian axis degrees of any
effective summand are bounded by those of \(T\), and their product bounds
its squared mixed norm. Thus only finitely many isogeny curves contribute.
At a labelled generic point the two reduced vertices satisfy
\(B_1+B_2=T\) and have mixed parts \(sz,-sz\) with \(s\ne0\).
The period intersection has multiplicity one, with both node flags
retained. Here only the fibre Betti scalar is needed modulo the
marked error ideal; the Chow factorization for the product pencil
is used in the unmarked relations. Terms over point pairs act
by zero. We must determine the fibre degrees, since the marked
pencil can lose nefness.

Choose a real decomposable mixed test \(v\perp\Lambda\), \(v^2=0\),
\((v,z)=1\): decomposable mixed tensors span the mixed space.
Equivalently, choose a decomposable mixed tensor on \(E\times E'\)
pairing to one with the mixed graph of \(\eta\) and pull the test
back by the adjoint correspondence. The mixed entry of \(z\) fixes
this normalization up to the common sign of the Poincar\'e divisor.
The axis and exceptional classes have zero mixed entry. The
\(m\)-fold tensor of the test pairs to one with the signed graph,
since the signs from permuting degree-one factors cancel those of
the projectors and \(m\) is even. This normalization works for
every primitive degree and every lift of a framing.

Put \(D_0=\langle F,v,B_1,B_2\rangle\) and
\(K=\ker((-,B_1):\mathbb C v+U_0\to\mathbb C)\).
The isotropic pair \(F,v\) couples to \(B_1,B_2\) by the matrix
with rows \((FB_1,FB_2)\) and \((s,-s)\). Its determinant is
\(-sF(B_1+B_2)=-2s\ne0\), so \(D_0\) is nondegenerate.
Every element of \(K\) is orthogonal to \(F,v,B_1\) and to
\(B_1+B_2=T\), hence \(K\subset D_0^\perp\). This isotropic
seven-plane has a dual in the nondegenerate orthogonal complement.
Subtract half the Gram matrix along \(K\) to make the dual isotropic.
Scaling this pair of dual spaces by \(\lambda,\lambda^{-1}\)
fixes \(D_0\). Each fixed-surface descendant tensor is invariant
under this torus, also for divisible effective classes
\cite[\S2.2, (7)]{OberdieckNL}, since cotangent line bundles extend in the same deformation.
For a positive-degree rational vertex with fewer than three flags
and markings this invariance follows from the divisor equation as
well. If \(F B_i\ne0\), add the fixed divisor \(F\). If \(F B_i=0\)
and the tensor can contribute, the root bound above gives
\(B_i=n_i\delta_i\) with \(\delta_i^2=-2\).
Use \(B_i\) itself as the added divisor, since \(B_i^2=-2n_i^2\ne0\).
In either case the torus fixes the divisor. The divisor equation,
inducting on the number of marks missing from stability and then
on the total cotangent exponent as in
Proposition~\ref{family-unmarked-rows}, preserves the invariance.
Taking \(\lambda\to0\) replaces each \(u_j\) by \((u_j B_1/s)v\).
By the root bound, \(P\) kills every fibre-degree pair other than
\((1,1)\).

For the surviving pairs put
\[
 B_{1,2}=O+x_{1,2}F\ \pm(a+sz),\qquad
 g_j=1+B_j^2/2,\qquad e=\|a\|^2+s^2i/h .
\]
Here \(e\) is positive and even by integrality of \(a+sz\).
Squaring gives \(x_j=g_j+e/2,\ t=e+g_1+g_2\) and
\(B_1B_2=2t-2-g_1-g_2\). Both classes are primitive.
A negative \(g_j\) gives a zero tensor, and the missing noneffective
pairs with \(g_j\ge0\) have already been killed by \(P\).
The Gram matrix of \(F,v/s,B_1,B_2\) agrees with the four-vector
matrix in Proposition~\ref{family-unmarked-rows}: the first two
vectors pair against the two vertices as \((1,1),(1,-1)\).
To realize this matrix on the unmarked product, take a primitive
cyclic isogeny of degree \(e\) with all four unselected branch labels
matched. The classes \(M_i+g_iF\) have the prescribed squares and
intersection \(2e-2+g_1+g_2\). The resulting isometry between
nondegenerate four-dimensional spans extends to real K3 cohomology.
The primitive vertex comparison \cite[Corollary~3]{BuellesK3Classes}
therefore gives exactly its contraction multiplied by \(s^mP(a)\).
Indeed the all-\(v\) test supplies \(s^b\) and the \(l\) replacements
divide by \(s^l\). Interchanging vertices negates \(a,s\) and
exchanges their arithmetic genera, as in the ordered unmarked sum.

\smallskip\noindent\emph{4. A recurrence for the scalar graph series.}
For each normalized isogeny curve \(D\), let \(i_D\) be the degree
of its primitive isogeny and \(G_D\) its signed primitive graph.
After transport by \(\eta\), the realization of the relative class
on an oriented lift \(D'\) commutes with finite-index modular
monodromy on \(\operatorname{Sym}^m H^1\). This monodromy group is
Zariski dense in \(\mathrm{SL}_2\), and the symmetric power is
irreducible. Its scalar therefore determines the relative class
modulo the error ideal. Any two extensions to the normalized curve
differ by terms over point pairs, and proper trace multiplies the
graph by \([D':D]\). For a functional \(\ell\) on the quotient, set
\[
 \begin{aligned}
 f_\alpha(\tau)&=
 \sum_D\sum_{D'\mid D}[D':D]\!
 \sum_{\substack{0\ne s\in d_D\mathbb Z\\\alpha_{D'}(s)=\alpha}}
s^m\widetilde\ell(G_D)q^{s^2i_D/(2h)},\\
 \theta_{\alpha,P}(\tau)&=\sum_{a\in W+\alpha}P(a)q^{\|a\|^2/2}.
 \end{aligned}
\]
Each coefficient is a finite sum. Negating \(s\) negates its
discriminant label and preserves \(s^m\), so \(f_{-\alpha}=f_\alpha\).
Put
\[
 y_{2e}=[q^e]\sum_\alpha f_\alpha(2\tau)\theta_{\alpha,P}(2\tau),
 \qquad y_{2e+1}=0 .
\]
The marked Chow row reads
\[
 P_b(2t-2)y_{2t}
 +\sum_{u=1}^{t-1}C_b(u,2t-2-u)y_{2t-2u}=0.
\]
All odd-index equations vanish. The equation at \(t=1\) also holds,
since \(e\ge2\) gives \(y_2=0\). Hence the marked series satisfies
the whole recurrence \eqref{rd-m-rec} in degree \(b=2d'\).
By Proposition~\ref{framed-comparison-solutions} the entire
solution space is \(\Delta_{\rm mod}R_{d'-1}\) for either parity
of \(d'\). Its computation on parabolic cohomology applies because
a labelled lift and its graph with the computed fibre scalar act
identically there (Lemma~\ref{labelled-two-actions}).
Translation by one fixes the integral-exponent series
and sends \(B\) to \(-B\), while fixing \(A\) and \(\Delta_{\rm mod}\).
Consequently
\[
 \sum_e y_{2e}q^e\in\Delta_{\rm mod}\mathbb C[A,B^2]_{d'-1}.
\]

\smallskip\noindent\emph{5. Separating the discriminant labels.}
The frame projection of the ample cone is open. Apply the
\hyperlink{frame-interpolation}{frame interpolation lemma} to this
open set and the exceptional directions
\(\mathcal B\). Lift the resulting \(a_H\) to a rational ample class
\(H\), and use its plane \(U_0\) in all the marked relations.
Since \(f_{-\alpha}=f_\alpha\), the unknowns are indexed by
\(o\in(W^\vee/W)/\{\pm1\}\). The minor supplied by the lemma recovers
them from the corresponding marked rows.

To bound the level explicitly, let \(c=\operatorname{lcm}(2h,4)\).
The two abelian axes in \(A_\phi\) have intersection \(h\) and
self-intersection zero. Pull them back to the blowup at the two-torsion
points and push them to the Kummer surface. Their intersection becomes
\(2h\), and both classes are orthogonal to the sixteen exceptional
curves. Thus \(\Lambda\) contains the full-rank sublattice
\(U(2h)\oplus\langle-2\rangle^{16}\).
The quadratic form on its dual has denominator dividing \(c\).
The same holds on \(W^\vee\) by the orthogonal splitting of \(\Lambda\).
Consequently the rank-sixteen Weil representation is trivial on
\(\Gamma(c)\) \cite[\S2]{BruinierSteinWeil}. Poisson summation gives
weight \(8+\deg P_j\) for each theta row, and each row after dilation
by two is modular on \(\Gamma(2c)\). For the chosen rows,
the marked recurrence places the right-hand sides in the fixed
spaces \(\Delta_{\rm mod}R_{(m+\deg P_j)/2-1}\), of weight
\(m+\deg P_j+10\). In Cramer's rule every cofactor term in a
given numerator replaces one row weight by that weight plus \(m+2\).
The quotients \(f_\alpha(2\tau)\) therefore have weight \(m+2\)
and belong to a fixed finite-dimensional space of meromorphic
modular forms, independent of \(\ell\).

\smallskip\noindent\emph{6. The prime-to-level Hecke series.}
Sum the \(f_\alpha(2h\tau)\) over the labels and all connected
components of the complete frame cover. At index \(n=s^2i\)
coprime to \(h\), the condition \(s\eta L\subset L'\) forces
\(\eta L=L'\), and every framing contributes. The sum of
\([D':D]\) is now the degree of the complete cover, independent
of \(D\). Both signs of \(s\) give the same graph. Multiplication
by \(s\) acts by \(s^m\) on the signed power, so including all
positive \(s\) and primitive graphs counts exactly all kernels
of order \(n\), also when \(n\) is a square. At both endpoints we
sum over all lifts of the cyclic-level marking to the fine cover.
This finite flat pullback of the cyclic Hecke correspondence is
\(P^tT_nP=\delta I(T_n)\). Its degree over the correspondence is
\(\delta^2\), independently of \(n\); no isogeny has to transport
a fine frame. Thus the surviving series is a fixed nonzero multiple
of \eqref{md-t-series}. Before imposing \((n,h)=1\), the final
dilation gives level \(\Gamma(2ch)\). The exponents \(s^2i_D\)
are integral, so translation by one gives level \(\Gamma_1(2ch)\).
At each \(p\mid h\) apply \(1-V_pU_p\) to delete the multiples of
\(p\). The finite correspondence \(U_p\) preserves meromorphicity
at this level, and \(V_p\) raises the level by at most \(p\).
Thus \(M=2ch\operatorname{rad}(h)\) is an auxiliary level common
to all \(\ell\). Its part supported on \(h\) is the stated \(H_h\).
The kernel count of
Lemma~\ref{full-graph-hecke} applies unchanged for \(p\nmid h\):
prime-to-\(h\) isogenies transport the cyclic subgroup.
Proposition~\ref{md-t-spectrum} gives the asserted inverses for \(h\ge2\).

For the filtration, let \(\mathcal H\subset\mathcal A/J\) be the span
of the \(T_n\bmod J\) for \((n,h)=1\). The recursions make it a unital
commutative algebra. Its full dual embeds by \eqref{md-t-series}
into the classical space at level \(M_h\), since every functional on
\(\mathcal H\) extends to \(\mathcal A/J\). This embedding intertwines
the good Hecke operators, which are simultaneously diagonal in a
newform and Eisenstein basis. The dual regular action is faithful,
so \(\mathcal H\) is a finite product of copies of \(\mathbb C\).
Its character idempotents split each \(J^iV_h/J^{i+1}V_h\), even if
an eigenspace has infinite multiplicity. The bound \(J^{m-1}V_h=0\)
gives the stated length. In particular the action on \(V_h\) is
already semisimple when \(m=2\).

At level one the unmarked graph \(H_m(n)\) acts as \(T_n\) on
\(V_1\). Point-pair correspondences kill every class pushed forward
from a smooth fibre by Lemma~\ref{labelled-two-actions}. Thus the
action factors through the cuspidal Hecke algebra identified after
Corollary~\ref{all-untwisted-modularity}. Its idempotents split
\(V_1\) itself, and Deligne's bound excludes both exceptional
eigenvalues. If the cusp space is zero, the identity graph acts by
zero and \(V_1=0\).
\end{proof}

For \(h\ge2\) the level bound gives a finite test for graph relations
modulo \(J\). Put
\(B=\lfloor(m+2)[{\rm SL}_2(\mathbb Z):\Gamma_0(H_h)]/12\rfloor\).
Then the \(T_n\bmod J\) with \((n,h)=1\) are spanned by those with
\(n\le B\). The series of any functional vanishing on these first
operators has zero constant term and order at infinity greater than
\(B\), so the valence formula makes it zero. At level one the same
argument applies to the full-level graphs modulo point pairs.

\subsection{Primes shared by the level and the twist}

Fix a shared prime \(p\). The shared-prime calculation in
Section~\ref{rd-cycles} gives \(p\Vert N\), including \(p=2\).
Put \(h=N/p\) and \(t=p^d\), and write \(D\) for
pullback along the forgetful map, and \(D^*\) for trace. Let \(W\)
be the normalized reversal of a marked order-\(p\) edge. Then
\[
          D^*D=(p+1)1,\qquad D^*WD=t^{-1}T_p^{(h)} .
\]
The first identity counts the \(p+1\) lines above each source.
In the second one, quotienting by the selected line, then
forgetting the line, gives the ordinary Hecke correspondence
with the normalization \(t^{-1}\).
In terms of directed edges
\[
 U_p=tW(DD^*-1),\qquad U_p^*=t(DD^*-1)W .
\]
One changes the outgoing line to any other and reverses the edge,
in the indicated order. Counting on a fine cover and then tracing
preserves the stabilizer factors by the orbit calculation above.
Put \(A=(1+W)D\) and \(c=2(p+1)\). On the \(W\)-fixed space at
level \(N\) and on the full space at level \(h\), respectively, we obtain
\begin{equation}\label{shared-edge-kernel-operators}
 B_p-2pt=t(AD^*-c),\qquad D^*A-c=t^{-1}T_p^{(h)}-(p+1).
\end{equation}
Thus trace identifies the first kernel with the second, with inverse
\(A/c\). Indeed \(AD^*\) and \(D^*A\) are each multiplication by \(c\)
on the corresponding kernel, and the two maps intertwine these products.

\begin{theorem}[Successive kernels on cycles of weakly holomorphic forms]
\label{actual-defect-kernels}
Write \(D=D_{\Delta,r}(N)\). Split each good-prime factor of
\eqref{md-t-product} into \(T_p-(p+1)t_p\) and \(T_p+(p+1)t_p\).
List these and the shared-prime factors \(B_p-2pt_p\) in any order
as \(A_1,\ldots,A_s\). Put \(P_0=1\), \(P_j=A_jP_{j-1}\), and
\(D_j=\ker(P_j|_D)\). There is a filtration
\(0=D_0\subseteq D_1\subseteq\cdots\subseteq D_s\subseteq D\) with
\begin{align}
 D_j/D_{j-1}&\ \simeq\ P_{j-1}D\cap\ker A_j,\qquad 1\le j\le s,
 \label{defect-kernel-layer}\\
 D/D_s&\ \simeq\ L_{\Delta,N}D\ \subseteq\ \sum_i C_iD_{1,1}(M_i).
 \label{defect-transfer-layer}
\end{align}
At a shared prime the trace to level \(N/p\) identifies the
intersection in \eqref{defect-kernel-layer} with its image in
\(\ker(T_p-(p+1)t_p)\) there. In particular, once the untwisted
relations hold at the auxiliary levels, Chow-valued modularity is
equivalent to the simultaneous vanishing of the intersections
in \eqref{defect-kernel-layer}.
\end{theorem}
\begin{proof}
Work on the smooth base using Lemma~\ref{scholl-open-chow}.
Each coefficient, hence each element of \(D\), is fixed by
the shared-prime reversals by the full-test calculation in
Proposition~\ref{md-transfer-proof}. All the factors preserve those
fixed spaces: prime-to-\(p\) isogenies transport the \(p\)-line,
and reversal interchanges the two summands of \(B_p\).
The map \(P_{j-1}:D_j\to P_{j-1}D\cap\ker A_j\) is onto by
definition and has kernel \(D_{j-1}\). The last quotient follows
from \(P_s=L_{\Delta,N}\) and Proposition~\ref{md-transfer-proof}.
At a shared prime \eqref{shared-edge-kernel-operators} identifies each
intersection with its trace image at the lower level. Theorem~\ref{rd-main}
identifies \(D=0\) with modularity.
\end{proof}

This argument tests the images \(P_{j-1}D\) without assuming that any
individual Hecke operator preserves \(D\). To identify the kernel on
the entire coefficient space we impose the following compatibilities.
\begin{proposition}[The kernel of transfer on the coefficient space]
\label{merged-coefficient-kernel}
At level \(m\) write \(W_m,D_m\) as in Theorem~\ref{rd-main} and set
\(Q_m=W_m/D_m\), retaining the twist and reducing the root at lower levels.
Put \(L=L_{\Delta,N}\) as in \eqref{md-t-product}.
Suppose that the factors of \(L\) commute on \(W_N\) and preserve \(D_N\).
At each shared prime assume that edge reversal \(W\) fixes \(W_N\), and
that trace \(D^*\) and the symmetric pullback \(A=(1+W)D\) carry both coefficient and defect
spaces to their counterparts at the other level. At each prime
\(p\mid\Delta\), take \(m=N\) if \(p\nmid N\), and \(m=N/p\) otherwise.
Assume that \(T_p\) acts on \(Q_m\) and that its eigenvalues there have
absolute value at most \(2p^{d+1/2}\). Then there is an exact sequence
\[
 0\longrightarrow\ker(L\mid W_N)\longrightarrow D_N
       \xrightarrow{\;L\;} L D_N\longrightarrow 0 .
\]
At a shared prime, trace identifies
\(\ker(B_p-2pt_p\mid W_N)\) with
\(\ker(T_p-(p+1)t_p\mid W_{N/p})\), with inverse \(A/(2(p+1))\).
\end{proposition}
\begin{proof}
By the assumed compatibilities, \eqref{shared-edge-kernel-operators}
restricts to the coefficient spaces and descends to their quotients.
It gives the trace isomorphism in both cases. Since
\((p+1)p^d>2p^{d+1/2}\), the spectral bound excludes both exceptional
roots on \(Q_N\) at unshared primes and the positive root on \(Q_{N/p}\)
at shared primes. The trace isomorphism thus excludes the latter kernel
upstairs as well.
Every factor of \(L\) acts injectively on \(Q_N\). It follows that
\(\ker(L\mid W_N)\subset D_N\), which proves exactness.
\end{proof}

Thus \(D_N=0\) is equivalent to \(L D_N=0\) together with the absence
of eigenvalues \(\pm(p+1)t_p\) on \(W_N\) at unshared primes and of
\(+(p+1)t_p\) on \(W_{N/p}\) at shared primes. The quotient
bound in the proposition can be checked on cusp forms when the surjection
\(S^\vee\to Q_m\) of Theorem~\ref{rd-main} intertwines the good Hecke
operator with the adjoint of an operator having the stated bound.

\begin{corollary}[Inverting the exceptional Hecke factors]
\label{merged-hecke-inverses}
Under the hypotheses of Proposition~\ref{md-t-spectrum}, with
\(d\ge1\) and \(p\nmid h\), each of \(T_p\pm(p+1)p^d\) has an inverse
on \(V\) given by a polynomial in \(T_p\). On the full CM module
\(V_h\) of Proposition~\ref{marked-cm-spectrum}, a single polynomial
works for every CM field and conductor at that level.
\end{corollary}
\begin{proof}
For \(h\ge2\) use \(P_p^{\,m-1}\) from
Proposition~\ref{md-t-spectrum} and the filtration in
Proposition~\ref{marked-cm-spectrum}. For \(h=1\) use the
characteristic polynomial on the full-level cusp space. Each
B\'ezout inverse is independent of the CM field and conductor.
\end{proof}

At a shared prime let \(S\) be the polynomial inverse of
\(D^*A-c=t^{-1}T_p^{(N/p)}-(p+1)\) on the lower CM space.
On the \(W\)-fixed space upstairs, \eqref{shared-edge-kernel-operators} gives
\[
            (B_p-2pt)^{-1}=(ct)^{-1}(ASD^*-1).
\]
Multiplying in either order verifies the formula using
\(S(D^*A-c)=(D^*A-c)S=1\).

\begin{proposition}[A prime-power test on CM cycles]\label{prime-power-cm-test}
Fix \(m=2d\ge2\), \(p\nmid h\), and write \(S_h=S_{m+2}(\Gamma_0(h))\).
Let \(P(X)\) be the product of \(X-\alpha\) over the distinct eigenvalues
of the classical \(T_p\) on \(S_h\), with empty product \(1\). Define
\[
 U_0=1,\qquad U_1=X,\qquad U_{v+1}=XU_v-p^{m+1}U_{v-1}\quad(v\ge1).
\]
Write \(P=\sum_{v=0}^{s}b_vU_v\), where \(s=\deg P\le\dim S_h\).
On the CM cycle space \(V_h\) of Proposition~\ref{marked-cm-spectrum} set
\[
                R_{h,p}=\sum_{v=0}^{s} b_v T_{p^v}.
\]
If for each \(z\in V_h\) some power of \(R_{h,p}\) kills \(z\), then
both \(T_p\pm(p+1)p^d\) are injective on \(V_h\).
At \(h=1\) the single graph identity
\begin{equation}\label{prime-power-graph}
                 \sum_{v=0}^s b_v H_m(p^v)=0\pmod{J_m}
\end{equation}
implies \(R_{1,p}=0\). In that case, for \(\lambda=\pm(p+1)p^d\),
the inverse of \(T_p-\lambda\) on \(V_1\) is \(Q_\lambda(T_p)\), where
\[
 Q_\lambda(X)=\frac{P(\lambda)-P(X)}{P(\lambda)(X-\lambda)} .
\]
\end{proposition}
\begin{proof}
The count in Lemma~\ref{full-graph-hecke} applies with the cyclic
marking since \(p\nmid h\). Thus \(T_{p^v}=U_v(T_p)\) on \(V_h\)
and \(R_{h,p}=P(T_p)\). At this good prime the cusp eigenvalues in
weight \(m+2\) have absolute value at most \(2p^{d+1/2}\)
\cite{DeligneWeil}, strictly less than \((p+1)p^d\).
For either exceptional value \(\lambda\), \(P(\lambda)\ne0\).
On a \(\lambda\)-eigenvector every power \(R_{h,p}^j\) acts by
\(P(\lambda)^j\), so local nilpotence forces that vector to vanish.

At level one \(H_m(p^v)\) acts on \(V_1\) by \(T_{p^v}\). Point-pair
correspondences act trivially on CM fibres by
Lemma~\ref{labelled-two-actions}, so \eqref{prime-power-graph} gives
\(P(T_p)=0\) on \(V_1\).
The asserted inverse follows by evaluating
\((X-\lambda)Q_\lambda(X)=1-P(X)/P(\lambda)\).
\end{proof}

In particular \eqref{prime-power-graph} follows from the full-level
graph relations in degree \(m\) alone. For \(g\in S_1\) the scalar
relation is
\(\sum_v b_v [q^{p^v}]g=[q]P(T_p)g=0\), so these relations apply to
exactly the coefficients in \eqref{prime-power-graph}. This replaces
the marked spectral input whenever the good level in the transfer
is one, including at a shared prime with \(N=p\). For larger good
levels the test concerns the action of \(R_{h,p}\) on cycles modulo
rational equivalence. Its vanishing on cusp forms was used only to
choose the coefficients \(b_v\).

\begin{proposition}[A criterion in ordinary Chow]\label{merged-chow-interface}
Fix \(N,\Delta,r\) and \(\kappa\ge2\), and put \(a=\kappa-1\).
At each of the auxiliary levels \(M_i\) of
Proposition~\ref{md-transfer-proof}, take the finitely many
even integers \(b\) occurring in the scalar criterion
of Theorem~\ref{med-u-scalar}. For each such \(b\), consider all
finitely supported families \((u_n)\) such that \(\sum_n u_n c_h(n)=0\)
for every \(h\in S_{b+2}(\mathrm{SL}_2(\mathbb Z))\). Suppose
that every such family satisfies
\[
                \sum_n u_n H_b(n)=0 \pmod {J_b},
\]
where \(J_b\) denotes the classes supported over finite sets of base
point pairs. For \(\Delta=1\) use the identity transfer and \(M_1=N\).
Finally suppose that \(L_{\Delta,N}\) of \eqref{md-t-product} is
injective on the restriction of \(D_{\Delta,r}(N)\) to the smooth
base. Then \(Z^{(N)}_{\Delta,r,\kappa}(F)=0\) in ordinary compact Chow for every weakly
holomorphic \(F\) of weight \(3/2-\kappa\) and type \(\rho_{\rm in}\).
\end{proposition}
\begin{proof}
Theorem~\ref{graph-to-moments} carries the graph relations to the
trace moments and kills \(J_b\). The scalar criterion and compact
restriction (Lemma~\ref{scholl-open-chow}) therefore give
\(D_{1,1}(M_i)=0\) at every auxiliary level. Now
\eqref{md-t-inclusion} gives \(L_{\Delta,N}Z^{(N)}_{\Delta,r,\kappa}(F)=0\).
Injectivity on the smooth base kills the restriction of \(Z_\kappa(F)\),
and the same compact restriction lemma kills its compact class.
\end{proof}

The unmarked calculation supplies the graph relations in Chow. On the
marked family the comparison of fibre realizations suffices because
the error ideal acts nilpotently on CM cycles. We now combine these
two results.
\begin{theorem}[Vanishing in compact Chow]\label{rd-m-vanishing}
For every cyclic level \(N\), codimension \(\kappa\ge2\), compatible
fundamental twist \((\Delta,r)\), and
\(f\in M^!_{3/2-\kappa,\overline{\widetilde\rho_L}}\), one has
\[
 \sum_{\mu,n>0}c_f(-n,\mu)Z^{(N)}_{\Delta,r,\kappa}(n,\mu)=0
 \qquad\text{in }\CH^\kappa(\mathcal Y_\kappa(N)_{\overline{\mathbb Q}})_{\mathbb C}.
\]
\end{theorem}
\begin{proof}
Write \(D=D_{\Delta,r}(N)\). Corollary~\ref{all-untwisted-modularity}
gives \(D_{1,1}(M_i)=0\) at every auxiliary level. Hence the last
quotient in \eqref{defect-transfer-layer} is zero. Each earlier
quotient in \eqref{defect-kernel-layer} embeds in a good-prime
exceptional eigenspace on CM cycles, at level \(N\) for \(p\nmid N\)
and at level \(N/p\) by trace for a shared prime.
All intermediate images lie in the full CM cycle space, since isogenies
carry CM fibres to CM fibres. At a shared prime the earlier factors
also preserve the subspace fixed by reversal. Proposition~\ref{marked-cm-spectrum}
makes each of the resulting eigenspaces zero.
The filtration of Theorem~\ref{actual-defect-kernels} now gives \(D=0\).

The Kummer comparisons used in these two inputs can be made over
\(\mathbb C\). Extension from \(\overline{\mathbb Q}\) to \(\mathbb C\)
is injective on rational Chow groups: a rational equivalence spreads
over a finite-type integral base, and specialization at a closed
point of a common open gives an equivalence over \(\overline{\mathbb Q}\).
The injection persists after tensoring with \(\mathbb C\).
\end{proof}

\begin{remark}[The family comparison in low and high codimension]\leavevmode
\begin{enumerate}
\item \emph{Codimension two.}
Here the unmarked recurrence can be taken modulo
\(J_b^{\,\mathrm{fib}}\). On each labelled support the numerical
pairings determine the generic fibre class. The error contracts
to a homologically trivial divisor on \(E^2\) odd in both factors
and hence vanishes in \(\Pic(E^2)_{\mathbb Q}\) by
Proposition~\ref{restricted-divisor-graph-errors}. This gives the
untwisted input at all auxiliary levels. The marked spectral argument
then applies as above.
\item \emph{Higher codimension.}
The unmarked recurrence uses the Chow factorization
\eqref{med-g-factor} of Proposition~\ref{fixed-evaluation}, obtained
by cutting the reduced map space and including colliding positions.
\end{enumerate}
\end{remark}

\subsection{Modularity of the Abel--Jacobi images}

Write \(P_N\) for the parabolic cohomology of
\(\operatorname{Sym}^{m}R^1\pi_*\mathbb Q\) on \(Y_0(N)\), with
\(m=2\kappa-2>0\). It is the Scholl summand in degree \(m+1\)
of the compact Kuga--Sato variety. Write \(P_{N,\ell}\) for its
\(\ell\)-adic realization. Fix an embedding
\(\overline{\mathbb Q}\hookrightarrow\mathbb C\). The rational
Abel--Jacobi targets, before extending scalars, are
\[
 \begin{split}
 J_{\rm H}(N)&=\operatorname{Ext}^1_{\rm MHS}
                      (\mathbb Q(0),P_N(\kappa)),\\
 J_\ell(N)&=\varinjlim_K H^1_{\rm cont}(G_K,P_{N,\ell}(\kappa)).
 \end{split}
\]
The limit uses restriction over number fields in \(\overline{\mathbb Q}\).
We view each target as a rational vector space. Write \(J_*(N)_{\mathbb C}\)
for its tensor product over \(\mathbb Q\) with \(\mathbb C\), and
\({\rm AJ}_*\) for the corresponding complex-linear extension of
the Abel--Jacobi map.

\begin{theorem}[Abel--Jacobi modularity]\label{general-aj-modularity}
For every cyclic level \(N\), \(\kappa\ge2\), compatible fundamental
discriminant \(\Delta\) and root \(r\), the series
\[
       \sum_{\mu,n>0} q^n e_\mu\otimes
              {\rm AJ}_*\bigl(Z_\kappa(n,\mu)\bigr)
                \ \in\ S\otimes J_*(N)_{\mathbb C}
\]
is a cusp form, for the Hodge realization and for each prime \(\ell\).
In particular, every weakly holomorphic input satisfies
\({\rm AJ}_*(Z_\kappa(f))=0\) in each of these targets.
\end{theorem}
\begin{proof}
By the principal-part criterion (Theorem~\ref{rd-main}), it suffices
to show that \({\rm AJ}_*\) kills \(D_{\Delta,r}(N)\).

\smallskip\noindent\emph{The untwisted diagonal.}
For \(\Delta=r=1\), Proposition~\ref{diagonal-hecke-defect} expresses
the relations on the smooth base as \(R_i(TH_b(1))\), where
\(T\in I_b^{\rm cusp}\) acts on the target of the signed diagonal.
On a fine cover \(Y\), let \(U_b\) be the \(b\)-fold power of the
universal elliptic curve. Taking the odd part in each coordinate
and alternating the coordinates leaves only
\(\operatorname{Sym}^b R^1\pi_*\) in relative degree \(b\).
Since \(b>0\), finite-index \(\mathrm{SL}_2(\mathbb Z)\) monodromy
excludes global sections. Affineness of \(Y\) excludes cohomology
above degree one. Thus the signed cohomology of \(U_b\) occurs only
in degree \(b+1\), and that of its square, with both projectors,
only in degree \(2b+2\).

Take level-group invariants on both factors. By Poincar\'e duality
the diagonal represents \(H_c^{b+1}\to H^{b+1}\), whose image is
parabolic cohomology. The multiplier \(A^{b+1}/D\) for \(AD=n\)
gives the usual \(T_n\) on cusp forms of weight \(b+2\).
Eichler--Shimura therefore gives \(\operatorname{cl}(TH_b(1))=0\):
\(T\) annihilates the holomorphic summand, and \(I_b^{\rm cusp}\)
is the complexification of the annihilator on the rational Fourier
structure. This kills the conjugate summand as well and, by
comparison, the \(\ell\)-adic realization.

\smallskip\noindent\emph{The compact Abel--Jacobi class.}
Work over a number field of definition of the correspondences and
cycles involved. In total degree \(2b+2\), concentration of signed
cohomology leaves only the column of arithmetic degree zero in the
absolute Hodge and Hochschild--Serre spectral sequences
\cite{BeilinsonAbsoluteHodge,JannsenContinuous}. The geometric
vanishing thus kills the absolute class of \(TH_b(1)\). Pullback,
diagonal intersection, contraction and proper pushforward commute
with the cycle class maps, so \(R_i(TH_b(1))\) has zero absolute
class on the smooth locus.

On the compact Scholl summand cohomology occurs only in degree
\(m+1\), so its absolute class in degree \(m+2\) with twist
\(\kappa\) is the Abel--Jacobi class. Restriction to the smooth
locus is injective on both Abel--Jacobi targets. Indeed, after
splitting the cusps, the quotient of open by compact cohomology
in degree \(m+1\) embeds, after twist, in a sum of copies of
\(\mathbb Q(1-\kappa)\), also at cyclic level. Before the twist
the boundary lies in a sum of copies of \(\mathbb Q(-m-1)\)
\cite[\S1.3.4]{SchollMotives}. This boundary has no Hodge classes
of type \((0,0)\) and no \(\ell\)-adic Galois invariants.
The extension and Galois cohomology exact sequences give the
injections. In the \(\ell\)-adic case corestriction gives injectivity
under the finite field extension as well. This proves the untwisted
vanishing at every level.

\smallskip\noindent\emph{Transfer of a twist.}
Functoriality and untwisted vanishing kill the right side of
\eqref{md-t-inclusion} after applying \({\rm AJ}_*\).
For \(p\nmid h\) and \(d=\kappa-1\), the Hecke bound on \(P_h\)
in weight \(2d+2\) gives
\(\lvert\lambda_p\rvert\le 2p^{d+1/2}<(p+1)p^d\)
\cite{DeligneWeil}. The rational characteristic polynomial of
\(T_p\) therefore gives polynomial inverses to
\(T_p\pm(p+1)p^d\) on both realizations and their Abel--Jacobi
targets, including \(\ell=p\). At a shared prime the coefficients
are \(W\)-fixed (Proposition~\ref{md-transfer-proof}), and trace
embeds the exceptional kernel in the good-prime kernel at \(N/p\)
by \eqref{shared-edge-kernel-operators}. Each factor preserves the
required fixed spaces. Their product is therefore injective on
the images in question, proving \({\rm AJ}_*(D_{\Delta,r}(N))=0\).
\end{proof}

\section{Heegner divisors and the arithmetic inner product}\label{rd-h-lift}

For \(\kappa=1\), the coefficients of the weight-\(3/2\) series
are divisors on the modular curve. Subtracting a multiple of a
rational cusp places each coefficient in the Jacobian. Write
\(J=J_0(N)\), \(K_\Delta=\mathbb Q(\sqrt{\Delta})\), and let
\(P_\infty\) denote the cusp at infinity. Put
\begin{align*}
 y_{\Delta,r}(n,\mu)
   &=[Z_1(n,\mu)-\deg Z_1(n,\mu)P_\infty]
     \in J(K_\Delta)\otimes_{\mathbb Z}\mathbb C,\\
 A_{\Delta,r}(\tau)
   &=\sum_{\mu,\,n>0}y_{\Delta,r}(n,\mu)q^n e_\mu .
\end{align*}
For \(\Delta\ne1\) the genus character already gives each divisor
\(Z_1(n,\mu)\) degree zero. On the coarse curve we count each
oriented vector orbit in \eqref{rd-coeff} with multiplicity
\(2/|\Gamma_x|\). Bruinier and Ono count each such orbit with
multiplicity \(1/|\Gamma_x|\) in \cite[(5.1)]{BO}, so our series is twice
theirs. We use the corresponding height convention of
\cite[(4.13), Proposition~7.2]{BruinierYangCM}
and \cite[\S2]{DY2}.

\begin{theorem}[The divisor series]\label{med-divisors}
The series \(A_{\Delta,r}\) belongs to
\(S_{3/2,\widetilde\rho_L}\otimes_{\mathbb C}
 (J(K_\Delta)\otimes_{\mathbb Z}\mathbb C)\).
For a harmonic input \(F\) of weight \(1/2\) with cuspidal shadow,
\begin{equation}\label{merge-div-shadow}
 y_{\Delta,r}(F):=
   \sum_{\mu,\,n>0}c_F^+(-n,\mu)y_{\Delta,r}(n,\mu)
   =(A_{\Delta,r},\xi_{1/2}F)_{\rm Pet}.
\end{equation}
In particular every weakly holomorphic input gives zero in the
complexified Jacobian.
\end{theorem}
\begin{proof}
The modularity assertion follows from \cite[Corollary~6.7]{BO}
after reversing the sign of the quadratic
form and doubling the coefficients. For the second assertion,
apply a linear functional to the Jacobian factor and call the
resulting cusp form \(a\). By Stokes' theorem on a truncated
fundamental domain, \((a,\xi_{1/2}F)_{\rm Pet}\) equals the limit
of the integral of the componentwise contraction of \(a\) and \(F\)
along the upper edge. The nonholomorphic part of \(F\) tends to
zero in this integral. Integrating horizontally then extracts the
constant term of the contraction of \(a\) with \(F^+\). Since
\(a\) is cuspidal, only the negative powers of \(q\) in \(F^+\)
contribute. This proves \eqref{merge-div-shadow} for every linear
functional and hence in the complexified Jacobian. For weakly holomorphic inputs
with coefficients in a number field, twisted Borcherds products
give the divisor relation directly \cite[Corollary~6.6]{BO}.
\end{proof}

\subsection{The admissible lift}
At every level the Petersson product
\(\theta_{\Delta,r}(g)=(A_{\Delta,r},g)_{\rm Pet}\) defines a class in
the complexified Jacobian. To lift it to the arithmetic Chow group, we must
specify an integral model and Green currents. Take first
\(\Delta=r=1\) and \(N\) squarefree, and use the regular modular
stack \(\mathcal X_0(N)\) of \cite{DY1}. We use the effective
stack degree \(2/|\operatorname{Aut}(x)|\) of \cite[\S6.2]{DY1}
and the absolute arithmetic degree with factor \(1/2\) at infinity.
The arithmetic pairing extends the Gillet--Soul\'e intersection to
the logarithmic Green currents at the cusps \cite{GS,BKK,Kuhn2}. In this
normalization the Faltings--Hriljac formula identifies the
N\'eron--Tate pairing on degree-zero classes with the negative of
the arithmetic intersection of their admissible extensions
\cite[\S5]{BruinierYangCM}.

Let \(\widehat\omega_N\) be the metrized Hodge bundle of
\cite{DY1}. For \(y\in J(\mathbb Q)\otimes\mathbb R\), let
\(i_{\omega}(y)\) denote its extension of zero curvature that is
orthogonal to every vertical component and to \(\widehat\omega_N\).
The intersection matrix of the components determines the vertical
correction up to whole fibres, and the Green equation determines
the archimedean correction up to a constant. In the arithmetic Chow
group a whole fibre is equivalent to a constant Green function.
Since \(\deg\omega_N>0\), orthogonality to \(\widehat\omega_N\)
fixes the remaining constant. We extend \(i_{\omega}\)
complex-linearly.

Write \(a(f)=(0,f)\) for a pure Green class and put \(v=\Im\tau\).
For \(z\in\mathbb H\), let \(x_-\) be the projection of
\(x\in L\otimes\mathbb R\) onto the negative plane perpendicular
to the positive line represented by \(z\). Put
\(R(x,z)=-(x_-,x_-)\) and
\(\beta_1(t)=\int_1^\infty e^{-t u}\,du/u\).
The Kudla Green function for the untwisted coefficient is
\[
 \Xi(n,\mu,v)(z)=
     \sum_{\substack{x\in L+\mu,\ Q(x)=n\\x\ne0}}
                  \beta_1\bigl(2\pi v R(x,z)\bigr).
\]
For \(n>0\) its singularities are the CM divisor \(Z_1(n,\mu)\).
When \(-4Nn>0\) is a square, its logarithmic terms at the cusps
give the cusp divisor. In the zero coefficient one must also
correct the log-log term. More precisely,
\begin{equation}\label{p19a-zero-divisor}
 \widehat{\mathcal Z}(0,0,v)
 =\widehat{\mathcal Z}(0,0,v)^{\rm naive}
       -2\widehat\omega_N-\sum_{p\mid N}\mathcal X_p^0
       -a\bigl(\log(v/N)\bigr).
\end{equation}
Here the naive class comprises the cusp divisor prescribed by the
Green function and \(\Xi(0,0,v)\), and \(\mathcal X_p^0\)
is the component of the fibre at \(p\) containing the cusp \(0\).
The two components at \(p\mid N\) are denoted by
\(\mathcal X_p^0,\mathcal X_p^\infty\).
The resulting weight-\(3/2\) arithmetic theta series is
\[
 \widehat\phi(\tau)=
 \sum_{\mu\in L^\sharp/L}\,\sum_{n\in Q(\mu)+\mathbb Z}
                 \widehat{\mathcal Z}(n,\mu,v)q^n e_\mu .
\]
Its curvature, including the contribution of the zero vector,
is the Kudla--Millson theta form \(\Theta^{\rm KM}(\tau,z)\)
\cite[Theorem~1.1, (6.11) and (8.4)]{DY1}.

\begin{theorem}[Projection of the arithmetic theta series]
\label{p19a-lift}
Suppose \(N\) is squarefree and \(\Delta=r=1\).
For \(g\in S_{3/2,\rho_L}\) one has
\begin{equation}\label{rd-h-proj}
 \widehat\theta(g):=(\widehat\phi,g)_{\rm Pet}
          =i_\omega\bigl(\theta_{1,1}(g)\bigr).
\end{equation}
In particular its arithmetic self-intersection is
\(-\|\theta_{1,1}(g)\|_{\rm NT}^{2}\) for the Hermitian extensions
of the pairings.
\end{theorem}
\begin{proof}
To identify \(\widehat\theta(g)\) with the normalized admissible
extension, we determine its generic class, its pairings with the
fibre components and \(\widehat\omega_N\), and its curvature.

\smallskip\noindent\emph{The fibre and Hodge pairings.}
Let \(\mathcal E_L(\tau,s)\) be the normalized vector-valued
Eisenstein series of \cite[(1.6)]{DY1}, and put
\(\mathcal E=\mathcal E_L(\tau,1)\),
\(\mathcal E'=\partial_s\mathcal E_L(\tau,1)\).
With \(\varphi\) denoting Euler's totient, the arithmetic
Siegel--Weil identities \cite[Theorem~1.3]{DY1} give
\begin{align}
  \deg\widehat\phi&=\frac{2}{\varphi(N)}\mathcal E,\notag\\
  \langle\widehat\phi,\mathcal X_p^0\rangle_{\rm GS}
    =\langle\widehat\phi,\mathcal X_p^\infty\rangle_{\rm GS}
        &=\frac{\log p}{\varphi(N)}\mathcal E
           \qquad(p\mid N),\label{p19a-vertical}\\
  \langle\widehat\phi,\widehat\omega_N\rangle_{\rm GS}
    &=\frac{1}{\varphi(N)}
       \left(\mathcal E'-\sum_{p\mid N}
                  \frac{p\log p}{p-1}\mathcal E\right).\notag
\end{align}
Unfold \((\mathcal E_L(\cdot,s),g)_{\rm Pet}\) in the domain of
absolute convergence. It vanishes because the constant term of
\(g\) is zero. Analytic continuation gives vanishing at \(s=1\)
and for the derivative there. Thus \(\widehat\theta(g)\) has
degree zero and pairs to zero with \(\widehat\omega_N\) and both
components at every \(p\mid N\). At a good prime its degree also
gives zero intersection with the fibre.

\smallskip\noindent\emph{The Green function after projection.}
Fix a smooth positive probability measure \(\nu\) on
\(X=X_0(N)(\mathbb C)\) invariant under complex conjugation
\(c\). One obtains such a measure by averaging any smooth positive
probability measure with its pushforward by \(c\). Its symmetric
Green kernel satisfies
\[
  dd_z^c[\mathcal G_\nu(w,z)]+\delta_w=\nu,\qquad
  \int_X\mathcal G_\nu(w,z)\nu(z)=0.
\]
Subtract the \(\nu\)-admissible Green function of each horizontal
divisor from its completed Green function, and let \(s(\tau,z)\)
be the mean-zero part of the difference. It satisfies
\[
       dd_z^c s=\Theta^{\rm KM}(\tau,z)-\deg(\widehat\phi)\nu,
       \qquad \int_Xs(\tau,z)\nu(z)=0 .
\]
The curvature has integral \(\deg\widehat\phi\). The Green kernel
therefore gives the unique mean-zero solution
\begin{equation}\label{p19a-green-inversion}
       s(\tau,z)=-\int_X\mathcal G_\nu(w,z)\,
                              \Theta^{\rm KM}(\tau,w).
\end{equation}
Indeed the right-hand side has the required \(dd_z^c\) and zero
mean. Equivalently, against a mean-zero real eigenfunction \(f\)
with \(dd^cf=-(\lambda/2)f\nu\), integration by parts gives
\[
         \int_X s(\tau,z)f(z)\nu(z)=
               -\frac{2}{\lambda}\int_X f(z)\Theta^{\rm KM}(\tau,z).
\]
These current identities hold on the compact curve. The cusp and
index-zero corrections of \eqref{p19a-zero-divisor} supply the
boundary terms \cite[\S8]{DY1}.

The curvature itself has zero Petersson pairing with \(g\) at
each \(z\). Divide \((\Theta^{\rm KM}(\cdot,z),g)_{\rm Pet}\) by
the invariant area form. The Casimir identity and integration by
parts on a truncated fundamental domain show that the quotient
is harmonic in \(z\), since the cusp decay of \(g\) kills the
boundary integrals. The kernel estimate gives limit zero at
each cusp, so the harmonic function vanishes
\cite[proof of Theorem~5.1, \(k=0\)]{AlfesKM}. Consequently the
projected mean-zero Green correction in \eqref{p19a-green-inversion}
vanishes. The pairing with \(\widehat\omega_N\) in \eqref{p19a-vertical}
fixes the remaining Green constant.

\smallskip\noindent\emph{The admissible extension.}
On the generic fibre the divisor of the modular discriminant
section used in \cite[Lemma~6.4, (6.6)]{DY1} is supported at
\(P_\infty\). Thus the Hodge term of \eqref{p19a-zero-divisor}
has zero degree-zero projection, as do all cusp divisors by
Manin--Drinfeld. The remaining coefficients are
\(Z_1(n,\mu)-\deg Z_1(n,\mu)[P_\infty]\) for \(n>0\).
This is \(A_{1,1}\) in the complexified Jacobian
\cite[proof of Theorem~8.4]{DY1}. Its pairing with
\(g\) is \(\theta_{1,1}(g)\). The component pairings, Hodge pairing
and curvature calculated above characterize
\(i_\omega(\theta_{1,1}(g))\), proving \eqref{rd-h-proj}.
The self-intersection formula follows from Faltings--Hriljac.
\end{proof}

For any \(N\) and either sign of \(\Delta\), the Mordell--Weil
component of the arithmetic lift has a direct construction. Put
\(K=\mathbb Q(\sqrt\Delta)\). On a regular proper model over
\(\mathcal O_K\), choose a metrized line bundle of positive generic
degree and use it to normalize the admissible extension \(i_\omega\).
Then \(i_\omega(A_{\Delta,r})\) is a cusp form and its Petersson
projection is \(i_\omega(\theta_{\Delta,r}(g))\).
Indeed \(i_\omega\) is linear on degree-zero divisor classes and
Theorem~\ref{med-divisors} gives the assertion on the generic fibre.
Faltings--Hriljac computes the normalized self-intersection of this
lift at every level.

Let \(N\) be squarefree and let \(\Delta\ne1\) be fundamental, with
\(r^2\equiv\Delta\bmod 4N\). Over \(K=\mathbb Q(\sqrt\Delta)\),
form the regular model \(\mathcal X_K\) by blowing up the singular
points of \(\mathcal X_0(N)\otimes\mathcal O_K\) above the common
primes of \(N\) and \(\Delta\), and pull back \(\widehat\omega_N\).
For \(n>0\) take the horizontal closure of the signed divisor on
this model. The horizontal part is zero for \(n\le0\). At one
embedding of \(K\) use the twisted Kudla Green function
\(\Xi_{\Delta,r}(n,\mu,v)\) of \cite[\S4]{DY2}, in the orbit
convention of \eqref{rd-coeff}, and at the conjugate embedding use its
negative. Write \(\widehat\phi_{\Delta,r}^K\) for the resulting series.

The nontrivial automorphism of \(K/\mathbb Q\) negates the divisor
by \cite[Lemma~5.1]{BO}. For imaginary \(K\) the Green functions
also satisfy the required complex conjugation law. Reflection
\(z\mapsto-\bar z\) carries \([A,B,C]\) to \([A,-B,C]\), which
preserves the genus character and changes \(\mu\) to \(-\mu\).
Negation of the vector then multiplies the genus character by
\(\operatorname{sgn}(\Delta)\) and preserves the majorant.
Thus \(\Xi_{\Delta,r}(n,\mu,v)(-\bar z)
=\operatorname{sgn}(\Delta)\Xi_{\Delta,r}(n,\mu,v)(z)\).

\begin{proposition}\label{prop:resolved-signed-lift}
Suppose that \(\Delta>1\), or that \(\Delta<0\) and every \(p\mid N\)
ramifies in \(K\). For every \(g\in S_{3/2,\widetilde\rho_L}\),
\[
  (\widehat\phi^K_{\Delta,r},g)_{\rm Pet}
        =i_\omega\bigl(\theta_{\Delta,r}(g)\bigr).
\]
\end{proposition}
\begin{proof}
Away from \(N\) the model is smooth with geometrically irreducible
fibres, so the degree-zero horizontal divisor has zero intersection
with each fibre. If \(p\mid N\) is unramified in \(K\), the square
congruence implies that \(p\) splits. This case occurs here only for
\(\Delta>1\). Above each of the two places, Fricke interchanges the
two components and fixes the signed divisor
\cite[Lemma~5.1(i),(iv)]{BO}. The component intersections are equal
and their sum is zero.

Now let \(p\mid(N,\Delta)\). On an \'etale chart with auxiliary
prime-to-\(p\) level, the two reduced branches of the original fibre
meet transversely at each supersingular point
\cite[Theorem~13.4.7, \(n=1\)]{KatzMazurModuli}. Regularity gives the
completed local equation \(xy=p\) after strict henselization.
Over the ramified quadratic base it becomes \(xy=u\pi^2\) for a
unit \(u\) and a uniformizer \(\pi\). Blowing up the closed point
resolves it. The exceptional curve is the smooth conic
\(XY=\bar u T^2\), also in characteristic two. Inertia fixes the
original components and every geometric node, hence the unique
exceptional component over each node. As it negates the horizontal
divisor, its intersection with each of these components vanishes.

At each embedding the curvature is the signed Kudla--Millson kernel,
with its current equation on the compact curve
\cite[Theorem~4.4 and Corollary~4.5]{DY2}. Its Petersson pairing
with \(g\) vanishes pointwise \cite[proof of Theorem~5.1,
\(k=0\)]{AlfesKM}. Subtract \(i_\omega(A_{\Delta,r})\) coefficientwise.
Modulo pullbacks from \(\operatorname{Spec}\mathcal O_K\), the
component intersections just computed identify the difference with
the mean-zero Green inverses of the curvatures as in
\eqref{p19a-green-inversion}. Use the same conjugation-invariant
measure at both embeddings. By uniqueness of the mean-zero
normalization \(G_\nu(cw,cz)=G_\nu(w,z)\), so the inverses respect
the real structure. They transform modularly and their Petersson
pairings vanish. In the arithmetic Chow group with real coefficients
the base class is determined by its degree, by finiteness of the ideal
class group and Dirichlet's unit theorem. The signed coefficients,
the normalized admissible extensions and the mean-zero Green functions
are anti-invariant, whereas Galois fixes the degree. The remaining
base class is therefore zero.
\end{proof}

For comparison, Du and Yang use a different metric convention for
the twisted arithmetic series. At squarefree level they prove
weight-\(3/2\) modularity for both signs and compute the Hodge pairing
in their convention \cite[Theorems~1.2 and~1.6]{DY2}:
\[
 \left\langle\widehat\phi_{\Delta,r},\widehat\omega_N\right\rangle_{\rm GS}
 =\begin{cases}
 \dfrac{h(\Delta)\log\varepsilon_\Delta}{\varphi(N)}
                  \mathcal E_L(\tau,1),&\Delta>1,\\[4pt]
 0,&\Delta<0 .
 \end{cases}
\]
Here \(h(\Delta)\) and \(\varepsilon_\Delta>1\) are the class
number and the fundamental unit of the real quadratic field in the
first case. These Hodge pairings have zero cuspidal projection.
For the signed series over \(K\) the absolute Hodge pairing
vanishes already coefficientwise: Galois negates the series and
fixes the pulled-back Hodge bundle. The same argument makes the
sum above a rational prime vanish for a component defined over the
rational base. At a fixed split place the two intersections can
nevertheless be nonzero. We compute them on the resolved model over \(K\).

\begin{proposition}[The component series for a negative twist]
\label{unary-split-components}
Let \(\Delta<0\), let \(N\) be squarefree and let \(p\mid N\)
split in \(K=\mathbb Q(\sqrt\Delta)\). For \(Q\mid N\), let
\(\sigma_Q\) act on residues by \(-1\bmod 2Q\) and \(+1\bmod 2N/Q\).
With \(h_K=|\operatorname{Pic}(\mathcal O_K)|\) and
\(w_K=|\mathcal O_K^\times|\), put
\[
 \vartheta'_{N,\mu}(\tau)
      =\sum_{x\equiv\mu\ (2N)}x\,q^{x^2/(4N)},\qquad
 V_{p,\mu}=\sum_{Q\mid N}(-1)^{\mathbf 1_{p\mid Q}}\,
                              \vartheta'_{N,\sigma_Q\mu}.
\]
At either place \(v\mid p\), the two series of geometric
intersection degrees of the horizontal divisors with the fibre
components are \((h_K/w_K)V_p\) and \(-(h_K/w_K)V_p\).
The signs interchange at the conjugate place. Degrees use the orbit
weights \(2/|\Gamma_x|\), before multiplication by \(\log p\).
At \((n,\mu)=(1/(4N),1)\) the degrees are \(2h_K/w_K\) and its negative.
\end{proposition}
\begin{proof}
\smallskip\noindent\emph{The CM field.}
The horizontal divisor is zero for \(n\le0\). For \(n>0\) put
\(K'=\mathbb Q(\sqrt{-4N|\Delta|n})\). If \(K'\ne K\), choose
a prime \(\ell\) split in \(K'\), inert in \(K\) and prime to
\(2N\Delta\) and to every CM conductor in the coefficient.
On each oriented CM order a prime ideal over \(\ell\) permutes
the marked CM pairs by isogenies of degree \(\ell\). This preserves
the residue and stabilizer and changes the genus sign by
\(\chi_\Delta(\ell)=-1\), by CM reciprocity and the determinant
law of Section~\ref{p19-translated}. Compare the endpoints at the
same place of a common field of definition. Both degeneracy maps
of this prime-to-\(p\) correspondence pull back each vertical
component to the same divisor: on the dense ordinary locus the
isogeny preserves the connected or \'etale type of the \(p\)-subgroup.
The maps are finite and generically \'etale there, so the Cartier
pullbacks agree with multiplicities, including at a node. The
permutation therefore negates the intersection with either component
\cite[proof of Proposition~4.4.2]{ZhangShimuraPoints}.

\smallskip\noindent\emph{The square indices.}
If \(K'=K\), then \(4Nn=t^2\) for an integer \(t>0\) with
\(t^2\equiv\mu^2\bmod4N\). For \(p\mid t\), one has
\(\sigma_p\mu=\mu\). The involution \(W_p\) fixes the coefficient
and exchanges the two components. Their sum has degree zero, so
both degrees are zero, as are the corresponding terms of \(V_{p,\mu}\).
Assume \(p\nmid t\). The reduction of a CM elliptic curve is then
ordinary, and the two conjugate \(p\)-level ideals give the
connected and the \'etale lines. Fix a component and a place.
Relative to \(t=\mu=1\), the positive definite branch has the same
local orientation for \(\mu\equiv t\bmod p\) and the opposite
orientation for \(\mu\equiv-t\bmod p\); at \(p=2\) use modulus \(4\).
Negation reverses the genus sign and exchanges the two branches.

Write the positive quadratic form of the source elliptic curve as
\([Na,b,c]=\delta[A,B,C]\) with \([A,B,C]\) primitive. The letters
\(a,b,c\) in this count are coefficients of the form: a matrix
with entries \(b_x/(2N),-a_x/N,c_x,-b_x/(2N)\) has fixed-point
equation \(Nc_xz^2-b_xz+a_x=0\).
The primitive discriminant is \(\Delta f^2\), where \(\delta f=t\).
At level one the mass of the proper classes is \(2h_f/w_f\).
Here \(h_f,w_f\) are the class number and number of units of the
order of conductor \(f\). The ring class number formula gives
\cite[Theorem~7.24]{CoxPrimes}
\[
 \frac{h_f}{w_f}=\frac{h_K}{w_K}\,
         f\prod_{\ell\mid f}(1-\chi_\Delta(\ell)/\ell).
\]
Multiplying by the number of allowed cyclic lines gives the mass
with stabilizers, since the units act on these lines. The factor
\(2/|\Gamma_x|\) enters the degree once. Equivalently, after
rigidification by \(\{\pm1\}\), the coarse point of multiplicity
\(2/|\Gamma_x|\) pulls back to the residual gerbe of multiplicity one.

\smallskip\noindent\emph{Counting the local lines.}
If \(q\mid N\) and \(q\nmid\delta\), the middle residue selects
one line. Indeed, with \(e\in\{0,1\}\) of the parity of \(B\), an
order generator on the underlying lattice has matrix
\[
 T=\begin{pmatrix}(e-B)/2&-C\\A&(e+B)/2\end{pmatrix}.
\]
Primitivity makes \(T\bmod q\) nonscalar. The allowed line is its
eigenspace for the specified root: in an adapted basis
\(A\equiv0\bmod q\), and its eigenvalue is \((e-B)/2\bmod q\).
At odd \(q\), first determine \(B\bmod q\) by dividing the prescribed
middle residue by \(\delta\), then use the parity of \(\Delta f^2\).
At \(q=2\), \(\delta\) is odd and division is modulo \(4\).
The discriminant congruence gives the required root in both cases.
If \(q\mid\delta\) and \(q\nmid\Delta\), all \(q+1\) lines are allowed.
For \(q=2\) this uses \(\delta B\equiv t\equiv r\mu\bmod4\),
from \(B\equiv f\bmod2\) and \(t^2\equiv\mu^2\bmod8\).

At a shared prime \(q\mid(N,\Delta,t)\) a nonzero genus weight
requires exactly one of \(a,c\) to be a unit. Indeed, \(q\mid b\)
and \(q\mid ac\); at \(2\), \(4\mid b\) and \(16\mid\Delta t^2\)
give \(2\mid ac\). The involution \(W_q\) interchanges the two
unit conditions, exchanging \(a,c\bmod q\) up to units. It
preserves the genus character \cite[\S4]{BO}, fixes the residue
and preserves the component at \(p\ne q\) by its prime-to-\(p\)
isogeny. At the other shared primes its matrix belongs to the local
level group and preserves the two unit conditions. Thus we may
count the branch \(q\nmid c\) twice at each
such \(q\). There \(q\nmid\delta\), so there is a single line.
For \(q\mid(N,\Delta)\) and \(q\nmid t\) the ramified maximal order
gives a single line without this doubling.

On the branches thus chosen, \((\delta,\Delta)=1\) for nonzero
weights. The primitive positive form represents an integer prime to
\(f\Delta\) that is an ideal norm in \(K\), so the genus weight is
\(\chi_\Delta(\delta)\). For \(q^a\Vert t\) outside \(N\), summing
the local class factors over \(\delta f=t\) with this character gives
\(q^a\). For \(q\mid N\), \(q\nmid\Delta\), compatibility gives
\(\chi_\Delta(q)=1\). When \(a>0\), terms with \(q\mid\delta\)
give \((q+1)q^{a-1}\); terms with \(q\nmid\delta\) give
\(q^a-q^{a-1}\). At a shared prime with \(a>0\), the class
factor is \(q^a\) and the two branches again give \(2q^a\).
The full positive mass is therefore \(2(h_K/w_K)t\,2^{\omega((t,N))}\).

For each of \(x=t,-t\), precisely \(2^{\omega((t,N))}\) divisors
\(Q\mid N\) satisfy \(x\equiv\sigma_Q\mu\bmod2N\). Complementary
divisors exchange the two choices and reverse the \(p\)-sign.
The two contributions to \(V_{p,\mu}\) therefore add, with the
orientation sign found above. This gives the stated series.
\end{proof}

If \(\Delta<0\) and every \(p\mid N\) ramifies, squarefreeness gives
\(N\mid r\). Each residue \(r\mu\) then equals its negative, and the
two vectors \(x,-x\) cancel in both the divisor and the Green function.

\begin{theorem}[The signed arithmetic theta lift at squarefree level]
\label{squarefree-signed-theta}
Let \(N\) be squarefree and let \(\Delta\ne1\) be a compatible
fundamental discriminant. The series \(\widehat\phi_{\Delta,r}^K\)
transforms with weight \(3/2\) and type \(\widetilde\rho_L\), in
the sense of currents at the infinite places. Write
\(\widehat\theta_{\Delta,r}^K(g)=(\widehat\phi_{\Delta,r}^K,g)_{\rm Pet}\).
For every \(g\) in a good Hecke eigenspace with cuspidal weight-two
Shimura correspondent,
\[
       \widehat\theta_{\Delta,r}^K(g)=
                      i_\omega\bigl(\theta_{\Delta,r}(g)\bigr).
\]
For \(\Delta<0\) and general \(g\in S_{3/2,\bar\rho_L}\) the
correction is explicit. For each split \(p\mid N\) choose a place
\(v\) and label the components \(C_v^+,C_v^-\) so that
\(C_v^+\) has positive intersection degree at \((1/(4N),1)\).
Transport these labels by Galois to \(\bar v\), so \(C_{\bar v}^+\)
has negative degree there. Put \(D_v=C_v^--C_v^+\)
and \(s_p=\deg_v(C_v^+\cdot C_v^-)\). Then
\[
 \widehat\theta_{\Delta,r}^K(g)=i_\omega(\theta_{\Delta,r}(g))
       +\frac{h_K}{w_K}\sum_{\substack{p\mid N\\p\ {\rm split}}}
              \frac{(V_p,g)_{\rm Pet}}{2s_p}\,(D_v-D_{\bar v}).
\]
In particular, equality with the admissible lift holds precisely
when \((V_p,g)_{\rm Pet}=0\) at every such \(p\).
\end{theorem}
\begin{proof}
Subtract the normalized admissible series \(i_\omega(A_{\Delta,r})\)
and the mean-zero Green inverse \(a(s)\) of its curvature difference,
as in \eqref{p19a-green-inversion}. The curvature is modular,
and the Green inverse is a continuous linear operator on mean-zero
currents on the compact curve. Consequently \(s\) transforms in
the same way and its Petersson projection against a cusp form
vanishes by the pointwise kernel calculation of
Proposition~\ref{prop:resolved-signed-lift}.

What remains has zero generic class and curvature. Subtracting a
principal divisor leaves a vertical divisor and a constant Green
function on each connected compact fibre. Modulo pullbacks from
the arithmetic base, the component intersections of
Proposition~\ref{prop:resolved-signed-lift} determine it except at
split \(p\mid N\) for negative \(\Delta\). At such a place
\(D_v\cdot C_v^+=2s_p\) and \(D_v^2=-4s_p\), since the whole
fibre has intersection zero with each component. By
Proposition~\ref{unary-split-components} the vertical term there
is \((h_K/w_K)V_pD_v/(2s_p)\), and the conjugate term has the
opposite sign. In \(\widehat{\operatorname{Pic}}(\mathcal O_K)_{\mathbb R}\)
arithmetic degree is an isomorphism by the unit theorem and
finiteness of the class group. Its Galois-anti-invariant part is
zero, which disposes of the remaining base class.
Thus the asserted decomposition holds coefficientwise before
projection, with \(a(s)\) added on the right.

The vector \(\vartheta'_N\) is the theta series of the positive
rank-one lattice of quadratic form \(Nx^2\), with the linear
polynomial insertion. Poisson summation gives its weight \(3/2\)
and its type \(\bar\rho_L\). Its constant coefficient vanishes.
The \(\sigma_Q\) commute with the Weil action, so the \(V_p\) are
cuspidal of this type. This proves the transformation law for
the signed arithmetic series. It also explains why the component
correction disappears on an entire cuspidal Shimura eigenspace.
For every prime \(\ell\equiv1\bmod 4N\), the good operator
\(T(\ell^2)\) acts on \(\vartheta'_N\) and its residue permutations
by \(1+\ell\), by the Fourier coefficient formula. On the eigenspace
of a weight-two cusp form \(G\) its eigenvalue satisfies
\(|a_G(\ell)|\le2\sqrt\ell<1+\ell\). Self-adjointness gives
\((V_p,g)_{\rm Pet}=0\). Conversely the differences \(D_v-D_{\bar v}\)
at distinct split primes are independent, as seen by intersection
with \(C_v^+\). The case \(\Delta>1\) follows directly from
Proposition~\ref{prop:resolved-signed-lift}.
\end{proof}

Two consequences record the normalizations. For a prime
\(\ell\nmid 2N\Delta\), put \(D_\ell=T_\ell-(\ell+1)\) on the
arithmetic curve. On the eigenspace of a weight-two newform \(G\)
one also has
\begin{equation*}\tag{\ref*{rd-h-proj}a}\label{eq:signed-hecke-lift}
 i_\omega(\theta_{\Delta,r}(g))=
 \frac{(D_\ell\widehat\phi_{\Delta,r}^K,g)_{\rm Pet}}
 {a_G(\ell)-(\ell+1)} .
\end{equation*}
Indeed the isogeny correspondence preserves zero curvature and
vertical orthogonality by pullback and trace. At a common ramified
prime, its two degeneracy maps are \'etale near the supersingular
nodes and have the same inverse image of the node locus, so both
extend to the chosen resolution by blowing up that inverse image.
The normalized extensions commute with \(T_\ell\) because the
remaining anti-invariant base class is zero. Good Hecke equivariance
on the Jacobian \cite[(7.3)]{BO} now gives the formula.

For negative \(\Delta\), using the Hermitian extensions and the
absolute arithmetic degree, the full inner product for \(g,g'\) is
\begin{equation}\label{vertical-arithmetic-pairing}
\begin{split}
 -\langle\widehat\theta_{\Delta,r}^K(g),\widehat\theta_{\Delta,r}^K(g')\rangle_{\rm GS}
  &=\langle\theta_{\Delta,r}(g),\theta_{\Delta,r}(g')\rangle_{\rm NT}\\
  &\quad+\left(\frac{h_K}{w_K}\right)^2
       \sum_{\substack{p\mid N\\p\ {\rm split}}}
       \frac{\log p}{s_p}\,(V_p,g)_{\rm Pet}\,
                          \overline{(V_p,g')_{\rm Pet}} .
\end{split}
\end{equation}
Indeed, distinct fibres are disjoint and the admissible extension
pairs to zero with each component. The two conjugate places contribute
equally, cancelling the divisor \([K:\mathbb Q]=2\) in the absolute
degree. The \(V_p\) for different \(p\) lie in different characters
of the residue involutions, negative at \(p\) and positive at the
other level primes. Each is nonzero by its first coefficient.
Thus the rank of the additional Hermitian form is exactly the
number of split primes dividing \(N\).

The size of each vertical direction can be computed without a CM
height. Let \(g_p=V_p/\|V_p\|_{\rm Pet}\). For any good
\(\ell\equiv1\bmod4N\), Hecke equivariance gives
\(T_\ell\theta_{\Delta,r}(g_p)=(1+\ell)\theta_{\Delta,r}(g_p)\).
All the eigenvalues on the Jacobian have absolute value at most
\(2\sqrt\ell\). Thus \(\theta_{\Delta,r}(g_p)=0\), and
\begin{equation}\label{unary-vertical-unit-norm}
 -\langle\widehat\theta_{\Delta,r}^K(g_p),\widehat\theta_{\Delta,r}^K(g_p)
   \rangle_{\rm GS}
       =2^{\omega(N)+2}\sqrt N\,(h_K/w_K)^2\log p .
\end{equation}
To verify the constant, sum the squares of the entries of \(V_p\)
at \(t^2/(4N)\), \(t>0\). The coefficient count above gives
\(2^{\omega(N)+2}t^2\,2^{\omega((t,N))}\) for \(p\nmid t\), and
zero otherwise. Unfolding the Petersson integral against the
full-level real analytic Eisenstein series at parameter \(z\)
therefore gives the Dirichlet factor
\[
  \zeta(2z-1)(1-p^{1-2z})
       \prod_{q\mid N,\ q\ne p}(1+q^{1-2z}) .
\]
Its residue at \(z=1\), multiplied by
\(2^{\omega(N)+2}(4N)^{3/2}\Gamma(3/2)/(4\pi)^{3/2}\)
and divided by the Eisenstein residue \(3/\pi\), is
\[
 \|V_p\|_{\rm Pet}^2
   =\frac{2^{\omega(N)}N^{3/2}}{3}
     (1-p^{-1})\prod_{q\mid N,\ q\ne p}(1+q^{-1}).
\]
In the effective stack convention,
\(s_p=(p-1)[{\rm SL}_2(\mathbb Z):\Gamma_0(N/p)]/12\)
\cite[\S7, (7.1) and the following computation]{DY1}.
Their ratio is \(2^{\omega(N)+2}\sqrt N\), proving
\eqref{unary-vertical-unit-norm}.

\subsection{The height on a Shimura newform line}
Let \(G=\sum_{j\ge1}a_G(j)q^j\) be a newform of weight two,
level \(N\) and trivial character, normalized by \(a_G(1)=1\).
Let \(g\in S_{3/2,\widetilde\rho_L}\) be the Shimura newform
with real coefficients chosen in Theorem~\ref{med-inner}.
The sign of \(\Delta\) determines the Weil representation of \(g\).
Both Petersson norms use the usual weight factor and the measure
\(du\,dv/v^2\), without division by the volume. We integrate \(G\)
over \(\Gamma_0(N)\backslash\mathbb H\) and \(g\) over
\(\mathrm{SL}_2(\mathbb Z)\backslash\mathbb H\), using the
standard Hermitian product on the components of \(g\).
Choose \(F_g\) with
\(\xi_{1/2}F_g=g/\|g\|_{\rm Pet}^2\).
By \eqref{merge-div-shadow},
\(\theta_{\Delta,r}(g)=\|g\|_{\rm Pet}^2 y_{\Delta,r}(F_g)\).
On this line the \(G\)-isotypical projection of
\(A_{\Delta,r}\) equals \(g\otimes y_{\Delta,r}(F_g)\)
\cite[Theorem~7.6]{BruinierYangCM}, \cite[Theorem~7.7]{BO}. Thus if \(b=c_g(n,\mu)\), then
\begin{equation}\label{med-proj}
 y_{\Delta,r}^{G}(n,\mu)
       =\frac{b}{\|g\|_{\rm Pet}^2}\,\theta_{\Delta,r}(g).
\end{equation}

For an integer \(j\), let \(\chi_\Delta(j)=(\frac{\Delta}{j})\)
be the Kronecker character, with \(\chi_1=1\). We use the
classical normalization
\(L(G,\chi_\Delta,s)=\sum_{j\ge1}a_G(j)\chi_\Delta(j)j^{-s}\)
in the half-plane of absolute convergence. The prime in
\(L'(G,\chi_\Delta,1)\) denotes the derivative of its analytic
continuation with respect to \(s\) at the centre \(s=1\).

\begin{proof}[Proof of Theorem~\ref{med-inner}]
Extend \(\theta_{\Delta,r}(g)\) admissibly to a regular proper
model over the integers of a field of definition. After division
by the field degree, Faltings--Hriljac identifies its negative
arithmetic self-intersection with its N\'eron--Tate height. This
proves the first equality in \eqref{rd-h-nt}. At squarefree level
and \(\Delta=r=1\), \eqref{rd-h-proj} identifies the extension
defined by the arithmetic theta series.

\smallskip\noindent\emph{The untwisted height.}
For \(\Delta=r=1\), Bruinier--Yang give
\[
 \|y_{1,1}(F_g)\|_{\rm NT}^2
   =\frac{2\sqrt N}{\pi\|g\|_{\rm Pet}^2}\,L'(G,1).
\]
In \cite[Theorems~7.6 and~7.7]{BruinierYangCM}, \(g\) and the
principal part of its preimage are normalized over the Hecke field.
By \eqref{merge-div-shadow} the class depends only on the shadow.
Since the Petersson lift is \(\|g\|_{\rm Pet}^2 y_{1,1}(F_g)\),
its squared height has the additional factor \(\|g\|_{\rm Pet}^4\).
This proves \eqref{rd-h-nt} for every real normalization of \(g\).
For any other root of \(r^2\equiv1\bmod4N\), the isometry
\(\sigma_r:e_\mu\mapsto e_{r\mu}\) gives
\(A_{1,r}=\sigma_r A_{1,1}\) by the orbit definition.
It is a Petersson isometry commuting with the good Hecke operators,
so \(\sigma_rg\) belongs to the same exact-level newform line.
Pairing \(A_{1,1}\) with \(\sigma_rg\) therefore gives the same
self-height.

\smallskip\noindent\emph{An auxiliary coefficient for the twist.}
Put \(\sigma=\operatorname{sgn}(\Delta)\). The congruence
\(r^2\equiv\Delta\bmod4N\) and \((N,\Delta)=1\) give
\(\chi_\Delta(N)=1\). The hypothesis on the root number then gives
\(\varepsilon(G)\chi_\Delta(-N)=-1\), hence \(\varepsilon(G)=-\sigma\).
By Bump, Friedberg and Hoffstein, choose a fundamental discriminant
\(d\) of sign \(-\sigma\) such that \(|\Delta d|>4\),
every prime dividing \(2N|\Delta|\) splits in \(\mathbb Q(\sqrt d)\),
and \(L(G,\chi_d,1)\ne0\)
\cite[Theorem~(ii)]{BumpFriedbergHoffstein}. For a root
\(\mu^2\equiv d\bmod4N\), put \(n=|d|/(4N)\) and \(b=c_g(n,\mu)\).
Then \eqref{med-wald} gives \(b\ne0\), also when \(\Delta=1\).

For the twisted height take \(\Delta\ne1\). The coprime fundamental
discriminants \(d,\Delta\) give a fundamental discriminant \(D=\Delta d<-4\),
and every level prime splits in \(K=\mathbb Q(\sqrt D)\).
The height and Fourier coefficient at this index satisfy
\begin{align}
 \|y_{\Delta,r}^{G}(n,\mu)\|_{\rm NT}^2
  &=\frac{\sqrt{|\Delta d|}}{4\pi^2\|G\|_{\rm Pet}^2}
       L'(G,\chi_\Delta,1)L(G,\chi_d,1),\label{med-gz}\\
 |b|^2
  &=\frac{\|g\|_{\rm Pet}^2\sqrt{|d|}}
              {8\pi\sqrt N\|G\|_{\rm Pet}^2}\,L(G,\chi_d,1).
       \label{med-wald}
\end{align}
The Waldspurger formula \eqref{med-wald} in these Petersson norms
is stated for both signs of \(d\) in the published proof of
\cite[Theorem~7.8]{BO}, with the factor \(1/(8\pi\sqrt N)\) for
the vector-valued Petersson norm.
To compute the height, let \(\eta\) be the genus character of
\(\operatorname{Pic}(\mathcal O_K)\) attached to \(D=\Delta d\).
Write \(P_\beta\) for the \(\eta\)-weighted sum of CM points
represented by positive definite forms \([A,B,C]\) of discriminant
\(D\) with \(N\mid A\) and \(B\equiv\beta\bmod 2N\).
Each ideal class occurs once, including when \(D\) is even
\cite[I, \S1, p.~505]{GKZ}. The character is nontrivial, so
\(P_\beta\) has degree zero.

For a vector \(x\) with entries \(b_x/(2N),-a_x/N,c_x,-b_x/(2N)\),
the fixed point satisfies \(Nc_xz^2-b_xz+a_x=0\). The terms with
\(c_x>0\) therefore give \(P_\beta\) for \(\beta=-r\mu\bmod 2N\).
Their stabilizer weight in \eqref{rd-coeff} is one since \(D<-4\).
Negating \(x\) contributes the other branch with weight
\(\sigma=\operatorname{sgn}(\Delta)\). The Fricke involution
\(w_N\) sends \(P_\beta\) to \(P_{-\beta}\), preserving the genus
character \cite[I, \S2, Proposition~1, P2]{GKZ}. Its eigenvalue
on the \(G\)-component is \(-\varepsilon(G)=\sigma\). Hence
\[
 y_{\Delta,r}^G(n,\mu)=(P_\beta+\sigma P_{-\beta})^G
                     =2P_\beta^G .
\]

The explicit Gross--Zagier formula of Cai, Shu and Tian applies to
\(P_\beta^G\) for either parity of \(D\)
\cite[Theorem~1.1]{CaiShuTian}. The ring class conductor is one,
\((D,N)=1\), and \(\mathcal O_K^\times=\{\pm1\}\). It gives
\[
 \widehat h_K(P_\beta^G)=
 \frac{\sqrt{|D|}}{8\pi^2\|G\|_{\rm Pet}^2}\,
 L'(G/K,\eta,1).
\]
Inducing \(\eta\) gives \(\chi_\Delta\oplus\chi_d\), so
\(L(G/K,\eta,s)=L(G,\chi_\Delta,s)L(G,\chi_d,s)\). Only the
derivative of the first factor contributes at \(s=1\) because its
central value is zero. Passing to absolute height divides by two
\cite[I, (6.4)]{GrossZagier}. Passing from \(P_\beta^G\) to
\(y_{\Delta,r}^G(n,\mu)=2P_\beta^G\) multiplies the height by
four. This proves \eqref{med-gz}. The common nonzero factor
\(L(G,\chi_d,1)\) cancels on dividing by \eqref{med-wald}.
Equation~\eqref{med-proj} then gives \eqref{rd-h-nt}. The same
equalities recover \eqref{med-gz} from \eqref{rd-h-nt}.

For any \(a\) with \(b_g(a)\ne0\), take the squared norm of
\eqref{med-proj} at \(a\) and substitute \eqref{rd-h-nt}.
Solving for the derivative gives \eqref{eq:curve-single-derivative}.
Using \eqref{med-proj} again at \(a\) gives the divisor identity
even when the derivative vanishes.
\end{proof}

\smallskip\noindent\emph{Proof of the shared-prime formula
\eqref{eq:curve-shared-height}.}
Choose \(d\) as above by \cite[Theorem~(ii)]{BumpFriedbergHoffstein},
using \(\varepsilon(G)=-\sigma\). Thus \(d\) is fundamental and
prime to \(2N\Delta\), every prime dividing \(2N|\Delta|\) splits in
\(\mathbb Q(\sqrt d)\), and \(L(G,\chi_d,1)\ne0\). At an index
\(n=|d|/(4N)\), \(\mu^2\equiv d\bmod 4N\), formula \eqref{med-wald}
gives \(b_g(n,\mu)\ne0\). The coefficient identity \eqref{med-proj}
still applies. Indeed \cite[(5.1), Corollary~6.7 and Theorem~7.7]{BO}
uses the full vector-valued series with the orbit sum in \(r\mu+L\)
for each original component \(\mu\). Its \(G\)-projection lies on the
unique Shimura line with these Hecke eigenvalues, since \(G\) is new
of exact level \(N\). The factor of two between the orbit weights
in (5.1) of \cite{BO} and in \eqref{rd-coeff} multiplies both sides
of \eqref{med-proj} equally.

Put \(D=\Delta d\), \(K=\mathbb Q(\sqrt D)\), and
\(\beta=-r\mu\bmod 2N\). Each positive definite branch \(P_\beta\)
still contains each ideal class once. In fact the integer
\(\gcd(N,\beta,(\beta^2-D)/(4N))\) in the orbit classification
of \cite[I, \S1, p.~505]{GKZ} is one. A prime dividing both
\(N\) and \(\beta\) belongs to \(\mathcal R\) and divides \(N\) exactly
once, by fundamentality and the square congruence. For an odd such
prime \(p\), \(v_p(D)=1\), so \(p\nmid(\beta^2-D)/(4N)\).
At a shared \(2\), if \(v_2(D)=2\), then \(D/4\equiv3\bmod4\)
and \(\beta\equiv2\bmod4\). If \(v_2(D)=3\), then \(4\mid\beta\).
In both cases \(v_2(\beta^2-D)=3=v_2(4N)\).
On these primitive forms the generalized genus character is the
ordinary ideal class character \(\eta\), by taking \(N_1=N\) in
\cite[I, \S2, Proposition~1]{GKZ}. Fricke invariance and
\(w_N|_G=\sigma\) again give \(y_{\Delta,r}^G(n,\mu)=2P_\beta^G\).

Every level prime outside \(\mathcal R\) splits in \(K\).
At \(p\in\mathcal R\) the unique prime \(\mathfrak p\) above \(p\)
ramifies. Evaluating the genus character with the factor \(d\) of
\(D\) gives \(\eta([\mathfrak p])=\chi_d(\mathrm N\mathfrak p)=\chi_d(p)=1\).

If \(a_p(G)=-1\) at every shared prime, the local condition
\(\eta([\mathfrak p])\ne a_p(G)\) holds throughout \(\mathcal R\).
The Gross--Zagier formula then acquires a factor \(2^t\) and the
Rankin product has root number \(-1\) \cite[Theorem~1.1]{CaiShuTian}.
Since \(L(G,\chi_d,1)\ne0\), the factor \(G\otimes\chi_\Delta\)
has root number \(-1\). Normalizing the height as above and applying
\eqref{med-wald} and \eqref{med-proj} at the same component \(\mu\)
gives \eqref{eq:curve-shared-height}.

Suppose instead that \(a_p(G)=+1\) at some shared prime, still
with \(\varepsilon(G)=-\sigma\). Use the same \(d\) and nonzero
coefficient of \(g\). The \(p\)-part of the level
ideal is the unique ramified prime \(\mathfrak p\) of norm \(p\).
Quotienting by this subgroup replaces an ideal class
\([\mathfrak a]\) by \([\mathfrak p^{-1}\mathfrak a]\) and
preserves the oriented level condition, since
\(\bar{\mathfrak p}=\mathfrak p\). As \(\eta([\mathfrak p])=1\),
\(W_p\) fixes \(P_\beta\). Its eigenvalue on \(G\) is
\(-a_p(G)=-1\), so \(P_\beta^G=0\). Equation~\eqref{med-proj}
and \(b_g(n,\mu)\ne0\) give \(\theta_{\Delta,r}(g)=0\).

In \eqref{rd-h-nt} and under the local signs of
\eqref{eq:curve-shared-height}, the factor multiplying the central
derivative is positive. The lift of the real form \(g\) lies in the
real span of the divisor classes, where the N\'eron--Tate pairing
is positive definite modulo torsion. Hence
\(L'(G,\chi_\Delta,1)=0\) if and only if
\(\theta_{\Delta,r}(g)=0\) in the complexified Jacobian. For squarefree \(N\)
and \(\Delta=r=1\), \eqref{rd-h-proj} gives the same criterion
for the arithmetic lift.

For \(\Delta=1\) and \(\kappa=1\) one can also keep the constant
term of the input. Let \(U\) be an auxiliary CM plane whose
negative fundamental discriminant is prime to \(2N\). In the
normalization of \cite[Theorem~7.14]{BruinierYangCM}, the Faltings
intersection of the compactified arithmetic divisor
\(\widehat Z^{\,c}(F)\) with the CM zero-cycle \(Z(U)\) is
\[
 \frac{\deg Z(U)}2
       \left(c_F^+(0,0)k_U(0,0)-
               L'(\xi_{1/2}F,U,0)\right).
\]
Here \(k_U(0,0)\) is the constant coefficient of the holomorphic
part of the incoherent Eisenstein derivative, and
\(L'(\xi_{1/2}F,U,0)\) uses the Rankin--Selberg variable of that
formula. The finite intersection subtracts the contributions of
the polar part from the constant-term expression at infinity,
leaving the displayed \(c_F^+(0,0)\) term. In the Jacobian
identity \eqref{merge-div-shadow}, the cusp form has zero constant
coefficient.

\subsection{The mixed height on the modular curve}
Work under Theorem~\ref{med-inner}, or under the sign conditions of
\eqref{eq:curve-shared-height} when \(t=\#\{p:p\mid(N,\Delta)\}>0\).
Choose \(F_g\) as in \eqref{med-proj}. Write \(P=y_{\Delta,r}(F_g)\) and
\(y(a)=y_{\Delta,r}(m,\mu)\) for \(a=(m,\mu)\). The self-adjoint
Hecke projector onto the \(G\)-component fixes \(P\) and sends \(y(a)\)
to \(b_g(a)P\). Thus \(\langle P,y(a)\rangle_{\rm NT}
=b_g(a)\|P\|_{\rm NT}^2\). With \(\theta_{\Delta,r}(g)=\|g\|_{\rm Pet}^2P\),
equations~\eqref{rd-h-nt} and~\eqref{eq:curve-shared-height} give
\begin{equation}\label{p19-divisor-mixed}
 \langle P,y(a)\rangle_{\rm NT}
   =2^t\,\frac{2\sqrt{N|\Delta|}}{\pi}\,
       \frac{b_g(a)}{\|g\|_{\rm Pet}^{2}}\,
           L'(G,\chi_\Delta,1).
\end{equation}
Conversely, dividing one such mixed height by a nonzero \(b_g(a)\)
recovers \(\|P\|_{\rm NT}^2\) and gives the derivative in
\eqref{eq:curve-single-derivative}. If a shared prime has \(a_p(G)=+1\)
and \(\varepsilon(G)=-\operatorname{sgn}(\Delta)\), the preceding
Atkin--Lehner calculation gives \(P=0\), so every such mixed height
is zero. For a calculation by regularized integration, the cusp subtraction
in \(y(a)\), the boundary term and the finite intersections must all
refer to this same pair (Section~\ref{rd-h-height}).

\subsection{A height determinant over a multiquadratic field}
The absolute N\'eron--Tate pairing allows the lifts for several
quadratic twists to be compared on one Jacobian. Over a multiquadratic
field these lifts belong to distinct Galois character spaces. We take
\(K\) totally real so that all its characters, including the trivial
one, use the same Shimura line. Their height determinant gives a
higher derivative of the base-change \(L\)-function.

\begin{theorem}\label{thm:multiquadratic-height}
Let \(K/\mathbb Q\) be a totally real multiquadratic extension of
degree \(d\), with Galois group \(\mathcal G\). Keep the normalized
weight-two newform \(G\) on \(\Gamma_0(N)\) and its real Shimura
correspondent \(g\) from Theorem~\ref{med-inner}, with a fixed real
embedding of the Hecke field. For each character \(\chi\) of
\(\mathcal G\), let \(\Delta_\chi\) be the positive fundamental
discriminant giving \(\chi\) as a Dirichlet character, and set
\(\Delta_{\mathbf1}=1\). For each \(p\mid N\) let \(I_p\) be
its inertia group in \(K/\mathbb Q\), and put
\(a=\sum_{p\mid N}d(1-1/|I_p|)\). Assume \(a_p(G)=-1\) at each
\(p\mid N\) ramified in \(K\). Choose
\(r_\chi\bmod 2N\) with \(r_\chi^2\equiv\Delta_\chi\bmod 4N\)
for each character, taking \(r_{\mathbf1}=1\).
Assume that \(G\) has root number \(-1\). Define the vectors
\[
       P_\chi=
       \frac{\theta_{\Delta_\chi,r_\chi}(g)}
            {\|g\|_{\rm Pet}\,\Delta_\chi^{1/4}}
       \ \in\ J_0(N)(K)\otimes_{\mathbb Z}\mathbb R .
\]
Write \(L(G/K,s)\) for the finite \(L\)-function of \(G\) after base
change from \(\mathbb Q\) to \(K\), in the classical variable with
centre \(s=1\). Then
\begin{equation}\label{eq:multiquadratic-leading}
  \det\!\bigl(\langle P_\chi,P_\psi\rangle_{\rm NT}\bigr)_{\chi,\psi}
       =2^a\left(\frac{2\sqrt N}{\pi}\right)^d
         \frac{L^{(d)}(G/K,1)}{d!}.
\end{equation}
The determinant is positive precisely when
\(\operatorname{ord}_{s=1} L(G/K,s)=d\). If \(e\) denotes the
nullity of the height matrix, then
\[
  e=\#\{\chi:L'(G,\chi,1)=0\},\qquad
  \operatorname{ord}_{s=1}L(G/K,s)\ \ge\ d+2e .
\]
\end{theorem}
\begin{proof}
For a nontrivial \(\chi\), the coefficients of the twisted divisor
series are defined over \(\mathbb Q(\sqrt{\Delta_\chi})\), and its
nontrivial automorphism negates them \cite[Lemma~5.1(iii),(v)]{BO}.
For the trivial character the degree-zero coefficients are rational.
The genus character therefore gives
\(\sigma P_\chi=\chi(\sigma)P_\chi\) for every
\(\sigma\in\mathcal G\). These are real vectors: in
\eqref{merge-div-shadow} one may take a preimage of the normalized
shadow with real principal part \cite[Lemma~7.4]{BruinierYangCM}.
The canonical height is absolute and Galois invariant. For
\(\chi\ne\psi\) take \(\sigma\) with \(\chi(\sigma)\psi(\sigma)=-1\).
Invariance forces \(\langle P_\chi,P_\psi\rangle_{\rm NT}=0\).
Write \(t_\chi=\#\{p:p\mid (N,\Delta_\chi)\}\). If \(t_\chi=0\),
the character is even and its level congruence gives \(\chi(N)=1\),
so the twist has root number \(-1\). If \(t_\chi>0\), the sign
follows from the shared-prime calculation in \eqref{eq:curve-shared-height}.
By \eqref{rd-h-nt} and \eqref{eq:curve-shared-height},
\[
       \|P_\chi\|_{\rm NT}^2
          =2^{t_\chi}\frac{2\sqrt N}{\pi}L'(G,\chi,1).
\]
A character is unramified at \(p\) exactly when it is trivial on
\(I_p\). There are \(d/|I_p|\) such characters. Counting the
remaining characters at each \(p\mid N\) gives
\(\sum_\chi t_\chi=a\), also at \(p=2\).

The regular representation of \(\mathcal G\) is the direct sum of
its characters. Artin formalism gives
\(L(G/K,s)=\prod_\chi L(G,\chi,s)\). At a shared prime the square
congruence forces \(p\Vert N\), and the local representation is an
unramified twist of Steinberg. Its ramified quadratic twist has
Euler factor one, agreeing with the Dirichlet series. All \(d\) factors vanish at \(1\), so the coefficient of \((s-1)^d\)
is \(\prod_\chi L'(G,\chi,1)\). Taking the product of the diagonal
heights proves \eqref{eq:multiquadratic-leading}. Root number \(-1\)
makes the completed twist odd in \(s-1\). Its completion factor is
holomorphic and nonzero at \(1\), so the finite twist has the same
odd vanishing order. Each zero first derivative thus adds at least two to the
lower bound \(d\) for the order of the product.
\end{proof}

Positivity of the N\'eron--Tate pairing on
\(J_0(N)(K)\otimes\mathbb R\) identifies \(d-e\) with the
dimension of the span of these lifts. Formula~\eqref{eq:multiquadratic-leading}
is the Gram determinant of the specified real vectors. Computing the
Mordell--Weil regulator would also require an integral basis of the
full group modulo torsion.

\section{CM orbits and regularized theta lifts}\label{p19-cm-section}

\subsection{Translating a twisted Heegner divisor}\label{p19-translated}
Throughout this section assume \((\Delta,2N)=1\) unless specified
otherwise. To evaluate the theta lift on a CM orbit, we first
determine how translation changes the genus weight in \eqref{rd-coeff}.
Write \(K_0(N)=\prod_p K_p\), where \(K_p\) consists of the matrices
in \(\mathrm{GL}_2(\mathbb Z_p)\) with lower left entry in
\(N\mathbb Z_p\). Identify \(\mathrm{GSpin}(V)\) with
\(\mathrm{GL}_2\), acting on trace-zero matrices by conjugation.

For \(x\in L^\sharp\) with entries \(b/(2N),-a/N,c,-b/(2N)\),
write \(D_x=b^2-4Nac\). The genus character is supported on
\(\Delta\mid D_x\), \(D_x/\Delta\) a square modulo \(4N\), and
\(\gcd(a,b,c,\Delta)=1\). On this support \(\chi_\Delta(x)\)
is the Kronecker symbol \((\frac{\Delta}{t})\) for any integer
\(t\) prime to \(\Delta\) represented by an admissible form
\([aN_1,b,cN_2]\), \(N_1N_2=N\) \cite[\S4]{BO}.
For \((N,\Delta)=1\) the form \([a,b,Nc]\) gives the level-one
character law applied to \(Nx\) \cite[Lemma~7.3]{BEY}. The same
law on units also holds at primes shared by \(N\) and \(\Delta\).
Indeed, for \(j\mid N\) and \(v=(X,Y)^t\), put
\[
 M_j=\begin{pmatrix}0&1\\-j&0\end{pmatrix},\qquad
 Q_{x,j}(v)=\frac{N}{j}\det(M_jv,xM_jv)
            =jaX^2+bXY+(N/j)cY^2 .
\]
If \(k_p\in K_p\), then \(M_j^{-1}k_pM_j\in\mathrm{GL}_2(\mathbb Z_p)\)
because the lower left entry of \(k_p\) is divisible by \(N\). Thus
\(Q_{k_pxk_p^{-1},j}(v)=\det(k_p)Q_{x,j}(M_j^{-1}k_p^{-1}M_jv)\).
The represented units at \(p\) are multiplied by \(\det k_p\) for
every \(j\), and conjugation preserves \(D_x\) and \(b\bmod 2N\).
The represented-unit definition above, applied at the primes of
\(\Delta\) by the Chinese remainder theorem, gives
\[
    \chi_{\Delta,p}(k_p x k_p^{-1})
       =(\Delta,\det k_p)_p\,\chi_{\Delta,p}(x).
\]
This also covers zero weights and \(p=2\). At \(p\nmid\Delta\)
the character is trivial on units.

Put \(\eta_\Delta(h)=\prod_{p<\infty}(\Delta,\det h_p)_p\) and
\(K_0(N)^{\chi_\Delta}=\ker(\eta_\Delta|_{K_0(N)})\). On the upper half-plane
component write a finite adelic element \(h=\gamma k\), with
\(\gamma\in\mathrm{GL}_2(\mathbb Q)^+\) and \(k\in K_0(N)\).
Such a decomposition exists by strong approximation and the
surjectivity of the determinant of \(K_0(N)\) onto
\(\widehat{\mathbb Z}^{\,\times}\). It identifies \((z,h)\) with
\(\gamma^{-1}z\) at level \(K_0(N)\). To translate a divisor at
level \(K_0(N)^{\chi_\Delta}\), take the same orbit sum as in \eqref{rd-coeff}
in codimension one, using rational \(x\) such that \(h^{-1}xh\)
belongs to the finite adelic coset \(\widehat L+r\mu\). Evaluate
the character at \(h^{-1}xh\), retaining the orbit stabilizers.
Reindexing by \(\gamma^{-1}x\gamma\) gives the original divisor
multiplied by \(\eta_\Delta(k)=\eta_\Delta(h)\). The last equality
is Hilbert reciprocity for the positive rational number \(\det\gamma\).
The same reindexing in the theta kernel, and hence in its
regularized integral, gives
\begin{equation}\label{p19-translate-lift}
 \Phi_{\Delta,r}^{\,j}(z,h,F)
    =\eta_\Delta(h)\,
      \Phi_{\Delta,r}^{\,j}(\gamma^{-1}z,1,F).
\end{equation}
This formula fixes the sign of the translation independently of the
choice of \(\gamma\). The same calculation on the lower half-plane
uses a positive determinant there as well, so each orientation
retains its own kernel value.

Let \(U\subset V\) be a negative definite rational plane. Write
\(K\) for its imaginary quadratic field and identify
\(\mathrm{GSpin}(U)(\mathbb A_f)\) with
\(\mathbb A_{K,f}^{\times}\). For an optimal embedding of the
maximal order with \(K_0(N)\cap\mathbb A_{K,f}^{\times}
=\widehat{\mathcal O}_{K}^{\,\times}\), the points of the CM
cycle \(Z(U)\) are indexed by the two orientations \(z_U^\pm\)
and the ideal classes of \(K\). If \(w_{K}\) and \(h_{K}\)
are the number of roots of unity and the class number, respectively,
the orbit sum is
\begin{equation}\label{p19-auxiliary-cycle}
   Z(U)=\frac{2}{w_{K}}
       \sum_{\epsilon=\pm}\sum_{[h]\in\operatorname{Cl}(K)}
              [(z_U^\epsilon,h)],\qquad
   \deg Z(U)=\frac{4h_{K}}{w_{K}} .
\end{equation}

\begin{proposition}\label{p19-genus-translation}
Suppose the fundamental discriminant of \(K\) factors as
\(D_0=\Delta d_0\) with coprime fundamental discriminants. Let
\(\chi\) be the genus character of this factorization. For every
\(h\in\mathbb A_{K,f}^{\times}\) one has
\(\eta_\Delta(h)=\chi([h])\).
\end{proposition}
\begin{proof}
Under the torus embedding, \(\det h\) is the norm of \(h\) to the
finite ideles of \(\mathbb Q\). The quadratic character associated
with \(\mathbb Q(\sqrt\Delta)\), pulled back by this norm, factors
through \(\operatorname{Cl}(K)\). Indeed its class field over
\(K\) is the genus extension \(K(\sqrt{\Delta})/K\),
unramified at the finite primes by the fundamental discriminant
factorization. The character is trivial if this extension is
trivial. On an ideal prime to \(D_0\) the resulting value is
\((\frac{\Delta}{\mathrm N\mathfrak a})\), which defines the genus
character. Such ideals represent every class. The local norm
description of the quadratic character gives the Hilbert symbols
in \(\eta_\Delta\), proving the identity.
\end{proof}

Thus evaluating the translated lift on \eqref{p19-auxiliary-cycle}
introduces one factor \(\chi([h])\) for the auxiliary CM orbit.
The vector defining the first divisor already carries
\(\chi_\Delta(x)\). Both factors enter the finite intersection
calculation in Section~\ref{rd-h-height}.

\subsection{The Green kernel and the polarized theta lift}\label{p19g-main}
Write \(Q_\Delta=Q/|\Delta|\), \(e(t)=\exp(2\pi i t)\) and
\(\tau=u+iv\). The point \(z\) of the upper half-plane determines
a negative plane and its positive orthogonal line. Write \(x_-\) for the projection to the negative plane and
\(x_+\) for that to the positive line. The Siegel theta
kernel is obtained by summing
\[
 v\,e\bigl(Q_\Delta(x_+)\tau+
              Q_\Delta(x_-)\bar\tau\bigr)
\]
over the lattice cosets in each component, with the residue
condition and genus weight of \eqref{rd-coeff}. The Millson kernel
has the additional factor \(p_z(x)=-(x,\lambda(z))_\Delta/\sqrt2\),
where \((\,\cdot\,,\,\cdot\,)_\Delta\) is the bilinear form of
\(Q_\Delta\). The vector \(\lambda(z)\) has squared length one
and we choose its ray so that, for \(z=t+iy\) and
\(x=\left(\begin{smallmatrix} b/(2N)&-d/N\\ c&-b/(2N)\end{smallmatrix}\right)\),
\[
 p_z(x)=-\frac{Nc|z|^2-bt+d}{2\sqrt{N|\Delta|}\,y}.
\]
In particular
\(p_z(x)^2=Q_\Delta(x_+)\). Denote the regularized integrals
of \(F\) against these kernels by \(\Phi_{\Delta,r}(z,F)\) and
\(\Phi^{\mathcal M}_{\Delta,r}(z,F)\). The regularization takes
the constant term at \(t=0\) of the continued truncated integral
with the factor \(v^{-t}\) and measure \(du\,dv/v^2\).

For \(n>0\) in component \(\mu\), let \(F_{n,\mu}(\tau,s,k)\) be
the Poincar\'e series obtained by summing the weight \(k\) translates
of the seed
\[
 \frac{(4\pi nv)^{-k/2}M_{-k/2,s-1/2}(4\pi nv)}
      {2\Gamma(2s)}\,e(-nu)e_\mu
\]
over the translation subgroup in the metaplectic group. Here
\(M_{\alpha,\beta}\) is the Whittaker function. The harmonic
specialization is at \(s=1-k/2\), with both reflected polar terms
as in the proof of Theorem~\ref{merged-height-theorem}. For
\(z\ne z'\) put
\[
 g_t(z,z')=-2Q_{t-1}\!\left(1+
              \frac{|z-z'|^2}{2\Im z\,\Im z'}\right),
\]
where \(Q_{t-1}\) is the Legendre function of the second kind.
To specify the orbit multiplicities, write
\[
 \mathcal G_t^\epsilon(z;n,\mu)
   =\sum_{x\in L_{r\mu}[|\Delta|n]}
      \chi_\Delta(x)\epsilon(x)g_t(z,z(x)).
\]
This is the automorphic Green function evaluated on the weighted
CM divisor: summing over the effective modular group gives each
orbit the factor \(2/|\Gamma_x|\) of \eqref{rd-coeff}. We use
\(\epsilon=1\) for the Siegel kernel and
\(\epsilon(x)=-\operatorname{sgn}(p_z(x))\) for the Millson kernel.
For \(Q_\Delta(x)>0\) the latter sign is constant as \(z\) varies
on the chosen half-plane.

\begin{proposition}[The resolvent of the two theta kernels]
\label{p19g-green-kernel}
Away from the CM divisor, for \(\Re(s)>1\), the two lifts satisfy
\begin{align*}
 \Phi_{\Delta,r}(z,F_{n,\mu}(\cdot,s,\tfrac12))
   &=-\frac{2}{\Gamma(s+\tfrac14)}
        \mathcal G_{2s-1/2}^{1}(z;n,\mu),\\
 \Phi^{\mathcal M}_{\Delta,r}(z,F_{n,\mu}(\cdot,s,-\tfrac12))
   &=\frac{2\sqrt n}{\Gamma(s-\tfrac14)}
        \mathcal G_{2s-1/2}^{\epsilon}(z;n,\mu).
\end{align*}
The identities extend by meromorphic continuation.
\end{proposition}
\begin{proof}
Unfold the Poincar\'e series on the strip \(0\le u\le1\).
The integral in \(u\) imposes \(Q_\Delta(x)=n\). Put
\(X=n/Q_\Delta(x_+)\), so that \(0<X<1\) away from
\(z(x)\). In weight \(1/2\) the Laplace transform of the Whittaker
seed gives the radial factor
\[
 \frac{2\Gamma(s-\tfrac14)}{\Gamma(2s)}X^{s-1/4}
 {}_2F_1(s-\tfrac14,s+\tfrac14;2s;X).
\]
The hypergeometric expression for
\(Q_{2s-3/2}(\cosh d(z,z(x)))\) and the duplication formula for
\(\Gamma\) identify this with \(4Q_{2s-3/2}/\Gamma(s+1/4)\)
\cite[Proposition~6.2]{BEY}. In weight \(-1/2\) the transform
interchanges the two upper hypergeometric parameters, increases the
power of \(X\) by \(1/2\), and replaces
\(\Gamma(s-1/4)\) by \(\Gamma(s+1/4)\). The extra factor is
\(p_z(x)X^{1/2}=-\epsilon(x)\sqrt n\).
Substitution and summation give both formulas.
\end{proof}

Put \(w=3/2-\kappa\) and \(j=\lfloor(\kappa-1)/2\rfloor\).
Write \(\Phi^j_{\Delta,r}\) for the lift of \((4\pi)^{-j}R_w^jF\),
using the Siegel kernel for odd \(\kappa\) and the Millson kernel
for even \(\kappa\), where \(R_k=2i\partial_\tau+k/v\).
At the harmonic specialization the Whittaker raising identity
multiplies the raised seed by \(j!(4\pi n)^j\) relative to the
seed in its final weight. Proposition~\ref{p19g-green-kernel}
therefore gives the factor \((-1)^\kappa2n^{(\kappa-1)/2}\)
in front of the corresponding \(\mathcal G_\kappa\). The linear
insertion in even codimension retains the orientation of the vector
at its CM point.

\subsection{The CM value and its boundary term}
Let \(D_0=\Delta d_0=-4N|\Delta|m_0<0\) be fundamental,
with \(d_0\) fundamental and \((\Delta,2N)=1\).
Fix a primitive CM vector \(x_0\) in component \(r\mu_0\) with
\(Q_\Delta(x_0)=m_0\), \(\chi_\Delta(x_0)\ne0\), and put
\(U=x_0^\perp\). Normalize the genus weight on the orbit by its
value at \(x_0\), equivalently multiply the translated average by
\(\chi_\Delta(x_0)\). This gives the weights of the CM cycle
itself by the genus translation law. The mass in the formulas is
the unsigned mass \(\deg Z(U)=4h(D_0)/|\mathcal O_{D_0}^{\times}|\).
For the linear kernel we choose the orientation with
\(p_{z_U^+}(x_0)=-\sqrt{m_0}\). In both parities the notation
\(\Phi^j_{\Delta,r}(Z(U),F)\) means the following CM average,
with \(a=\kappa-1\):
\[
 \frac{2\chi_\Delta(x_0)}{w_K}\sum_{[h]\in\operatorname{Cl}(K)}
 \left(\Phi^j_{\Delta,r}(z_U^+,h,F)
       +(-1)^a\Phi^j_{\Delta,r}(z_U^-,h,F)\right).
\]
Reversing the orientation fixes both orthogonal projections and
reverses the linear factor \(p_z\). The two summands therefore
contribute equally with these signs. Division by \(\deg Z(U)\)
gives the average over one oriented branch.

Set \(\Lambda=(|\Delta|L,Q_\Delta)\),
\(P=\Lambda\cap\mathbb Qx_0\) and \(B=\Lambda\cap U\).
At the CM plane the Gaussian splits on \(P\oplus B\). To perform
this splitting on the original theta series, let \(\iota_U\) be
the following transfer on coefficient vectors. Replace \(e_\mu\)
by the sum of the admissible \(e_\delta\) for
\(\delta\in\Lambda^\sharp/\Lambda\) above \(r\mu\), weighted by
\(\chi_\Delta(x_0)\chi_\Delta(\delta)\). Replace each \(e_\delta\)
in turn by the sum of the \(e_h\otimes e_\beta\) for
\[
 (h,\beta)\in\Lambda^\sharp/(P\oplus B),\qquad
 h+\beta\equiv\delta\pmod\Lambda .
\]
Here \(h\) and \(\beta\) lie in the discriminant modules of \(P\)
and \(B\), respectively. In particular, their joint residue
condition is retained. Write \(\vartheta\) for the unary theta
series of \(P\), with insertion \(1\) or \(p_{z_U}\) according
to the parity above. Its weight is \(\alpha=1/2\) or \(3/2\).

Let \(E_B(\tau,s;1)\) be the weight-one incoherent Eisenstein
series with finite Schwartz functions the characteristic functions
of the cosets of \(B\), in the normalization of
\cite[\S2]{BruinierYangCM}. It vanishes at \(s=0\). Write
\(\mathcal E_B=\sum_{\beta,t\ge0}k(t,\beta)q^te_\beta\)
for the holomorphic Fourier part of its derivative there, including
the constant left after subtracting the logarithm in the zero
component. Define
\begin{align*}
 B_U(F)&=\operatorname{CT}\langle\iota_U F^+,
                                  [\vartheta,\mathcal E_B]_j\rangle,\\
 \mathcal L_U(g,s)&=\langle[\vartheta,E_B(\cdot,s;1)]_j,
                                \iota_U g\rangle_{\rm Pet}.
\end{align*}
The bracket is the Rankin--Cohen bracket for weights \((\alpha,1)\),
with derivatives \((2\pi i)^{-1}\partial_\tau\) and the order of
the factors as displayed \cite[(3.6)]{BEY}.

\begin{proposition}[The CM value before finite intersections]
\label{p19g-cm-derivative}
For a harmonic input with cuspidal shadow whose polar divisor is
disjoint from the CM orbit, the normalized CM value is
\begin{equation}\label{p19g-stokes}
 \frac{\Phi^j_{\Delta,r}(Z(U),F)}{\deg Z(U)}
      = B_U(F)-\mathcal L_U'(\xi_w F,0).
\end{equation}
For even \(\kappa\) the left side uses the oriented unary insertion.
\end{proposition}
\begin{proof}
Apply \(\iota_U\) to \(F\) and factor the theta kernel. By
Siegel--Weil, the binary theta series has CM average
\(\deg Z(U)E_B(\tau,0;-1)/2\) \cite[\S4]{BruinierYangCM}.
On the truncated domain move the \(j\) raising operators from
\(F\) to the kernel. The identity \(L_1E'_B(\tau,0;1)=E_B(\tau,0;-1)/2\)
and the Rankin--Cohen lowering identity give the raised kernel
\(2(-4\pi)^j L[\vartheta,E'_B(\cdot,0;1)]_j\), with
\(L=-2iv^2\partial_{\bar\tau}\) \cite[Proposition~3.6]{BEY}.

The product rule yields an interior term involving \(L\iota_U F\)
and a total derivative. Since \(\iota_U\) is real, the interior
term is the negative Petersson pairing with \(\iota_U\xi_w F\).
The top-edge integral of the total derivative takes the constant
Fourier coefficient. Regularization removes its logarithmic term,
and the decaying terms tend to zero. The remaining constant is
\(B_U(F)\) \cite[Theorem~5.4 and the proof of Theorem~7.13]{BEY}.
\end{proof}

There is a finite test for the boundary in \eqref{p19g-stokes}.
Write \(c_U^+(-n;h,\beta)\) for the coefficients of
\(\iota_U F^+\) and \(r(t,h)\) for those of \(\vartheta\). With
\[
 P_j(t,m)=\sum_{b=0}^j(-1)^b
      \binom{\alpha+j-1}{b}\binom{j}{j-b}t^{j-b}m^b,
\]
the test is
\begin{equation}\label{p19g-boundary}
 B_U(F)=\sum_{\substack{h,\beta\\t,m\ge0}}
      c_U^+(-t-m;h,\beta)r(t,h)k(m,\beta)P_j(t,m).
\end{equation}
This is the constant term of the product of the three Fourier expansions.
The term \(t=m=0\) contributes only for \(j=0\); in that case
the zero-component constant of \(\mathcal E_B\) must be retained.
At \(m>0\) the Eisenstein coefficient is a product of local
Whittaker functions before differentiation. Two vanishing factors
at the centre give \(k(m,\beta)=0\). Indeed, after extracting
the vanishing factors, bundle the remaining holomorphic factors
into one function. If only the factor at \(p\) vanishes, the
product rule gives its derivative times the value of that function.
With two vanishing factors every first-order term is zero. For instance the
normalized unramified local density at \(p\nmid 2ND_0\), for an
integral index and the unshifted local lattice, is
\[
 (1-\chi_{D_0}(p)p^{-1}X)
       \sum_{i=0}^{v_p(mN)}\chi_{D_0}(p)^iX^i,\qquad X=p^{-s}.
\]
For inert \(p\) and odd \(v_p(mN)\) this vanishes at \(s=0\)
and its derivative is
\((1+p^{-1})(v_p(mN)+1)\log(p)/2\).

Two terms in \eqref{p19g-boundary} cancel if their genus weights
are opposite and every other factor agrees. For instance, negating
a unit at \(q\mid\Delta\) with \(q\equiv3\bmod4\) reverses the local
character. To pair the two residues one must also preserve the
transferred multiplicity, the unary coefficient and the binary
Whittaker coefficient.

When all pole orders are less than \(m_0\), the line \(\mathbb Zx_0\)
contributes no positive unary norm. The test still includes the other cosets
in \(P^\sharp/P\). The unary projection of
\(y\in\Lambda^\sharp\) is \((y,x_0)_\Delta x_0/(2m_0)\), of norm
\((y,x_0)_\Delta^2/(4m_0)\), and the lattice condition
\(h+\beta\in\Lambda^\sharp\) couples the unary and binary residues.

For an odd negative fundamental \(\Delta\), one always has a prime
\(q\mid\Delta\) with \(q\equiv3\bmod4\), since
\(|\Delta|\equiv3\bmod4\). In the positive-twist case with even
\(\kappa\), reflection of the input gives \(c_F^+(-n,0)=0\). These
observations explain the distinction between the signs in the
boundary calculation. For the other residue classes the joint
sum in \eqref{p19g-boundary} still has to be evaluated.

\subsection{The Dirichlet series of the CM pairing}
\begin{proposition}[The shift to the centre of the Shimura lift]
\label{p19g-shimura}
Assume in addition that \(\mathbb Qx_0\cap L^\sharp=\mathbb Zx_0\).
Let \(g\) have weight \(\kappa+1/2\) and real Fourier coefficients
in the cusp space, on the cuspidal Shimura component if \(\kappa=1\).
Put \(S_g=\operatorname{Sh}_{m_0,\mu_0}(g)\) as in
\eqref{merged-shimura-coeff}. Then
\begin{equation}\label{p19g-shift}
 \mathcal L_U(g,s)=A_\kappa(s)
          \frac{L(S_g,\chi_\Delta,s+\kappa)}
               {L(\chi_{D_0},s+1)},
\end{equation}
where
\[
 A_\kappa(s)=
 \frac{(-1)^{\kappa-1}\Gamma(\kappa-\tfrac12)\Gamma(s+\kappa)}
 {2^{2s+2\kappa-2}\pi^{(s-1)/2+\kappa}
             m_0^{(s+\kappa)/2}\Gamma(\kappa)\Gamma(1+s/2)} .
\]
Consequently
\begin{equation}\label{p19g-central}
 \mathcal L_U'(g,0)=
 \frac{(-1)^{\kappa-1}2^{4-2\kappa}\Gamma(\kappa-\tfrac12)
                     \sqrt{N|\Delta|}}
 {m_0^{(\kappa-1)/2}\deg Z(U)\,\pi^{\kappa+1/2}}\,
                 L'(S_g,\chi_\Delta,\kappa).
\end{equation}
\end{proposition}
\begin{proof}
Against a cusp form the bracket can be replaced by
\(\vartheta\otimes E_B(\tau,s;1+2j)\) multiplied by
\[
 (4\pi)^{-j}\binom{\alpha+2j-1}{j}
       \frac{\Gamma(1+s/2+j)}{\Gamma(1+s/2)}.
\]
Indeed, the other term in the bracket identity is a total raising
operator, orthogonal to holomorphic cusp forms \cite[Lemmas~3.2--3.3]{BEY}.
Unfold the Eisenstein series for \(\Re(s)\) large. Its seed selects
the zero binary coset, leaving the vectors in
\(\mathbb Qx_0\cap\Lambda^\sharp=\mathbb Zx_0\), since
\(\Lambda^\sharp=L^\sharp\). The vector \(ax_0\) has
\(Q_\Delta(ax_0)=a^2m_0\) and residue \(ar\mu_0\) in \(L^\sharp/L\).
It therefore selects input component \(a\mu_0\), with normalized
genus weight \(\chi_\Delta(a)\). Pairing the two signs of \(a\)
doubles the sum, by the component symmetry of \(g\) and the parity
of the unary insertion. The integral in \(v\) is
\(\Gamma(s/2+\alpha+j)/(4\pi a^2m_0)^{s/2+\alpha+j}\).
In even codimension the unary coefficient at \(ax_0\), for
\(a>0\), has the factor \(-a\sqrt{m_0}\).
The duplication formula now gives
\[
 \mathcal L_U(g,s)=A_\kappa(s)
   \sum_{a\ge1}\chi_\Delta(a)b(a^2m_0,a\mu_0)a^{-s-\kappa}.
\]
In the Dirichlet series for \(S_g\), put \(\ell=ab\) in the
divisor sum \eqref{merged-shimura-coeff}. The \(b\)-sum at
\(s+\kappa\) is \(L(\chi_\Delta\chi_{d_0},1+s)
=L(\chi_{D_0},1+s)\); the \(a\)-sum is the one just displayed.
This proves \eqref{p19g-shift} by continuation. Incoherence gives
\(\mathcal L_U(g,0)=0\), and \(L(\chi_{D_0},1)\ne0\), so
the Shimura central value vanishes. Differentiating and substituting
the class number formula
\(\deg Z(U)=2\sqrt{|D_0|}L(\chi_{D_0},1)/\pi\)
gives \eqref{p19g-central}.
\end{proof}

\section{Heights of higher Heegner cycles}\label{rd-h-height}

For \(\kappa\ge2\), we pair the polar cycle \(Z_\kappa(F)\) with an
auxiliary CM cycle \(Z_\kappa(U)\). The discriminant of \(U\) determines
which Shimura lift enters the height formula.
Write \(H(A,B)=(-1)^\kappa
\langle\widehat A,\widehat B\rangle_{\rm GS}\), using admissible
extensions and dividing the arithmetic degree by the degree of
the field of definition \cite{Zhang,BEY}. For the span of CM cycles
on a fine-level model, Qiu constructs this pairing on Chow classes,
including self-pairings, and proves the Hecke projection formula
\cite[\S\S2.2.4--2.2.5, (2.12)--(2.14)]{QiuCM}.
At cyclic level one pulls the cycles to a common fine cover and
divides the pairing by the degree of the cover. The projection formula
on a common refinement makes this independent of the cover.
We first work with real Fourier coefficients.

For an auxiliary fundamental discriminant \(d_0\) of sign opposite
to \(\Delta\), choose \(\mu_0^2\equiv d_0\bmod 4N\), and put
\[
 m_0=\frac{|d_0|}{4N},\qquad D_0=\Delta d_0=-4N|\Delta|m_0 .
\]
Thus \(d_0\) indexes the Shimura lift and \(D_0\) is the
discriminant of the CM datum. We use a fundamental \(D_0\) in the
height comparisons below. Write \(Z_\kappa(U)=Z_\kappa(m_0,\mu_0)\) for
the corresponding cycle with the orientations and stabilizers of
\eqref{rd-coeff}. For \(g=\sum b(n,\mu)q^n e_\mu\), our
normalization of the Shimura lift is
\begin{equation}\label{merged-shimura-coeff}
 [q^\ell]\operatorname{Sh}_{m_0,\mu_0}(g)=
 \sum_{u\mid\ell}\chi_{d_0}(u)u^{\kappa-1}
       b\big(m_0(\ell/u)^2,(\ell/u)\mu_0\big).
\end{equation}
For a cusp form \(S=\sum_{\ell\ge1}a_S(\ell)q^\ell\) of weight
\(2\kappa\), we use the classical normalization
\[
 L(S,\chi_\Delta,s)=\sum_{\ell\ge1}
       \chi_\Delta(\ell)a_S(\ell)\ell^{-s},\qquad
 L'(S,\chi_\Delta,\kappa)=
      \left.\frac{d}{ds}L(S,\chi_\Delta,s)\right|_{s=\kappa}.
\]
The series defines the function for \(\Re(s)\) large, and the
derivative is taken after analytic continuation. We omit
\(\chi_\Delta\) when \(\Delta=1\) and use the constant
\(C_{\kappa,N,|\Delta|}\) of Section~\ref{rd-intro}.

\subsection{The fundamental untwisted coefficient}

For the untwisted identity \eqref{med-untwistedheight}, the polar
indices are prime to \(D_0\) in the sense stated in
Theorem~\ref{merged-height-theorem}. The finite and
archimedean terms can both be expressed using the CM plane \(U\).
Their constant terms cancel in the global height.

\begin{proof}[Proof of the untwisted height identity \eqref{med-untwistedheight}]
For odd \(\kappa>1\), the formula follows from Theorem~6.5 and
Corollary~6.6 of \cite{BEY}. For \(\kappa=2j+2\), put \(a=\kappa-1\).
We express the boundary and the finite height in terms of the same
binary Eisenstein coefficients, following \cite[\S4]{DuPeng}.

It suffices to treat the Poincar\'e series \(F_1=F_{m_1,\mu_1}\)
at each admissible pole. Subtracting a real linear combination of these
series with the polar part of \(F\) leaves a form with zero shadow
by Proposition~\ref{med-exact}. It is holomorphic of negative weight
and hence zero. The auxiliary vector \(x_0\) is primitive as required
by Proposition~\ref{p19g-shimura}. A factorization \(x_0=b x'\) with
\(x'\in L^\sharp\) and integral \(b>1\) would give
\(b^2\mid 4NQ(x_0)=|D_0|\), but \(D_0\) is odd and fundamental.

Write \(D_i=-4Nm_i\), choose positive square roots of \(|D_i|\),
and let \(r_i\) represent \(\mu_i\bmod 2N\). Sum over
\(n\equiv r_0r_1\bmod 2N\) with \(n^2<D_0D_1\). The endpoints
are excluded by \((D_0,2ND_1)=1\). Set
\begin{align*}
 t_n&=\frac{n}{\sqrt{|D_0||D_1|}},&
 v_n&=m_1(1-t_n^2),\\
 \delta_n&=\frac{n+r_1\sqrt{D_0}}{2\sqrt{D_0}},&
 S&=\sum_n k(v_n,\delta_n)P_a(t_n).
\end{align*}
Here \(k(v_n,\delta_n)\) is the binary Eisenstein derivative
coefficient at the indicated index and residue.

\smallskip\noindent\emph{The boundary at infinity.}
To determine the analytic residue, take a lift \(x_1\in\mu_1+L\)
and put \(u_i=2Nx_i\), \(n=\operatorname{tr}(u_0u_1)/2\).
Here \(m_1\) is the total Fourier index. The positive complementary
index \(v_n\) comes from the incoherent Eisenstein series.
In the split matrix model the diagonal entries of
\(u_i\) are \(b_i,-b_i\), and the off-diagonal entries are
\(-2a_i,2Nc_i\). Thus \(n=b_0b_1-2N(a_0c_1+a_1c_0)\), giving
\(n\equiv r_0r_1\bmod 2N\). For the standard primitive representative
\(a_0=-1\), varying \(c_1\) gives every such \(n\).
Polarizing \(Q=N\det\) gives
\((x_1,x_0)=-n/(2N)\). Since \((x_0,x_0)=-D_0/(2N)\), the unary
projection is \((n/D_0)x_0\). The complementary binary residue is
\(\delta_n\) \cite[Lemma~7.13]{BruinierYangCM},
\cite[proof of Proposition~6.7]{BEY}.
Conjugation by \(\operatorname{diag}(-1,1)\) fixes the diagonal
entries and the signed trace. The two rays have the same CM point
and enter the Millson kernel with their own orientations in
\eqref{arch-signed-double-sum}. Our choice \(p_{z_U^+}(x_0)=-\sqrt{m_0}\)
therefore gives the unary insertion
\[
  -\frac{n\sqrt{m_0}}{D_0}=\sqrt{m_1}\,t_n .
\]
The two indices in the Rankin--Cohen bracket are \(m_1t_n^2\)
and \(v_n\). With weights \((3/2,1)\) and the bracket convention
of \cite[(3.6)]{BEY}, the Leibniz expansion and Rodrigues' formula give,
after removing \(k(v_n,\delta_n)\),
\[
 m_1^{a/2}\sum_{s=0}^j
     \binom{j+1/2}{s}\binom js
     t_n^{2j+1-2s}(t_n^2-1)^s
   = m_1^{a/2}P_a(t_n).
\]
Simultaneous reflection of \(r_1,n,\delta_n\) preserves the binary
coefficient. Both \(P_a\) and the reflected polar term change sign,
so the two polar terms contribute equally. The odd unary series
has zero constant term. Thus
\[
                  B_U(F_{m_1,\mu_1})=2m_1^{a/2}S .
\]

\smallskip\noindent\emph{The finite height.}
Each prime dividing \(N\) splits in \(\mathbb Q(\sqrt{D_0})\), so
the auxiliary elliptic curve has ordinary reduction there. Its
endomorphism algebra cannot contain both distinct CM fields. Hence
the horizontal supports are disjoint at these primes
\cite[after (6.23)]{BEY}. The standard closures used in the untwisted
calculation are admissible \cite[(6.14)--(6.16)]{BEY}, so these
primes contribute zero to the finite height. At a common reduction, necessarily over
\(p\nmid N\), let \(u_i\) be the image of \(\sqrt{D_i}\) under
the marked CM embedding. The relation \(u_0u_1+u_1u_0=2n\)
fixes the sign of the same index \(n\) \cite[Lemma~6.11]{BEY}.
Order both symplectic eigenbases by \(i\sqrt{|D_i|},-i\sqrt{|D_i|}\)
with the same square root \(i\) of \(-1\). For the change of basis
\(\left(\begin{smallmatrix}A&B\\ C&D\end{smallmatrix}\right)\)
the trace gives
\[
 AD+BC=\frac{\operatorname{tr}(u_0u_1)}
                  {2(i\sqrt{|D_0|})(i\sqrt{|D_1|})}=-t_n.
\]
The normalized fibre pairing is therefore
\((-1)^aP_a(AD+BC)=P_a(t_n)\)
\cite[Proposition~3.3.3]{Zhang}, \cite[(3.7)]{XueGKZ}.
The degree of the base intersection is
\(-\deg Z(U)\,k(v_n,\delta_n)/2\), with stabilizers and deformation
lengths included \cite[Theorem~6.4]{BruinierYangCM},
\cite[Lemma~6.11]{BEY}. This base degree is independent of \(\kappa\).
The fibre pairing above supplies the parity-dependent factor.
Multiplying by the fibre pairing, the cycle normalizations, and
the height sign \((-1)^\kappa\) gives
\[
 H_{\rm fin}\bigl(Z_\kappa(m_1,\mu_1),Z_\kappa(U)\bigr)
  =(-1)^{\kappa-1}(m_0m_1)^{a/2}\frac{\deg Z(U)}2\,S .
\]

\smallskip\noindent\emph{Cancellation in the global height.}
Put \(A_U=m_0^{a/2}\deg Z(U)/2\). Since
\(Z_\kappa(F_1)=2Z_\kappa(m_1,\mu_1)\), the two evaluations of \(S\)
give
\[
 H_{\rm fin}(Z_\kappa(F_1),Z_\kappa(U))
      =(-1)^{\kappa-1} A_U B_U(F_1).
\]
The oriented Green identity \eqref{arch-signed-comparison} identifies
the archimedean integral with the CM value in the Stokes formula
\eqref{p19g-stokes}. Thus
\[
 H_\infty=(-1)^\kappa A_U
             \bigl(B_U(F_1)-\mathcal L'_U(\xi_wF_1,0)\bigr).
\]
Adding the finite height cancels \(B_U(F_1)\).
Equation \eqref{p19g-central} gives the constant in
Theorem~\ref{merged-height-theorem}.
Linearity in the polar terms gives \eqref{med-untwistedheight} for \(F\).
\end{proof}

The polewise coprimality in \cite[Corollary~6.6]{BEY}, together with
\((D_0,2N)=1\), implies the hypothesis \((D_0,2ND(F))=1\) of
\cite[Theorem~0.6]{DuPeng}. Indeed, removing content \(c\) divides a
quadratic form's discriminant by \(c^2\), so each absolute CM
discriminant in the least common multiple \(D(F)\) divides \(4Nn\)
at its pole. A prime dividing the content can therefore disappear
from \(D(F)\). Improper intersections can occur at the excluded
indices \cite[\S6]{BEY}. Odd codimension \(\kappa>1\) imposes no
condition on \(c_F^+(0,0)\). In even codimension, untwisted component
symmetry forces \(c_F^+(0,0)=0\).

For complex \(F\) we extend the left side linearly in \(F\).
Under this convention the \(L'\)-value on the right is conjugated,
because \(\xi_w\) is conjugate-linear. If \(F\) is weakly
holomorphic, its shadow vanishes. Thus \(H(Z_\kappa(F),Z_\kappa(U))=0\)
for every auxiliary cycle satisfying the hypotheses of
\eqref{med-untwistedheight}.
For fixed \(D_0\) the coprime heights are Fourier coefficients
of a cusp form of weight \(k\). In odd codimension this is
\cite[Theorem~6.12]{BEY}. Corollary~\ref{merged-height-representative}
gives the assertion in both parities.

\subsection{The local height comparison}

Suppose \((\Delta,2N)=1\), taking the genus character to be trivial
when \(\Delta=1\). We use the CM orbit of \eqref{p19-auxiliary-cycle}
with its two orientations.

\begin{proposition}[The two orientations in the archimedean height]
\label{arch-oriented-packets}
Let \(\kappa\ge2\), \(a=\kappa-1\), and take two coefficients
\(Z_\kappa(n,\mu)\), \(Z_\kappa(m,\nu)\) of \eqref{rd-coeff}
with disjoint CM supports on \(X_0(N)\). At a fixed complex embedding
write \(z_x\in\mathfrak H\) for the point defined by \(x\),
\(A_x=2\chi_\Delta(x)/|\Gamma_x|\), and
\(\epsilon_x=\operatorname{sgn}(c)\) when \(c\) is the lower left
entry of \(x\). With \(G_{N,\kappa}\) formed by summing
\(g_\kappa\) over the effective modular group,
\begin{equation}\label{arch-signed-double-sum}
\begin{aligned}
 H_\infty^{\mathbb C}\bigl(Z_\kappa(n,\mu),Z_\kappa(m,\nu)\bigr)
   =\frac{(nm)^{a/2}}{2}
     \sum_{[x],[y]} A_x A_y(\epsilon_x\epsilon_y)^a
                     G_{N,\kappa}(z_x,z_y).
\end{aligned}
\end{equation}
The two sums run over the orbits defining the respective coefficients.
For cycles defined over a common number field \(E\), the absolute
archimedean contribution is the average of \eqref{arch-signed-double-sum}
for the conjugate pairs over all embeddings \(E\hookrightarrow\mathbb C\).
\end{proposition}
\begin{proof}
At \(z_x\), the trace-zero endomorphism acts on the tangent line by
\(cz_x-b/(2N)\), up to a positive factor clearing denominators.
Its imaginary part has the sign of \(c\). The form
\(Nc|z|^2-b\Re z+d\), where the upper right entry of \(x\) is
\(-d/N\), has this same sign on \(\mathfrak H\), since its
discriminant \(b^2-4Ndc\) is negative. Thus
\(\epsilon_x=-\operatorname{sgn}p_z(x)\) in the convention of
Proposition~\ref{p19g-green-kernel}.
With the positive normalization of the projected graph, reversing
the CM embedding multiplies its fibre class by \((-1)^a\)
\cite[Proposition~2.4.1]{Zhang}. For the reference embeddings at two
distinct points, the local height on \(X_0(N)\), in the convention
\(H=(-1)^\kappa\langle\, ,\,\rangle_{\rm GS}\), is
\(\tfrac12G_{N,\kappa}(z_x,z_y)\)
\cite[Proposition~4.1.2, equation~(4.1.1)]{Zhang}. Applying these
two identities to \eqref{rd-coeff} gives \eqref{arch-signed-double-sum}.
The fine-level descent in Zhang's formula gives the height of one
pair of CM cycles. The factors \(A_x,A_y\) are the multiplicities
of this pair in the two coefficients. Averaging over embeddings
uses the absolute arithmetic degree.
\end{proof}

For a fundamental auxiliary discriminant \(D_0\) prime to \(2N\),
with \((\Delta,2N)=1\), we can describe the second sum explicitly.
Choose the reference vector \(x_0\) in the specified residue class
with \(\epsilon_{x_0}=1\), and let \(w_K=|\mathcal O_K^\times|\).
On the upper half-plane put
\[
  \Phi^j_{\rm or}(U,F)
   =\frac{4\chi_\Delta(x_0)}{w_K}
       \sum_{[h]\in\operatorname{Cl}(K)}
          \eta_\Delta(h)\,
          \Phi^j_{\Delta,r}(\gamma_h^{-1}z_{x_0},1,F),
  \qquad h=\gamma_h k_h
\]
with the positive-determinant representatives of
\eqref{p19-translate-lift}. The factor \(4/w_K\) includes the two
auxiliary orientations. For a harmonic input of weight
\(3/2-\kappa\) with real Fourier data and disjoint polar support,
the fixed-component formula becomes
\begin{equation}\label{arch-positive-orbit-height}
 H_\infty^{\mathbb C}\bigl(Z_\kappa(F),Z_\kappa(U)\bigr)
     =\frac{(-1)^\kappa}{2}\,m_0^{a/2}\Phi^j_{\rm or}(U,F).
\end{equation}

Indeed, the involution \(x\mapsto x^\star=JxJ^{-1}\), with
\(J=\operatorname{diag}(-1,1)\), preserves \(b\bmod 2N\).
It sends \(z_x\) to \(-\overline{z_x}\) and
\(\epsilon_x\) to \(-\epsilon_x\). The associated quadratic
form changes from \(Q(X,Y)\) to \(-Q(X,-Y)\), so
\(\chi_\Delta(x^\star)=\operatorname{sgn}(\Delta)\chi_\Delta(x)\).
Simultaneous reflection of both arguments preserves the Green
function, since \(J\) normalizes \(\Gamma_0(N)\). Changing both
summation variables in \eqref{arch-signed-double-sum} therefore
identifies the negative auxiliary contribution with the positive one.

The positive auxiliary vectors form a single ideal-class orbit.
Indeed \(D_0\) is fundamental, and the prescribed root modulo
\(2N\) selects the cyclic level ideal of norm \(N\) in
\(\mathcal O_K\). The CM pairs in this orientation are
\((\mathbb C/\mathfrak a,\mathfrak n^{-1}\mathfrak a/\mathfrak a)\),
one for each ideal class \([\mathfrak a]\), with stabilizer
\(\mathcal O_K^\times\) \cite[I, \S6]{GrossZagier}.
The translate has vector \(\gamma_h^{-1}x_0\gamma_h
=k_h x_0 k_h^{-1}\). Its genus weight is
\(\chi_\Delta(x_0)\eta_\Delta(h)\) by the local character law
preceding \eqref{p19-translate-lift}. For \(F_{n,\mu}\), the
harmonic specialization of Proposition~\ref{p19g-green-kernel}
gives \((-1)^\kappa 2n^{a/2}\) times the signed source Green sum,
and its two reflected polar terms give \(2Z_\kappa(n,\mu)\).
Substitution proves \eqref{arch-positive-orbit-height} for this
series. The harmonic Poincar\'e expansion gives the result for \(F\).

In this prime-to-level setting, \eqref{arch-positive-orbit-height}
also gives the absolute archimedean contribution. Indeed both
weighted orbit sums, with their oriented graph representatives,
acquire the same factor \(\chi_\Delta(\sigma)\) under transport,
as in the proof of Proposition~\ref{coef-galois-character}. Reindex
both sums before evaluating the Green function; the two factors
multiply to one. Together with the simultaneous reflection above,
this identifies the value for the two complete packets at each
embedding. Individual pairs of CM points need not have equal
Green values at conjugate embeddings.

Let \(F\) be a harmonic input of weight \(w=3/2-\kappa\) with
real Fourier data. Suppose that for the chosen \(U\) the theta
lift computes the archimedean height, its Stokes boundary
vanishes, and the total finite contribution is zero. Then the
comparison below gives the mixed height identity
\begin{equation}\label{merged-twisted-height}
 H(Z_\kappa(F),Z_\kappa(U))=C_{\kappa,N,|\Delta|}
 L'\big(\operatorname{Sh}_{m_0,\mu_0}(\xi_wF),\chi_\Delta,\kappa\big).
\end{equation}
We describe when the two corrections can vanish separately.
Section~\ref{sec:newform-heights} uses a comparison of their sum
to choose an auxiliary cycle for a fixed \(F\).
One can choose \(m_0\) larger than every pole order, with \(d_0\)
prime to \(2N\Delta D(F)\) and all primes dividing \(N\) split in
\(\mathbb Q(\sqrt{d_0})\). For a positive nontrivial twist and
odd \(\kappa\), one possible local condition on the twist is
\[
   q\mid\Delta,\qquad q\equiv3\bmod4 \quad\text{for some prime }q .
\]
For negative odd fundamental \(\Delta\) this condition is automatic.
For positive \(\Delta\) and even \(\kappa\) the zero residue of the
polar part vanishes by reflection. In each case the remaining
residues in \eqref{p19g-boundary} and the finite heights must give
the required total cancellation for the same \(U\).
For one Poincar\'e series, when its boundary and finite terms
vanish and the oriented archimedean identity applies, the formula gives
\[
 H(Z_\kappa(n,\mu),Z_\kappa(U))
   =\frac{C_{\kappa,N,|\Delta|}}{2}\,
    L'\bigl(\operatorname{Sh}_{m_0,\mu_0}(\xi_w F_{n,\mu}),
                       \chi_\Delta,\kappa\bigr).
\]
The factor \(1/2\) accounts for the two reflected polar terms.
Choosing a separate \(U\) for each such coefficient need not give
a common \(U\) for their linear combination \(Z_\kappa(F)\).

Consider a polar term of index \(n_1\), with \(D_1=-4N|\Delta|n_1\).
Assume \((D_0,2N)=1\) and \(D_0D_1\) is not a square.
The auxiliary order is maximal because \(D_0\) is fundamental,
and the two CM fields are distinct. At a prime of good reduction
split in \(K=\mathbb Q(\sqrt{D_0})\), the auxiliary elliptic
curve has ordinary reduction with rational endomorphism algebra
\(K\). A common reduction would embed the other CM field in
\(K\), so the two horizontal closures are disjoint there.
At a nonsplit prime one instead counts special endomorphisms. For a
fixed summand choose residues \(r_i^2\equiv D_i\bmod 4N\)
for the two CM embeddings. Put
\(\mathfrak a=(N,(r_0+\sqrt{D_0})/2)\), of norm \(N\),
and write \(\mathfrak d_K\) for the different. A class
\(\delta\in\mathfrak d_K^{-1}\mathfrak a/\mathfrak a\)
specifies the residue attached to the pair of embeddings.
The local count is over special endomorphisms \(\beta\)
of norm \(t\,\mathrm N\mathfrak a\) satisfying
\[
 \beta\in V(E,\iota)\mathfrak d_K^{-1}\mathfrak a,\qquad
 \delta+\beta\in\operatorname{End}(E)\mathfrak a .
\]
Indeed, the CM stack maps to the modular curve by
\((E,\iota)\mapsto(E\to E/E[\mathfrak a])\). A point of its
fibre product with the other CM divisor carries a second embedding
\(\phi:\mathcal O_{D_1}\hookrightarrow\operatorname{End}(E\to E/E[\mathfrak a])\)
respecting the level ideal. In the supersingular quaternion algebra,
decompose \(\phi((r_1+\sqrt{D_1})/2)=\delta+\beta\), taking the
\(K\)-component first. Its trace and the level ideal give
\(2\delta-r_1=n/\sqrt{D_0}\) with \(n\equiv r_0r_1\bmod 2N\).
Squaring \(\phi(\sqrt{D_1})\) gives
\(D_1=n^2/D_0-4\mathrm N(\beta)\), since \(\beta\) acts
antilinearly on \(K\). Conversely, the norm and residue conditions
reconstruct the second embedding
\cite[Lemma~7.10 and proof of Proposition~7.11]{BruinierYangCM}.
For each \(n\equiv r_0r_1\bmod 2N\) with \(n^2<D_0D_1\), the
corresponding norm parameter and residue are
\[
       t(n)=\frac{D_0D_1-n^2}{4N|D_0|},\qquad
       \delta(n)=\frac{n+r_1\sqrt{D_0}}{2\sqrt{D_0}}
                         \pmod{\mathfrak a}.
\]
If no special endomorphism satisfies both conditions for any
summand at \(p\), the horizontal closures are disjoint there.
Here \(V(E,\iota)\) consists of the endomorphisms that intertwine
the chosen CM action with its conjugate. Each object carries its
deformation length divided by its stabilizer order. The higher
cycle contributes \(P_{\kappa-1}(n/\sqrt{D_0D_1})\), where \(n\) is their mixed
pairing in the normalization of \cite[(4.3)]{DuPeng}.
The denominator normalizes both oriented directions. In the untwisted
case with \((D_0,2ND_1)=1\), \cite[Proposition~4.2]{DuPeng}
sums these insertions for both parities of \(\kappa\).
For a twisted pair the auxiliary genus character is the character \(\chi\) of
\(\operatorname{Cl}(K)\) attached to \(D_0=\Delta d_0\). For a finite
idele \(h\) of \(K\), the torus action on the theta kernel contributes
the product of finite Hilbert symbols
\((\Delta,\det h)_{\mathbb A_f}=\chi([h])\)
by Proposition~\ref{p19-genus-translation}. The first cycle also carries its own
genus weight.

\begin{proposition}[Possible primes for a fixed local intersection]
\label{p19-local-support}
Consider the horizontal summands described by the preceding
special-endomorphism model. For a fixed \(n\) put
\(M_n=(D_0D_1-n^2)/(4N)\). A nonempty summand at an inert
prime \(p\) requires \(M_n/p\) to be the norm of an integral ideal
of \(K\). At a ramified \(p\) the same condition holds with \(M_n\)
in place of \(M_n/p\). In particular, if two distinct inert primes
have odd valuation in \(M_n\), the summand has empty support
at every prime in this model. If exactly one inert prime has
odd valuation, it is the only possible prime in that support.
\end{proposition}
\begin{proof}
The congruence on \(n\) gives \(M_n\in\mathbb Z_{>0}\).
For nonsplit \(p\), present the endomorphism algebra as
\(K\oplus Kj\), where \(j\alpha=\bar\alpha j\).
Choose an auxiliary split prime \(p_0=\mathfrak p_0\bar{\mathfrak p}_0\)
with \(\mathrm N(j)=\eta_p\), where \(\eta_p=pp_0\) for an inert
\(p\) and \(\eta_p=p_0\) for a ramified \(p\)
\cite[\S6 and proof of Proposition~7.11]{BruinierYangCM}.
In ideal coordinates on a class-group orbit one has
\[
 V(E_J,\iota_J)=J\bar J^{-1}\mathfrak p_0^{-1}j ,
\]
where \(J\) is a fractional ideal. Translating the CM point by
the ideal \(I\) replaces \(J\) by \(IJ\)
\cite[\S5]{KRYIMRN}. Write \(\beta=\alpha j\). Its membership in
\(V(E_J,\iota_J)\mathfrak d_K^{-1}\mathfrak a\) makes
\[
   \mathfrak c=(\alpha)J^{-1}\bar J\,
        \mathfrak p_0\mathfrak d_K\bar{\mathfrak a}^{-1}
\]
an integral ideal. Since \(\mathrm N(\beta)=Nt(n)\), its norm is
\(\mathrm N(\mathfrak c)=p_0 |D_0| t(n)/\eta_p\).
This is \(M_n/p\) or \(M_n\), respectively. Inert primes have
even valuation in every ideal norm. At a split prime the two
fields cannot embed in the ordinary endomorphism algebra. These
facts prove the support assertions for the horizontal summand.
\end{proof}

The same computation gives the unrestricted norm count in each
ideal class. Write \(\rho(s,C)\) for the number of integral ideals
of norm \(s\) and ideal class \(C\), with value zero if \(s\)
is not a positive integer. The ideal \(\mathfrak c\) above has
class \([J]^{-2}[\mathfrak p_0\mathfrak d_K
\bar{\mathfrak a}^{-1}]\). Conversely each such ideal gives
\(|\mathcal O_K^\times|\) generators \(\alpha\) of the required
norm after removing the fixed ideal factors. Thus the number
of \(\beta\) before imposing the residue condition is
\[
 |\mathcal O_K^\times|\,
 \rho\!\left(\frac{p_0M_n}{\eta_p},
    [J]^{-2}[\mathfrak p_0\mathfrak d_K\bar{\mathfrak a}^{-1}]\right).
\]
Imposing \(\delta(n)+\beta\in\operatorname{End}(E_J)\mathfrak a\)
restricts this count. The actual intersection further weights each
resulting point by its two genus factors and its deformation length.

\begin{proposition}[A sufficient condition for local cancellation]
\label{finite-genus-cancellation}
Keep the preceding CM data, including the two distinct CM fields.
At a nonsplit prime \(p\), let \(\mathcal I_p\) consist of the local
labels in the intersection, including the norm and residue.
Write \(A=\operatorname{Cl}(K)\). For \(\lambda\in\mathcal I_p\)
and \(h\in A\), let \(a_{p,\lambda}(h)\) be the full local weight
after removing the auxiliary genus character \(\chi(h)\). It
includes the first cycle's genus weight, the stabilizer and
deformation factors, and the higher-cycle insertion if present.
Suppose these weights account for the local term \(S_p\) under
consideration, and put \(b_p(h)=\sum_{\lambda\in\mathcal I_p}
a_{p,\lambda}(h)\). Then \(S_p=0\) if and only if
\[
   \sum_{\chi(h)=1} b_p(h)=\sum_{\chi(h)=-1} b_p(h).
\]
A sufficient condition is a permutation of the full set of pairs
\((h,\lambda)\) that preserves \(a_{p,\lambda}(h)\) and reverses
\(\chi(h)\). In particular, one may take
\((h,\lambda)\mapsto(hu,\sigma(\lambda))\) if \(\chi(u)=-1\) and
\(a_{p,\sigma(\lambda)}(hu)=a_{p,\lambda}(h)\) for all pairs.
When \(\sigma\) fixes every label, the signed sum vanishes for
each \(\lambda\).
\end{proposition}
\begin{proof}
Summing over \(\lambda\) gives \(S_p=\sum_{h\in A}\chi(h)b_p(h)\),
the difference of the two displayed sums. The permutation identifies
the positive and negative lists of full weights with multiplicities.
Their sums agree. When it preserves a given label, this argument
already takes place within that label.
\end{proof}

For the two cycles with the same twist, translation by \(u\in A[2]\)
preserves the unrestricted norm count, which depends only on \(h^2\).
To apply the proposition with this translation one also needs
\(\chi(u)=-1\) and invariance of the residue and full weight.
The character may be trivial on \(A[2]\).
An odd \(q\mid\Delta\) inert in \(\mathbb Q(\sqrt{d_0})\) supplies
a class with the required sign. The prime above \(q\) in \(K\)
squares to \((q)\) and has genus value \((d_0/q)=-1\), since
\((d_0,q)=1\). Its class has order two. Invariance of the residue
and weights still requires a local calculation.

For different twists, fix a place \(v\) of a common finite Galois
field of definition containing both quadratic fields. Suppose the
local pairing is bilinear and invariant under the decomposition
group at \(v\), and that the representatives, including all vertical
corrections, transform there by \(\chi_{\Delta_1},\chi_{\Delta_2}\).
If their product character is nontrivial on that group, an element
fixing \(v\) negates the pairing, so the local contribution vanishes.
Proposition~\ref{coef-galois-character} gives the characters of the
generic orbit sums. The transformation of the vertical corrections
must be verified on the chosen model.

\begin{remark}[The role of the residue condition]
\label{nonzero-local-genus-height}\leavevmode
\begin{enumerate}
\item \emph{A single intersection.}
A nontrivial twist and distinct CM fields allow a nonzero local
intersection. On the smooth part of the
\(j\)-line at a good prime, two CM
sections with distinct generic \(j\)-invariants \(a,b\) meet with
length \(v(a-b)\). If exactly one pair of sections meets and both
its genus weights are nonzero, the signed intersection is nonzero.
Stack multiplicities at that point are positive.

\item \emph{The two places above \(71\).}
Take \(N=1,\ \Delta=5,\ r=1,\ \kappa=1\) and the fundamental
discriminants \(D_0=-15,\ D_1=-20\). The corresponding indices
are \((m_0,\mu_0)=(3/4,1)\) and \((m_1,\mu_1)=(1,0)\).
Each discriminant has two CM points, giving the two conjugate
sections over \(\mathbb Q(\sqrt5)\), with nonzero opposite genus
weights. The singular moduli in \cite[\S1]{DukeCMValues} give
\begin{align*}
 H_{-15}(X)&=X^2+191025X-121287375,\\
 H_{-20}(X)&=X^2-1264000X-681472000 .
\end{align*}
The prime \(71\) splits in \(\mathbb Q(\sqrt5)\) and is inert in
\(\mathbb Q(\sqrt{-15})\). At either place above \(71\) the two
root sets reduce to \(\{41,66\}\) and \(\{17,41\}\).
The only common root is \(41\), away from \(0\) and \(1728\).
The resultant of the two displayed polynomials has \(71\)-adic
valuation one. The unique meeting pair therefore has length one.
Conjugation interchanges both pairs of genus weights, so the
two places give the same sign in the absolute local intersection.

\item \emph{Scope of the calculation.}
Here \(m_0=3/4<m_1=1\), and \(\Delta=5\) has no divisor
congruent to \(3\bmod 4\). The calculation concerns horizontal
divisors at a specified prime. In the height formula the auxiliary
cycle is chosen subject to additional conditions, and the total
height includes its archimedean term.
\end{enumerate}
\end{remark}

If disjointness at the good split primes and the signed cancellation
of Proposition~\ref{finite-genus-cancellation} at the nonsplit primes
account for all horizontal terms, and any remaining bad-prime and
vertical contributions vanish or cancel, then
\(H_{\rm fin}(Z_\kappa(F),Z_\kappa(U))=0\).
In the Petersson integral \(\mathcal L_U(g,s)\) of
Proposition~\ref{p19g-cm-derivative} the centre is \(s=0\). The Dirichlet
series in Proposition~\ref{p19g-shimura} is evaluated at \(s+\kappa\).

\begin{proposition}[Global height with its boundary and finite terms]
\label{p19-height-comparison}
Use the hypotheses of Propositions~\ref{p19g-cm-derivative}
and~\ref{p19g-shimura}, including \((\Delta,2N)=1\) and
\(\mathbb Q x_0\cap L^\sharp=\mathbb Zx_0\). Let \(F\) have real Fourier
data and polar support disjoint from the auxiliary CM cycle.
For \(D_0\) prime to \(2N\), \eqref{arch-positive-orbit-height} gives
the following archimedean comparison for both parities. For the
other auxiliary data assume that
\[
 H_\infty(Z_\kappa(F),Z_\kappa(U))=
    \frac{(-1)^\kappa}{2}\,m_0^{(\kappa-1)/2}
             \Phi_{\Delta,r}^{\,j}(Z(U),F).
\]
Take \(B_U(F)\) with the normalization per unit of CM mass in
\eqref{p19g-stokes}.
Let \(\mathcal P\) contain the finite primes with potentially
nonzero contributions on the arithmetic model. Then
\begin{equation}\label{height-residual-finite}
\begin{aligned}
 H(Z_\kappa(F),Z_\kappa(U))
   &=C_{\kappa,N,|\Delta|}\,
       L'\bigl(\operatorname{Sh}_{m_0,\mu_0}(\xi_wF),\chi_\Delta,\kappa\bigr)\\
   &\quad+\frac{(-1)^\kappa m_0^{(\kappa-1)/2}\deg Z(U)}{2}\,B_U(F)
             +\sum_{p\in\mathcal P}H_p(Z_\kappa(F),Z_\kappa(U)).
\end{aligned}
\end{equation}
\end{proposition}
\begin{proof}
Substitute \eqref{p19g-stokes} in the archimedean height. The boundary
term is the displayed multiple of \(B_U(F)\). For the derivative term
use \eqref{p19g-central} under its primitive-line hypothesis. The unsigned mass \(\deg Z(U)\) and the factor
\(m_0^{(\kappa-1)/2}\) cancel. The signs multiply to
\(-(-1)^\kappa(-1)^{\kappa-1}=1\), leaving the coefficient
\(C_{\kappa,N,|\Delta|}\) of the central derivative.
Adding the finite heights gives \eqref{height-residual-finite}.
\end{proof}
Each \(H_p\) includes the local logarithm and is divided by the
field degree. Disjoint generic supports make the sum finite.
For the chosen \(U\), equation~\eqref{merged-twisted-height} holds
precisely when this sum cancels the displayed boundary term.
Contributions at different places can cancel. In the untwisted
case \eqref{med-untwistedheight} gives this cancellation under its
hypotheses. For overlapping supports one would also need the
regularized diagonal Green value and the arithmetic self-intersection.

\subsection{Auxiliary cycles and the twisted mixed height}\label{sec:newform-heights}

\noindent\emph{The local comparison used for twists.}
Write \(h=(\kappa-1)/2\), \(D_0=\Delta d_0\), and
\(\delta=4h_K/w_K\) for the mass of the two oriented CM orbits of
\(K=\mathbb Q(\sqrt{D_0})\). We use a fundamental \(D_0\), a
primitive auxiliary line, and polar CM fields distinct from \(K\).
The arithmetic height uses the field-normalized intersection of
admissible extensions of the two genus-weighted cycles. For CM
cycles and the behavior of these extensions under change of field
and level, see \cite[\S\S2.2.3--2.2.5]{QiuCM}.

Fix a positive perpendicular norm in \eqref{p19g-boundary}.
Write \(\operatorname{Diff}\) for the places where the incoherent
binary space does not represent it. Two vanishing local factors
force the derivative of their product to vanish. If
\(\operatorname{Diff}=\{p\}\), the only possible contribution is
\(W'_p(0)\prod_{v\ne p}W_v(0)\). The geometric term counts
representations in the coherent space obtained by switching at \(p\).
It uses the local intersection weight at \(p\) and the original
residue conditions elsewhere. Decomposing the joint residue
condition into a finite sum of product functions gives this
description term by term.

At a supersingular point let \(R\) be the maximal order of the
quaternion division algebra over \(\mathbb Q_p\). Fix the maximal
order \(O\) of a quadratic field \(K_p\) embedded in \(R\).
For an integral endomorphism \(\alpha=A+Yj\), write \(A,Y\in K_p\),
\(jz=\bar z j\), and suppose \(Y\ne0\).
Put \(b=\mathrm{Nrd}(Yj)\), and let \(\varpi\) be a uniformizer of \(O\).
On the canonical lift of this \(O\)-action, the largest \(s\)
for which \(\alpha\) lifts modulo \(\varpi^s\) is
\[
 s=
 \begin{cases}
   (v_p(b)+1)/2,& K_p/\mathbb Q_p\ \text{unramified},\\
   v_p(b)+v_p(\operatorname{disc}O),& K_p/\mathbb Q_p\ \text{ramified}.
 \end{cases}
\]
Indeed, Gross's lifting theorem gives the endomorphism order
\(O+\varpi^{s-1}R\) at this thickness \cite[Theorem~1.4]{WewersCanonical}.
In the unramified case choose \(j\) with \(v_p(\mathrm{Nrd}(j))=1\).
Then \(R=O\oplus Oj\) and \(v_p(b)=2v_{K_p}(Y)+1\).
In the ramified case choose \(j\) with \(v_p(\mathrm{Nrd}(j))=0\)
and put \(f=v_p(\operatorname{disc}O)\). The two blocks of the trace
pairing on \(O\oplus Oj\) each have discriminant valuation \(f\).
The trace discriminant of \(R\) has valuation two, so
\([R:O\oplus Oj]=p^{f-1}\). Since \(R\cap K_p=O\), projection onto
\(K_pj\) gives \(\operatorname{pr}_-(R)=\varpi^{1-f}Oj\).
Thus \(\alpha\) lifts exactly when \(v_{K_p}(Y)\ge s-f\), and
\(v_{K_p}(Y)=v_p(b)\). This calculation includes the ramified
quadratic extensions of \(\mathbb Q_2\). At primes dividing the
level one must also impose the cyclic subgroup condition. The
signed intersection uses both genus characters and the
multiplicities of the CM orbits.

For the comparison take \(F\) with real Fourier data and a fundamental
CM datum \(U\) admissible for \(F\). The congruences
\(D_0=\Delta d_0\equiv(r\mu_0)^2\bmod 4N\) and \((d_0,2N\Delta)=1\)
determine the local type at every \(p\mid N\). If \(p\nmid\Delta\),
this is a unit square modulo \(p\), or modulo \(8\) at \(p=2\), so
\(p\) splits in \(K\). If \(p\mid\Delta\), it ramifies in \(K\)
and \(p\Vert N\). Thus among the level primes only these ramified
primes can support a supersingular intersection of the two CM
closures. The full arithmetic comparison includes vertical terms
and uses both genus characters and every residue class in the
theta kernel. Choose the reference vector \(x_0\) with leading
coefficient \(-1\), of genus value \(\operatorname{sgn}(\Delta)\).
With the same unary orientation as in the CM integral, define
\(\iota_U^{\rm ref}\) by the transfer of Section~\ref{p19-cm-section}
with weight \(\chi_\Delta(\delta')\) at each discriminant label
\(\delta'\), retaining every joint residue. Set
\[
 B^{\rm ref}=\operatorname{CT}
      \langle\iota_U^{\rm ref}F^+,[\vartheta,\mathcal E_B]_j\rangle .
\]
The boundary is normalized per unit of CM mass. The archimedean
comparison uses the Green integral and the same oriented transfer in
the Stokes formula. In our normalization these give
\begin{equation*}\tag{AC}\label{arch-signed-comparison}
\begin{aligned}
 H_\infty(Z_\kappa(F),Z_\kappa(U))
    &=\tfrac12(-1)^\kappa m_0^h\Phi_{\Delta,r}^{\,j}(Z(U),F),\\
 B_U(F)&=\chi_\Delta(x_0) B^{\rm ref}.
\end{aligned}
\end{equation*}
We verify both equalities for every admissible fundamental \(U\).
The maximal order \(\mathcal O_K\) has a unique level ideal for the
given root and orientation. At a split prime the root specifies the
ideal. At a shared prime \(p\Vert N\), including \(p=2\), this is the
unique ramified prime of norm \(p\). Each orientation consequently
has one orbit of the ideal class group with stabilizer
\(\mathcal O_K^\times\), giving the mass \(\delta=4h_K/w_K\).

\smallskip\noindent\emph{The two Green values.}
The reference \(x_0\) has negative lower left entry. At its upper
half-plane fixed point \(p_z(x_0)>0\), so \(z_U^+\), chosen by
\(p_{z_U^+}(x_0)=-\sqrt{m_0}\), is the other orientation.
Conjugation by \(J=\operatorname{diag}(-1,1)\) takes this point to
the upper half-plane together with its conjugate level ideal and its
positive Lie orientation. The determinant law of
Section~\ref{p19-translated} gives
\(\eta_\Delta(J)=\prod_{p<\infty}(\Delta,-1)_p=\operatorname{sgn}(\Delta)\).
Multiplication by \(\chi_\Delta(x_0)\) thus gives weight one to the
transported branch. The branch at the original upper half-plane
point has weight \((-1)^{\kappa-1}\chi_\Delta(x_0)\).
These are the two weights in \eqref{arch-signed-double-sum}.
The Poincar\'e series contributes \((-1)^\kappa 2n^h\) to its
Green value and \(2Z_\kappa(n,\mu)\) to its cycle. With the factor
\(1/2\) for a pair of fibre cycles, this gives the first equality
in \eqref{arch-signed-comparison}. Proposition~\ref{coef-galois-character}
and the reindexing following \eqref{arch-positive-orbit-height} give
the absolute archimedean contribution for the complete packets.

\smallskip\noindent\emph{The Stokes boundary.}
The transfer to \(\Lambda\) sums over \(\delta'\mapsto r\mu\) with
\(Q_\Delta(\delta')\equiv\operatorname{sgn}(\Delta)Q(\mu)\bmod\mathbb Z\).
Weighted by \(\chi_\Delta(\delta')\), it intertwines the two Weil
representations for every \(\Delta\equiv r^2\bmod4N\)
\cite[Proposition~3.2]{AETraces}. Restriction to \(P\oplus B\)
retains all joint cosets above each \(\delta'\). The torus fixes the
unary line. On the binary space average first over a compact open
subgroup of \(\widehat{\mathcal O}_K^\times\) fixing all the individual
cosets. The full weighted sum is invariant under
\(\widehat{\mathcal O}_K^\times\). On the torus the determinant
character is the genus character of the coprime fundamental
factorization \(D_0=\Delta d_0\), hence trivial on units.
The normalized average is therefore independent of this refinement.
Applying binary Siegel--Weil to each Schwartz function gives the same
mass \(\delta\) \cite[Proposition~2.2 and Lemma~4.6]{BruinierYangCM}.
The Stokes argument of Proposition~\ref{p19g-cm-derivative} now
applies to the full transfer. Its factor \(\chi_\Delta(x_0)\)
multiplies the entire sum and gives the second equality of
\eqref{arch-signed-comparison}. We unfold the corresponding Rankin
term below with the same residue labels.

Put \(\varepsilon=\operatorname{sgn}(\Delta)(-1)^{\kappa-1}\).
For the maximal auxiliary order, the comparison with the Stokes
boundary takes the following form. Since \((d_0,2N\Delta)=1\),
the auxiliary order is maximal and each shared prime has two source
labels. We sum over these labels before comparing the finite height
of the whole Fourier coefficient with the Stokes boundary. The
identity to be established is
\begin{equation*}\tag{FC}\label{finite-signed-comparison}
 \sum_{p<\infty} H_p(Z_\kappa(F),Z_\kappa(U))
     =\varepsilon\,\frac{m_0^h\delta}{2}\, B^{\rm ref}.
\end{equation*}
We describe the normalization of the local terms required in this identity.
Fix the joint residues, the trace of the two endomorphisms and a
positive perpendicular norm. Write \(k_p\) for its contribution
to the binary Eisenstein derivative. If the difference set is
\(\{p\}\), let \(E^{(p)}\) be the Fourier coefficient of the coherent
series with the local space switched at \(p\). In the coherent
Siegel--Weil formula the ideal-class average is
\(E^{(p)}=(2/h_K)\sum_{[\mathfrak a]}\theta_{[\mathfrak a]}\),
with the prescribed local test functions. A fixed residue embedding
has effective automorphism weight \(2/w_K\). Thus its geometric
representation mass is \(\delta E^{(p)}/4\).

The remaining local input is the comparison of the signed
intersection multiplicity, including the residue degree, with
the derivative of the local Whittaker function. With the above mass,
a termwise comparison requires the fixed-trace intersection on the base
to equal \(-\operatorname{sgn}(\Delta)\delta k_p/2\).

For instance, if \(2\mid(N,\Delta)\), the comparison uses the ramified
quadratic field \(K_2\) and the unique auxiliary subgroup
\(E[\mathfrak p_2]\). Suppose a shell of depth \(l\) and sign \(s_b\)
has binary coefficient \(-s_b(l+1)E^{(2)}\log 2/2\).
The required field-normalized signed length for both CM types together
is then \(\chi_{\delta_2}(-1)s_b(l+1)\), where \(\delta_2\) is the
dyadic prime discriminant in \(\Delta\). The sign follows from
\(\prod_{q\mid\Delta}\chi_{\delta_q}(-1)=\operatorname{sgn}(\Delta)\).
If the two strict sections meet only the central exceptional line of
the resolved node and have the same projected graph \(S_\kappa\),
their signed specialization there is
\(([P_+]-[P_-])\times S_\kappa\). This class vanishes because the
points have the same class on \(\mathbb P^1\). The global argument
below proves \eqref{finite-signed-comparison} after summation.

For a polar index \(n\), the oriented fibre intersection supplies
\(-P_{\kappa-1}(t_{\rm or})\), with \(t_{\rm or}\) the normalized
oriented trace. The unary insertion gives
\((-1)^{\kappa-1}n^h P_{\kappa-1}(t_{\rm or})\). Multiplication by
\((m_0n)^h\) on the geometric side and summation over the polar
terms give \eqref{finite-signed-comparison} whenever all the local
comparisons hold with the stated multiplicities. A positive polar index of perpendicular norm zero
would be \(a^2m_0\) for \(a\in\mathbb Q\) and would have auxiliary
CM field \(K\), excluded by admissibility. At total index zero the
positive bracket order in odd codimension, or the odd unary theta
series in even codimension, kills the constant. By \eqref{arch-signed-comparison} the archimedean
height contains the opposite boundary term.

To unfold the remaining integral, take the primitive vector \(x_0\)
with leading coefficient \(-1\) in this convention. The pure unary
vectors are \(a x_0\). When \((a,\Delta)=1\), multiplication by \(a\)
identifies the residue labels over \(x_0\) with those over \(a x_0\):
at a shared prime it permutes the possible roots, and at each
nonshared prime the specified residue determines the root. When
\((a,\Delta)>1\), the genus weight is zero. Thus, apart from the
reference genus value, the coefficient of \(a x_0\) is
\(\chi_\Delta(a)\sum_\nu c_{\xi_wF}(a^2m_0,a\nu)\), summed over the
original labels above \(r\mu_0\). Apply the divisor convolution in
\eqref{merged-shimura-coeff} to this finite sum. Its denominator is
\(L(\chi_\Delta\chi_{d_0},1+s)=L(\chi_{D_0},1+s)\), with the zero
Euler coefficients retained. The two unary orientations and the
class number formula \(\delta=2\sqrt{|D_0|}L(\chi_{D_0},1)/\pi\)
then give the constant of \eqref{p19g-central}, with the reference
sign included in the integral. For a real input write \(S_{\rm tr}\)
for the sum of the Shimura lifts at the residues
\(r\nu\equiv r\mu_0\bmod 2N\), \(\nu^2\equiv d_0\bmod4N\).
Substitution in the oriented Stokes evaluation gives the exact difference
\begin{equation}\label{total-height-discrepancy}
\begin{aligned}
 H(Z_\kappa(F),Z_\kappa(U))-
 C_{\kappa,N,|\Delta|}L'(S_{\rm tr}(\xi_wF),\chi_\Delta,\kappa)
   &=\sum_{p<\infty}H_p(Z_\kappa(F),Z_\kappa(U))\\
   &\quad-\frac{\varepsilon m_0^h\delta}{2}B^{\rm ref}.
\end{aligned}
\end{equation}
The same unfolding before differentiation replaces the numerator of
\eqref{p19g-shift} by
\(L(\sum_\nu\operatorname{Sh}_{m_0,\nu}(g),\chi_\Delta,s+\kappa)\)
for real \(g\), with the same residue fibre.
Both \(A_\kappa(0)\) and \(L(\chi_{D_0},1)\) are nonzero.
The incoherent Eisenstein series vanishes at zero for each finite
Schwartz function, so the numerator vanishes there as well.

Thus the mixed-height identity at an admissible \(U\) is equivalent
to the equality in \eqref{finite-signed-comparison} once the trace
is identified with the residue projection. Vanishing of each local
height would be a stronger condition.

\begin{auxlemma}[The geometric projector on a newform line]
Let \(\kappa\ge2\). At every prime \(p\nmid 2N\Delta\), the geometric
Hecke operator acts on the series \(\Phi_\kappa\) as the Hecke operator on
its modular variable. If the simultaneous eigenspace in \(S\) for the
good eigenvalues of \(G\) is \(\mathbb Cg\), there is a polynomial
\(e_G\) in the geometric good Hecke operators which is a projector on
\(W_N\) and satisfies
\[
       e_G X_g=X_g,\qquad e_G Z_\kappa(n,\mu)=c_g(n,\mu) X_g .
\]
It is self-adjoint for a height on \(W_N\) satisfying the
projection formula for these correspondences.
\end{auxlemma}
\begin{proof}
Put \(a=\kappa-1\), \(d=-\operatorname{sgn}(\Delta)4Nn\), and
\(D=\Delta d\). We check the oriented version of the divisor count
in \cite[(7.3)]{BO}. A cyclic \(p\)-isogeny corresponds to a
sublattice of index \(p\). Fix a target quadratic form \(Q'\) of
discriminant \(p^2D\). Its occurrences in the correspondence are
indexed by the superlattices of index \(p\) on which \(Q'\) is
integral. If \(Q'\) is primitive at \(p\), its local diagonal form
is \(\langle u,p^2v\rangle\) with \(u\) a unit. There is one such
superlattice. For \(Q'=pR\) with \(R\) primitive at \(p\), the
choices are the isotropic lines of \(R\bmod p\), numbering
\(1+(\frac{D}{p})\). If \(p^2\) divides all coefficients of \(Q'\),
all \(p+1\) lines occur.

The triangular substitutions of determinant \(p\) preserve the sign of
a definite form and multiply its middle residue modulo \(2N\) by \(p\).
They also preserve the represented units defining the genus character
at every prime dividing \(\Delta\). Indeed, for a source form
\([NA,B,C]\) and \(N=N_1N_2\), a substitution
\((x,y)\mapsto(hx+by,jy)\), \(hj=p\), acts on the auxiliary form
\([N_1A,B,N_2C]\) by \((x,y)\mapsto(hx+N_2by,jy)\),
which is invertible at these primes. Since
\(\chi_\Delta(pR)=(\frac{\Delta}{p})\chi_\Delta(R)\), the extra
middle contribution has symbol \((\frac{d}{p})\).
These counts retain the stabilizer weights in \eqref{rd-coeff} by
orbit--stabilizer, separately for both orientations.

Transport of the normalized fibre class along each isogeny contributes
\(p^a\) \cite[(2.4) and Lemma~2.1]{XueGKZ}. Inserting \(n^{a/2}\)
now gives
\[
 T_pZ_\kappa(n,\mu)=Z_\kappa(p^2n,p\mu)
 +p^a\Big(\frac{d}{p}\Big)Z_\kappa(n,\mu)
 +p^{2a+1}Z_\kappa(n/p^2,\mu/p),
\]
with an inadmissible index read as zero. This is the Fourier
coefficient formula \cite[(7.1)]{BO}, in the sign convention of
\eqref{rd-coeff}. The coefficient map \(S^\vee\to W_N\) from
Theorem~\ref{med-target} therefore intertwines the two actions.
The good operators on \(S\) commute and are self-adjoint.
Their joint cuspidal eigensystems come from classical newforms in
both Weil types. To see this, descend a nonzero eigenvector by
the adjoint of an isotropic lift whenever its image is nonzero.
This decreases the index and preserves the eigenvalues away from
\(2N\). At the last level the vector belongs to \(S_*\) of
Lemma~\ref{fundamental-separation}. A nonzero unit fundamental
coefficient is the first coefficient of its Shimura lift, a
classical cusp form with those good eigenvalues
\cite[Proposition~1]{SkoruppaBinary}. This argument suffices even
when the multiplicity in the original space is greater than one.
A real polynomial in finitely many of them selects their joint
\(G\)-eigenspace and descends to \(W_N\). The normalized shadow
expansion \eqref{rd-basis} gives the two assertions about its image.
Transposition reverses each isogeny, giving the dual cyclic isogeny
with the same fibre factor. The height projection formula thus gives
self-adjointness for \(T_p\), and hence for \(e_G\). If the good
eigenvalues have multiplicity greater than one, this argument selects
their full joint eigenspace.
\end{proof}

\begin{auxlemma}[Admissible indices and the finite comparison]
Every \(F\) admits infinitely many fundamental auxiliary indices in
the sense of the introduction. At each admissible index the
mixed-height identity \eqref{intro-newform-height} holds for
\(\Delta=1\). For a general compatible \(\Delta\) and real \(F\),
it is equivalent to \eqref{finite-signed-comparison}.
\end{auxlemma}
\begin{proof}
\smallskip\noindent\emph{One choice for all polar terms.}
Fix \(F\) and put \(\epsilon=-\operatorname{sgn}(\Delta)\).
Choose a prime \(q\equiv\epsilon\bmod 4N\) outside the finite set
of primes dividing \(2N\Delta\), a polar \(4Nn\), or the discriminant
of a CM field in the polar support. Dirichlet's theorem gives
infinitely many choices. Set
\[
 d_0=\epsilon q,\qquad D_0=\Delta d_0,\qquad
 \mu_0=1\bmod 2N,\qquad m_0=q/(4N).
\]
The odd prime discriminant \(d_0\) is \(1\bmod 4N\).
The product \(D_0<0\) of the coprime fundamental discriminants
\(\Delta,d_0\) is fundamental, also when \(\Delta\) is even.
Put \(c=(r^2-D_0)/(4N)\in\mathbb Z_{>0}\). The vector
\[
 x_0=\begin{pmatrix}r/(2N)&1/N\\ -c&-r/(2N)\end{pmatrix}
 \quad\text{satisfies}\quad Q(x_0)=-D_0/(4N)=|\Delta|m_0 .
\]
Its middle label is \(r\mu_0=r\bmod 2N\), and its upper-right entry
makes it primitive in \(L^\sharp\). Its fixed point corresponds to
the positive definite form \([Nc,r,1]\). In the orientation of
the finite comparison the auxiliary form \([-1,r,-Nc]\) represents
\(-1\), so its genus value is \(\operatorname{sgn}(\Delta)\). The CM field ramifies at \(q\)
and differs from every field in the polar support. Distinct \(q\)
give distinct fields. Excluding the polar data of a basis gives one
infinite family for any fixed finite-dimensional space of inputs.

For \((\Delta,2N)=1\) these indices satisfy the hypotheses on the
auxiliary cycle in Proposition~\ref{p19-height-comparison}. One can
also require \(2\) to split in \(K\) by imposing \(d_0\equiv\Delta\bmod8\).
This agrees with \(d_0\equiv1\bmod4N\) modulo four for odd \(N\),
and modulo eight for even \(N\). The level primes already split by
the congruence for \(D_0\).

\smallskip\noindent\emph{The untwisted height at every admissible index.}
If \(\Delta=r=1\), admissibility gives \((D_0,2N)=1\) and
\((D_0,4Nn)=1\) for every polar index \(n\). For real \(F\),
Theorem~\ref{merged-height-theorem} now gives the height directly.
For another root \(r^2\equiv1\bmod4N\), multiplication by \(r\)
is an involutive isometry of \(L^\sharp/L\). Relabel \(F\) by
this isometry and \(U\) by \((m_0,r\mu_0)\). Both orbit sums become
the ones for \(r=1\), and the coefficient at \(ar\mu_0\) of the
relabeled shadow is the original coefficient at \(a\mu_0\).
Thus the Shimura convolution \eqref{merged-shimura-coeff} is unchanged.

\smallskip\noindent\emph{The trace at a shared prime.}
For a fundamental index with \((d_0,2N\Delta)=1\), the original
residue labels above \(r\mu_0\bmod2N\) form the set
\[
 \mathcal T(U)=\{\nu\bmod 2N:
       r\nu\equiv r\mu_0\bmod 2N,\quad \nu^2\equiv d_0\bmod 4N\}.
\]
At a nonshared odd prime dividing \(N\), \(r\) is a unit and the
first congruence fixes the label. At an odd prime
\(p\mid(N,\Delta)\), compatibility gives \(p\Vert N\) and \(p\mid r\).
The two roots are \(\mu_0,-\mu_0\bmod p\). If \(2\mid(N,\Delta)\)
then \(2\Vert N\), and the two odd classes modulo four give the
two labels. If \(2\mid N\) and \(2\nmid\Delta\), then \(r\) is odd and fixes the label.
For odd \(N\), the parity of the label follows from the square
condition modulo four. The Chinese remainder theorem thus
identifies \(\mathcal T(U)\) with the orbit of \(\mu_0\) under the
\(t\) independent involutions of Section~\ref{rd-cycles}.
By the coefficient formula \eqref{merged-shimura-coeff},
\[
 \sum_{\nu\in\mathcal T(U)}\operatorname{Sh}_{m_0,\nu}(\xi_wF)
     =2^t\operatorname{Sh}_{m_0,\mu_0}(\Pi_{\Delta,N}\xi_wF)
     =2^t S_U(\Pi_{\Delta,N}F).
\]
This counts the inputs to the Shimura lift. Since \(D_0\) is
fundamental, the auxiliary order is maximal and its conductor
introduces no further multiplicity. Substituting this trace into
\eqref{total-height-discrepancy} shows that \eqref{intro-newform-height}
is equivalent to \eqref{finite-signed-comparison} for real inputs.
For complex \(F\) put \(JF(\tau)=\overline{F(-\bar\tau)}\) in
each component. Fourier conjugation preserves the input type and
commutes with \(\xi_w\) (with the same conjugation on its target).
In \(F=F_1+iF_2\), take \(F_1=(F+JF)/2\) and
\(F_2=(F-JF)/(2i)\). Their polar coefficients are the real and
imaginary parts of those of \(F\), so the same \(U\) is admissible
for both. The height is complex linear in \(F\), and the shadow is
conjugate linear. Applying the real identity twice therefore gives
the complex conjugate in \eqref{intro-newform-height}.
In the coprime range, comparison with \eqref{height-residual-finite}
identifies the difference from the single derivative with \(R_U(F)\).
\end{proof}

If \(g\) is a primitive Jacobi newform, one may also require
\(b_g(U)\ne0\). In either signature \(g\) is orthogonal to the
images of \(\phi(\tau,z)\mapsto\phi(\tau,pz)\) from index \(N/p^2\)
\cite[\S5]{ABS}. Under theta decomposition these are the isotropic
lifts in Lemma~\ref{fundamental-separation}. Applying that lemma
with the primes of the polar data excluded gives infinitely many
admissible \(U\) with the required coefficient. At squarefree
level the same reasoning applies to every nonzero \(g\).
Section~\ref{sec:primitive-pairings} computes the height first on
primitive lines. Theorem~\ref{positive-full-mixed} treats the
remaining multiplicities by comparison with the untwisted formula.

\subsection{The projected self-height at exact level}\label{sec:projected-exact-height}
We compute the self-height of a fundamental coefficient on the
\(G\)-component. The split case fixes the global constant. A local
period at each common prime then determines the ramified case.

\begin{lemma}[Comparing CM periods on a fixed Hecke component]
\label{finite-toric-tests}
Fix an imaginary quadratic field \(K\), a quadratic ring class
character \(\Omega\) trivial on the ideles of \(\mathbb Q\), and the
unitary representation \(\pi\) of a normalized newform \(G\) of
weight \(2\kappa\) and trivial character. We allow \(\Omega=1\).
Assume \(L(1/2,\pi_K\otimes\Omega)=0\). Consider finitely many
pairs of \(\Omega\)-weighted CM orbit sums, pulled back to a common
fine level with the level degree normalization. Assume the Hecke
modules used for each pair, in Chow or modulo the radical of its
pairing, are finite sums of good eigenspaces attached to weight
\(2\kappa\) cuspidal automorphic representations. Let
\(e_G\) select their common good eigenvalues belonging to \(G\).
Then the projected heights are proportional, with the same scalar,
to the normalized toric periods of their convolution tests. They
are all zero if \(L'(1/2,\pi_K\otimes\Omega)=0\).
\end{lemma}
\begin{proof}
Work in the quotient of the CM orbit module by the radical of its
height pairing. Passage to this quotient preserves the height and
commutes with the good Hecke operators. Choose a finite set \(S\)
containing the primes of the levels, the ramified primes of \(K\),
the nonspherical places of the tests, and \(2\), enlarged so that
\(|S|\ge3\). The unramified Hecke decomposition
\cite[Theorem~2.3.4]{QiuCM} is a direct sum on this quotient.
The finitely many orbit sums therefore involve only finitely many
eigensystems. Strong multiplicity one separates them away from \(S\).
A polynomial in good Hecke operators selecting \(G\) on the given
modules acts on their images as Qiu's \(\pi\)-projection
\cite[(2.18)--(2.20)]{QiuCM}.

Let \(f,h\) be the two orbit tests, with \(h^\vee(g)=h(g^{-1})\).
The height of their translates is the height distribution at
\(F=f*h^\vee\) by the projection formula \cite[(2.14)]{QiuCM}.
The assumption \(L(1/2,\pi_K\otimes\Omega)=0\) gives precisely the
\(\pi\)-module identity required by Qiu's projector formula
\cite[Theorem~2.3.7]{QiuCM}. The smooth operator
\(\widetilde\pi^\infty(F^\vee)\) has finite rank: it is supported on
the finite-dimensional invariant spaces for sufficiently small compact
open subgroups on each side. Decompose it into rank-one tensors.
Formula~(2.22) of \cite[Theorem~2.4.3]{QiuCM} gives the height of
each tensor as its normalized toric period times a scalar depending
only on \(K,\Omega,\pi\) and the measures. By linearity the same
holds for their sum. The scalar contains
\(L'(1/2,\pi_K\otimes\Omega)\), proving the last assertion.
\end{proof}

\begin{lemma}[Split-level normalization of the primitive height]
\label{split-primitive-absolute}
Let \(G\) be a normalized newform of exact level \(N\) and weight
\(2\kappa\), with \(\kappa\ge2\). Let \((N,\Delta)=1\), and take a
compatible fundamental index \(a_U=(m,\mu)\) with \(d\) prime to
\(2N\Delta\), \(D=\Delta d<-4\) fundamental and \(m=|d|/(4N)\).
Suppose \(\varepsilon(G)=-\operatorname{sgn}(\Delta)\).
For the associated genus character \(\Omega\) of
\(K=\mathbb Q(\sqrt D)\),
\[
 H(e_GZ_\kappa(U),e_GZ_\kappa(U))
  =\frac{(2\kappa-2)!\sqrt{|D|}\,m^{\kappa-1}}
         {2^{4\kappa-2}\pi^{2\kappa}\|G\|_{\rm Pet}^2}
       L'_K(G,\Omega,\kappa).
\]
This includes both signs of the twist. For this height identity,
\(e_G\) may be taken on the CM orbit module modulo the radical of
its height pairing.
\end{lemma}
\begin{proof}
Every prime dividing \(N\) splits in \(K\): for an odd level prime
use \(D\equiv(r\mu)^2\bmod 4N\) and the unit condition, and for
\(2\mid N\) use \(D\equiv1\bmod8\). Put \(h=|\operatorname{Pic}
(\mathcal O_K)|\) and \(u_K=|\mathcal O_K^\times|/2\). For one
choice of orientation let \(C\) be the sum of the normalized fibre
cycles over this class group with weights \(\Omega\). Each ideal
class occurs once on the chosen branch by
\cite[I, \S1]{GKZ}. The level ideal has genus character one because
\(\Omega_p=\chi_{d,p}\circ\mathrm N\) is trivial at each \(p\mid N\).
The Fricke involution \(W_N\) exchanges the two roots of the
oriented level ideal. Fibre transport by the isogeny preserves the
normalized CM class for the transported complex embedding
\cite[Proposition~2.4.1]{Zhang}. Negating the CM vector multiplies
the signed fibre class by
\(\epsilon=\operatorname{sgn}(\Delta)(-1)^{\kappa-1}\).
Up to the sign of the initial genus representative, the full
coefficient in \eqref{rd-coeff} is therefore
\(m^{(\kappa-1)/2}(1+\epsilon W_N)C/u_K\).

The central value of the base change is zero: at coprime level
\(\varepsilon(G\otimes\chi_\Delta)
=\varepsilon(G)\chi_\Delta(-N)=-1\). To compare the two
orientations in the height, use Lemma~\ref{finite-toric-tests}.
Both \(C\) and \((1+\epsilon W_N)C\) are linear combinations of
translates of the same \(\Omega\)-weighted orbit. In the quotient by
the radical of the height, each has finite cuspidal Hecke support
\cite[Theorem~2.3.4 and \S2.4.2]{QiuCM}. We can thus compare their
\(G\)-projections in this quotient, without a Chow projector. At exact level
\(\pi^{K_0(N)}\) is the newvector line. The normalized Fricke
operator acts there by
\(w_N=(-1)^\kappa\varepsilon(G)=\epsilon\).
If \(P\) denotes averaging over \(K_0(N)\), the toric kernels of
the two self-pairings differ by replacing \(P\) by
\((1+\epsilon W_N)P(1+\epsilon W_N)^\vee\). In \(\pi\) the latter
is \(4P\). Translation from the chosen CM embedding conjugates both
kernels by the same element. The lemma therefore gives a factor four
also for the heights, including when the derivative is zero.

For the Hilbert class field of \(K\), the signed arithmetic
intersection of the single branch with its \(G\)-projection is
\[
 \frac{(2\kappa-2)!\,u_K^2 h\sqrt{|D|}}
 {2^{4\kappa-1}\pi^{2\kappa}\|G\|_{\rm Pet}^{2}}\,
 L'_K(G,\Omega,\kappa)
 \qquad\text{\cite[Corollary~0.3.2]{Zhang}}.
\]
Zhang takes the pullback to a fine cover divided by the square
root of its degree; the auxiliary fine level can be chosen for any
cyclic \(N\) \cite[\S4.1]{Zhang}. His arithmetic intersection
is over the Hilbert class field \cite[\S4.2]{Zhang}. Division by
its degree \(2h\) gives our absolute height. Multiplication by
\(4m^{\kappa-1}/u_K^2\) gives the formula.
\end{proof}

\begin{lemma}[A fundamental coefficient in either Jacobi signature]
\label{signed-fundamental-norm}
Let \(G\) be a normalized newform of weight \(2\kappa\), trivial
character and exact level \(N\), with \(\kappa\ge2\). Suppose its good
Hecke eigenspace in the full Jacobi cusp space of index \(N\), weight
\(\kappa+1\) and signature \(\varepsilon(G)\) is a line. Write \(g\)
for the vector-valued form of a generator. Let \(d_0\) be a fundamental
discriminant of that signature, prime to \(2N\), and fix one root
\(\mu_0^2\equiv d_0\bmod 4N\). With \(b_0=c_g(|d_0|/(4N),\mu_0)\),
\begin{equation}\label{eq:gzz-waldspurger}
 \frac{|b_0|^2}{\|g\|_{\rm Pet}^2}
 =\frac{(\kappa-1)!\,|d_0|^{\kappa-1/2}}
 {2^{2\kappa+1}\pi^\kappa N^{\kappa-1/2}\|G\|_{\rm Pet}^2}
 L(G,\chi_{d_0},\kappa).
\end{equation}
\end{lemma}
\begin{proof}
Write \(\phi\) for the Jacobi form attached to \(g\). Its theta
components are \(g_\mu\) in the holomorphic case and \(\overline{g_\mu}\)
in the skew-holomorphic case. We may choose \(g\) with real coefficients.
For \(d_0<0\) apply \cite[Ch.~II, \S4, Corollary~1]{GKZ} to the
chosen residue. For \(d=d_0>0\) we prove the skew-holomorphic formula
of \cite{SkoruppaBonn} with the same normalization.

Let \(S\) be the Shimura map at \((d,\mu_0)\). The elliptic
\(G\)-eigenspace at its exact level is a line, so the first
Fourier coefficient gives \(S\phi=b_0G\). By Hecke equivariance
and the Jacobi multiplicity hypothesis, \(S^*G\) is a multiple
of \(\phi\). Its coefficient at \((d,\mu_0)\) is therefore
\(|b_0|^2\|G\|_{\rm Pet}^2/\|\phi\|_{\rm Jac}^2\).

Put \(B=\binom{2\kappa-2}{\kappa-1}\) and let \(f_d\) be the
binary quadratic form kernel of discriminant \(d^2\), middle
residue \(\mu_0^2\bmod 2N\) and genus character \(\chi_d\).
Skoruppa uses the Jacobi weight \(\kappa+1\) and evaluates his
kernel at \(-\bar t\) in the elliptic variable. With
\(f_d^\#(t)=\overline{f_d(-\bar t)}\), his adjoint formula gives
the same coefficient as
\[
 c\,d^{2\kappa-1}(G,f_d^\#)_{\rm Pet},\qquad
 c=\frac{i^\kappa(2/N)^{\kappa-1}}{\pi B}
 \qquad\text{\cite[(3), Proposition~1]{SkoruppaBinary}}.
\]
Reflection is antiunitary and \(G^\#=G\). Hence
\((G,f_d^\#)_{\rm Pet}=\overline{(G,f_d)_{\rm Pet}}\).
The cycle integral formula \cite[Ch.~II, \S1, Proposition~2]{GKZ}
gives \((G,f_d)_{\rm Pet}=\pi B\,2^{2-2\kappa}d^{1-2\kappa}R_d\),
where \(R_d\) is the sum of the oriented periods with genus weights.
Consequently
\[
 \frac{|b_0|^2\|G\|_{\rm Pet}^2}{\|\phi\|_{\rm Jac}^2}
       =i^\kappa(2N)^{1-\kappa}\overline{R_d}.
\]

The contributing orbits have representatives \([0,d,a]\) for
\(a\bmod d\) with \((a,d)=1\). Indeed, a form of discriminant
\(d^2\) with nonzero genus character is primitive. Factor it as
\((AX+BY)(CX+EY)\) with primitive factors ordered by \(AE-BC=d\).
The level and middle residue give \(N\mid AC\) and \(N\mid BC\),
so \(N\mid C\). A matrix in \(\Gamma_0(N)\) carries the cusp of
the second factor to infinity, and the remaining translations change
\(a\) by multiples of \(d\). For this representative the character
is \(\chi_d(a)\) and the path runs from \(-a/d\) to \(i\infty\).
Putting \(z=-a/d+iy\) gives
\[
 R_d=i^\kappa d^{\kappa-1}\int_0^\infty
       \sum_{a\bmod d}\chi_d(a)G(-a/d+iy)y^{\kappa-1}\,dy .
\]
The primitive even Gauss sum turns the sum over \(a\) into
\(\sqrt d\,(G\otimes\chi_d)(iy)\). Its Mellin integral is
\(\Gamma(\kappa)(2\pi)^{-\kappa}L(G,\chi_d,\kappa)\).
Before imposing any Jacobi multiplicity condition, the same adjoint
kernel calculation gives, for every real elliptic cusp form \(R\),
\begin{equation}\label{eq:adjoint-kernel-coefficient}
 c_{S_{d,\mu_0}^*R}(d/(4N),\mu_0)
       =\gamma_{\kappa,N,d}L(R,\chi_d,\kappa),
 \qquad \gamma_{\kappa,N,d}=(2N)^{1-\kappa}d^{\kappa-1/2}\frac{\Gamma(\kappa)}{(2\pi)^\kappa}.
\end{equation}
The factor \(i^\kappa\) in \(R_d\) cancels the phase of the adjoint
after conjugation, since the twisted central value is real.
For the normalized eigenform \(G\), the nonnegativity of the twisted central value also follows from \cite[Theorem~1]{LapidRallis}.

For either signature, write \(z=x+\tau y\) in the Jacobi norm.
Integration over \(x\) separates the theta components. The Gaussian
integral over \(y\) contributes \(1/(2\sqrt{N\operatorname{Im}\tau})\),
and the quotient by \(z\mapsto-z\) contributes \(1/2\).
Thus \(\|\phi\|_{\rm Jac}^2=\|g\|_{\rm Pet}^2/(4\sqrt N)\).
Substitution gives \eqref{eq:gzz-waldspurger}.
\end{proof}

\begin{lemma}[The factor at the primes common to the twist and the level]
\label{ramified-primitive-period}
Let \(G\) be a normalized newform of weight \(2\kappa\ge4\),
exact level \(N\), trivial character and root number
\(-\operatorname{sgn}(\Delta)\). Suppose \(W_pG=G\) for every
\(p\mid(N,\Delta)\). For a fundamental index \(v=(m,\mu)\) with
\((d,2N\Delta)=1\) and \(D=\Delta d<-4\) fundamental, its projected
self-height satisfies \eqref{eq:intro-gzz} with \(D_0=D\), \(m_0=m\).
These sign conditions hold for the correspondent of a retained real
primitive Jacobi newform.
\end{lemma}
\begin{proof}
Write \(K=\mathbb Q(\sqrt D)\), let \(\Omega\) be its genus character,
and put \(J=\{p:p\mid (N,\Delta)\}\). At the primes of \(N\) the
discriminant \(d\) is a unit square. Thus \(\Omega_p=\chi_{d,p}
\circ {\rm N}_{K_p/\mathbb Q_p}=1\). The primes in \(J\) ramify in
\(K\) and have exponent one in \(N\). All other level primes split.

\smallskip\noindent\emph{The sign of the newvector.}
When \(G\) is the correspondent of a retained primitive \(g\), its
signature gives \(\varepsilon(G)=-\operatorname{sgn}(\Delta)\).
To check \(W_pG=G\) for \(p\in J\), choose a unit fundamental coefficient
\(b_g(v')\ne0\) away from \(2N\Delta\), with discriminant \(d'\).
In the kernel of the Shimura lift \(S_{v'}\), a determinant-\(p\)
Atkin--Lehner matrix \(W\) acts on the quadratic form \(Q\) by
\(Q\mapsto p^{-1}Q\circ W\). This permutes the forms in the kernel
and changes their middle residue by \(\sigma_p\)
\cite[Ch.~I, \S1, (6)--(7)]{GKZ}. The genus character of discriminant
\(d'\) is unchanged, since \(W\) is invertible modulo \(|d'|\)
and \(\chi_{d'}(p)=1\). The elliptic slash on the kernel therefore
acts by the residue involution on the Jacobi variable, in either
signature \cite[Proposition~1]{SkoruppaBinary}. This gives
\((S_{v'}g)|W_p=S_{v'}(\sigma_p g)=b_g(v')G\).
The equality \(\sigma_pg=g\) follows from \(\Pi_{\Delta,N}g=g\).
The argument includes \(p=2\), where \(d'\equiv1\bmod8\).

We can now compute the root number before using any height formula.
At \(p\mid N\) outside \(J\) the twisting character is locally trivial.
At \(p\mid\Delta\) the two characters in the unramified parameter,
or in the semisimplification of the Steinberg parameter at \(p\in J\),
acquire the same ramified quadratic twist. Their unramified factors
cancel by triviality of the central character. The product of the two
normalized Gauss sums is \(\chi_{\Delta,p}(-1)\). In the Steinberg
case the ramified parameter has no inertia invariants, so the monodromy
correction is one. Thus this local sign replaces \(w_p(G)\) at
\(p\in J\) and \(+1\) elsewhere on the conductor of the twist.
The archimedean discrete series is invariant under the sign twist.
Using \(\varepsilon(G)=-\operatorname{sgn}(\Delta)\) gives
\[
 \varepsilon(G\otimes\chi_\Delta)
 =\varepsilon(G)\operatorname{sgn}(\Delta)\prod_{p\in J}w_p(G)=-1 .
\]

\smallskip\noindent\emph{Normalization of the orbit.}
The root number gives \(L(G,\chi_\Delta,\kappa)=0\). Since
\(L_K(G,\Omega,s)=L(G,\chi_\Delta,s)L(G,\chi_d,s)\), the central
vanishing required by Qiu's
rank-one projector formula follows \cite[Theorems~2.3.7 and~2.4.3, (2.22)]{QiuCM}.
On the finite sum of constituents of the orbit, Qiu's projector and
\(e_G\) select the same good Hecke eigenvalues and hence give the
same height \cite[(2.18)--(2.20)]{QiuCM}.

In the representation \(\pi\) of \(G\), averaging both arguments
over \(K_0(N)\) projects onto the newvector line. The kernel of this
projector is a unit vector tensored with its dual. Translating the
oriented CM isogeny to the standard level structure conjugates the
local torus in this matrix coefficient. A common fine cover gives the
same level-volume factor for every \(K\) \cite[(2.12)--(2.14)]{QiuCM}.

Qiu averages a single oriented class group orbit with total mass one.
With \(h=|\Pic(\mathcal O_K)|\) and \(u_K=|\mathcal O_K^\times|/2\),
our sum over that orbit has multiplier \(h/u_K\). Indeed the ideal
\(\mathfrak n=(N,(r\mu+\sqrt D)/2)\) specifies the cyclic isogeny
on each ideal class \cite[Introduction, \S1]{GKZ}. At a ramified
level prime its local factor is the unique prime of norm \(p\).
Passage to the dual isogeny conjugates \(\mathfrak n\) and translates
the ideal class by \(\mathfrak n\) or its inverse. Its genus value is
one at every level prime. Consequently the two signed vector
orientations act on the single-branch sum by \(1+\epsilon W_N\),
as in Lemma~\ref{split-primitive-absolute}, even when a local ideal
equals its conjugate. Here \(\epsilon=\operatorname{sgn}(\Delta)
(-1)^{\kappa-1}\) is the Fricke eigenvalue of \(G\) by the
assumed root number.
Thus the multiplier in the self-height is
\(4m^{\kappa-1}(h/u_K)^2\). The normalization of the fibre class
is fixed by its self-intersection \cite[\S2.2.1]{QiuCM}.

Write \(\beta_p(K)\) for Qiu's normalized toric integral with
newvector matrix coefficient equal to one at the identity.
For the finite Euler product \(L_{\rm fin}\), the class number
formula gives \((h/u_K)^2/L_{\rm fin}(1,\chi_D)^2=|D|/\pi^2\).
The archimedean factors in (2.22) are the same for all these
imaginary quadratic fields and genus characters. Because the
central value vanishes, only the derivative of the finite Euler
product contributes when differentiating the completed function.
Thus the projected height is
\[
 c(G,N)\,m^{\kappa-1}|D|\,L'_K(G,\Omega,\kappa)
                  \prod_{p<\infty}\beta_p(K),
\]
where \(c(G,N)\) is independent of \(K,\Omega\).

Lemma~\ref{finite-toric-tests} applies directly to this orbit test,
including at the shared primes, since the base-change central value
vanishes. It remains to compute the local factors \(\beta_p(K)\).

\smallskip\noindent\emph{The local integrals.}
In the notation of \cite[\S2.4.1]{QiuCM} the normalizing multiplier is
\[
 \frac{L_p(1,\chi_D)L_p(1,\pi,\mathrm{ad})}
      {\zeta_p(2)L_p(1/2,\pi_K\otimes\Omega)} .
\]
Away from \(ND\) the spherical value is one. At \(p\mid N\),
\(p\notin J\), the oriented level segment lies in the apartment
of the split torus with trivial character. Its period
\(\beta_p^{\rm sp}\) depends only on \(\pi_p\) and the level:
such segments are related by a translation of that apartment or
by interchanging the two eigenlines.

At \(p\mid D\) put \(\delta_p=|D|_p^{1/2}\). The quotient
\(K_p^\times/\mathbb Q_p^\times\) consists of two unit cosets
of measure \(\delta_p\) each. For \(p\nmid N\), let
\(\alpha,\alpha^{-1}\) be the unitary-normalized Satake parameters,
\(y=p^{-1/2}\), and \(s=\Omega_p(\varpi_K)=\pm1\).
The units fix the lattice vertex and a uniformizer takes it to a
neighbor. The spherical Hecke eigenvalue divided by the \(p+1\)
neighbors gives matrix coefficient \((\alpha+\alpha^{-1})y/(1+y^2)\).
The integral and its multiplier are respectively
\[
 \delta_p\left(1+\frac{s(\alpha+\alpha^{-1})y}{1+y^2}\right),
 \qquad \frac{1+y^2}{(1+s\alpha y)(1+s\alpha^{-1}y)} .
\]
Their product is \(\delta_p\).

At \(p\in J\) write \(\pi_p=\mathrm{St}\otimes\nu\) and
\(\lambda=\nu(p)=\pm1\), so \(w_p=-\lambda\) is its
Atkin--Lehner eigenvalue. At a split torus the matrix coefficient
on the shell \(\operatorname{diag}(p^n,1)\) is
\((\lambda/p)^{|n|}\): the Iwahori double coset for \(n\ge0\)
has \(p^n\) right cosets and its Hecke eigenvalue is \(\lambda^n\).
Each shell has measure one. Multiplying their sum by
\((1-\lambda/p)^2/(1-1/p)\) gives \(\beta_p^{\rm sp}=1+1/p\).
For the ramified torus, units fix the edge
\([\mathcal O_{K_p},\varpi_K\mathcal O_{K_p}]\) pointwise and
\(\varpi_K\) interchanges its endpoints. The two matrix coefficients
are \(1,w_p\). Since \(\Omega_p=1\), the normalizing multiplier is
\(1-\lambda/p\). We have proved \(w_p=1\), hence
\[
 \beta_p(K)=\delta_p(1+w_p)(1-\lambda/p)
            =2\delta_p\,\beta_p^{\rm sp}.
\]
For \(p=2\), \(\delta_2\) is \(2^{-v_2(D)/2}\), with
\(v_2(D)=2\) or \(3\) in the ramified case. The extra factor two
is the ratio after removing this discriminant factor.
These local ratios also occur in
\cite[Proposition~3.12 and Lemma~3.13]{CaiShuTian}.

\smallskip\noindent\emph{The global constant.}
Multiplication gives \(|D|\prod_p\beta_p(K)=
2^t\sqrt{|D|}\prod_{p\mid N}\beta_p^{\rm sp}\).
To evaluate the common constant, choose coprime fundamental
discriminants \(a,b\) with every prime of \(2N\) split in both
quadratic fields, \(\operatorname{sgn}(a)=-\varepsilon(G)\),
\(\operatorname{sgn}(b)=\varepsilon(G)\), and
\(L'(G,\chi_a,\kappa)L(G,\chi_b,\kappa)\ne0\).
Apply \cite[Theorem, parts (i) and (ii)]{BumpFriedbergHoffstein}
successively, requiring the second field to split also at the primes
of \(a\). We may take \(ab<-4\). For \(D_*=ab\) the derivative of
the genus \(L\)-function is nonzero and \(t_*=0\).
Lemma~\ref{split-primitive-absolute} at this discriminant gives
\[
 c(G,N)\prod_{p\mid N}\beta_p^{\rm sp}
   =\frac{(2\kappa-2)!}{2^{4\kappa-2}\pi^{2\kappa}\|G\|_{\rm Pet}^2}.
\]
Substitution gives \eqref{eq:intro-gzz}, also for target
discriminants at which the derivative vanishes.
\end{proof}

\begin{proof}[Proof of Corollary~\ref{p19-higher-gzz}]
For \(D_0<-4\), Lemmas~\ref{split-primitive-absolute}
and~\ref{ramified-primitive-period} give the self-height directly.
For every \(U\) in the statement we may instead deduce it from the
mixed height. That identity is calibrated at a separate fundamental
index of discriminant less than \(-4\) (Theorem~\ref{positive-full-mixed}).
At the chosen \(U\) we use the Fourier coefficient with its
stabilizer weight \(2/|\Gamma_x|\) from \eqref{rd-coeff}.
The calculation also
applies on any retained line satisfying the Shimura, central vanishing
and Waldspurger hypotheses in part~(2) of the remark after
Corollary~\ref{p19-higher-gzz}.
Write \(X_g=Z_\kappa(F_g)\), \(Y=e_G Z_\kappa(U)=b_0 X_g\), and
\(A_g=C_{\kappa,N,|\Delta|}L'(G,\chi_\Delta,\kappa)/\|g\|_{\rm Pet}^2\).
Self-adjointness and \(e_GX_g=X_g\) give
\[
 H(X_g,Z_\kappa(U))=H(X_g,Y)=b_0H(X_g,X_g),\qquad
 H(Y,Y)=b_0^2 H(X_g,X_g).
\]
Thus \(H(Y,Y)=b_0 H(X_g,Z_\kappa(U))\), including when \(b_0=0\).
Since \(\operatorname{Sh}_{m_0,\mu_0}(g)=b_0G\), multiplying
\eqref{intro-newform-height} by \(b_0\) gives \(H(Y,Y)=2^t b_0^2 A_g\).
Use \eqref{eq:gzz-waldspurger} at the stated index \(a_U\), as supplied
by Lemma~\ref{signed-fundamental-norm} on an exact-level line of
multiplicity one in the full Jacobi space.
The coefficient \(b_0\) is real. The genus character is the
restriction of \(\chi_{d_0}\) to \(K\). Its induction to \(\mathbb Q\)
is \(\chi_{d_0}\oplus\chi_\Delta\), since \(\chi_K=\chi_{d_0}\chi_\Delta\).
The corresponding \(L\)-functions therefore satisfy
\[
  L_K(G,\chi,s)=L(G,\chi_{d_0},s)L(G,\chi_\Delta,s).
\]
This uses all finite local factors. At a prime shared by \(\Delta\)
and \(N\), compatibility gives \(p\Vert N\). The representation of
a newform of exact level \(N\) and trivial character there is an
unramified twist of Steinberg. The ramified quadratic twist has
local \(L\)-factor one, as in the Dirichlet series used in
\eqref{intro-newform-height}. At primes dividing \(d_0\), the
original newform is unramified. On a retained primitive line the local root numbers computed in
Lemma~\ref{ramified-primitive-period} give
\(L(G,\chi_\Delta,\kappa)=0\), independently of \(U\). On the other lines we have assumed
central vanishing. Differentiating at \(s=\kappa\) gives \(L'_K(G,\chi,\kappa)=
L(G,\chi_{d_0},\kappa)L'(G,\chi_\Delta,\kappa)\).
Multiply \eqref{eq:gzz-waldspurger} by
\(2^t C_{\kappa,N,|\Delta|}L'(G,\chi_\Delta,\kappa)\).
Substitution of \(|d_0|=4Nm_0\) and
\(\Gamma(\kappa-\tfrac12)=(2\kappa-2)!\sqrt\pi/
(4^{\kappa-1}(\kappa-1)!)\) gives precisely the constant in
\eqref{eq:intro-gzz}. Its right side is therefore \(2^t b_0^2 A_g\),
as required, with \(2^t=1\) when \((N,\Delta)=1\).
Conversely, a single \(U\) with \(b_0\ne0\) determines
\(H(X_g,X_g)\) by division by \(b_0^2\). An index with \(b_0=0\)
gives zero on both sides for every value of \(H(X_g,X_g)\).
\end{proof}

\begin{corollary}[Infinite rank at one Fourier index]\label{fixed-index-infinite-rank}
Let \(G\) be a normalized newform of weight \(2\kappa\ge4\),
trivial character and exact level \(N\). There is a single Fourier
index \((m,\mu)\) and infinitely many pairwise coprime fundamental
discriminants \(\Delta\), each with a compatible root \(r_\Delta\),
such that the classes
\[
     Z^{(N)}_{\Delta,r_\Delta,\kappa}(m,\mu)
\]
are linearly independent over \(\mathbb C\) in \(A_N\) and pairwise
orthogonal for the absolute height. Let \(L\) be the compositum of
their quadratic fields \(\mathbb Q(\sqrt\Delta)\). Dividing each class
by its scalar \(c_{a,N,\Delta}\) from \eqref{eq:common-chow-scalar}
gives linearly independent classes in
\(\CH^\kappa_{\rm hom}(\mathcal Y_\kappa(N)_L)_{\mathbb Q}\).
For any given number field \(B\), the discriminants can be chosen
so that \(L\) is linearly disjoint from the Galois closure of \(B\)
over \(\mathbb Q\), with any prescribed finite set of primes split
completely in \(L\). If \(\varepsilon(G)=-1\), one can take \(L\)
totally real.
\end{corollary}
\begin{proof}
First apply \cite[Theorem, part (ii)]{BumpFriedbergHoffstein}
to \(G\), requiring splitting at \(2N\). It gives a fundamental
discriminant \(d\) of sign \(\varepsilon(G)\) with
\(L(G,\chi_d,\kappa)\ne0\). Fix a root \(\mu^2\equiv d\bmod4N\)
and put \(m=|d|/(4N)\). Part (i) of the same theorem gives a
fundamental \(\Delta\) of the opposite sign, split at the primes of
\(2Nd\), for which \(L(G,\chi_\Delta,s)\) has a simple central zero.
Adding the ramified primes of the previous choices to this finite
splitting set produces infinitely many pairwise coprime \(\Delta\).
Each has a root \(r_\Delta\bmod 2N\). Discard any choices with
\(\Delta d\ge-4\).

Both discriminants are odd, fundamental and coprime, so their product
is fundamental. Since \(d\) and \(\Delta\) have opposite signs, the
fixed pair \((m,\mu)\) is an allowed index for every chosen \(\Delta\).
For \(K=\mathbb Q(\sqrt{\Delta d})\) and the associated genus
character \(\Omega\), factorization gives
\[
 L'_K(G,\Omega,\kappa)
     =L'(G,\chi_\Delta,\kappa)L(G,\chi_d,\kappa)\ne0 .
\]
Lemma~\ref{split-primitive-absolute} now gives a nonzero projected
height of the coefficient at the fixed index. This projection
can be taken in the height quotient: a zero Chow class has zero
image there and hence zero projected height. The original coefficient
is therefore nonzero in ordinary Chow. Proposition~\ref{coef-galois-character}
separates the classes, gives orthogonality and descends each
normalized class to its quadratic field.

To obtain disjointness, include all primes ramified in the Galois
closure of \(B\), as well as the prescribed splitting primes, at
every choice of \(\Delta\). Take \(|\Delta|>1\). Each nonempty finite
product of the chosen discriminants has an odd prime dividing exactly
one factor. The corresponding quadratic field ramifies there, whereas
the Galois closure of \(B\) is unramified. Every nontrivial
subextension of a finite multiquadratic extension contains a quadratic
field, proving the disjointness at each stage and hence in the union.
The sign of \(\Delta\) gives total reality when \(\varepsilon(G)=-1\).
The proof uses only the height quotient and the character law.
\end{proof}

\subsection{Heights on primitive newform lines}\label{sec:primitive-pairings}
The self-height determines a single scalar on each one-dimensional
good Hecke eigenspace. We first identify that scalar on a primitive
line and then pair arbitrary linear combinations of the resulting
cycles.

\begin{proposition}[Auxiliary cycles on a primitive newform line]
\label{primitive-auxiliary-height}
Let \(g\) be a real primitive Jacobi newform with normalized
correspondent \(G\) of exact level \(N\). Suppose the good Hecke
eigenspace in the retained sector is \(\mathbb Cg\).
There are infinitely many fundamental \(U\) admissible for \(F_g\)
with \(b_g(U)\ne0\). For either sign of \(\Delta\), every admissible
\(U\) satisfies \eqref{intro-newform-height}.
\end{proposition}
\begin{proof}
Lemma~\ref{fundamental-separation} supplies infinitely many \(v\) with
\(b_g(v)\ne0\) after excluding all primes in the polar data of \(F_g\).
Fix one with \(D_v<-4\). Self-adjointness of \(e_G\) gives
\[
 H(X_g,Z_\kappa(U))-2^t b_g(U) A_g
   = b_g(U)\bigl(H(X_g,X_g)-2^t A_g\bigr).
\]
It remains to determine the scalar in parentheses. At every shared
prime the kernel identity \(S_v\sigma_p=W_p S_v\) proved in
Lemma~\ref{ramified-primitive-period} holds on the full Jacobi space.
As \(W_pG=G\), its adjoint gives \(\sigma_p S_v^*G=S_v^*G\).
Thus \(S_v^*G\) belongs to the retained \(G\)-eigenspace
\(\mathbb Cg\), with the usual identification by theta decomposition.
This is precisely the multiplicity condition used in the kernel
calculation of Lemma~\ref{signed-fundamental-norm}. For the
holomorphic signature the same adjoint calculation is
\cite[Ch.~II, \S\S3--4]{GKZ}. Hence \eqref{eq:gzz-waldspurger}
holds at \(v\). Multiplying that identity by the derivative, as in
the genus factorization above, identifies the right side of
\eqref{eq:intro-gzz} with \(2^t b_g(v)^2 A_g\).
Lemma~\ref{ramified-primitive-period} identifies its left side with
\(b_g(v)^2 H(X_g,X_g)\). Since \(b_g(v)\ne0\), the required scalar
identity follows, also when \(A_g=0\).
\end{proof}

\noindent\emph{Self-opposite auxiliary indices.}
Let \(N\) be squarefree and let \(\Delta\ne1\) be a compatible twist
with \((\Delta,2N)=1\) and \(\operatorname{sgn}(\Delta)(-1)^{\kappa-1}=-1\).
Put \(\epsilon=-\operatorname{sgn}(\Delta)\). When \(N\) is odd, choose primes \(\ell\equiv\epsilon N\bmod 4\) away from
\(2N\Delta\) and the discriminants of the finitely many CM fields
in the polar support. There are infinitely many by Dirichlet's theorem.
Take
\[
 d_0=\epsilon N\ell,\qquad \mu_0=N\bmod 2N,\qquad m_0=\ell/4 .
\]
The integer \(d_0\) is squarefree, \(d_0\equiv1\bmod4\) and
\(d_0\equiv N^2\bmod4N\). Thus \(d_0\) is a fundamental discriminant prime to \(\Delta\);
their product \(D_0=\Delta d_0<0\) is also fundamental.
The form \([N,rN,((rN)^2-D_0)/(4N)]\) is integral and primitive of
discriminant \(D_0\). It represents \(N\), which is prime to
\(\Delta\), so its genus weight is nonzero. When \(N>1\), the primes
of \(N\) ramify in the CM field of discriminant \(D_0\). 
If \(N\) is even and squarefree, take instead
\(d_0=4\epsilon N\ell\), \(\mu_0=0\), \(m_0=\ell\), for any odd
prime \(\ell\) away from \(2N\Delta\) and the polar CM field
discriminants. Both \(d_0/4\) and \(D_0/4=\Delta d_0/4\) are
squarefree and congruent to \(2\bmod 4\), so the discriminants are
fundamental. The primitive form \([N,0,|\Delta|\ell]\) has
discriminant \(D_0\) and represents \(N\), which is prime to
\(\Delta\).

In both constructions \(\mu_0\) is its own negative. The negative component sign in \eqref{rd-coeff} gives
\(Z_\kappa(U)=0\) in Chow and the same symmetry gives \(b_0=0\).
In fact every residue in \eqref{merged-shimura-coeff} is a multiple
of \(\mu_0\) and is self-opposite. The Shimura lift at \(U\) therefore
vanishes for every cuspidal shadow of this type. Since the height
on the CM subspace descends to Chow, both sides vanish even at
these discriminants sharing primes with \(2N\).

\begin{proposition}[Testing the twisted mixed height]\label{twisted-height-test}
Assume \((\Delta,2N)=1\) and fix the newform data \(g,G,F_g\)
introduced after Theorem~\ref{intro-newform-mixed}. Suppose
\(H(X_g,\cdot)\) is real linear on the Chow coefficient space.
Extend it complex linearly and put
\[
 A_g=C_{\kappa,N,|\Delta|}\frac{L'(G,\chi_\Delta,\kappa)}
                                  {\|g\|_{\rm Pet}^2},\qquad
 E_g=(\mathrm{id}\otimes H(X_g,\cdot))\Phi_\kappa-A_g g .
\]
Then \(E_g\in S\) depends only on the normalized shadow. At any
admissible fundamental \(U\), \eqref{intro-newform-height} holds
exactly when \(c_{E_g}(a_U)=0\). Under
Proposition~\ref{p19-height-comparison} this coefficient is the total
boundary and finite correction in \eqref{height-residual-finite}.

For any nontrivial fundamental twist prime to \(2N\), if \(N\)
is squarefree and \(\kappa\ge2\), one can choose \(\dim S\) fundamental indices,
whose input discriminants \(d_0\) avoid any prescribed finite set of
primes and whose Fourier coefficients determine \(E_g\).
Equality at these indices implies \eqref{intro-newform-height}
at every admissible fundamental index. Infinitely many admissible
\(U\) have \(b_g(U)\ne0\), and their CM fields may all be chosen
distinct from those in the polar support of \(F_g\).
\end{proposition}
\begin{proof}
By Theorems~\ref{med-target} and~\ref{intro-chow-vanishing},
contracting \(\Phi_\kappa\) gives a cusp form independent of the lift
of its fixed shadow. The coefficient of \(E_g\) at \(a_U\) is
\(H(X_g,Z_\kappa(U))-A_g b_g(U)\), as asserted. On the retained
line of Proposition~\ref{primitive-auxiliary-height}, self-adjointness
gives \(H(X_g,Z_\kappa(a))=b_g(a)H(X_g,X_g)\) for every index \(a\).
Here \(t=0\) and that proposition gives \(H(X_g,X_g)=A_g\), so
\(E_g=0\) on the entire Fourier expansion.

For an independent coefficient test, Lemma~\ref{med-eta} gives the
bound \(0<n\le k/12\) at every level and in both parities.
At squarefree level Lemma~\ref{fundamental-separation} gives a basis
of fundamental evaluations after excluding \(2N\Delta\), the polar
primes and any further prescribed primes. For each such index
\(D_0=\Delta d_0<0\) is fundamental, and
\([N,r\mu_0,((r\mu_0)^2-D_0)/(4N)]\) represents a CM point.
Successively excluding the primes of the previously chosen \(d_0\)
gives infinitely many admissible fields with \(b_g(U)\ne0\).
\end{proof}

At squarefree level the same basis treats any real input \(F\) after
projection to the invariant sector. Choose real harmonic lifts \(F_i\)
of the projected monomial principal parts \(p_i\) at \(a_{U_i}\). Then
\(\beta(P_{F_i})=\ell_{U_i}\) and \(Z_\kappa(F_i)=Z_\kappa(U_i)\).
Write \(\beta(P_F)=\sum_i a_i\ell_{U_i}\) with \(a_i\) real.
Proposition~\ref{med-exact} makes \(F-\sum_i a_iF_i\) weakly holomorphic,
so Theorem~\ref{intro-chow-vanishing} gives
\[
 Z_\kappa(F)=\sum_i a_i Z_\kappa(U_i),\qquad
 \xi_w F=\sum_i a_i\,\xi_wF_i .
\]
At a common auxiliary \(U\) it therefore suffices to test these
fundamental pairs. Even with pairwise coprime \(d_i\), the CM
discriminants \(\Delta d_i\) share the primes of \(\Delta\), where
both genus weights enter the local multiplicities. Theorem~\ref{intro-newform-mixed}
allows the same \(U\) for all shadows. Its derivative term on a
normalized primitive line is \(2^t A_g b_g(U)\), so at \(b_g(U)\ne0\)
the vanishing of \(X_g\) forces \(L'(G,\chi_\Delta,\kappa)=0\).

For a fixed primitive constituent \(G\) one can count the admissible
indices with \(b_G(U)\ne0\) more precisely. Prescribe
\(\chi_d(p)=1\) at every \(p\mid2N\Delta\) and at every prime in the
polar data of the chosen lift. Ono and Skinner give
\(\gg X/\log X\) squarefree fundamental discriminants
\(|d|\le X\) with these signs and \(L(G,\chi_d,\kappa)\ne0\)
\cite[Corollary~3]{OS}. The root number is then
\(\varepsilon(G)\operatorname{sgn}(d)=+1\), so
\(\Delta d<0\) on a retained component. The prescribed signs give
a unit square root of \(d\bmod4N\), including at \(2\) since
\(d\equiv1\bmod8\). Omitting the finitely many polar CM fields
leaves the same lower bound. At the primitive conductor of \(G\),
\eqref{eq:gzz-waldspurger} gives \(b_G(U)\ne0\) at each of these
discriminants. The implied constant may depend on \(G\) and the
finite set of prescribed primes.

\begin{auxlemma}[Pairings on the sum of primitive newform lines]
Choose a real generator \(g\) on each of a nonempty collection
\(\mathcal B\) of distinct lines in \(S_{\rm prim}^{(1)}\), and put \(h_g=2^t A_g\).
Suppose the projected shadows of \(F,F'\) have expansions
\[
 \Pi_{\Delta,N}\xi_w F=\sum_{g\in\mathcal B}
       \alpha_g\,g/\|g\|_{\rm Pet}^2,\qquad
 \Pi_{\Delta,N}\xi_w F'=\sum_{g\in\mathcal B}
       \alpha'_g\,g/\|g\|_{\rm Pet}^2.
\]
For the complex-bilinear extension of the height and every Fourier
index \(v\),
\[
 \begin{aligned}
 H(Z_\kappa(F),Z_\kappa(v))&=\sum_{g\in\mathcal B}\overline{\alpha_g}\,b_g(v)h_g,\\
 H(Z_\kappa(F),Z_\kappa(F'))&=\sum_{g\in\mathcal B}
           \overline{\alpha_g\alpha'_g}\,h_g .
 \end{aligned}
\]
In particular, \eqref{intro-newform-height} holds at every fundamental
index admissible for \(F\).
\end{auxlemma}
\begin{proof}
Since \(\xi_w\) is antilinear, \(\Pi_{\Delta,N}F-
\sum_g\overline{\alpha_g}F_g\) is weakly holomorphic. By
Theorem~\ref{intro-chow-vanishing} and Corollary~\ref{med-antisector},
\(Z_\kappa(F)=\sum_g\overline{\alpha_g}X_g\).
The good Hecke projectors are self-adjoint and mutually orthogonal.
They give \(e_GZ_\kappa(v)=b_g(v)X_g\) at every index.
Proposition~\ref{primitive-auxiliary-height} computes
\(H(X_g,X_g)=h_g\) on each line. Bilinearity now gives both formulas.
At an admissible \(U\) the Shimura lift of \(g\) is \(b_g(U)G\).
Since \(b_g(U)\) and \(L'(G,\chi_\Delta,\kappa)\) are real, the first
formula is precisely \eqref{intro-newform-height}. Its validity at
\(U\) does not require a simultaneous nonzero coefficient of all the
forms in \(\mathcal B\).
\end{proof}

\begin{theorem}[Mixed heights on the primitive part at arbitrary level]
\label{primitive-height-pencil}
Fix \((N,\Delta,r)\) and the primitive lines \(\mathcal B\) above. Write
\(s=|\mathcal B|\) and \(e_{\mathcal B}=\sum_{g\in\mathcal B}e_G\).
Fix a finite set of primes containing those of \(2N\Delta\).
There are two ordered sets of \(s\) fundamental CM data \(\mathcal U=(U_i)\)
and \(\mathcal V=(V_j)\) for which the coefficient matrices
\(B_{\mathcal U}=(b_g(U_i))_{i,g}\) and \(B_{\mathcal V}=(b_g(V_j))_{j,g}\) are invertible.
All \(2s\) input discriminants can be taken pairwise coprime,
of absolute value greater than one and divisible by none of the
prescribed primes. For any such choice, set
\[
 H_{ij}=H(e_{\mathcal B}Z_\kappa(U_i),Z_\kappa(V_j)),\qquad
 K_{ij}=\sum_{g\in\mathcal B}
 \frac{b_g(U_i)b_g(V_j)}{\|g\|_{\rm Pet}^2}.
\]
The generalized eigenvalue belonging to \(g\) is
\(\|g\|_{\rm Pet}^2 h_g=2^t C_{\kappa,N,|\Delta|}L'(G,\chi_\Delta,\kappa)\).
In particular,
\[
 \frac{\det(zK-H)}{\det K}
 =\prod_{g\in\mathcal B}
 \bigl(z-2^t C_{\kappa,N,|\Delta|}L'(G,\chi_\Delta,\kappa)\bigr).
\]
The rank of \(H\) is the number of nonzero derivatives
in this product, and the corresponding \(X_g\) are linearly independent
in the complexified Chow group.
\end{theorem}
\begin{proof}
The real span of \(\mathcal B\) lies in \(S_*\) of
Lemma~\ref{fundamental-separation}. Select rows successively. If the
rows already chosen fail to span the dual, take a nonzero form they
all annihilate. The lemma gives infinitely many nonzero fundamental
coefficients of that form away from the prescribed primes. Choose
one with \(|d|>1\), append its row and exclude its prime divisors
at every later choice. Each row raises the rank by one. Repeat for
\(\mathcal V\), retaining all exclusions from \(\mathcal U\). This proves the asserted
existence with both matrices invertible.

Put \(Q=\operatorname{diag}(\|g\|_{\rm Pet}^{-2})\) and
\(\Lambda=\operatorname{diag}(\|g\|_{\rm Pet}^2h_g)\).
Self-adjointness gives
\(H(e_{\mathcal B}x,y)=H(e_{\mathcal B}x,e_{\mathcal B}y)\)
on the coefficient space. The preceding lemma therefore gives
\(H=B_{\mathcal U}Q\Lambda B_{\mathcal V}^{\mathsf T}\), and by definition
\(K=B_{\mathcal U}QB_{\mathcal V}^{\mathsf T}\). Thus \(HK^{-1}=B_{\mathcal U}\Lambda B_{\mathcal U}^{-1}\).
In particular the column \((b_g(U_i))_i\) identifies the eigenvalue
belonging to \(g\), even when two eigenvalues coincide. Substituting \(h_g=2^t A_g\)
gives the stated determinant and rank.
The \(X_g\) are mutually orthogonal.
Pairing a Chow relation with each \(X_g\) of nonzero height gives
the asserted independence. On \(\sum_{g\in\mathcal B}\mathbb C X_g\)
the height radical is precisely the span of the \(X_g\) of zero height.

Each good Hecke correspondence in \(e_{\mathcal B}\) preserves the
CM field up to isomorphism, since conjugation by an isogeny identifies
the rational endomorphism algebras. Its image may occupy several
fibres with that field. The fields of discriminants \(\Delta d_{U_i}\)
and \(\Delta d_{V_j}\) are distinct, so every entry admits
representatives with disjoint generic support. Primes of \(\Delta\) can divide both CM discriminants and
contribute to the finite height.
\end{proof}

One harmonic input also recovers all the scalars \(h_g\).
Take a real \(F_*\) with normalized shadow coefficients \(\alpha_g\ne0\),
and exclude its polar primes when choosing \(U\). Then
\[
 \bigl(H(Z_\kappa(F_*),Z_\kappa(U_i))\bigr)_i
       =B_U(\alpha_g h_g)_g .
\]
Inverting \(B_U\) determines every \(h_g\) from these total heights,
even when some derivatives vanish.

\begin{corollary}[Detection by fundamental discriminants]
\label{primitive-height-count}
Fix a finite set of primes \(\mathcal P\) containing those of
\(2N\Delta\). Let \(\lambda\) be a complex-linear functional on the
Chow coefficient span \(W_N\), and put
\[
 h_\lambda=(\mathrm{id}\otimes\lambda)\Phi_\kappa,\qquad
 h_\lambda^*=\operatorname{pr}_{S_*}h_\lambda ,
\]
where \(\operatorname{pr}_{S_*}\) denotes Petersson projection onto
the space of Lemma~\ref{fundamental-separation}.
\begin{enumerate}
\item If \(h_\lambda^*=0\), then \(\lambda(Z_\kappa(U))=0\) at
every fundamental index whose input discriminant is prime to
\(\mathcal P\). Otherwise, for each \(0<\eta<1\) there are
\(\gg X^\eta\) distinct such discriminants \(|d_U|<X\) with
\[
 \bigl|\lambda(Z_\kappa(U))\bigr|\ge c\,|d_U|^{\kappa/2-1/4}
\]
at some index \(a_U\), for a fixed \(c>0\). In particular these
fundamental indices represent infinitely many distinct Chow classes.
At squarefree level the vanishing alternative is equivalent to
\(\lambda=0\) on \(W_N\).
\item Under the hypotheses of
Theorem~\ref{primitive-height-pencil}, take
\(\lambda=H(Z_\kappa(F),\,\cdot\,)\) and a real harmonic input
\(F\) with \(\xi_wF=\sum_{g\in\mathcal B}\alpha_g g/\|g\|_{\rm Pet}^2\).
The vanishing alternative holds precisely when
\(\alpha_g L'(G,\chi_\Delta,\kappa)=0\) for every \(g\), and then
all Fourier indices have zero mixed height. Otherwise the lower bound
counts discriminants supporting a nonzero mixed height at an index
admissible for \(F\).
\end{enumerate}
\end{corollary}
\begin{proof}
Chow-valued modularity gives \(h_\lambda\in S\) with coefficient
\(\lambda(Z_\kappa(U))\) at \(U\).
Within \(S\), the orthogonal complement of \(S_*\) is spanned by the
invariant projections of the lifts from levels \(N/p^2\).
Each such lift is supported on residues divisible by \(p\), and
the residue symmetries preserve this divisibility. A fundamental
index prime to \(2N\) has unit residue. Its coefficient therefore
equals that of \(h_\lambda^*\), to which
Lemma~\ref{fundamental-separation} applies. If \(N\) is squarefree
then \(S_*=S\), and the Fourier coefficients of \(\Phi_\kappa\) span \(W_N\).
More generally, write \(O=S_*^\perp\subset S\), and let
\(T:S^\vee\to W_N\) be the coefficient map. The unit fundamental
evaluations annihilate \(O\) and span \(S_*^\vee\) on restriction,
so they span \(\operatorname{Ann}(O)\) in \(S^\vee\). The corresponding
Chow classes therefore span \(T(\operatorname{Ann}(O))\), independently
of the finite set of excluded primes.

In the height specialization the formula for all primitive
coefficients, proved above, gives
\(h_\lambda=2^t\sum_{g\in\mathcal B}\alpha_g A_g g\in S_*\).
The newforms are linearly independent. To obtain indices admissible
for \(F\), enlarge \(\mathcal P\) by the primes dividing its polar
indices \(4Nn\) and the discriminants of its polar CM fields, and
discard \(|d_U|=1\). The first part now gives the asserted count.
At every index counted, the nonzero height detects a nonzero
Chow coefficient.
\end{proof}

On a single primitive line, there is also a conclusion that uses
only Chow-valued modularity and Hecke equivariance. If \(X_g\ne0\),
then \(e_GZ_\kappa(U)=b_g(U)X_g\ne0\) for every index counted by
Lemma~\ref{fundamental-separation} for \(g\). These projected classes
all belong to the same line, although their CM fields are distinct.

\begin{conjecture}[A Fourier coefficient of the harmonic lift]
\label{merged-three-way}
On a sign-\(-1\) newform line as above, normalize the coefficients
of \(g\) and the principal part of \(F_g\) in the Hecke field \(K_G\).
The following conditions are equivalent.
\begin{enumerate}
\item The central derivative \(L'(G,\chi_\Delta,\kappa)\) vanishes.
\item The polar cycle \(Z_\kappa(F_g)\) is rationally equivalent to zero.
\item The coefficient \(c^+_{F_g}(|\Delta|/(4N),r)\) belongs to \(K_G\).
\end{enumerate}
\end{conjecture}
For the exact-level newform and sign conditions of \cite[\S5]{ABS},
one can choose \(g\) over \(K_G\) and \(F_g\) with principal part
over \(K_G\) and shadow \(g/\|g\|_{\rm Pet}^2\). Under the hypotheses
of their Theorem~5.1, the third condition is equivalent to all Fourier
coefficients of the canonical meromorphic modular form for the Heegner
divisor belonging to \(i\pi^\kappa\sqrt{\Delta}\,K_G\).
Here \(|\Delta|/(4N)\) is the exponent of \(q\) in the \(r\)-component
of the input Weil representation.

Both the cycle and the coset of this positive coefficient in
\(\mathbb C/K_G\) are independent of the choice of lift with the
prescribed shadow and \(K_G\)-rational principal part. The difference
of two such lifts is weakly holomorphic. Its cycle vanishes by
Theorem~\ref{intro-chow-vanishing}, and its Fourier coefficients
belong to \(K_G\) by the rationality argument in the proof of
Theorem~\ref{med-complete-basis}. The conjecture asks whether the
resulting coset detects the Chow class of a lift with nonzero shadow.

\begin{corollary}[Nonvanishing detected by a mixed height]
\label{cor-mixed-detects-nonzero}
Let \(g\) be as in Proposition~\ref{primitive-auxiliary-height}.
If \(L'(G,\chi_\Delta,\kappa)\ne0\), then
\(Z_\kappa(F_g)\ne0\) in the complexified Chow group. Moreover,
\(Z_\kappa(U)\ne0\) for infinitely many distinct fundamental CM
discriminants.
\end{corollary}
\begin{proof}
Apply Lemma~\ref{fundamental-separation} away from the primes
in the polar data of \(F_g\). This gives infinitely many admissible
fundamental indices with distinct discriminants and \(b_g(U)\ne0\).
At each such index \eqref{intro-newform-height} gives a nonzero
pairing, since \(2^t C_{\kappa,N,|\Delta|}\ne0\). A zero Chow
class pairs to zero with every class, proving both assertions.
\end{proof}

\section{Reconstructing heights from Heegner cycles}
\label{merged-conclusion}

Let \(H\) be the signed Beilinson--Bloch pairing of
Section~\ref{rd-h-height}, normalized by absolute arithmetic degree.
It factors through Chow on a common field and fine level, including
for cycles from different CM fields \cite[\S\S2.2.4--2.2.5]{QiuCM}.
For a CM cycle \(Y\), extend \(\lambda_Y(Z)=H(Z,Y)\) complex
linearly from real cycles and put
\[
        h_Y(\tau)=\sum_{\mu,n>0} H(Z_\kappa(n,\mu),Y)\,q^n e_\mu .
\]
Applying \(\lambda_Y\) to Theorem~\ref{med-target} gives
\(h_Y=\sum_j H(Z_\kappa(F_j),Y)\,g_j\in S\).
More generally, the quotient criterion of Theorem~\ref{rd-main}
gives the same expansion for any linear functional annihilating
\(D_N\). For every harmonic input \(F\) the principal-part pairing
then gives
\begin{equation}\label{merged-height-shadow}
      H(Z_\kappa(F),Y)=\beta(P_F)(h_Y)
                 =\langle h_Y,\xi_w F\rangle_{\rm Pet}.
\end{equation}
The Petersson product is linear in its first variable. Conjugate-linearity
in its second variable cancels that of \(\xi_w\).

\begin{corollary}[A cusp form for the coprime heights]
\label{merged-height-representative}
Take \(\Delta=r=1\), an integer \(\kappa\ge2\), and a fundamental CM
discriminant \(D_0\) as in Theorem~\ref{merged-height-theorem},
with \((D_0,2N)=1\). Write \((m_0,\mu_0)\) for the index of
the auxiliary CM cycle \(Z_\kappa(U)\).
For an orthogonal basis \((g_j)\) of \(S\), define
\[
 h_U^{\rm an} =
 C_{\kappa,N,1}\sum_j
 \frac{\overline{L'(\operatorname{Sh}_{m_0,\mu_0}(g_j),\kappa)}}
      {\|g_j\|_{\rm Pet}^{2}}\,g_j .
\]
It is the unique element of \(S\) whose Fourier coefficients at admissible
positions \((n,\mu)\) with \((4Nn,D_0)=1\) satisfy
\[
       c_{h_U^{\rm an}}(n,\mu)= H(Z_\kappa(n,\mu), Z_\kappa(U)).
\]
\end{corollary}
\begin{proof}
The displayed Riesz representative of
\(g\mapsto C_{\kappa,N,1}\overline{L'(\operatorname{Sh}_{m_0,\mu_0}(g),\kappa)}\)
is independent of the orthogonal basis. At a self-opposite residue
of negative reflection sign both coefficients vanish. Otherwise
take the real harmonic Poincar\'e series \(F_i\) with signed polar
monomial \(p_i\) at \(i=(n,\mu)\) as in Section~\ref{med-basis}.
Then \(\beta(p_i)(g)=c_g(n,\mu)\) and
\(Z_\kappa(F_i)=Z_\kappa(n,\mu)\). Removing the content of a polar
quadratic form divides its discriminant by a square. Thus the polar
field discriminants divide \(4Nn\) in absolute value, and
\((4Nn,D_0)=1\) permits Theorem~\ref{merged-height-theorem} at \(F_i\).
It gives
\[
 H(Z_\kappa(n,\mu),Z_\kappa(U))
 =\langle h_U^{\rm an},\xi_w F_i\rangle_{\rm Pet}
 =\beta(p_i)(h_U^{\rm an})=c_{h_U^{\rm an}}(n,\mu).
\]
For uniqueness, let \(R\) be prime to \(2N\). We claim that the
coefficients with \((4Nn,R)=1\) determine every \(g\in S\), for either
reflection sign. For \(g\ne0\) the span of all \(g_\mu(4N\tau)\) is
nonzero and stable under the diamond operators at level \((4N)^2\),
by the transformation law in the proof of Lemma~\ref{fundamental-separation}.
Take a nonzero character eigenvector in this span. At each prime \(p\mid R\)
delete the terms with exponent divisible by \(p\).
The level increases by a factor of \(p^2\), and the form stays nonzero
because its previous level was prime to \(p\)
\cite[Lemmas~3.1--3.2]{SahaSquarefree}. A surviving coefficient comes
from a coefficient of \(g\) with \((4Nn,R)=1\). Taking \(R=|D_0|\)
gives uniqueness. This recovers the criterion of
\cite[Lemma~6.15]{BEY} in their weight range.
\end{proof}

Deleting from \(h_U^{\rm an}\) the terms with \((4Nn,D_0)>1\)
gives the series of coprime geometric heights on \(\Gamma_0(D_0^2)\)
by the Jacobi discriminant projection
\cite[Lemma~2.4]{SchwagJacobi}. The projection preserves the elliptic
translation law since it depends only on the discriminant.
For odd \(\kappa\) this is \cite[Corollary~6.13]{BEY}.
For the Beilinson--Bloch pairing the full series \(h_{Z_\kappa(U)}\)
belongs to \(S\), and agrees with \(h_U^{\rm an}\) on the coprime
positions. Uniqueness gives \(h_{Z_\kappa(U)}=h_U^{\rm an}\).
Thus at any cyclic level the complete column of global heights
against this fundamental \(U\) is determined, including its
self-height.

For a finite collection of auxiliary cycles \(Z_\kappa(U_i)\) one can make the
same choice of polar support for every lift. Put \(R=\prod_i|D_i|\).
The coefficient criterion just proved shows that the functionals
\(\ell_a:g\mapsto c_g(a)\) at positions \(a=(n,\mu)\) with
\((4Nn,R)=1\) span \(S^\vee\). Choose a basis
from these positions over the real coefficient space. For a real
harmonic input \(F\), write \(\beta(P_F)\) in the
basis and lift the corresponding polar polynomial to \(F'\).
Equation~\eqref{rd-principal} gives \(\xi_wF'=\xi_wF\). Thus
\(Z_\kappa(F)\) and \(Z_\kappa(F')\) agree in Chow, and \(F'\) can be
paired using \eqref{med-untwistedheight} with every one of the \(U_i\).

\subsection{The contribution of a primitive newform}
Continue to assume \(\Delta=r=1\) and \(\kappa\ge2\). Let \(a_U=(m,\mu)\)
be a fundamental CM index with \(D_U=-4Nm\) prime to \(2N\). Write
\(Z_\kappa(U)=Z_\kappa(m,\mu)\), \(h_U=h_U^{\rm an}\) and \(C=C_{\kappa,N,1}\).
The form \(h_U\) is defined at every cyclic level by
Corollary~\ref{merged-height-representative}.

Fix a normalized primitive newform \(G\) of weight \(2\kappa\),
level \(M\mid N\), trivial character and root number \(-1\).
Let \(\phi_G\) be a corresponding real Jacobi newform of weight
\(\kappa+1\) and index \(M\). Denote its Jacobi coefficient of
discriminant \(D_U\) and residue \(\mu\bmod 2M\) by \(b_G(U)\).
For each positive \(d\) with \(d^2\mid N/M\) put
\(l=N/(Md^2)\), and write \(g_{G,d}\) for the vector form
corresponding to \(\phi_G|U_d V_l\). Here \(U_d\) replaces \(z\)
by \(dz\) and \(V_l\) raises the index by \(l\). These vectors,
as \(G\) and \(d\) vary, form a basis of \(S\). The Jacobi
newform theorem and the index-raising identities hold for every integral
Jacobi weight \(\kappa+1\) \cite[Theorem~5 and \S3, (2)--(5)]{SkoruppaZagier}.
For even \(\kappa\), \(b_G(m,-\mu)=-b_G(m,\mu)\), so we retain the
residue \(\mu\) in the coefficient. We keep the Petersson norm
at level \(N\) throughout.

\begin{proposition}[The height against a fixed fundamental cycle]
\label{height-old-component}
For every cyclic level \(N\) and every \(U\) as above,
\begin{equation}\label{height-old-line}
 \langle h_U,g_{G,d}\rangle_{\rm Pet}
       = C\, b_G(U)\,\frac{\tau_0(l)}{d^\kappa}\,L'(G,\kappa),
        \qquad \tau_0(l)=\sum_{t\mid l}1 .
\end{equation}
Thus one auxiliary \(U\) with \(b_G(U)\ne0\) determines the
dependence on every other fundamental auxiliary index in the
same primitive Hecke component.
\end{proposition}
\begin{proof}
The Shimura coefficient convention in \eqref{merged-shimura-coeff}
is the Jacobi lifting convention of Skoruppa and Zagier. Their
newform identity and the compatibility with \(U_d V_l\)
\cite[\S3, (iii) and (5)]{SkoruppaZagier} give
\[
  \operatorname{Sh}_{U}(g_{G,d})(z)
       = b_G(U)\sum_{t\mid l} t^\kappa G(dt z).
\]
Here \(z\) is the modular variable of the integral-weight form.
Taking the Dirichlet series multiplies \(b_G(U)L(G,s)\) by
\(d^{-s}\sum_{t\mid l} t^{\kappa-s}\). Since \(L(G,\kappa)=0\),
differentiating the finite factor contributes zero at the centre.
Its value there is \(d^{-\kappa}\tau_0(l)\).
The Riesz formula in Corollary~\ref{merged-height-representative}
now gives \eqref{height-old-line}, since the Fourier data are real.
\end{proof}

To use this formula with several test cycles, fix
a harmonic lift of \(g_{G,d}\) supported at positions prime to all
their discriminants, as above. Its actual mixed height against
each cycle is the corresponding left side of \eqref{height-old-line}.
In particular, even when several old vectors occur for the same \(G\),
the dependence on the auxiliary cycle is through the single
coefficient \(b_G(U)\).

\subsection{A symmetric generating series for mixed heights}
For a prime \(p\) with \(p^2\mid N\), let \(H_p\subset A_L\) be the
subgroup of order \(p\). It is isotropic, and
\(H_p^\perp/H_p\) is the discriminant module at level \(N/p^2\).
The isotropic lift sends \(e_{\bar x}\) to
\(\sum_{h\in H_p}e_{x+h}\) for \(x\in H_p^\perp\).
Let \(S_*\subset S\) consist of the forms Petersson-orthogonal
to all such lifts from the corresponding lower-level cusp spaces.
At squarefree level \(S_*=S\).

\begin{lemma}[Fundamental CM coefficients away from prescribed primes]
\label{fundamental-separation}
For every compatible \((N,\Delta,r)\), every \(\kappa\ge2\) and
every finite set of primes \(\mathcal P\), the coefficient functionals
\(\ell_U(g)=c_g(m,\mu)\) at fundamental indices \(a_U=(m,\mu)\) with
\(d_U=-\operatorname{sgn}(\Delta)4Nm\) prime to \(2N\Delta\) and
to \(\mathcal P\) span \(S_*^\vee\) on restriction to \(S_*\).
Fix \(0<\eta<1\). There are \(c,C,X_0>0\), depending only on
\(\kappa,N,\Delta,\mathcal P,\eta\), such that for every nonzero
\(g\in S_*\) and every \(X\ge X_0\), at least \(CX^\eta\)
distinct such \(d_U\) with \(|d_U|<X\) satisfy
\[
  |\ell_U(g)|\ge c\,\|g\|_{\rm Pet}\,|d_U|^{\kappa/2-1/4}
\]
at some index \(a_U\). Two bases of \(S_*^\vee\) may be chosen with
all input discriminants pairwise coprime.
\end{lemma}
\begin{proof}
Put \(\sigma=-\operatorname{sgn}(\Delta)\) and \(M=4N\). The
output quadratic module
\(A=(\mathbb Z/2N,\mu\mapsto\sigma\mu^2/M)\) is an orthogonal
sum of its primary parts \(A_p\). In \(\mathbb C[A_p]\) let \(W_p\)
be the span of the Weil translates of \(e_x\) for
\(x\notin pA_p\). These local operations can be performed
independently, up to nonzero scalars. Indeed the Weil action
factors projectively through the local levels, and the Chinese
remainder theorem and the surjection
\(\mathrm{SL}_2(\mathbb Z)\to\mathrm{SL}_2(\mathbb Z/M)\)
allow a local Fourier transform or a local \(T\) at any one prime
\cite[\S5]{StrombergWeil}.

If \(A_p\) has prime order, the unit vectors and the Fourier
transform of any one of them span \(\mathbb C[A_p]\). When \(2\Vert N\),
the dyadic part is \(\mathbb Z/4\) with \(Q(x)=vx^2/8\), \(v\) odd.
After its unit components are removed, the Fourier transform of
\(e_1\) is a nonzero multiple of \(e_0-e_2\). Applying \(T\)
gives \(e_0+e_2\), so this factor too is spanned.

It remains to treat \(p^2\mid N\). If \(v\in W_p^\perp\), both \(v\)
and its finite Fourier transform are supported on \(pA_p\).
Fourier inversion then makes \(v\) invariant under translation
by \(H_p=(pA_p)^\perp\). Thus \(v\) is the lift of a function on
\(pA_p/H_p\). Conversely the lifted space is invariant under the
Weil action: \(T\) acts on each coset by a single scalar, and
Fourier transformation of the coset sum is supported on
\(H_p^\perp=pA_p\) and constant on its \(H_p\)-cosets.
This proves that \(W_p^\perp\) is exactly the lifted space.
Globally \(H_p\) is generated by \(2N/p\bmod 2N\); its norm
is \(\sigma N/p^2\), and division by \(p\) identifies the
quotient quadratic module with that of level \(N/p^2\).

The translates of the unit basis vectors span
\(\bigotimes_p W_p\). Its orthogonal complement is the sum of
the lifted spaces with the other factors unrestricted.
Orthogonal projection onto each such space commutes with the
Weil action. It also commutes with \(\Pi_{\Delta,N}\): when
\(p^2\mid N\), compatibility gives \(p\nmid\Delta\), so the
shared-prime involutions act on other primary factors.
For a cusp form \(g\) its projection is therefore a
cusp form in the image of the lower-level lift. If \(g\in S_*\),
Petersson orthogonality forces every such projection to vanish.
A nonzero \(g\in S_*\) consequently has a nonzero unit component:
otherwise its transformation law would make all its values
orthogonal to \(\bigotimes_p W_p\).

The functions \(g_\mu(M\tau)\) for unit \(\mu\) span a nonzero
space of scalar cusp forms on \(\Gamma_1(M^2)\). Indeed, for
\(\gamma=\left(\begin{smallmatrix}a&b\\c&d\end{smallmatrix}\right)
\in\Gamma_0(M^2)\), the matrix acting on \(M\tau\) is
\(\left(\begin{smallmatrix}a&Mb\\c/M&d\end{smallmatrix}\right)\).
It acts on the basis by multiplication of the index by \(d\), up to
a common scalar \cite[\S5]{StrombergWeil}. On \(\Gamma_1(M^2)\)
this matrix belongs to \(\Gamma(M)\), and the scalar is the
Kronecker symbol \(\bigl(\frac{c/M}{d}\bigr)\).
Since \(d\equiv1\bmod M^2\), one has \(\bigl(\frac M d\bigr)=1\)
and \(d\equiv1\bmod4\). The scalar is therefore
\(\bigl(\frac c d\bigr)\), the usual theta multiplier to the
power \(2\kappa+1\) on this group. Rational dilation preserves
cuspidality. The span of the unit components is stable under the
diamond operators of \(\Gamma_0(M^2)/\Gamma_1(M^2)\), so it contains
a nonzero cusp form with a Dirichlet character on \(\Gamma_0(M^2)\).
This character is even since the scalar theta multiplier is \(1\) at \(-I\).
For \(\Delta>0\) and \(N\) odd the Jacobi theta decomposition gives
another construction of such a form
\cite[Theorem~3.3 and Propositions~5.5--5.6]{WashburnJacobi}.

Every exponent \(t\) of this scalar form satisfies
\(t\equiv\sigma\mu^2\bmod M\) for some unit \(\mu\), so \((t,M)=1\).
At each prime \(p\mid\Delta\) or \(p\in\mathcal P\) outside \(M\),
delete the terms at multiples of \(p\). This multiplies the level
by \(p^2\) and preserves nonvanishing, since a nonzero form supported
on multiples of \(p\) must already have level divisible by \(p\)
\cite[Lemmas~3.1--3.2]{SahaSquarefree}. Call the resulting form \(f\). Its exponents are prime to its level.
For an odd squarefree \(t\) put \(A(t)=t^{1/4-\kappa/2}[q^t]f\).
Mao adelizes scalar forms of this weight and even character on
\(\widetilde{\mathrm{SL}}_2\). For the standard additive character
\(\psi\), the images of their irreducible constituents under \(S_\psi\)
are cuspidal representations \(\pi_i\) of \(\mathrm{PGL}_2\), since
\(\kappa\ge2\) excludes unary theta series. They correspond to fixed
newforms of weight \(2\kappa\) and trivial nebentypus
\cite[Appendix~2, \S7.1]{MaoAppendix}. For each pure tensor,
the estimate preceding (7.15) in \emph{loc. cit.}, \S7.4, bounds
\(t|\rho_i(t)|^2\) by a constant times \(L(\tfrac12,\pi_i\otimes\chi_t)\),
where \(\chi_t\) is the quadratic character of the square class of \(t\).
The estimate includes the local factor at \(2\). The conductor of
\(\chi_t\) is \(t\) or \(4t\), whereas the fundamental discriminant
used for the CM test is \(\sigma t\), since \(\sigma t\equiv1\bmod4\).
The two characters can differ.
Since \(A(t)\) is a fixed linear combination of \(\sqrt t\,\rho_i(t)\),
\[
 |A(t)|^4\ll_f\sum_i |L(\tfrac12,\pi_i\otimes\chi_t)|^2 .
\]
For each fixed \(\pi_i\) the sum of these squared central values over
\(t<X\) is \(O_{\pi_i,\epsilon}(X^{1+\epsilon})\). The conductor is
\(O_{\pi_i}(t^2)\). In the approximate functional equation, Mellin
inversion on a dyadic interval for \(t\) separates the smooth weights
in sums of length \(O_{\pi_i,\epsilon}(X^{1+\epsilon})\).
Write the summation variable as \(b^2m\) with \(m\) squarefree and
apply the quadratic large sieve in \(m\). The quadratic character
on \(b^2\) restricts \(t\) to the twists coprime to \(b\).
The factor \(1/b\) and the normalized Hecke bound
\(\lvert\lambda_{\pi_i}(b^2m)\rvert\le\tau(b^2m)\) allow the sum
in \(b\) at a cost \(O_{\pi_i,\epsilon}(X^\epsilon)\);
see \cite[Theorem~2.4]{SoundYoung} and its full-level application in
Corollary~2.5 there. At fixed bad primes, first partition the twists
by their local square classes and separate their Euler factors.
This proves \(\sum_{t<X}|A(t)|^4\ll_{f,\epsilon}X^{1+\epsilon}\)
over odd squarefree \(t\). After removing the primes of the scalar
level, the proof of \cite[Proposition~2.3, \S2.5]{SahaSchmidtYoshida}
gives the complementary bound
\[
   \sum_{t\ {\rm odd\ squarefree}} |A(t)|^2e^{-t/Y}\gg_f Y
   .
\]

To pass to all \(g\) with uniform constants, put
\[
 B_g(t)^2=t^{1/2-\kappa}
       \sum_{(\mu,2N)=1}|c_g(t/M,\mu)|^2 ,
\]
with missing coefficients read as zero. Restrict \(t\) to odd squarefree
integers prime to \(M\Delta\) and to \(\mathcal P\).
For each nonzero \(g\), the scalar form constructed above gives
\(\liminf_{Y\to\infty}Y^{-1}\sum_t B_g(t)^2e^{-t/Y}>0\).
These quadratic forms are uniformly bounded on the Petersson unit sphere.
Indeed, for a fixed scalar component \(h=g_\mu(M\tau)\), the function
\(y^{\kappa/2+1/4}|h(x+iy)|\) is bounded. Parseval at \(y=1/T\)
therefore gives \(\sum_{T\le t<2T}B_g(t)^2\ll_g T\).
Applying this bound to a fixed basis gives the asserted uniform bound.
If the least eigenvalue tended to zero along a sequence \(Y\to\infty\),
a convergent subsequence of the corresponding unit eigenvectors would
contradict the positive lower limit for its limit vector.

The fourth-moment bound applies to every diamond component of
\(g_\mu(M\tau)\). Applying it to a fixed basis of \(S_*\) gives
\(\sum_{t<X} B_g(t)^4\ll_\epsilon X^{1+\epsilon}\) uniformly for
\(\|g\|_{\rm Pet}=1\). Take \(Y=X/(\log X)^2\).
The dyadic bound makes the smoothed second moment above \(X\)
negligible uniformly, while the terms with \(B_g(t)<c_0\) contribute
at most \(c_0^2Y\). For a sufficiently small fixed \(c_0>0\), the
remaining terms contribute \(\gg Y\). Cauchy--Schwarz gives
\(\gg_\epsilon Y^2/X^{1+\epsilon}\) such \(t<X\). Choosing
\(\epsilon<1-\eta\) and taking a largest component in \(B_g(t)\)
proves the claimed bound after rescaling \(g\).
Discard \(t=1\). Each remaining coefficient comes from some
\(c_g(t/M,\mu)\ne0\). The integer \(d_U=\sigma t\) is squarefree
and congruent to \(\mu^2\bmod M\), hence to \(1\bmod4\).
It is a fundamental discriminant prime to \(2N\Delta\) and to
\(\mathcal P\). Consequently \(D_U=\Delta d_U\) is negative
fundamental and \(D_U\equiv(r\mu)^2\bmod M\). The positive
definite form \([N,r\mu,((r\mu)^2-D_U)/(4N)]\) gives the CM index.
Thus the indicated functionals have zero common kernel on \(S_*\). To select
a basis, apply the nonvanishing just proved to a nonzero vector in
the current common kernel, excluding the primes of all input
discriminants already chosen. Each choice increases the rank.
Repeat for the second basis with the first basis's primes excluded.
\end{proof}

\begin{theorem}[A symmetric height series and its central derivatives]
\label{height-kernel}
Let \(N\) be squarefree, \(\Delta=r=1\), and \(\kappa\ge2\).
There is a unique tensor \(\mathcal H\in S\otimes S\) such
that
\[
   (\mathrm{id}\otimes\ell_U)\mathcal H=h_U
   \quad\text{for every fundamental }U\text{ with }(D_U,2N)=1 .
\]
The tensor is symmetric and defines a form of weight
\((\kappa+\tfrac12,\kappa+\tfrac12)\), of Weil type \(\rho_L\) in
each variable. At every pair of Fourier indices \(a,b\), its
coefficient is the global Beilinson--Bloch height
\[
         \mathcal H(a,b)=H(Z_\kappa(a),Z_\kappa(b)).
\]
Put \(g_G=g_{G,1}\) and write \(b_g(a)=c_g(a)\). Then
\begin{equation}\label{eq:spectral-heights}
 \mathcal H =
   C\sum_{\substack{M\mid N\\G\ {\rm at}\ M}}
       \frac{2^{\omega(N/M)}L'(G,\kappa)}
            {\|g_G\|_{\rm Pet}^2}\,g_G\otimes g_G ,
\end{equation}
where the sum runs over the normalized primitive forms specified
above and \(\omega(t)\) counts the distinct prime divisors of \(t\).
In particular, the rank of \(\mathcal H:S^\vee\to S\) equals
the number of these forms with \(L'(G,\kappa)\ne0\).
\end{theorem}
\begin{proof}
Use the real Skoruppa--Zagier basis of Proposition~\ref{height-old-component}.
The index-raising formula applies to every divisor of \(N\), including
\(2\) when \(N\) is even. Squarefreeness gives \(d=1\) in every term.
The first Fourier coefficient
of \(\operatorname{Sh}_U(g_G)\) in that proposition gives
\(\ell_U(g_G)=b_G(U)\), since \(G\) is normalized and every
\(G(tz)\) with \(t>1\) has zero coefficient at \(q\). The \(g_G\) are
orthogonal for the Petersson product. Indeed, multiplicity one distinguishes
their primitive forms by eigenvalues of self-adjoint Hecke operators away
from \(N\).

The right side of \eqref{eq:spectral-heights} has fundamental
columns \(h_U\) by \eqref{height-old-line}, since
\(\tau_0(N/M)=2^{\omega(N/M)}\). Pairing the Chow coefficients
of \(\Phi_\kappa\) produces a tensor in \(S\otimes S\) by bilinearity
and Theorem~\ref{med-target}. Corollary~\ref{merged-height-representative}
identifies the two tensors at pairs of fundamental indices with
coprime discriminants prime to \(2N\). By Lemma~\ref{fundamental-separation}
these pairs contain the Cartesian product of two bases of \(S^\vee\).
This proves equality at every pair of indices and uniqueness.
Symmetry and rank follow from the diagonal expansion.
\end{proof}

\begin{introremark}[The height columns at arbitrary level]
Keep \(\kappa\ge2\) and \(\Delta=r=1\). At every level \(N\), the
map \(\ell_U\mapsto h_U\) on the span of the fundamental evaluations
has rank equal to the number of primitive forms \(G\) in
Proposition~\ref{height-old-component} with \(L'(G,\kappa)\ne0\).
These columns determine a symmetric tensor up to an element of
\(\operatorname{Sym}^2 O\), where \(O=S_*^\perp\).

To see this, let \(Y_G\) be the real vector in the \(G\)-space
specified by
\[
 \langle Y_G,g_{G,d}\rangle_{\rm Pet}
     =d^{-\kappa}\tau_0\bigl(N/(Md^2)\bigr).
\]
This vector is nonzero. In the Jacobi basis \(O\) is precisely the
span of the \(g_{G,d}\) with \(d>1\): the index operators \(U_p\)
commute with \(V_l\), and the images of the lower-level bases under
\(U_p\) are exactly the basis vectors whose \(d\) is divisible by
\(p\). Thus the \(g_{G,1}\) give a basis of \(S/O\).
Each \(\ell_U\) kills \(O\), and Lemma~\ref{fundamental-separation}
shows that they span \(V=\operatorname{Ann}(O)\).
The first Shimura coefficient gives \(\ell_U(g_{G,1})=b_G(U)\).
By \eqref{height-old-line}, the map extends to \(V\) as
\[
        T\ell=C\sum_G L'(G,\kappa)\,\ell(g_{G,1})Y_G .
\]
Here the \(Y_G\) lie in distinct Hecke eigenspaces. This proves the
rank formula, including when several old vectors belong to one \(G\).

For \(\ell_U,\ell_v\) chosen from two bases with all discriminants
pairwise coprime, symmetry of the geometric height gives
\(\ell_v(T\ell_U)=\ell_U(T\ell_v)\). Thus the prescribed columns
are symmetric on \(V\times V\). To complete them to a symmetric
form on \(S^\vee\), choose its values on pairs in a complement
of \(V\). The ambiguity is exactly
\(\operatorname{Sym}^2((S^\vee/V)^\vee)=\operatorname{Sym}^2 O\).
If the form is also required to have self-adjoint good Hecke
operators, distinct primitive eigenspaces are orthogonal. Writing
\(O_G=O\cap S_G\), with \(S_G=\operatorname{span}\{g_{G,d}\}_d\),
the remaining freedom in this linear algebra problem is
\(\bigoplus_G\operatorname{Sym}^2 O_G\). The following result
determines the geometric pairing on each of these spaces as well.
\end{introremark}

The formula at a coincident pair gives its global height. Disjoint
representatives compute the sum over places; a local contribution for
the original representatives requires its own intersection or Green
function calculation.

\subsection{The height on the full oldspace}

\begin{proposition}[Two descriptions of a height block]\label{oldspace-height}
Let \(\Delta=r=1\), let \(N\) be arbitrary, and fix a constituent
\(G\) of primitive level \(M\) in the Jacobi decomposition above.
Put \(I_G=\{j\ge1:j^2\mid N/M\}\), and let \(e_G\) be the projector
onto the full good Hecke eigenspace. Choose a fundamental index
\(a_U=(m,\mu)\) prime to \(2N\) with \(b_G(U)\ne0\). For \(j\in I_G\)
put \(a_j=(j^2m,j\mu)\).
\begin{enumerate}
\item The classes \(e_G Z_\kappa(a_j)\) span the \(G\)-part of the
Chow coefficient space, and all have representatives in the one CM
field of discriminant \(d_U\). Two choices of \(U\) with distinct
fields give a matrix of mixed heights with disjoint generic supports
which determines every height pairing on that part.
\item The height pairing on this part vanishes identically if and
only if \(L'(G,\kappa)=0\). For \(L'(G,\kappa)\ne0\), one fundamental
height column determines the normalization of the toric period
formula \eqref{eq:oldspace-periods} for the entire part.
\end{enumerate}
\end{proposition}
\begin{proof}
\emph{The coefficient matrix.}
The evaluations at unit fundamental indices span
\(\operatorname{Ann}(O)\). Since \(g_{G,1}\) has nonzero image
in \(S/O\), one may choose \(U\) as stated, avoiding any fixed
finite set of primes. Put \(B_{jd}=c_{g_{G,d}}(a_j)\) for \(j,d\in I_G\)
in increasing order. In the Shimura convolution the coefficient
of \(q^j\) is \(B_{jd}\) plus a sum of \(B_{hd}\) over proper
divisors \(h\) of \(j\). By the proof of
Proposition~\ref{height-old-component}, the first nonzero term of
\(\operatorname{Sh}_U(g_{G,d})\) is \(b_G(U)q^d\).
Induction gives \(B_{jd}=0\) for \(j<d\) and \(B_{dd}=b_G(U)\).
Thus \(\det B=b_G(U)^{|I_G|}\). The evaluations \(\ell_{a_j}\)
form a basis on restriction to \(S_G\); applying the Chow coefficient
map gives the spanning assertion. The discriminants \(j^2d_U\)
have the same square class. Imprimitive vectors can have smaller
conductors, but their rational CM endomorphism algebra is the same.
Good isogenies preserve that algebra. Choose a second \(U'\) with a
different field and form its matrix \(B'\) in the same way. If \(D_G\)
is the height matrix in the basis dual to \((g_{G,d})_d\), the measured
matrix is \(BD_G(B')^{\mathsf T}\). Inverting \(B,B'\) recovers \(D_G\).

\smallskip\noindent\emph{Heights within one field.}
Fix the first \(U\), set \(K=\mathbb Q(\sqrt{d_U})\), and write
\(\pi\) for the unitary representation of \(G\). Pull back to the full adelic fine cover for a normal principal
level \(U'\) contained in the cyclic level group, retaining all
determinant components. Each pulled-back \(Z_\kappa(a_j)\) is a sum of
torus orbits for \(K\) with the trivial character. Indeed, the
untwisted weight and the chosen CM embedding are constant on each
ideal-class orbit. The graph of its integral trace-zero endomorphism,
with its axis corrections and the normalization of \eqref{rd-coeff},
is the normalized fibre cycle of \cite[\S2.2.1]{QiuCM}, with its
orientation sign. For an orbit \(\mathcal O\), write \(w_{\mathcal O}\)
for the coefficient of this normalized fibre cycle, including its
sign and the multiplicity of the pullback. Choose \(P_{g_{\mathcal O}}\)
in the orbit in the notation of that reference. In its probability
Haar normalization the test function
\[
 f_j=\sum_{\mathcal O}w_{\mathcal O}\,|\mathcal O|\,
                         \mathbf1_{g_{\mathcal O}U'}
\]
represents the entire coefficient by \cite[(2.5)--(2.7)]{QiuCM}:
the integral of a coset indicator gives the average over its orbit.
If two oriented terms have the same underlying point their
multiplicities add in this formula.

The good Hecke actions on the Chow coefficient spaces intertwine by
the geometric projector lemma of Section~\ref{sec:newform-heights}.
The classical base-change \(L\)-function is
\(L(G,s)L(G,\chi_{d_U},s)\). Its central value is zero since \(G\)
has root number \(-1\). If \(L'(G,\kappa)=0\), its central
derivative also vanishes. Lemma~\ref{finite-toric-tests} makes the
height zero on all pairs of the one-field spanning classes, hence
on the whole \(G\)-part. Different good Hecke components are
orthogonal by self-adjointness.

For the converse, project the fundamental height column onto \(S_G\).
The resulting form is
\(C b_G(U)L'(G,\kappa)Y_G\). Since \(Y_G\ne0\), the triangular basis
provides \(j_0\in I_G\) with \(c_{Y_G}(a_{j_0})\ne0\).
When \(L'(G,\kappa)\ne0\) the known entry
\(H_0=C b_G(U)L'(G,\kappa)c_{Y_G}(a_{j_0})\) is nonzero.
This uses the full fundamental column, determined from coprime
supports by Corollary~\ref{merged-height-representative}.

For the other entries, let \(\alpha_{\pi^\infty}^{\sharp}\) be
Qiu's normalized product of local toric matrix-coefficient integrals,
extended linearly to finite-rank operators. With
\(f^\vee(g)=f(g^{-1})\) and convolution for the same Haar measure, put
\[
 P_{jk}=\alpha_{\pi^\infty}^{\sharp}
       \bigl(\widetilde\pi^\infty((f_j*f_k^\vee)^\vee)\bigr).
\]
Lemma~\ref{finite-toric-tests} shows that each projected height is
\(P_{jk}\) times one scalar. Descent from \(U'\) uses the same
degree for every pair. The entry \(H_0\) makes \(P_{j_0,1}\ne0\), giving
\begin{equation}\label{eq:oldspace-periods}
 H(e_GZ_\kappa(a_j),Z_\kappa(a_k))
       = H_0\,\frac{P_{jk}}{P_{j_0,1}}\qquad(j,k\in I_G).
\end{equation}
This formula retains the local periods of the actual level vectors.
In particular, it requires no rank-one assumption on their matrix.
\end{proof}

\begin{theorem}[The rank at arbitrary cyclic level]\label{oldspace-full-rank}
Let \(\kappa\ge2\), \(N\ge1\) and \(\Delta=r=1\). Pair the Chow
coefficients of \(\Phi_\kappa\) in its two variables to obtain
\(\mathcal H_N\in S\otimes S\). Let \(G\) range over the normalized
weight-\(2\kappa\) newforms of trivial character and root number \(-1\)
in the Jacobi decomposition, and let \(M\) be the exact level of \(G\).
Then
\[
 \operatorname{rank}\mathcal H_N
   =\sum_{\substack{G\ {\rm at}\ M,\ M\mid N\\L'(G,\kappa)\ne0}}
                   \#\{j\ge1:j^2\mid N/M\}.
\]
For every \(G\) occurring in the sum, the Chow coefficient map
\(S_G^\vee\to W_N\) is injective. The heights on this part are
given by \eqref{eq:oldspace-local-height} below.
\end{theorem}

\begin{proof}
Fix \(G,M,U\) as in Proposition~\ref{oldspace-height}, taking
\(d_U<-4\), and write \(a=\kappa-1\).
If \(L'(G,\kappa)=0\), that proposition gives the zero block.
Otherwise it gives a nonzero entry against \(Z_\kappa(U)\) and
identifies the one-field matrix up to a common scalar with the
toric pairing. We compute its rank on the triangular basis.

\smallskip\noindent\emph{The toric functional on level vectors.}
At \(p^E\Vert N\) put \(c=v_p(M)\), \(n=E-c\), and
\(s=\lfloor n/2\rfloor\). Since \(d_U\equiv\mu^2\bmod4N\) is a
unit discriminant, \(K_p\) is split. Fix the two eigendirections
of the torus, whose apartment in the tree of \(\mathrm{PGL}_2(\mathbb Q_p)\)
is the line of diagonal lattices. Write \(\mathcal L_p\) for the
split torus invariant functional on \(\pi_p\). The local toric
matrix-coefficient pairing is proportional to
\(\mathcal L_p(v)\overline{\mathcal L_p(w)}\)
\cite[\S3D, proof of Proposition~3.12]{CaiShuTian}.
Casselman's basis of \(\pi_p^{K_0(p^E)}\) consists of translates
\(v_t\), \(0\le t\le n\), of a unit newvector carried by the
oriented subpath \([t,t+c]\) of the standard segment \([0,E]\)
\cite{CasselmanNew}. For \(c=0\) the subpath is one vertex. Let
\(R_p=((v_t,v_{t'}))_{t,t'}\) be its positive definite Gram matrix.
A flag represented by \(gK_0(p^E)\) gives period coordinates
\(\mathcal L_p(\pi_p(g)v_t)\). Indeed, acting by its coset
indicator on \(v_t\) gives \(\pi_p(g)v_t\) times the compact volume.
For a sum of flags, \eqref{eq:oldspace-periods} therefore pairs the
two lists of period coordinates using \(R_p^{-1}\).

\smallskip\noindent\emph{The flags at index \(j^2m\).}
Put \(r_p=v_p(j)\le s\). At \(p\), write the level structure as a
lattice inclusion \(L'\subset L\) with cyclic quotient of order \(p^E\).
For one choice of CM orientation the endomorphism
\(j(\mu+\sqrt{d_U})/2\) kills \(L/L'\); the middle residue gives
exactly this condition as in Section~\ref{rd-h-height}. Choose the
local square root congruent to \(\mu\bmod 2p^E\), including at \(2\)
where \(d_U\equiv\mu^2\bmod 2^{E+2}\). Up to a unit and an
endomorphism divisible by \(p^E\), this condition is
\(p^{r_p}\operatorname{pr}_1 L\subset L'\), where
\(\operatorname{pr}_1\) is projection onto the first eigendirection.

Translate the nearest apartment vertex of \(L\) to the origin.
At depth \(i\) write \(L=\langle v_1,v_2\rangle\) with
\(v_1=p^i e_1,\ v_2=z e_1+e_2\); for \(i>0\), \(z\) is a unit
modulo \(p^i\). At \(i=0\) there is one choice. The admissible
sublattices are exactly
\[
 L'=\langle p^h v_1,\ b v_1+p^{E-h}v_2\rangle,\qquad i+h\le r_p .
\]
Indeed \(L'\cap\mathbb Q_p e_1=\mathbb Z_p p^{h+i}e_1\), whereas
\(\operatorname{pr}_1 L=\mathbb Z_p e_1\). For \(h>0\), cyclicity
requires \(b\) to be a unit modulo \(p^h\); for \(h=0\) take \(b=0\).
The two endpoint depths are \(i,h\). The segment runs \(i\) edges
to the apartment, \(E-i-h\) edges along it and \(h\) edges away.
Indeed, for \(h>0\) the first coordinate of the second
generator of \(L'\) has valuation \(i\), since \(E-h>i\);
the intersection with its first axis has valuation \(i+h\).
This gives depth \(h\) and the stated direction along the apartment.
Writing \(w_0=1\), \(w_i=(p-1)p^{i-1}\) for \(i>0\), there are
\(w_iw_h\) such normalized flags.

These flags enter with their individual multiplicities. A torus-unit
orbit of flags of maximal endpoint depth \(b_0\) has cardinality
\(w_{b_0}=[\mathcal O_{K,p}^\times:
(\mathbb Z_p+p^{b_0}\mathcal O_{K,p})^\times]\).
Indeed, in the two diagonal coordinates a unit fixes an off-axis
lattice of depth \(i\) precisely when its entries agree modulo \(p^i\).
For a global packet the size of the ideal-class orbit divided by
its unit stabilizer changes by the product of these local indices.
This is also the ring-class number formula for the orders of \(K\).
Thus the orbit averages in \(f_j\) sum over the individual normalized
flags, and the choices at different \(p\mid N\) form a product.
The remaining scalar from the factor \((j^2m)^{a/2}\) is \(j^a\)
relative to \(j=1\).

On the opposite CM orientation the root changes sign at every
level prime, and the fibre class acquires \((-1)^a\).
The element interchanging the two eigendirections takes the flags
to those for the opposite root. It multiplies \(\mathcal L_p\) by
a scalar, since the toric functional is unique up to scale.
Thus the second orientation and the primes outside \(N\) give the
same overall factor for all \(j\). This factor is nonzero because
the fundamental height column on this block is nonzero. The period
on a directed subpath of the apartment is independent of its initial
vertex by translation. Its value is nonzero at each \(p\mid N\):
otherwise the fundamental orbit would vanish on the full oldvector
basis there. Normalize this value to one.

\smallskip\noindent\emph{Computing the local matrix.}
Let \(F_d\) be the normalized value on a length-\(c\) subpath
starting \(d\) edges off the apartment and directed towards and
then along it. Thus \(F_0=1\). For a subpath leaving the apartment
at its other end the same values apply. Indeed, path reversal acts
on the newvector by its Fricke sign, and exchange of the two axes
acts on the toric functional by a scalar. Their product fixes the
value on an apartment segment, so is one. Set \(x=a_G(p)/p^{\kappa-1}\)
when \(c=0\), and let \(w\in\{1,-1\}\) be the local Fricke sign when
\(c=1\). The values at positive depth are
\[
\begin{array}{c|c|c}
 &F_1&\text{at greater depth}\\ \hline
 c=0 &(x-2)/(p-1)&pF_{d+1}=xF_d-F_{d-1}\ (d\ge1)\\
 c=1 &-(1+w)/(p-1)&pF_d=-wF_{d-1}\ (d\ge2)\\
 c\ge2&-1/(p-1)&F_d=0\ (d\ge2).
\end{array}
\]
For a spherical vector the sum over adjacent vertices acts by \(x\).
On the apartment two neighbors remain on it; off it there is one
inward neighbor. This gives the first row. For \(c\ge1\), tracing
a newvector to level \(p^{c-1}\) gives zero. At \(c=1\) the two
incoming axial edges have values \(1,w\); at greater depth the
single inward neighbor reverses the orientation and introduces \(w\).
At \(c\ge2\) fix the remaining \(c-1\) edges. At an axis vertex
one backward extension follows the axis and \(p-1\) leave it. At a
vertex of positive depth all \(p\) permitted backward extensions
go farther off the axis. They have equal values, since the diagonal
units fixing that vertex act transitively on its outward branches.
This proves the last two rows, also for \(p=2\).

For \(0\le r\le s\) the vector \(V_r\) of period coordinates therefore
has entries
\[
 V_r(t)=\sum_{i+h\le r} w_iw_h\,
       F_{\max\{0,\,i-t,\,h-(n-t)\}},\qquad 0\le t\le n ,
\]
where \(i,h\ge0\). The two positive depths in the maximum cannot
occur together since \(i+h\le r\le n/2\). For completeness
\(R_p(t,t')=D_{|t-t'|}\) with \(D_0=1\). In the spherical case
\(D_1=x/(p+1)\) and \(pD_{b+1}=xD_b-D_{b-1}\) for \(b\ge1\).
For \(c=1\), \(D_b=(-w/p)^b\); for \(c\ge2\), \(D_b=0\) for
\(b>0\). These are the same adjacency and trace computations with
the first vector fixed. For the spherical inner product all \(p+1\)
neighbors of its vertex have the same value. For \(c\ge2\), the
stabilizer of the segment up to the penultimate vertex fixes the
first vector and permutes the \(p\) permitted terminal extensions.

Define \(K_p(r,r')=V_r^{\mathsf T}R_p^{-1}V_{r'}\).
In particular \(K_p(0,0)>0\) since \(V_0=\mathbf1\). The toric
pairing, the flag count and \eqref{eq:oldspace-periods} now yield
\begin{equation}\label{eq:oldspace-local-height}
 H(e_G Z_\kappa(a_j),Z_\kappa(a_k))
 = H(e_G Z_\kappa(U),Z_\kappa(U))\,(jk)^a
      \prod_{p\mid N}\frac{K_p(v_p(j),v_p(k))}{K_p(0,0)}.
\end{equation}
The fundamental column gives the reference value
\(C b_G(U)L'(G,\kappa)c_{Y_G}(a_U)\). It is nonzero on an active
block by the common scalar and the positivity of \(K_p(0,0)\).

\smallskip\noindent\emph{Independence and rank.}
In columns \(0\le t\le s\), the entries of \(V_r\) are constant
for \(t\ge r\). For \(1\le r\le s\) the difference between its
entries at \(r-1\) and \(r\) is \(w_r(F_1-1)\): only flags
with \(i=r,h=0\) contribute to the difference. This number is
nonzero. In the spherical case \(|x|\le2\sqrt p<p+1\); the other
two values of \(F_1\) are nonpositive. Successive differences of
columns therefore give independence. In fact the coordinates satisfy
\(V_r(t)=V_r(n-t)\), by interchanging \(i\) and \(h\) in the flag sum.
It follows that the period vectors span exactly the subspace fixed by
reversal of the \(n+1\) coordinates: its dimension is \(s+1\).

This gives the local determinant as well. Let
\(J:\mathbb R^{s+1}\to\mathbb R^{n+1}\) extend the coordinates
\(0,\ldots,s\) by reflection, counting the middle coordinate just once
when \(n\) is even. Form the square matrix with rows
\((V_r(0),\dots,V_r(s))\). Replace its first \(s\) columns by successive
differences and retain its last column. The zeroth row then has only
its final entry \(1\), and the remaining difference block is triangular
with diagonal \(w_r(F_1-1)\), \(1\le r\le s\). Taking the determinant
of the Gram matrix gives
\[
 \det\bigl(K_p(r,r')\bigr)_{0\le r,r'\le s}
 =\det(J^{\mathsf T}R_p^{-1}J)
    \prod_{r=1}^{s} w_r^2(F_1-1)^2 .
\]
Both sides include \(s=0\) with the empty product convention.
The determinant is positive, and the global tensor product
has rank \(\prod_{p\mid N}(s_p+1)=|I_G|\). The triangular change
of basis in Proposition~\ref{oldspace-height} preserves this rank.
A vector in the kernel of the Chow coefficient map pairs to zero
with every other vector, so full rank implies injectivity. Adding
the ranks of the mutually orthogonal blocks proves the theorem.
\end{proof}

\begin{proposition}[Absolute heights at cyclic level]\label{coprime-toric-block}
Let \(r^2\equiv\Delta\bmod4N\), with \(\Delta\) fundamental, and let
\(G\) be a normalized newform of weight \(2\kappa\ge4\), primitive
level \(M\mid N\), trivial character and root number
\(-\operatorname{sgn}(\Delta)\). Put
\[
 A=N/M,\qquad J_M=\{p:p\mid(M,\Delta)\},\qquad
 J_A=\{p:p\mid(A,\Delta)\}.
\]
Assume \(W_pG=G\) for \(p\in J_M\).
Fix a fundamental index
\(a_U=(m,\mu)\) with \(d\) prime to \(2N\Delta\), \(D=\Delta d<-4\)
fundamental, and \(m=|d|/(4N)\). Put \(a=\kappa-1\),
\(I_G=\{j:j^2\mid N/M\}\), and \(a_j=(j^2m,j\mu)\).
Write \(\mathcal K_p\) for the local matrices denoted by \(K_p\)
in Theorem~\ref{oldspace-full-rank}. Put \(K=\mathbb Q(\sqrt D)\)
and \(K_p=K\otimes\mathbb Q_p\), and let \(\Omega\) be the genus
character of \(K\).
Then, for \(j,k\in I_G\),
\begin{align}
 H(e_GZ_\kappa(a_j),Z_\kappa(a_k)) &=
 \frac{2^{|J_M|}(2\kappa-2)!\sqrt{|D|}\,m^a L'_K(G,\Omega,\kappa)}
 {2^{4\kappa-2}\pi^{2\kappa}\|G\|_{{\rm Pet},M}^{2}
                      [\Gamma_0(M):\Gamma_0(N)]}\,
 J_{jk}, \label{eq:split-absolute-matrix}\\[-2pt]
 J_{jk}&=(jk)^a\prod_{p\mid N}\mathcal K_p(v_p(j),v_p(k)). \notag
\end{align}
In particular these \(|I_G|\) projected Chow classes are linearly
independent whenever \(L'_K(G,\Omega,\kappa)\ne0\). For either
sign of \(\Delta\), if \(L'(G,\chi_\Delta,\kappa)\ne0\) then
infinitely many fundamental \(d\) give such an independent family.

\end{proposition}
\begin{proof}
\emph{Comparison with the primitive level.}
At every \(p\mid N\), the discriminant \(d\) is a unit square in
\(\mathbb Q_p\), so \(\Omega_p=\chi_{d,p}\circ\mathrm N=1\).
The field \(K_p\) splits away from \(J_M\cup J_A\) and is ramified
there. Each shared prime divides \(N\) exactly once, so the two sets
are disjoint and \(j,k\) are units at every shared prime. The root
number calculation in Lemma~\ref{ramified-primitive-period}, applied
at \(M\) with the assumed \(W_p\)-signs, gives
\(\varepsilon(G\otimes\chi_\Delta)=-1\). The same lemma's
newvector calculation gives the factor \(2^{|J_M|}\) in the primitive
self-height. Once the signs are specified its toric calculation uses
only \(G\) and the oriented CM orbit. For \(J_M=\varnothing\) use
Lemma~\ref{split-primitive-absolute}.

At a split prime the level flags and period coordinates are those of
Theorem~\ref{oldspace-full-rank}. At \(p\in J_M\) the level flag
is unchanged. At \(p\in J_A\) the representation is spherical at
level \(M\). Let \(O=\mathcal O_{K_p}\) and let \(\varpi\) be a
uniformizer. The unique \(O\)-stable cyclic edge joins \([O]\)
and \([\varpi O]\). Since \(\varpi^2 O=pO\), the torus element
\(\varpi\) exchanges the endpoints, also at \(p=2\).
Let \(v_0,v_1=\pi_p(\varpi)v_0\) be their unit spherical vectors.
The toric functional has equal values on them. This value is nonzero:
on the two unit cosets in \(K_p^\times/\mathbb Q_p^\times\)
the matrix coefficient of \(v_0\) is \(1\) and \(D_1=x/(p+1)\),
respectively, where \(x=a_G(p)/p^{\kappa-1}\).
Both cosets have the same measure and \(|x|\le2\sqrt p<p+1\).
Relative to the primitive value the period coordinates are therefore
\((1,1)\). Pairing them with \(R_p^{-1}\) gives \(\mathcal K_p(0,0)\),
just as at a split prime. More explicitly,
\[
 \mathcal K_p(0,0)
 = (1,1)\begin{pmatrix}1&D_1\\D_1&1\end{pmatrix}^{-1}
        \binom{1}{1}=\frac{2}{1+D_1}.
\]
The local level index is \(p+1\), so the net factor from this edge
and the change of level measure is \(2/(p+1+x)\), also for \(p=2\).
The local Haar measure and Euler factors cancel in the comparison
for the fixed field and character.

The genus character of each level ideal is one. At split primes,
reversing an apartment segment multiplies its newvector period by
the Fricke sign. The sign is one at a new spherical prime. At a
new ramified prime the exchange just computed fixes the period;
at the old ramified primes both the flag and its sign are unchanged.
Thus the two orientations produce the same multiplier at \(M\) and
\(N\). The scalar from rescaling a CM vector is \(j^a\), since
\(\chi_\Delta(j)=1\). All tests outside \(N\) agree.
Here is the comparison of level measures. At a compact cyclic level
\(J\) the test of one fundamental orientation is \(h\,\mathbf1_{gJ}\),
up to the common genus sign and unit weight, where \(h\) is the class
number. At a fixed fine normal level \(U'\subset J\) this formula sums
the cosets of \(U'\) in \(gJ\), giving the pullback of the whole orbit
by \cite[(2.5)--(2.7)]{QiuCM}. Its operator on \(\pi\) is \(h\)
times \(\operatorname{vol}(J)\) times the translate of the orthogonal
projection onto \(\pi^J\). The height at level \(J\) divides the
pairing of the two pulled-back cycles by the degree of the cover.
Consequently the two compact volumes in a self-pairing leave one
factor of \(\operatorname{vol}(J)\) in the level ratio. The central
stabilizer and determinant components cancel between the two levels
on the common full cover. More explicitly, if \(d_J\) is that cover degree and
\(J_N'=K_0(N)\subset J_M'=K_0(M)\), the factor is
\[
   \left(\frac{\operatorname{vol}(J_N')}
                   {\operatorname{vol}(J_M')}\right)^2
   \frac{d_{J_M'}}{d_{J_N'}}=[\Gamma_0(M):\Gamma_0(N)]^{-1}.
\]

At \(M\) the local vector on each newvector line has norm one.
Its toric value on the apartment is nonzero at split level primes:
in the Whittaker model the split toric integral of the newvector is
its standard local zeta integral \(L(1/2,\pi_p)\)
\cite[\S3D]{CaiShuTian}. Normalize that value to one. The pairing
of the lists of local period values at \(N\) is then
\(V_{v_p(j)}^{\mathsf T} R_p^{-1} V_{v_p(k)}\).
All lengths outside \(N\) and the normalized Hilbert class orbit
weights are the same at both levels. Lemma~\ref{finite-toric-tests}
therefore compares the level-\(N\) heights with the primitive
self-height, with the ratio just computed. If the derivative is zero, that lemma
gives zero for every entry. At primitive level the index is
\(m_M=|d|/(4M)=Am\). In one argument the fibre normalization changes
by \((j^2m/m_M)^{a/2}=j^a A^{-a/2}\). The two arguments therefore
replace \(m_M^a\) by \((jk)^a m^a\). Applying this to the primitive
self-height gives
\eqref{eq:split-absolute-matrix}. The positive definite matrix \(J\)
has full rank by the local independence proof in
Theorem~\ref{oldspace-full-rank}.

For either sign choose \(d\) by the nonvanishing theorem of
Bump, Friedberg and Hoffstein \cite[Theorem~(ii)]{BumpFriedbergHoffstein},
prescribing splitting at the primes dividing \(2N\Delta\).
It gives \(\varepsilon(G)d>0\) and \(L(G,\chi_d,\kappa)\ne0\).
Hence \(D=\Delta d<0\) is fundamental, \(d\) has a unit square
root modulo \(4N\), and \(L'_K(G,\Omega,\kappa)\ne0\) when
\(L'(G,\chi_\Delta,\kappa)\ne0\). Adding the primes of previous
choices to the prescribed splitting set gives infinitely many
distinct \(d\), also after any fixed finite exclusions.
This establishes the lower bound on the Chow dimension for the
negative as well as the positive twists.

\end{proof}

\subsection{Jacobi oldforms and their Petersson matrices}\label{sec:jacobi-petersson}
Skoruppa and Zagier give the holomorphic decomposition at arbitrary
index. The height comparison for negative twists also requires the
multiplicities of skew Jacobi forms at nonsquarefree index.
For the history of the skew decomposition see \cite[p.~207]{RSSTJacobi}.

\begin{lemma}[Skew Jacobi oldforms at arbitrary level]
\label{skew-index-raises}
Let \(k=\kappa+1\ge3\). Let \(G\) be a normalized newform of
weight \(2\kappa\), trivial character and exact conductor \(M\).
The skew-holomorphic Jacobi cusp forms of weight \(k\) and index
\(N\) with the good Hecke eigenvalues of \(G\) vanish unless
\(M\mid N\) and \(\varepsilon(G)=+1\). If these conditions hold,
put \(A=N/M\). This space has basis
\[
 \phi_G\big|U_iV_{A/i^2},\qquad i^2\mid A,
\]
where \(\phi_G\) spans the primitive line at index \(M\).
\end{lemma}
\begin{proof}
By Dirichlet's
theorem choose a prime \(q\equiv-1\bmod4N\), and put
\(\Delta=-q\), \(r=1\). Then \(\Delta\) is a negative fundamental
discriminant prime to \(2N\), \(r^2\equiv\Delta\bmod4N\), and
\(\Pi_{\Delta,N}=1\). Thus
Lemma~\ref{fundamental-separation} applies to the full skew Jacobi
space. For each \(p^2\mid N\), the isotropic lift from \(N/p^2\)
is the Jacobi oldform map \(U_p\). Let \(S_*\) be the orthogonal
complement of the sum of these lower-level images in the
\(G\)-eigenspace.
Lemma~\ref{fundamental-separation}, applied in the skew Weil
representation, shows that every nonzero vector in \(S_*\) has a
nonzero coefficient at a unit fundamental index.

Suppose the \(G\)-component occurs at \(N\). Choose the smallest
\(N_0=N/s^2\) at which it occurs. Its lower-level \(U_p\)-images
on this component vanish by minimality and good Hecke equivariance.
Applying the preceding separation argument at \(N_0\) gives a positive
fundamental \(d\) prime to \(2N_0\) and a residue \(r\) at which
coefficient evaluation is nonzero. The first coefficient of the
Shimura lift \(S_{d,r}\) equals this evaluation
\cite[Proposition~1]{SkoruppaBinary}. The lift is cuspidal at
\(N_0\) and preserves the good Hecke eigenvalues. Strong
multiplicity one therefore gives \(M\mid N_0\mid N\).

Let \(R\) represent the first Fourier coefficient on the elliptic
\(G\)-oldspace of level \(N_0\) under the Petersson pairing.
Its Fourier coefficients are real by conjugation symmetry.
Hecke equivariance implies that \(S_{d,r}^*R\) represents
evaluation at \((d,r)\) on the skew \(G\)-component and belongs to
that component. Evaluating this representative at \((d,r)\) gives \(\|S_{d,r}^*R\|^2>0\).
Equation~\eqref{eq:adjoint-kernel-coefficient} at \(N_0\) gives
\(L(R,\chi_d,\kappa)\ne0\). Writing \(R=\sum_{a\mid N_0/M}c_aG(az)\)
expresses this value as
\[
 L(G,\chi_d,\kappa)\sum_{a\mid N_0/M} c_a\chi_d(a)a^{-\kappa}.
\]
Hence the twist has nonzero central value. Since \(d>0\) and
\(\chi_d(p)=1\) at every \(p\mid M\) (at \(2\mid M\) use
\(d\equiv r^2\bmod8\)), its root number is \(\varepsilon(G)\).
This proves necessity of the sign condition.

We use the Jacobi packet with the normalization of
\cite[(6), Proposition~4.2]{SchmidtJacobiOld}. Its
\({\rm PGL}_2\)-label is the representation of \(G\). Indeed, at
almost every \(q\) the spherical Jacobi principal series of parameter
\(\chi\) has \(T_{\rm EZ}(q)\)-eigenvalue
\(q^{k-3/2}(\chi(q)+\chi(q)^{-1})\), and its image has Satake
parameters \(\chi(q),\chi(q)^{-1}\)
\cite[proofs of Proposition~4.2 and Theorem~5.1]{SchmidtJacobiOld}.
The good Shimura correspondence gives \(a_G(q)\) for this
eigenvalue. The finite double coset has upper triangular
representatives with positive lower diagonal entry, on which the
two classical Jacobi weight factors agree. This calculation
therefore applies also in skew signature. Strong multiplicity one
identifies the label.

The unary theta representations have holomorphic lowest weights
\(1/2\) and \(3/2\). Our theta components are holomorphic
lowest-weight vectors of weight \(\kappa+1/2\ge5/2\) after
conjugation in the skew signature. Thus no global unary summand
enters the cuspidal packet comparison below.
The Fourier--Jacobi expansion
identifies the central-character-\(\psi_N\) part of the adelic
Heisenberg quotient with one copy of its Schr\"odinger representation.
Locally the Jacobi representation factors as
\(\widetilde\pi_p\otimes\omega_{N,p}\), with \(\omega_{N,p}\) the
Schr\"odinger--Weil representation \cite[\S2, Propositions~2.1--2.2]{SchmidtJacobiOld}.
In the theta decomposition the Fourier coefficient pairs the
Whittaker functional on the first factor with the residue evaluation
on the second. At infinity the skew
theta expansion fixes a one-dimensional Gaussian and weight line.
At \(p\nmid N\) the local representation is unramified and the
integral Jacobi line is unique, including at \(p=2\)
\cite[\S8]{SchmidtJacobiOld}.

At \(p\mid N\), take two unit fundamental indices \((d,r)\) and
\((d',r')\), so \(d\equiv r^2\) and \(d'\equiv (r')^2\bmod4N\).
Choose \(\alpha,\alpha'\in\mathbb Z_p^\times\) with
\(\alpha^2=d\), \((\alpha')^2=d'\),
\(\alpha\equiv r\) and \(\alpha'\equiv r'\bmod 2N\mathbb Z_p\).
For odd \(p\) these lifts follow from Hensel's lemma. If \(p=2\)
and \(e=v_2(N)\), write
\(\alpha=r+2^{e+1}y\). The equation for \(y\) is
\[
 ry+2^e y^2=\frac{d-r^2}{2^{e+2}},
\]
whose derivative is an odd unit. Hensel's lemma gives the required
root; the same argument applies to \(d'\).
Set \(a=\alpha'/\alpha\). Then \(a^2=d'/d\) and
\(ar\equiv r'\bmod 2N\mathbb Z_p\). In both Jacobi signatures
the horizontal frequency is \((r^2-d)/(4N)\) \cite[(3.1)]{RSSTJacobi}.
Subtracting the theta frequency \(r^2/(4N)\) gives the Whittaker
frequency \(-d/(4N)\) with respect to the standard additive
character. The argument compares two such frequencies in the
same packet and uses only their quotient \(d'/d\).
The residue test is \(e_r\) in the finite Weil representation.
The element \(D_a=\operatorname{diag}(a,a^{-1})\) lies in the
\(\Gamma_0(N)\)-part of the local integral Jacobi subgroup because
\(a\) is a unit. It fixes invariant vectors.
Indeed, \(D_a n(x)D_a^{-1}=n(a^2x)\), so precomposing the
Whittaker functional with \(\widetilde\pi_p(D_a)\) multiplies its
parameter by \(a^2\). In the Schr\"odinger model
\(\omega_{N,p}(D_a)\varphi(x)\) is a nonzero scalar times
\(\varphi(ax)\). Evaluation at \(r/(2N)\) thus becomes evaluation
at \(ar/(2N)\), the required residue modulo \(\mathbb Z_p\).
These are the changes from \((d,r)\) to \((d',r')\).
Uniqueness of the local Whittaker functional therefore
makes the two coefficient functionals proportional. The local packet
criterion \cite[Theorem~3.1]{BaruchMao} now shows that every unit
fundamental test selects
the same local member at primes dividing \(N\). At all remaining
finite places \(G\) is unramified and the packet is a singleton.
The prescribed skew weight fixes the member at infinity. Global multiplicity one for the
metaplectic representation \cite{WaldspurgerShimura1980} shows that
all these tests lie on one automorphic summand. The local
proportionalities above make the global coefficient functionals
proportional on its full level-\(N\) invariant space. Since every
nonzero vector in \(S_*\) is detected by one of them,
\(\dim S_*\le1\).

Conversely, suppose \(\varepsilon(G)=+1\). At exact level \(N=M\), the \(G\)-eigenspace has no lower-level
\(U_p\)-old part, so it is \(S_*\) and has dimension at most one.
Choose a positive fundamental \(d\), split at every prime of \(2M\),
with \(L(G,\chi_d,\kappa)\ne0\). The nonvanishing theorem of
Bump--Friedberg--Hoffstein permits these local conditions \cite{BumpFriedbergHoffstein}.
Choose a root \(r^2\equiv d\bmod 4M\). The kernel \(S_{d,r}^*G\) is cuspidal by
\cite[Theorem~2]{SkoruppaBinary}. Equation
\eqref{eq:adjoint-kernel-coefficient} with \(R=G\) shows that its
coefficient at \((d/(4M),r)\) is nonzero. Hence it spans the
primitive line \(\mathbb C\phi_G\).

For general \(N\), choose such a \(d\) split at \(2N\). The
coefficient of \(\phi_G|V_A\) at its unit fundamental index is the
same nonzero primitive coefficient. Every lower-level \(U_p\)-image
has zero coefficient at unit residues. Thus the nonzero orthogonal projection of \(\phi_G|V_A\) spans
\(S_*\). Induction on \(N\) now gives the displayed family:
the lower-level bases contribute precisely the terms with \(i>1\),
and \(\phi_G|V_A\) contributes \(i=1\). To see independence, choose a root \(r\) of \(d\) modulo \(4N\).
For each \(i^2\mid A\), evaluate at discriminant \(i^2d\) and residue
\(ir\). The term indexed by \(j\) vanishes unless \(j\mid i\), since
\(r\) is a unit at \(N\). At \(j=i\) the coefficient formula for
\(V_{A/i^2}\) gives the nonzero primitive coefficient at
\((d,r\bmod 2M)\). The resulting triangular matrix is invertible.
\end{proof}

\begin{lemma}[Petersson matrices under a quadratic twist]
\label{jacobi-local-gram}
Fix \(M\mid N\). Let \(G\) have exact conductor \(M\) and root number
\(\varepsilon=\varepsilon(G)\). Let \(\phi_G\) be its primitive
Jacobi form of weight \(k=\kappa+1\), index \(M\), and signature
\(\varepsilon\), with the convention that signature \(-1\) is
holomorphic and signature \(+1\) is skew-holomorphic. Choose a fundamental
discriminant \(t\) with \(\operatorname{sgn}(t)=-\varepsilon\)
which belongs to \((\mathbb Z_p^\times)^2\) for every \(p\mid2N\).

Put \(G'=G\otimes\chi_t\), \(M'=Mt^2\), \(N'=Nt^2\),
\(A=N/M=N'/M'\), and \(P_i=U_iV_{A/i^2}\). Then \(G'\) has
conductor \(M'\) and root number \(-1\). Let \(\phi_{G'}\) be its
primitive holomorphic Jacobi form. Then
\[
 \frac{(P_i\phi_G,P_j\phi_G)_{{\rm Jac},N}}
      {\|\phi_G\|_{{\rm Jac},M}^{\,2}}
 =
 \frac{(P_i\phi_{G'},P_j\phi_{G'})_{{\rm Jac},N'}}
      {\|\phi_{G'}\|_{{\rm Jac},M'}^{\,2}},
 \qquad i^2,j^2\mid A.
\]
The same identity holds for the corresponding vector-valued
Petersson matrices.
\end{lemma}
\begin{proof}
For \(l=A/i^2\), the Jacobi operator \(P_i\) is \(l^{k/2-1}\)
times the sum of determinant-normalized slashes by the integral
matrices of determinant \(A\) divisible by \(i\), modulo left
\({\rm SL}_2(\mathbb Z)\). This includes imprimitive matrices
\cite[\S\S4--5]{SchmidtJacobiOld,SchmidtJacobiLocal}.
Let \(\Phi\) be the adelization of the primitive form. Use the rational
automorphism of the Jacobi group
\[
 \sigma_A(g,(\lambda,\mu,z))
   =(DgD^{-1},(\lambda,A\mu,Az)),\qquad D=\operatorname{diag}(A,1)
\]
of \cite[(3), Proposition~4.1]{SchmidtJacobiOld}. Pullback by
\(\sigma_A\) changes the central character from \(\psi_M\) to
\(\psi_{MA}\) at every place and preserves the quotient Haar measure
by the product formula. For a representative \(\alpha\), put
\(\gamma_\alpha=\alpha D^{-1}\), \(h_\infty=D/\sqrt A\), and
\(h_\alpha=(h_\infty,(\gamma_{\alpha,p}^{-1})_{p<\infty})\).
The adelic function corresponding to the individual slash is
\(\sigma_A^*R(h_\alpha)\Phi\), where \(R(h)\Phi(g)=\Phi(gh)\).

Indeed, at finite argument one, left multiplication by
\(\gamma_\alpha\in{\rm SL}_2(\mathbb Q)\) removes the finite components.
At infinity one has
\(\sigma_A(g)h_\infty=h_\infty U_{\sqrt A}(g)\), where \(U_s\)
fixes the linear group and takes the Heisenberg coordinates to
\((s\lambda,s\mu,s^2z)\). Since \(\gamma_\alpha h_\infty=\alpha/\sqrt A\),
this gives the determinant-normalized slash in either signature
\cite[\S3]{RSSTJacobi}. At \(D\) it is
\(A^{k/2}\phi(A\tau,Az)\). Individual summands may be paired at a
common finer compact open subgroup. After normalization by the primitive
norm, pullback by \(\sigma_A\) and the common unitary translation by
\(h_\infty\) cancel. The contribution of a pair of summands is thus
the product of their normalized finite matrix coefficients. The Haar
measure at finite level is normalized by the volume of the compact
open subgroup, giving the same ratio as the classical Petersson integral
for the full sums.

For each pair of representatives \(\alpha,\beta\) in the two sums,
both \(\alpha D^{-1}\) and \(\beta D^{-1}\) lie in
\({\rm SL}_2(\mathbb Z_p)\) at \(p\nmid A\). They fix the integral
Jacobi vector, so their normalized local matrix coefficient is one,
including at the extra primes dividing \(t\). At \(p\mid A\), the local
\({\rm GL}_2\)-representations of \(G\) and \(G'\) agree because
\(t\) is a local square. By the fundamental-coefficient formula
\eqref{eq:gzz-waldspurger} and the nonvanishing theorem of
Bump--Friedberg--Hoffstein, choose nonzero fundamental coefficients
of the permitted signs, with unit roots at every prime dividing \(2N\) for \(G\) and \(2N'\) for \(G'\).
First identify the indices \(M,M'=Mt^2\) by the local change of
Heisenberg coordinates \(U_t:(x,y,z)\mapsto(tx,ty,t^2z)\).
This identifies their Jacobi packets over the same \({\rm PGL}_2\)
representation \cite[Proposition~4.2 and Corollary~4.3]{SchmidtJacobiOld},
with the packet labels fixed in Lemma~\ref{skew-index-raises}.
The linear group is unchanged. Thus the absolute Whittaker
frequencies \(-d/(4M)\) and \(-d'/(4M')\) have square quotient
\(d'/(dt^2)\) at \(p\). Uniqueness of the packet member for this
Whittaker square class \cite[Theorem~3.1]{BaruchMao} selects the
same metaplectic representation, including at \(2\).
To identify the integral lines locally, use the
restricted tensor product realization of each occurring Jacobi
representation. The finite integral subgroup is a product of local
subgroups. Fix the prescribed archimedean weight vector and
nonzero integral vectors at the other finite places. Tensoring with
these vectors embeds the space of integral vectors at \(p\) into the
classical good Hecke eigenspace at the primitive index. That
eigenspace is a line by Lemma~\ref{skew-index-raises} in skew
signature and by \cite[Theorem~5]{SkoruppaZagier} in holomorphic
signature. The local integral space is therefore a line in each case.
Since \(t\) is a \(p\)-adic unit, the change \(U_t\) preserves
the integral subgroup, the normalized Haar measure and the
determinant-\(A\) correspondences. It therefore identifies these
two integral lines \cite[Propositions~2.1--2.2]{SchmidtJacobiOld}.

The local comparisons identify these products term by term. Summing over the
two finite sets of representatives proves the Jacobi Gram identity.
Passing to theta components multiplies both sides by
the same norm factor \(\sqrt{N/M}=\sqrt{N'/M'}\). In the skew case
this passage conjugates the matrix; the primitive Jacobi forms can be normalized to have real Fourier
coefficients, and the operators \(P_i\) preserve this real structure.
The matrix entries are therefore real.
\end{proof}

\subsection{The mixed height and the arithmetic inner product}\label{sec:mixed-inner-product}
\begin{theorem}[The fundamental height column at cyclic level]
\label{positive-full-mixed}
Let \(\kappa\ge2\), with \(\Delta\) fundamental and
\(r^2\equiv\Delta\bmod4N\). On a retained component of primitive
conductor \(M\), put \(A=N/M\) and \(J_M=\{p:p\mid(M,\Delta)\}\).
Take the real Jacobi basis \(g_i\) obtained from the primitive form
by \(U_iV_{A/i^2}\), \(i^2\mid A\). If \(U\) is a fundamental
CM datum with \((d_0,2N\Delta)=1\), pairing \(\Phi_\kappa\) with
\(Z_\kappa(U)\) gives a cusp form \(h_U\) satisfying
\[
 (h_U,g_i)_N=
 2^{|J_M|} C_{\kappa,N,|\Delta|}\,b_G(U)L'(G,\chi_\Delta,\kappa)
                 \frac{\tau_0(A/i^2)}{i^\kappa}.
\]
Consequently \eqref{intro-newform-height} holds simultaneously for
all harmonic inputs \(F\) at each such \(U\).
\end{theorem}
\begin{proof}
\smallskip\noindent\emph{Fourier evaluation on an oldspace.}
Fix a primitive constituent \(G\) of conductor \(M\mid N\) in the Jacobi signature \(\varepsilon(G)=-\operatorname{sgn}(\Delta)\).
Put \(A=N/M\), \(I=[\Gamma_0(M):\Gamma_0(N)]\), and use
\(J_M,J_A\) from Proposition~\ref{coprime-toric-block}.
The good eigenspace has a basis \(g_i\) corresponding to
\(\phi_G|U_iV_{A/i^2}\), \(i^2\mid A\), by
\cite[Theorem~5]{SkoruppaZagier} in the holomorphic signature and
Lemma~\ref{skew-index-raises} in the skew signature.
At a shared prime \(p\), one has \(p\Vert N\), so \(p\nmid i\).
The coefficient sums for \(V_{A/i^2}\) commute with the residue
reflection at \(p\). It descends to the reflection at \(M\) for
\(p\mid M\) and acts trivially on the primitive residues otherwise
\cite[\S2]{BSConverse}. Thus \(\Pi_{\Delta,N}\) retains the whole
block or kills it. On a retained block \(W_pG=G\) for \(p\in J_M\)
by Lemma~\ref{ramified-primitive-period}.

Choose a fundamental index \(v=(m,\mu_0)\) with
\(D=\Delta d<-4\), \(d\) split at \(2N\), prime to \(\Delta\),
and \(b=b_G(v)\ne0\), using \cite[Theorem~(ii)]{BumpFriedbergHoffstein}
and \eqref{eq:gzz-waldspurger} on the primitive line at \(M\).
Put \(B_{ji}=c_{g_i}(j^2m,j\mu_0)\) for \(j^2\mid A\) and
\(B_*=B/b\). We check that this matrix depends only on the local
Hecke data at primes dividing \(A\).

Write \(l=A/i^2\). For any positive integer \(j\), at discriminant \(j^2d\) a divisor \(v_1\mid l\)
in the \(V_l\) sum must divide the residue \(j\mu_0\), so \(v_1\mid j\).
The subsequent primitive coefficient occurs only if \(i\mid j/v_1\).
These divisibilities suffice, since the elliptic Fourier index scales
by \(j^2\). Its multiplier is \(v_1^\kappa\).
Applying the Shimura convolution \eqref{merged-shimura-coeff}, which
commutes with dilation of the coefficient sequence, gives
\[
 \operatorname{Sh}_v(g_i)(z)
       =b\sum_{v_1\mid l}v_1^\kappa G(iv_1z)
\]
in either signature. For a matrix row \(j^2\mid A\), we have \(\chi_d(h)=1\) for \(h\mid j\). Inversion
of the divisor sum gives
\[
 (B_*)_{ji}
 =\sum_{h\mid j}\mu_{\rm ar}(h)h^{\kappa-1}
    \sum_{\substack{v_1\mid l\\iv_1\mid j/h}}
          v_1^\kappa a_G\bigl(j/(hiv_1)\bigr),
\]
where \(\mu_{\rm ar}\) denotes the M\"obius function. This is lower
triangular, with diagonal one. For any fundamental \(U\) as in the
theorem, its
discriminant satisfies \(d_0\equiv\mu_0^2\bmod 4N\) with
\((\mu_0,2N)=1\). Thus \(\chi_{d_0}(p)=1\) at every \(p\mid N\),
which suffices for the divisors of \(j\). At odd level this argument
also allows \(d_0\equiv5\bmod8\). At any such \(U\), the vector
of fundamental evaluations is \(b_G(U)e_1\), with \(e_1\) selecting
\(i=1\). This includes \(b_G(U)=0\) and \(D_0=-3,-4\) whenever
the stated congruences hold. The calibration uses \(D<-4\) at the
reference index \(v\); evaluation at \(a_U\) already incorporates
its own stabilizer factor \(2/|\Gamma_x|\) from \eqref{rd-coeff}.

\smallskip\noindent\emph{The geometric height matrix.}
Let \(J\) be the matrix in \eqref{eq:split-absolute-matrix} and set
\(T=B_*^{-1}JB_*^{-\mathsf T}\). The proposition gives this absolute
height matrix for either sign of \(\Delta\). The root number on the
retained component gives
\(L'_K(G,\Omega,\kappa)=L(G,\chi_d,\kappa)L'(G,\chi_\Delta,\kappa)\).
At level \(M\), \eqref{eq:gzz-waldspurger} and
Lemma~\ref{split-primitive-absolute} identify the self-height
of the projected fundamental cycle with
\(2^{|J_M|}b^2 C_{\kappa,M,|\Delta|}L'(G,\chi_\Delta,\kappa)/\|g_M\|_M^2\),
using Lemma~\ref{ramified-primitive-period} at the shared primes.
Pull the height pairing back along the coefficient map
\(S_G^\vee\to W_N\) of Theorem~\ref{med-target}. Its matrix
in the basis dual to the \(g_i\) is
\[
 D_{G,\Delta}=
 \frac{2^{|J_M|} C_{\kappa,N,|\Delta|}\,L'(G,\chi_\Delta,\kappa)}
      {I A^{\kappa-1/2}\|g_M\|_M^2}\,T .
\]
Indeed the two geometric arguments give a factor \(A^{1-\kappa}\)
on passing from \(m_M=Am\) to \(m\), and
\(C_{\kappa,M,|\Delta|}/C_{\kappa,N,|\Delta|}=A^{-1/2}\).
The calculation divides by \(b\), and also applies when the derivative vanishes.

\smallskip\noindent\emph{Recovering the fundamental height column.}
Let \(\mathcal P_{ij}=(g_i,g_j)_N/\|g_M\|_M^2\).
Write \(h_U=(\mathrm{id}\otimes H(\,\cdot\,,Z_\kappa(U)))\Phi_\kappa\).
Since evaluation at \(U\) has vector \(b_G(U)e_1\), the coefficients
of the \(G\)-part of \(h_U\) in the real basis \((g_i)\) form the
column \(D_{G,\Delta}\,b_G(U)e_1\). Its Petersson products with
the basis forms are obtained by multiplying this column by
\(\|g_M\|_M^2\mathcal P\). We therefore need the identity
\begin{equation}\label{eq:fundamental-gram-column}
 \mathcal P T e_1=I A^{\kappa-1/2}
       \bigl(i^{-\kappa}\tau_0(A/i^2)\bigr)_{i^2\mid A}.
\end{equation}
Choose a fundamental discriminant \(t_0\) of sign \(-\varepsilon(G)\)
split at \(2N\), such that \(G'=G\otimes\chi_{t_0}\) has a simple
central zero \cite[Theorem~(i)]{BumpFriedbergHoffstein}.
Use \(M'=Mt_0^2\) and \(N'=Nt_0^2\). The ratio \(A\) and the level
index \(I\) are unchanged. The formula above gives the same \(B_*\),
since \(a_{G'}(n)=a_G(n)\) on integers supported at primes dividing \(A\).
These primes have the same local representations and level exponents,
hence the same \(J\), including at a new ramified level prime by the
equal endpoint periods of Proposition~\ref{coprime-toric-block}.
Outside \(A\) the normalized local factors of \(J\) equal one.
Lemma~\ref{jacobi-local-gram} gives the same \(\mathcal P\) even
when \(t_0<0\) and the Jacobi signature changes.

For \(G'\) at \(N'\) take the trivial twist and a fundamental
coefficient \(b'\ne0\). Equating its geometric column with
\eqref{height-old-line} and cancelling \(b'L'(G',\kappa)\ne0\)
gives \eqref{eq:fundamental-gram-column}. The matrices and the
right side of this identity agree for \(G\) and \(G'\), so it
also applies when \(L'(G,\chi_\Delta,\kappa)=0\).
For \(i^2\mid A\) put \(l=A/i^2\). Every prime of \(J_A\) has
exponent one in \(l\). At such a prime the divisor count has local
factor \(2\), while the character sum has factor
\(1+\chi_\Delta(p)=1\). At every other prime \(p\mid l\) both
factors are \(1+v_p(l)\), and \(\chi_\Delta(i)=1\). Hence
\[
       2^{|J_M|}\tau_0(l)
       =2^{|J_M|+|J_A|}\sum_{v\mid l}\chi_\Delta(v).
\]
The vector identity now gives each real Petersson coordinate:
\[
 (h_U,g_i)_N=
 2^{|J_M|+|J_A|} C_{\kappa,N,|\Delta|}\, b_G(U)\,
 i^{-\kappa}\sum_{v_1\mid l}\chi_\Delta(v_1)
 L'(G,\chi_\Delta,\kappa).
\]
Indeed the Shimura convolution at \(U\) has \(b_G(U)\) in place of
\(b\). Its twisted Dirichlet series has finite factor
\(i^{-s}\sum_{v_1\mid l}\chi_\Delta(v_1)v_1^{\kappa-s}\).
The derivative of this factor is multiplied by \(L(G,\chi_\Delta,\kappa)=0\).
The formula therefore includes \(b_G(U)=0\) as well. On a discarded block the
Chow projection vanishes and the full Shimura residue sum cancels
under the same reflection. The oldforms span the cusp space, so
the principal-part pairing \eqref{merged-height-shadow} gives
\eqref{intro-newform-height} for every \(F\), with its indicated
conjugation. This last equality pairs Chow classes over a common
field and fine level. The fundamental column was fixed before the
choice of \(F\), so its polar support may meet the support of
\(Z_\kappa(U)\).
\end{proof}

\begin{proof}[Proof of Theorem~\ref{intro-newform-mixed}]
The identity for a fixed \(U\) is Theorem~\ref{positive-full-mixed}.
For simultaneous admissibility, let \(I\) be the union of the polar
supports of a basis of the given space. Every input has polar
support contained in this finite set. Put
\(\sigma=-\operatorname{sgn}(\Delta)\). By Dirichlet's theorem
there are infinitely many primes \(q\equiv\sigma\bmod 4N\).
Omit the primes dividing \(2N\Delta\) or any \(4Nn\) with
\((n,\mu)\in I\). Then \(d_0=\sigma q\) is a fundamental
discriminant congruent to \(1\bmod 4N\), so we can take
\(\mu_0=1\). The products \(\Delta d_0\) are distinct negative
fundamental discriminants. Only finitely many give one of the CM
fields occurring in \(I\). Every remaining \(U\) is admissible
for the whole space.
\end{proof}

For real \(F\), substitution in \eqref{total-height-discrepancy}
now gives \eqref{finite-signed-comparison} at every admissible \(U\),
with the orientations and measures of that comparison. This determines
the sum of the finite contributions. To deduce vanishing at a specified
prime one still needs a local argument, such as the disjointness or
sign cancellation in Proposition~\ref{finite-genus-cancellation}.

\begin{corollary}[Twisted oldspaces]\label{twisted-oldspace-rank}
Let \((N,\Delta,r)\) be compatible, with \(\Delta\) fundamental, allowing \(\Delta=1\).
In either signature the rank and height formulas of
Theorem~\ref{oldspace-full-rank} hold on the subspace fixed by
\(\Pi_{\Delta,N}\), with \(L'(G,\chi_\Delta,\kappa)\) in place of
\(L'(G,\kappa)\). On a retained component of primitive conductor \(M\)
the rank of the height pairing is \(\#\{i:i^2\mid N/M\}\) if the derivative is nonzero
and zero otherwise. If the derivative is nonzero, the shadows on
this block whose heights vanish against every admissible fundamental
cycle form a hyperplane independent of the compatible twist
retaining the block. Every shadow outside that hyperplane has a
nonzero mixed height for \(\gg X/\log X\) admissible discriminants
\(|d_0|\le X\), for the fixed component and polar data.

Let \(\mathbf b(a)=(c_{g_i}(a))_i\) and choose a fundamental
\(v\) with \((d_v,2N\Delta)=1\) and \(b_G(v)\ne0\).
With \(T=B_*^{-1}JB_*^{-\mathsf T}\) as above, one projected
self-height determines every global height on this component:
\begin{equation}\label{eq:calibrated-projected-heights}
 H(e_G Z_\kappa(a),e_G Z_\kappa(b))=
 \frac{\mathbf b(a)^{\mathsf T}T\mathbf b(b)}
      {b_G(v)^2 T_{11}}\,
 H(e_G Z_\kappa(v),e_G Z_\kappa(v)).
\end{equation}
The indices \(a,b\) are arbitrary, including \(a=b\).
\end{corollary}
\begin{proof}
Theorem~\ref{positive-full-mixed} gives the height matrix
\[
 D_{G,\Delta}=2^{|J_M|}\sqrt{|\Delta|}\,
     L'(G,\chi_\Delta,\kappa)\,Q_{G,N},\qquad
 Q_{G,N}=\frac{C_{\kappa,N,1}}{I A^{\kappa-1/2}\|g_M\|_M^2}\,T .
\]
The local period vectors are independent by Theorem~\ref{oldspace-full-rank},
so \(J\) and \(Q_{G,N}\) are positive definite. This gives the stated rank.
For fixed \(G,N\) and primitive normalization, \(B_*\) uses only the
Hecke data at \(A\). The normalized endpoint periods agree at split
and ramified level edges. In the ratios in \eqref{eq:oldspace-local-height}
the tests outside \(jk\) cancel, including at every \(p\mid(N,\Delta)\),
where \(p\Vert N\). Thus \(Q_{G,N}\) is independent of the compatible
twist retaining the block.
For two compatible twists \(\Delta_j\) of the same sign retaining
\(G\), fix the same oldform basis and write
\(c_j=2^{\omega((M,\Delta_j))}\sqrt{|\Delta_j|}\) and
\(L'_j=L'(G,\chi_{\Delta_j},\kappa)\).
Let \(H_j^G(a,b)\) be the height of the two projected coefficients
for the \(j\)-th twist. Evaluating the common matrix \(Q_{G,N}\)
on \(\mathbf b(a),\mathbf b(b)\) gives
\[
 c_2L'_2 H_1^G(a,b)=c_1L'_1 H_2^G(a,b).
\]
The identity holds also for coincident or imprimitive indices and
when either derivative vanishes.

By Chow-valued modularity the height at \(a,b\) is
\(\mathbf b(a)^{\mathsf T}D_{G,\Delta}\mathbf b(b)\).
At \(v\) one has \(\mathbf b(v)=b_G(v)e_1\). Since \(T_{11}>0\),
\eqref{eq:calibrated-projected-heights} follows even at zero derivative.
The two-variable series is the sum of these Hecke-orthogonal blocks.

For the hyperplane assertion suppose the projected shadow belongs
to this block and write \(\Pi_{\Delta,N}\xi_wF=\sum_i x_i g_i\). Put
\[
 \lambda_{G,N}(x)=\sum_{i^2\mid A}
                x_i\,i^{-\kappa}\tau_0(A/i^2).
\]
At each admissible \(U\), Theorem~\ref{positive-full-mixed}
expresses the mixed height as
\[
 2^{|J_M|} C_{\kappa,N,|\Delta|}\,L'(G,\chi_\Delta,\kappa)\,
 b_G(U)\,\overline{\lambda_{G,N}(x)} .
\]
The coefficient of \(x_1\) is \(\tau_0(A)>0\). For any fixed lift,
the prescribed-sign count in Section~\ref{sec:primitive-pairings}
gives \(\gg X/\log X\) admissible discriminants with \(b_G(U)\ne0\).
Thus the common kernel is \(\ker\lambda_{G,N}\), which depends only
on the oldspace in the indicated basis. The full height
matrix is nonsingular when the derivative is nonzero, so the
imprimitive rows of \(B\) detect the additional directions.

For a shadow with components \(x_G\) in several retained blocks, the
same argument gives an exact nonvanishing criterion. Write \(g_{G,1}\)
for the first oldform in each block. At the unit fundamental indices,
the Fourier coefficients of \(g_{G,1}\) and
\(\operatorname{pr}_{S_*}g_{G,1}\) agree, since the other summand is
a sum of \(U_p\)-images. Each projection is nonzero by the choice of
a calibration index, and projections from distinct Hecke eigenspaces
are orthogonal. Consequently the mixed heights vanish at all admissible
fundamental indices if and only if
\[
 L'(G,\chi_\Delta,\kappa)\lambda_{G,N}(x_G)=0
 \quad\text{for every retained }G.
\]
Indeed, on these indices the heights are the coefficients of the sum
of the projected \(g_{G,1}\), with weights
\(2^{|J_M|} C_{\kappa,N,|\Delta|}L'(G,\chi_\Delta,\kappa)
 \overline{\lambda_{G,N}(x_G)}\) on the block of conductor \(M\).
If this form is nonzero,
Lemma~\ref{fundamental-separation} gives \(\gg X^\eta\) distinct
admissible discriminants \(|d_U|<X\) with nonzero mixed height,
for each \(0<\eta<1\), after the polar primes have been excluded.
\end{proof}
\begin{proof}[Proof of Theorem~\ref{ks-arithmetic-inner-product}]
\emph{Arithmetic representatives.}
By Proposition~\ref{ks-admissible-restriction}, contraction of the
arithmetic series gives the admissible lift of
\(\theta_{\Delta,r,\kappa}(f)=Z_\kappa(F)\) for any \(\xi_wF=f\).
A finite realization on each retained component uses one CM field.
Take \(v=(m,\mu)\) with \(b_G(v)\ne0\) as in
Proposition~\ref{coprime-toric-block} and put
\(C_j=e_G Z_\kappa(j^2m,j\mu)\), \(j^2\mid N/M\).
If the coefficient classes dual to the real Jacobi basis form
the column \(X\), modularity gives \(C=BX\) with
\(B_{ji}=c_{g_i}(j^2m,j\mu)\). The invertible matrix \(B\) of
Theorem~\ref{positive-full-mixed} therefore gives the arithmetic
representatives \(B^{-1}\widehat C\).
By \cite[\S\S2.2.4--2.2.5, (2.12)]{QiuCM} the signed absolute
arithmetic product of two such lifts of \(X,Y\) equals the
bilinear height \(H(X,\overline Y)\). The bar conjugates only
the complex coefficients. Trace and pullback make the number
independent of the common field and fine level.

\emph{The full eigenspace.}
In the real Jacobi basis of Theorem~\ref{positive-full-mixed} write
\((\Phi_\kappa^{(N)})^G=\sum_i g_i\otimes X_i\), and let
\(P=((g_i,g_j)_N)_{ij}\). Put \(A=N/M\),
\(I=[\Gamma_0(M):\Gamma_0(N)]\), and let \(g_M\) be the
primitive form at index \(M\). With \(T=B_*^{-1}JB_*^{-\mathsf T}\)
from that theorem, the matrix of the operator is
\begin{equation}\label{eq:ks-old-operator}
       E_{G,N}=\frac{T P}{I A^{\kappa-1/2}\|g_M\|_M^2}.
\end{equation}
By Corollary~\ref{twisted-oldspace-rank}, \(T\) is positive definite
and independent of the compatible twist retaining the block.
Thus \(P E_{G,N}\) is symmetric positive definite, so \(E_{G,N}\)
is positive and self-adjoint for the Petersson product.

If \(f=\sum x_i g_i\) and \(h=\sum y_i g_i\), the coefficients of
their lifts in the expansion above are \(P\bar x\) and \(P\bar y\).
Write \(D_{G,\Delta}=(H(X_i,X_j))_{ij}\), the geometric matrix
computed in Theorem~\ref{positive-full-mixed}. Contracting it
with the two Petersson evaluations gives
\[
  \bar x^{\mathsf T}P D_{G,\Delta}P y
   = 2^{\omega((M,\Delta))} C_{\kappa,N,|\Delta|}
          L'(G,\chi_\Delta,\kappa)\,(E_{G,N}h,f)_N.
\]
Distinct eigenspaces are orthogonal by self-adjointness of the good
Hecke operators. At \(M=N\), \(E_{G,M}=1\) by
\eqref{eq:fundamental-gram-column}. Here \((\Phi_\kappa^{(M)})^G=g\otimes X_g\)
and \(\theta_{\Delta,r,\kappa}(g)=\|g\|_M^2 X_g\), giving the
stated primitive norm.

\emph{The factors at the primes dividing \(A\).}
In the basis indexed by the tuples \((v_p(i))_{p\mid A}\), one has
\begin{equation}\label{eq:ks-local-tensor}
 E_{G,N}=\bigotimes_{p^n\parallel A} E_{G,Mp^n}.
\end{equation}
Each factor is computed from \eqref{eq:ks-old-operator} with the
same primitive form and index ratio \(p^n\). The divisor sum for
\(B_*\) in Theorem~\ref{positive-full-mixed} factors by multiplicativity
of the Hecke coefficients. The matrix \(J\) factors by
\eqref{eq:split-absolute-matrix}. For the normalized Petersson matrix
\(P/\|g_M\|_M^2\), use the determinant-\(A\) slashes in the proof
of Lemma~\ref{jacobi-local-gram}. Their representatives modulo
left \({\rm SL}_2(\mathbb Z)\) correspond to index-\(A\) row
lattices in \(\mathbb Z^2\), hence to independent local row lattices
of index \(p^n\). Divisibility of a representative by \(i\) imposes
divisibility by \(p^{v_p(i)}\) on its local rows. At \(p^n\parallel A\)
right multiplication of both representatives \(\alpha,\beta\) by
\(\operatorname{diag}(p^n/A,1)\in{\rm GL}_2(\mathbb Z_p)\)
gives determinant \(p^n\) and preserves \(\alpha\beta^{-1}\).
The local matrix coefficient is therefore the one at index ratio
\(p^n\). The two sums of representatives, their slash factors and
the factor \(\sqrt A\) from theta decomposition all separate over
the primes. This proves the tensor identity, since the level index
\(I\) and \(A^{\kappa-1/2}\) also factor.

At \(p\parallel A\) the factor is one-dimensional and
\eqref{eq:fundamental-gram-column} at level \(Mp\) gives \(E_{G,Mp}=2\).
Thus primes dividing \(A\) exactly once contribute the scalar
\(2^{\#\{p:p\parallel A\}}\), including when other primes divide
\(A\) to higher order. For squarefree \(A\) this is the entire
operator. If \(p^2\mid A\) for some \(p\), the same column identity
at level \(N\) gives a positive off-diagonal entry of \(E_{G,N}^{\mathsf T}\).
In all cases replacing both inputs in \eqref{eq:ks-two-lifts} by their
images under \(E_{G,N}^{-1/2}\) gives the ordinary Petersson product,
including when the derivative vanishes.
\end{proof}

In particular, let \(\operatorname{rad}H\) denote the radical of
the pairing on the full coefficient span \(W_N\). The map
\(Y\mapsto h_Y\) gives a good Hecke-equivariant isomorphism
\[
 W_N/\operatorname{rad}H
   \ \xrightarrow{\ \sim\ }\
   \bigoplus_{\substack{G\ {\rm retained}\\
                         L'(G,\chi_\Delta,\kappa)\ne0}} S_{G,N},
\]
where \(S_{G,N}\subset S\) is the level-\(N\) eigenspace of \(G\).
Indeed its kernel is the radical because the Fourier coefficients of
\(\Phi_\kappa\) span \(W_N\). On the dual of each retained block its
composite with the Chow coefficient map has matrix \(D_{G,\Delta}\).
This matrix is invertible on the indicated blocks and zero on the
others. Thus any remaining kernel of \(S^\vee\to W_N\) lies on
the blocks with zero central derivative.

\smallskip\noindent\emph{Detecting weak holomorphy by heights.}
Fix a retained \(G\)-component with \(L'(G,\chi_\Delta,\kappa)\ne0\),
and a fundamental \(v=(m,\mu)\) as in Proposition~\ref{coprime-toric-block}
with \(b_G(v)\ne0\). Put \(I_G=\{j:j^2\mid N/M\}\) and
\(v_j=(j^2m,j\mu)\). For every invariant harmonic input \(F\) with
shadow in this component,
\begin{equation}\label{eq:height-detects-weak}
\begin{aligned}
 F\in M^!_{w,\rho_{\rm in}}
 &\quad\Longleftrightarrow\quad Z_\kappa(F)=0\ \text{in Chow}\\
 &\quad\Longleftrightarrow\quad
 H(Z_\kappa(F),Z_\kappa(v_j))=0\quad\text{for all }j\in I_G .
\end{aligned}
\end{equation}
Indeed, the fundamental coefficient \(b_G(v)\ne0\) gives
\(L(G,\chi_{d_v},\kappa)\ne0\) by Waldspurger at the primitive
level. Thus \(L'_K(G,\Omega,\kappa)\ne0\) by genus factorization.
The projected \(Z_\kappa(v_j)\) have nonsingular Gram matrix
by Proposition~\ref{coprime-toric-block} and span the Chow image of
the block by Corollary~\ref{twisted-oldspace-rank}. Self-adjointness
of \(e_G\) permits the unprojected tests displayed above. Their
simultaneous vanishing forces the coefficient functional of \(F\)
to vanish on the block, by injectivity there. It already vanishes
on the other blocks by its shadow, so \(\xi_wF=0\) by the
principal-part pairing. The forward implications follow from
Theorem~\ref{intro-chow-vanishing} and bilinearity. The elliptic
curves supporting these tests have endomorphism orders in the same
imaginary quadratic field. The count \(|I_G|\) is minimal among
families of complex-linear scalar height tests on the input. The
Chow image has dimension \(|I_G|\), and fewer linear functionals
have a nonzero common kernel there.

\begin{corollary}[A higher-weight multiquadratic height determinant]
\label{higher-multiquadratic-determinant}
Let \(K/\mathbb Q\) be totally real and multiquadratic of degree
\(d\). For each character \(\chi\) of its Galois group let
\(\Delta_\chi\) be the corresponding fundamental discriminant,
with \(\Delta_{\mathbf1}=1\). Fix compatible roots at level \(N\)
for all these discriminants. Suppose the good Hecke block of a
normalized newform \(G\) of trivial character, primitive level \(M\mid N\) and weight
\(2\kappa\ge4\) occurs and is retained for every \(\chi\).
Put \(s=\#\{i:i^2\mid N/M\}\). Use the same real basis \((g_i)\)
as in Theorem~\ref{positive-full-mixed}. For each twist, let
\(X_{\chi,i}\) be the image of the \(i\)-th dual basis vector under
its Chow coefficient map on this block. Write
\(Y_{\chi,i}=\Delta_\chi^{-1/4}X_{\chi,i}\). These are classes in
\(\CH^\kappa_{\rm hom}(\mathcal Y_\kappa(N)_K)_{\mathbb R}\) for the
chosen real embedding of the coefficient field.
With \(Q_{G,N}\) as in Corollary~\ref{twisted-oldspace-rank}, their
\(ds\)-by-\(ds\) height matrix satisfies
\[
 \frac{\det\bigl(H(Y_{\chi,i},Y_{\psi,j})\bigr)}
      {(\det Q_{G,N})^d}
 =\left(2^{\,b_M}\frac{L^{(d)}(G/K,\kappa)}{d!}\right)^s,
\qquad b_M=\sum_{p\mid M}\left(d-\frac d{|I_p|}\right).
\]
Here \(I_p\) is the inertia group at \(p\), and \(L(G/K,u)\) is
the finite base-change \(L\)-function in classical normalization,
with centre \(u=\kappa\). Both indices of the matrix use the same
ordering. If \(e\) characters have zero central derivative, the
matrix has nullity \(se\) and \(\operatorname{ord}_{u=\kappa}L(G/K,u)
\ge d+2e\). In particular, order \(d+2\) gives rank \(s(d-1)\).
\end{corollary}
\begin{proof}
Proposition~\ref{coef-galois-character} makes the matrix block
diagonal in \(\chi\). The rescaling removes \(\sqrt{\Delta_\chi}\)
from the matrix in Corollary~\ref{twisted-oldspace-rank}. The
remaining block equals \(Q_{G,N}\) times
\(2^{\omega((M,\Delta_\chi))} L'(G,\chi,\kappa)\). At \(p\mid M\),
exactly \(d/|I_p|\) characters are unramified, giving the exponent
\(b_M\) on multiplying the block scalars. Each retained twist has
root number \(-1\), by the sign calculation at the primitive level.
Artin factorization gives \(L(G/K,u)=\prod_\chi L(G,\chi,u)\).
At a prime ramified in a character and dividing \(M\), compatibility
forces \(p\Vert M\). The ramified twist of the special representation
there has Euler factor one, in agreement with the Dirichlet twist.
Thus its coefficient of \((u-\kappa)^d\) is
\(\prod_\chi L'(G,\chi,\kappa)\), also when a factor has a higher
order zero. Taking determinants proves the formula. Since \(Q_{G,N}\)
is nonsingular, each zero derivative removes one block of rank \(s\).
As in Theorem~\ref{thm:multiquadratic-height}, the completion factors
are nonzero at the centre. Write the orders as \(1+2m_\chi\).
Then the product has order \(d+2\sum_\chi m_\chi\), and a first
derivative vanishes precisely when \(m_\chi\ge1\). A total excess
of two means that one twist has order three and all others have order one.
\end{proof}

For a single line in the modular variable, Chow-valued modularity
itself gives ratios of cycles before taking heights.
\begin{corollary}[Chow classes on a projected line]
\label{one-cycle-height}
Suppose \(\Phi_\kappa\in S\otimes A_N\), with the level, weight and twist
of Theorem~\ref{med-target}. Choose \(0\ne g\in S\). Let \(Z_\kappa^{[g]}(a)\) denote the coefficient at \(a\) of
the Petersson projection of \(\Phi_\kappa\) onto \(\mathbb Cg\) in the
modular variable. Put \(b_g(a)=c_g(a)\). For any \(a_0\) with
\(b_g(a_0)\ne0\) we have
\begin{equation}\label{eq:one-cycle-ratio}
             Z_\kappa^{[g]}(a)=\frac{b_g(a)}{b_g(a_0)}Z_\kappa^{[g]}(a_0).
\end{equation}
In particular, for real \(g\) and a bilinear height on these Chow
classes, all heights on this line follow from one reference self-height,
whenever the pairing is defined:
\[
 H(Z_\kappa^{[g]}(a),Z_\kappa^{[g]}(b))=
   \frac{b_g(a)b_g(b)}{b_g(a_0)^2}\,H(Z_\kappa^{[g]}(a_0),Z_\kappa^{[g]}(a_0)).
\]

Under the hypotheses of Theorem~\ref{height-kernel},
take \(g=g_G\). Let \(\mathcal H_g\) be the summand indexed by \(G\)
in \eqref{eq:spectral-heights}. Two reference indices \(a_0,v\)
with \(b_g(a_0)b_g(v)\ne0\) give a formula using a mixed height:
\begin{equation}\label{eq:one-height-ratio}
   \mathcal H_g(a,b)=
       \frac{b_g(a)b_g(b)}{b_g(a_0)b_g(v)}\,\mathcal H_g(a_0,v).
\end{equation}
\end{corollary}
\begin{proof}
In an orthogonal basis containing \(g\), Theorem~\ref{med-target}
gives \(Z_\kappa^{[g]}(a)=b_g(a)Z_\kappa(F_g)\) for
\(\xi_w F_g=g/\|g\|_{\rm Pet}^2\). Evaluating at \(a_0\) gives
\eqref{eq:one-cycle-ratio}. Apply bilinearity for the self-height
and the factorization \(\mathcal H_g\in\mathbb C(g\otimes g)\)
for \eqref{eq:one-height-ratio}. Both hold at zero height as well.
\end{proof}

Under the untwisted squarefree-level hypotheses of
Theorem~\ref{height-kernel}, let \(g=g_G\) be the real level-\(N\)
generator attached to a newform of conductor \(M\mid N\). Expand the functional
\(\ell_g=(\,\cdot\,,g)_{\rm Pet}/\|g\|_{\rm Pet}^2\) in each of two
fundamental evaluation bases with pairwise coprime discriminants
(Lemma~\ref{fundamental-separation}). Their images under the Chow
coefficient map give representatives \(A,B\) of \(X_g\) in disjoint
CM fibres. Evaluating \eqref{eq:spectral-heights} on this pair gives
\[
 H(A,B)=(\ell_g\otimes\ell_g)\mathcal H
 =C_{\kappa,N,1}\,2^{\omega(N/M)}
       \frac{L'(G,\kappa)}{\|g\|_{\rm Pet}^2}.
\]
At \(M=N\), multiplication of both representatives by \(b_g(U)\)
recovers \eqref{eq:intro-gzz} for \(\Delta=1\) by \eqref{eq:gzz-waldspurger}.
For fundamental \(U,V\) of coprime discriminants with
\(b_g(U)b_g(V)\ne0\), choose \(F_g\) with normalized shadow and
polar support prime to \(D_V\). Then \eqref{med-untwistedheight}
and \eqref{height-old-line} give
\[
      \mathcal H_g(a_U,a_V)=b_g(U)\,H(Z_\kappa(F_g),Z_\kappa(V)).
\]
For several eigenlines one can choose a common reference class
as an integral linear combination of fundamental cycles. Take a basis
of fundamental evaluations \(\ell_{U_1},\ldots,\ell_{U_s}\) as in
Lemma~\ref{fundamental-separation}. For a positive integer \(m\) set
\[
 Y_m=\sum_{i=1}^s m^{i-1}Z_\kappa(U_i),\qquad
 c_G(m)=\sum_{i=1}^s m^{i-1}b_{g_G}(U_i).
\]
Each of the \(s\) polynomials \(c_G\) is nonzero by the basis property
and has degree at most \(s-1\). Some integer \(1\le m\le s(s-1)+1\)
therefore avoids all their roots. By multiplicity one, a
real linear combination \(T\) of finitely many good Hecke operators
has distinct eigenvalues \(\lambda_G\) on the lines in the
squarefree-level expansion. The projector
\(e_G=\prod_{J\ne G}(T-\lambda_J)/(\lambda_G-\lambda_J)\)
satisfies \(e_GY_m=c_G(m)X_{g_G}\) by geometric Hecke equivariance.
Summing the projections gives
\begin{equation}\label{eq:hecke-one-cycle}
 Z_\kappa(a)=\sum_G\frac{b_{g_G}(a)}{c_G(m)}\,e_GY_m .
\end{equation}
If a single fundamental \(U\) has nonzero coefficient on every line,
the same formula holds with \(Y_m=Z_\kappa(U)\) and
\(c_G(m)=b_{g_G}(U)\).

The rank can be read from a matrix of disjoint geometric heights.
Choose CM data \((U_i),(V_j)\) whose Fourier evaluations give the two bases in
Lemma~\ref{fundamental-separation}. Form the height matrix and the
matrix of the dual Petersson pairing in these bases of \(S^\vee\),
\[
    H_{ij}=H(Z_\kappa(U_i),Z_\kappa(V_j)),\qquad
    K_{ij}=\sum_G
       \frac{b_{g_G}(U_i)b_{g_G}(V_j)}{\|g_G\|_{\rm Pet}^2}.
\]
The Petersson matrix \(K\) is invertible. The diagonal expansion
\eqref{eq:spectral-heights} gives the characteristic polynomial
\begin{equation}\label{eq:height-pencil}
   \frac{\det(tK-H)}{\det K}
    = \prod_{\substack{M\mid N\\G\ {\rm at}\ M}}
       \left(t-C\,2^{\omega(N/M)}L'(G,\kappa)\right).
\end{equation}
As in Theorem~\ref{primitive-height-pencil}, the two evaluation
matrices cancel in the determinant quotient. The nonzero roots,
counted with multiplicity, give the rank using heights between
distinct CM fields.
The cycle-vanishing questions for a nonzero shadow and the positive
coefficients of its harmonic lift are discussed in
Section~\ref{rd-h-height}.
\appendix

\section{Scalar tests for both Weil representations}\label{med-all-signs}

The trace polynomials in Section~\ref{rd-m-untwisted} have a geometric
interpretation. Here a matrix of unary theta functions gives scalar
identities directly in the Chow coefficient space for both Weil
representations and every compatible twist \((\Delta,r)\).
Let \(\rho_+=\bar\rho_L\) be the Weil representation attached to
\(\mu\mapsto\mu^2/(4N)\), and put \(a=\kappa-1\). Write
\[
 (-1)^\epsilon=\operatorname{sgn}(\Delta)(-1)^a,\qquad
 \epsilon\in\{0,1\},\qquad v=\tfrac12+\epsilon .
\]
Take one representative of each independent reflection orbit: for
even parity use \(I=\{0,\ldots,N\}\), and for odd parity
\(I=\{1,\ldots,N-1\}\). Let \(h=|I|\). The additional shared-prime
symmetries can be kept in the coefficients; they need not be
incorporated into this choice of coordinates. Denote by \(G_I\)
the diagonal matrix of the reflection orbit sizes. In particular
its entries at \(0,N\), when present, are one. For a formal cycle
series write \(X=(X_\mu)_{\mu\in I}\) as a column.

Use the unary theta vector and its Serre derivatives of weights
\(v_j=v+2j\), \(0\le j<h\):
\[
 \theta_\mu=\sum_{t\equiv\mu\ (2N)}t^\epsilon q^{t^2/(4N)},
 \qquad U_0=\theta,\qquad U_{j+1}=\mathcal D_{v_j}U_j .
\]
Use \(t^0=1\), including at zero. The square matrix \(M\) has
column \(j\) equal to \((U_j(\mu))_{\mu\in I}\).
The determinant computation of Lemma~\ref{med-u-wronskian} applies
to either choice of parity: its proof uses only the initial exponents
and the transformation of the theta vector. Consequently
\[
       \det M=c\,\eta^{2W},\qquad W=h(v+h-1),\quad c\ne0 .
\]
All entries of \(M\) are holomorphic, and \(M\) is invertible throughout
the upper half-plane. If \(h=0\) the series is identically zero by the
reflection symmetry. In what follows assume \(h>0\), let
\(\Delta_{\rm mod}=\eta^{24}\), and choose an integer \(\ell_0>W/12\)
large enough that \(k-v_{h-1}+12\ell_0>2\).

\begin{theorem}[Simultaneous scalar test at arbitrary twist]\label{med-sign-test}
Fix a complex vector space \(V\), and let \(X\) be any sequence in
\(V\) with strictly positive powers, the allowed exponents for
\(\widetilde\rho_L\), and the indicated component symmetries. Define
scalar-type formal series \(Y_j\) and integral weights \(K_j\) by
\begin{align}
 Y&=M^tG_I X, & K_j&=k+v_j &&(\Delta>0),\label{med-dual-contraction}\\
 Y&=\Delta_{\rm mod}^{\ell_0}M^{-1}X, & K_j&=k-v_j+12\ell_0
        &&(\Delta<0).\label{med-inverse-contraction}
\end{align}
Each \(Y_j\) has strictly positive integral exponents and \(K_j\) is
even. The original vector series is in
\(S_{k,\widetilde\rho_L}^{\Pi_{\Delta,N}}\otimes V\) if and only if
\[
     Y_j\in S_{K_j}(\mathrm{SL}_2(\mathbb Z))\otimes V
     \qquad(0\le j<h).
\]
In particular all levels can be treated in their original Chow groups.
\end{theorem}
\begin{proof}
For positive twist \(M^tG_I X\) pairs each \(U_j\) with its dual
type, weighting the representatives by their orbit sizes. The exponents
in each pairing add to integers and are strictly positive. For negative
twist the translation phases of \(M\) and \(X\) cancel in \(M^{-1}X\).
The adjugate of \(M\) has nonnegative orders and
\(\operatorname{ord}_\infty(\det M)=W/12<\ell_0\), so multiplication by
\(\Delta_{\rm mod}^{\ell_0}\) gives strictly positive integral exponents
here as well. The parity defining \(\epsilon\) gives even \(K_j\).

On the upper half-plane \(M\) is an invertible frame of columns of
weights \(v_j\) and type \(\rho_+\). A modular \(X\) of the dual type
pairs with these columns to give scalar forms of weights \(k+v_j\).
For \(X\) of the same type, its coordinates \(M^{-1}X\) have weights
\(k-v_j\). This gives the claimed transformed weights in both formulas.
The positive expansions give cusp forms, since the full modular curve
has one cusp.

Conversely, the coefficients of the finitely many scalar cusp tensors
\(Y_j\) span a finite-dimensional subspace \(V_0\subset V\). Reconstruct
the column on the upper half-plane as
\[
  X=G_I^{-1}M^{-t}Y\quad(\Delta>0),\qquad
  X=\Delta_{\rm mod}^{-\ell_0}MY\quad(\Delta<0).
\]
The frame transformation law gives the required weight and the dual
or original type, respectively. Both expressions are holomorphic there.
Inversion of the formal Laurent series recovers the prescribed expansion,
including its positivity and its shared-prime symmetries. Translation by
the metaplectic group then gives vanishing at all cusps on a cover
trivializing the type. In a basis of \(V_0\) each coordinate is therefore
a cusp form of the stated type.
\end{proof}

\subsection{Fourier coefficient tests}
The theorem gives coefficient identities involving the actual cycles,
with known scalar coefficients. For even \(K\), put \(d_K=\dim S_K(\mathrm{SL}_2(\mathbb Z))\).
By Lemma~\ref{med-eta} for the trivial representation, the Fourier
coefficients through \(\lfloor K/12\rfloor\) determine the form.
Choose \(d_K\) independent positions \(A_K\) in this range and the
dual basis \(s_r\), so that \([q^t]s_r=\delta_{rt}\) for \(r,t\in A_K\).
For any scalar-type
series \(Y\) with coefficients in \(V\), set
\[
       \mathcal E_{K,n}(Y)=[q^n]Y-
             \sum_{r\in A_K}([q^n]s_r)\,([q^r]Y).
\]
If \(d_K=0\) the sum is empty. Existence and uniqueness of a modular
interpolant show that \(Y\) is cuspidal of weight \(K\) exactly when
all these differences are zero. Every coefficient in
\eqref{med-dual-contraction} and \eqref{med-inverse-contraction} is
a finite combination of original coefficients, since the multiplying
series have a finite lower bound on their exponents.

\begin{corollary}[Scalar generators of the Chow obstruction]\label{med-generated-rows}
For the actual cycles and either sign of \(\Delta\),
\[
 D_N=\operatorname{span}_{\mathbb C}
  \{\mathcal E_{K_j,n}(Y_j):0\le j<h,\ n\ge1\}\quad\text{in }W_N.
\]
\end{corollary}
\begin{proof}
For a subspace \(B\subset W_N\), all the displayed rows vanish
modulo \(B\) exactly when every \(Y_j\bmod B\) is a scalar cusp form,
by interpolation. By Theorems~\ref{med-sign-test} and~\ref{rd-main}
this is equivalent to \(D_N\subset B\). Taking \(B\) in turn to be
the span of the rows and \(D_N\) gives equality.
\end{proof}

These scalar rows describe the obstruction using unary theta coefficients,
level-one cusp forms and cycles at the original level. Proposition~\ref{med-boundary}
gives a geometric criterion for their vanishing.

\subsection{A finite module of weakly holomorphic generators}
At fixed weight the weakly holomorphic inputs form a finite free
module over the polynomial ring in the modular invariant. The
theta matrix gives explicit generators for either sign of the twist.

Write \(E_4,E_6\) for the normalized Eisenstein series and
\(j_{\rm mod}=E_4^3/\Delta_{\rm mod}\). For an even integer \(l\), choose the
unique \(\alpha\in\{0,1,2\}\), \(\beta\in\{0,1\}\) with
\(4\alpha+6\beta\equiv l\bmod12\), and put
\[
   b_l=E_4^\alpha E_6^\beta
               \Delta_{\rm mod}^{(l-4\alpha-6\beta)/12}.
\]
Then \(M^!_l(\mathrm{SL}_2(\mathbb Z))=b_l\,\mathbb C[j_{\rm mod}]\).
Indeed, clearing poles by a power of \(\Delta_{\rm mod}\) gives
the graded ring \(\mathbb C[E_4,E_6,\Delta_{\rm mod}^{-1}]\).
In a monomial of weight \(l\), the \(E_4\)- and \(E_6\)-exponents
are respectively \(\alpha\bmod3\) and \(\beta\bmod2\). After
dividing by \(b_l\), only powers of \(E_4^3/\Delta_{\rm mod}=j_{\rm mod}\)
and \(E_6^2/\Delta_{\rm mod}=j_{\rm mod}-1728\) remain.
The resulting polynomial is unique since \(j_{\rm mod}\) is nonconstant.

\begin{proposition}[Generators for both input types]\label{med-ring}
For the full reflection sector let \(U_j^\vee\) denote the \(j\)-th
column of \(G_I^{-1}M^{-t}\) (extend the representatives by their
prescribed signs). At fixed weight \(w=2-k\), a free basis over
\(\mathbb C[j_{\rm mod}]\) of the weakly holomorphic inputs is
\[
 \begin{cases}
   B_j=b_{\,w-v_j}U_j,&\Delta>0,\\
   B_j=b_{\,w+v_j}U_j^\vee,&\Delta<0,
 \end{cases}
 \qquad 0\le j<h .
\]
For the invariant sector the \(\Pi_{\Delta,N} B_j\) generate as a module.
In particular, in ordinary Chow one has the equivalence
\begin{equation}\label{med-all-powers}
 D_N=0\quad\Longleftrightarrow\quad
 Z_\kappa(j_{\rm mod}^t\,\Pi_{\Delta,N}B_j)=0
          \ \ \text{for all }t\ge0,\ 0\le j<h .
\end{equation}
\end{proposition}
\begin{proof}
For \(\Delta>0\) write \(F=MH\). The transformation law gives
each \(H_j\) scalar weight \(w-v_j\). For \(\Delta<0\) use the
dual matrix \(F=G_I^{-1}M^{-t}H\), so \(H=M^tG_IF\) and \(H_j\)
has weight \(w+v_j\). Both changes of basis and their inverses are
holomorphic on the upper half-plane and meromorphic at the cusp,
by the determinant formula. The preceding scalar calculation gives
the indicated free basis in either case. Applying \(\Pi_{\Delta,N}\),
which commutes with scalar multiplication, gives the invariant
generators. Linearity of \(Z_\kappa\) on their polynomial multiples
proves \eqref{med-all-powers}.
\end{proof}

The projected columns above give generators directly from theta
functions. A basis with the smallest possible pole bound can be
obtained using just the Fourier coefficients of cusp forms.

\begin{proposition}[A rational basis adapted to pole order]
\label{filtered-weak-basis}
Let \(V_B\) be the space over \(\mathbb Q\) of invariant weakly
holomorphic inputs with rational Fourier coefficients and pole order
at most the integer \(B\). Put \(h_B=\dim V_B\) and \(h_B=0\) for
\(B\le0\). Let \(d_\Pi\) be the dimension of the residue space
with the central and shared-prime symmetries, and
\(s_j=\dim S_{j,\widetilde\rho_L}^{\Pi}\). Then
\[
                 h_B=Bd_\Pi-s_k+s_{k-12B}\qquad (B\ge0).
\]
All the invariant inputs over \(\mathbb Q\) form a free module over
\(\mathbb Q[j_{\rm mod}]\) of rank \(d_\Pi\). It has a basis with
exactly \(h_B-2h_{B-1}+h_{B-2}\) generators of pole order in
\((B-1,B]\) for each \(B\ge1\). If \(d_\Pi>0\), the least integer
pole bound for a set of module generators is
\[
              1+\min\{B\ge0:s_{k-12B}=0\}.
\]
\end{proposition}
\begin{proof}
Each independent residue contributes one allowed exponent to
\((B-1,B]\), so the polar polynomials through \(B\) have dimension
\(Bd_\Pi\). A cusp form annihilates their image under \(\beta\)
precisely when its exponents all exceed \(B\). Dividing by
\(\Delta_{\rm mod}^B\) identifies this annihilator with
\(S_{k-12B,\widetilde\rho_L}^{\Pi}\). Principal-part duality
(Proposition~\ref{med-exact}) and its rational structure
(Theorem~\ref{med-complete-basis}) now give the formula for \(h_B\).

Write \(G_B=V_B/V_{B-1}\) for \(B\ge1\) and \(G_0=0\).
Multiplication by \(j_{\rm mod}=q^{-1}+O(1)\) embeds \(G_{B-1}\)
in \(G_B\) for \(B\ge2\). At each \(B\ge1\) lift a rational basis of
a complement of \(j_{\rm mod}G_{B-1}\) in \(G_B\) to \(V_B\).
The lifts and the powers of \(j_{\rm mod}\) times the earlier lifts
give a basis on this graded piece. Induction on the pole bound
therefore gives a basis of \(V_B\) consisting of the shifts within
that bound. Taking their union proves freeness and gives the stated
number of generators at each step. The lifts can be chosen by
elimination on the principal parts in \eqref{med-pivot}.

For \(B\) large, \(\dim G_B=d_\Pi\). Let \(B_0\) be the least
integer with \(s_{k-12B_0}=0\). If \(B_0>0\), the last nonempty
complement occurs at \(B_0+1\) and has dimension
\(s_{k-12(B_0-1)}>0\). If \(B_0=0\), all \(d_\Pi\) generators
occur at \(B=1\). A generating set bounded by \(B_0\) would lie in
the submodule spanned by the earlier lifts. Freeness makes that
submodule proper when \(d_\Pi>0\), proving minimality.
\end{proof}

If \(b_1,\ldots,b_{d_\Pi}\) are the integer pole bounds of this
basis, then \(\sum_i(b_i-1)=s_k\). Indeed, for \(B\) large its shifts
give \(\sum_i(B-b_i+1)=Bd_\Pi-s_k\) independent forms. Applied to
all their principal parts, Theorem~\ref{intro-chow-vanishing} gives
the cycle relations simultaneously. The polynomial action here is
on the input forms.
\section{Detection and restriction of fibre classes}\label{app-fibres}
\subsection{Detection on the separate fibres}
Fix an embedding of the geometric ground field into \(\mathbb C\).
Write \(a=\kappa-1>0\). Let \(b(P)\) be the tuple of fibre classes
defined by \eqref{rd-coeff} before pushforward, extended by zero at
the other smooth CM points. We can pull this tuple to a neat cover
to take ordinary fibre cohomology.

\begin{theorem}[Detection before pushforward]\label{fibre-detection}
On the space \(\mathcal P\) of principal parts with the central and
shared-prime symmetries, the map
\[
     P\longmapsto \bigl(\operatorname{cl}(b_y(P))\bigr)_y
             \ \in\ \bigoplus_y H^{2a}(E_y^{2a},\mathbb C)
\]
is injective for every \(N,\Delta,r\) as in Section~\ref{rd-intro}.
\end{theorem}
\begin{proof}
Choose the largest \(n\) occurring in a nonzero \(P\) and a residue
\(\mu\) with \(c_P(-n,\mu)\ne0\). For an integer \(B\) representing
\(r\mu\bmod 2N\), take the following element of \(L^\sharp\):
\[
   x=\begin{pmatrix} B/(2N)&-1/N\\ C&-B/(2N)\end{pmatrix},
   \qquad C=n|\Delta|+B^2/(4N).
\]
The allowed congruence on \(n\) gives \(C\in\mathbb Z\).
The vector has norm \(n|\Delta|\), and
\((B^2-4NC)/\Delta=-4Nn\operatorname{sgn}(\Delta)
\equiv\mu^2\bmod 4N\). Its genus character is one, since
\([1,B,NC]\) represents \(1\) in the defining formula
\cite[\S4]{BO}. The upper right entry makes \(x\) primitive in
\(L^\sharp\).

Choose a point over the CM point of \(x\). After conjugation by
\(\Gamma_0(N)\), which fixes each discriminant coset, every vector
contributing there lies on \(\mathbb R x\): trace-zero real matrices
fixing this point infinitesimally form the one-dimensional Lie algebra
of its stabilizer. The lattice points on this line are \(\mathbb Z x\).
By maximality of \(n\), only \(x,-x\) can contribute.

At this exponent the allowed residues \(\nu\) for \(x\) satisfy
\(r\nu=r\mu\bmod 2N\) and \(\nu^2=\mu^2\bmod4N\).
They form exactly the orbit of \(\mu\) under the shared-prime
involutions. Indeed at every odd prime power in \(N\) prime to \(r\)
the first condition forces equality. At an odd prime of \(N\)
dividing \(r\), fundamentality gives \(p\Vert N\), and the square
condition gives the two local signs. If \(N\) is odd the square
condition fixes parity. If \(2\mid N\) and \(r\) is odd it is
invertible at \(2\). In the remaining case \(2\Vert N\), and the
squares modulo \(8\) distinguish the two even residues modulo \(4\)
and allow the two signs of an odd residue. The Chinese remainder
theorem now gives exactly the involutions defined in
Section~\ref{rd-cycles}.

All coefficients in this orbit equal \(c_P(-n,\mu)\). For \(-x\)
the central symmetry multiplies the corresponding coefficient by
\(\operatorname{sgn}(\Delta)(-1)^a\). The product of genus character
and oriented fibre class acquires the same sign. The two contributions
therefore add, with positive stabilizer multiplicities. The resulting
fibre class has nonzero Betti class by the self-intersection
normalization of \(S_\kappa(x)\) before pushforward. This detects \(P\).
\end{proof}

In particular, if \(N\mid r\) and \(\operatorname{sgn}(\Delta)(-1)^a=-1\),
every coset \(r\mu\) is self-opposite. Pairing \(x\) with \(-x\)
in \eqref{rd-coeff} gives \(W_N=0\) for this whole family.

\begin{proposition}[Residues of correspondences of lower weight]\label{weight-residues}
Let \(K\) be the function field of a connected fine modular curve
and \(E/K\) its universal elliptic curve. Fix \(0\leq c<a\).
Even after finite extensions of \(K\) and trace, correspondences
from \(E^m\) to \(E^{2a}\) carrying
\(\CH^{c+1}(E^m,1)\) to \(\CH^{a+1}(E^{2a},1)\) have zero
projected Betti residue at every smooth fibre. In particular their
images cannot have boundary \(b(P)\) for any nonzero \(P\in\mathcal P\).
\end{proposition}
\begin{proof}
At a place of good reduction take the closure of a correspondence
of codimension \(m+a-c\) in the smooth proper relative product.
Flat pullback, intersection with this extending class, and proper
pushforward commute with the localization boundary, using the
refined restriction of the class at the place. The residue of the
image is thus the specialized correspondence applied to the residue.
In Betti cohomology this specializes a morphism of local systems
from \(H^{2c}(E^m)\) to \(\operatorname{Sym}^{2a}H^1(E)\), up to
Tate twist, after applying the signed projector. The correspondence
extends over the whole smooth base by closure in the smooth
relative product, so its fibre class is monodromy invariant.
Over a finite cover the monodromy is still Zariski dense in
\(\mathrm{SL}_2\): after removing branch points its image has finite
index in the original monodromy. By K\"unneth, the source is a sum
of tensors with at most \(2c\) standard factors and trivial factors.
Its highest weights are at most \(2c\), while the target is irreducible
of highest weight \(2a\). The morphism is zero. Residues commute
with finite pushforward, so trace preserves this vanishing.
The last claim follows from Theorem~\ref{fibre-detection}.
\end{proof}

In particular this applies to all products of a generic cycle of
codimension \(a\) with a base function: take \(m=c=0\) and the
cycle itself as correspondence. Thus any lift in
Proposition~\ref{med-boundary} for a nonzero invariant weakly
holomorphic input remains nonzero on the generic fibre modulo the
span of such products and their finite norms. The boundary detects
linear independence there as well. In codimension two an explicit
construction on the two-parameter family of products of elliptic
curves has an isogeny-graph difference as boundary
\cite[Theorem~4.6 of v1]{SreekantanProducts}. The prescribed
boundary here is the entire tuple on a single modular curve.

\subsection{Restriction at the cusps and the level-one calculation}
\begin{lemma}[Restriction on the Scholl summand]\label{scholl-open-chow}
Let \(\widetilde X\) be the full fine-level Scholl compactification of
the \(m=2\kappa-2\) fold elliptic power, \(U\) its smooth-base locus,
and \(e\) the Scholl projector. For \(\kappa\ge2\) and every \(n\ge0\),
restriction gives an isomorphism
\[
 e\CH^\kappa(\widetilde X,n)_{\mathbb Q}
       \xrightarrow{\;\sim\;} e\CH^\kappa(U,n)_{\mathbb Q}.
\]
For correspondences, restriction also gives
\[
 (e\boxtimes e)\CH^{m+1}(\widetilde X\times\widetilde X,n)_{\mathbb Q}
 \xrightarrow{\;\sim\;} (e\boxtimes e)\CH^{m+1}(U\times U,n)_{\mathbb Q}.
\]
Both assertions hold on the cyclic-level summands.
\end{lemma}
\begin{proof}
Write \(C\) for the finite cusp scheme of the full level, taken
at least three. Scholl's character projector is trivial on level
translations, odd in each
elliptic coordinate and alternating under permutations
\cite[\S1.1.2]{SchollMotives}. Scholl's localization sequence
\cite[Theorem~1.3.3 and \S1.4.0]{SchollMotives} gives the exact segment
\[
 \begin{aligned}
 H_{\mathrm M}^{i-m-2}(C,\mathbb Q(j-m-1))
 &\longrightarrow eH_{\mathrm M}^{i}(\widetilde X,\mathbb Q(j))
 \longrightarrow eH_{\mathrm M}^{i}(U,\mathbb Q(j))\\
 &\longrightarrow H_{\mathrm M}^{i-m-1}(C,\mathbb Q(j-m-1)).
 \end{aligned}
\]
Indeed, in Scholl's homological indexing the term preceding
restriction in bidegree \((a,b)\) is \(H_{a-m}(C,b)\).
Poincar\'e duality on \(\widetilde X\) of dimension \(m+1\) gives
\(a=2(m+1)-i\) and \(b=m+1-j\). Duality on the zero-dimensional
scheme \(C\) then gives the two displayed cusp terms in consecutive degrees. Scholl proves
the sequence over \(\mathbb Q\). The same boundary stratification
and character calculation apply after extension to the geometric
ground field \cite[\S\S2--3]{SchollMotives}.

Both cusp terms adjacent to restriction have twist \(j-m-1\).
For \(j=\kappa\) this is \(1-\kappa<0\). On the smooth
zero-dimensional cusp scheme the motivic groups of this weight
vanish in every degree: Bloch's cycle complex in negative
codimension is zero. Thus restriction is an isomorphism for all
\(i\). Taking \(i=2\kappa-n\) gives the assertion for
\(\CH^\kappa(-,n)_{\mathbb Q}\), including \(n=0\) for ordinary
Chow groups.

For the second assertion, take the product of the boundary
stratification with a smooth variety \(T\). Scholl's calculation
of the character part applies with this factor throughout
\cite[Lemma~1.3.1, Proposition~2.4.1 and \S3]{SchollMotives}.
In the sequence for restriction from \(\widetilde X\times T\) to
\(U\times T\), the two boundary groups are now motivic cohomology
of \(C\times T\) with the same twist \(j-m-1\). At \(j=m+1\)
this twist is zero. For a smooth variety, the normalized cycle
complex in codimension zero is concentrated in degree zero, with
one generator per connected component. Take \(T=\widetilde X\)
and then \(T=U\). In either case inversion in one elliptic
coordinate fixes the connected components, so the projector on the
second factor kills both boundary groups. Restricting first in one
factor and then in the other proves the correspondence assertion.
Averaging over the finite group preserving the cyclic marking
commutes with both restrictions.
\end{proof}

On the smooth full-level modular curve, level translations already
act as the identity in relative rational Chow. Indeed, for a
nonzero level section \(t\) of the universal elliptic curve
\(\pi:E\to Y\), a positive power of \(\mathcal O(t-0)\) is pulled
back from \(Y\), since \(t\) is torsion in the relative Picard scheme.
Pulling back along the zero section identifies this line bundle
on \(Y\) with a power of the Hodge bundle: the two sections are
disjoint and the normal bundle of zero is the dual Hodge bundle.
The discriminant trivializes its twelfth power on \(Y\).
Hence \([t]=[0]\) in \(\Pic(E)_{\mathbb Q}\). Pullback by the
difference map \(E\times_Y E\to E\) identifies the translation
graph with the diagonal in Chow. Taking products shows that the
signed CM classes and isogeny graphs on the open fine level already
belong to the full Scholl summands. The lemma identifies their
unique lifts in the corresponding compact summands.

\begin{theorem}[Detection on pointed elliptic curves]\label{pointed-vanishing}
Put \(m=2\kappa-2\), and let \(\mathcal M_{1,m+1}(\Gamma_0(N))\)
parametrize smooth elliptic curves with a cyclic subgroup of order
\(N\), their origin and \(m\) additional distinct markings.
Restriction to this stack detects vanishing of every compact
projected Heegner class. For \(N=1\) and \(2\leq\kappa\leq5\),
all coefficients of \(\Phi_\kappa^{(1)}\) vanish in ordinary compact Chow for
every compatible \((\Delta,r)\). In particular
this gives a separate proof of Theorem~\ref{med-target} in this range.
\end{theorem}
\begin{proof}
On the smooth modular stack at level \(N\) let \(Y\) be the
\(m\)-fold power of the universal elliptic curve. Remove the loci
\(z_i=0\) and \(z_i=z_j\) for \(i\ne j\). The resulting open
\(V\) is exactly the displayed moduli stack. Restriction to \(V\)
has trivial kernel on the signed part of \(\CH^\kappa(Y)_{\mathbb Q}\).
Indeed, by localization the full kernel is generated by pushforwards
of cycles lying on the excluded loci. Inversion in coordinate \(i\)
fixes the first locus pointwise, and transposition of \(i,j\) fixes
the second. Both have character \(-1\), so the signed projector
kills every such pushforward on \(Y\). This kernel argument takes
place on \(Y\) and uses no action on \(V\). A compact projected class
vanishing on \(V\) thus vanishes on \(Y\), hence on the full
fine-level interior by pullback. Lemma~\ref{scholl-open-chow} and
finite-level descent give compact vanishing.

At level one the computations of Belorousski and Bishop give
\(\CH^{>0}(\mathcal M_{1,n})_{\mathbb Q}=0\) for \(3\le n\le10\)
\cite{BishopPointed}. Here \(3\le m+1\le9\) in the asserted
range, so all the restricted coefficients vanish.
\end{proof}

There is an explicit description of the surviving cycles on this
open stack at arbitrary level. For a vector \(x\) choose its
integral oriented endomorphism \(u_x\) as in Section~\ref{rd-cycles}.
In the fibre with \(2a\) coordinates take the locus
\[
             z_{2j}=u_x(z_{2j-1}),\qquad 1\leq j\leq a .
\]
Apply the signed projector in the full elliptic power, restrict
to distinct nonzero markings, and include the scalar used to
normalize the graph power. Denote the result pushed into the
pointed stack by \(G_x\). This is exactly the restriction of
\(S_\kappa(x)\): every axis term removed from a graph divisor
has some coordinate zero, also after signed permutation. Thus the
remaining equality in ordinary Chow for a weakly holomorphic \(f\)
can be tested entirely on pointed smooth curves as
\[
 \sum_{\mu,n>0} c_f(-n,\mu)n^{a/2}
 \sum_{x\in\Gamma_0(N)\backslash L_{r\mu}[n|\Delta|]}
       \frac{2\chi_\Delta(x)}{|\Gamma_x|}\,G_x=0 .
\]
Both sums are finite for each \(f\). The order of projection and
restriction in the definition of \(G_x\) keeps the signs meaningful
even when a coordinate inversion carries a configuration out of
the open set.

\subsection{How many independent fibre boundaries occur?}
Let \(\mathcal P_B\) consist of the invariant polar polynomials of
pole order at most a nonnegative integer \(B\). Let \(d_\Pi\) be the
dimension of the component space with the central and shared-prime
symmetries, and write \(s=\dim S\). Define
\[
      K_B=\mathcal P_B\cap\ker\beta,\qquad
      D_N^{\le B}=Z_\kappa(K_B)\subset W_N .
\]
By Proposition~\ref{med-exact}, \(K_B\) is precisely the space of
principal parts of invariant weakly holomorphic forms with this pole
bound. Use a common fine cover, and let \(T\) contain the supporting
fibres for every element of \(\mathcal P_B\). This is a finite set
because there are finitely many allowed indices below the bound and
finitely many CM orbits at each index.

The integer \(d_\Pi\) can be read off from the level and the twist.
Put \(\epsilon=\operatorname{sgn}(\Delta)(-1)^{\kappa-1}\),
\(J=\{p:p\mid(N,\Delta)\}\), \(h_p=p^{v_p(2N)}\) for \(p\in J\),
and \(h=2N/\prod_{p\in J}h_p\). Then
\[
 d_\Pi=\frac{h+\epsilon\gcd(2,h)}{2}
                \prod_{p\in J}\frac{h_p+\gcd(2,h_p)}{2}.
\]
Indeed, each shared-prime involution acts by negation on the factor
\(\mathbb Z/h_p\mathbb Z\) alone. Its orbits give the corresponding
factor in the product. On the remaining cyclic group of order \(h\),
negation has \(\gcd(2,h)\) fixed points, so its \(\epsilon\)-eigenspace
has the dimension in the first factor.

\begin{corollary}[The dimension of the space of fibre boundaries]
\label{bounded-fibre-boundaries}
Put \(s_j=\dim S_{j,\widetilde\rho_L}^{\Pi}\), so \(s_k=s\).
For every integer \(B\ge0\), the space of tuples \(b_T(P)\),
\(P\in K_B\), which are boundaries from \(e\CH^\kappa(U_T,1)\)
has dimension
\[
       B d_\Pi-s+s_{k-12B}-\dim D_N^{\le B}.
\]
Theorem~\ref{rd-m-vanishing} makes the last term zero. In particular,
for \(B\ge k/12\) the dimension is \(B d_\Pi-s\).
\end{corollary}
\begin{proof}
By the principal-part map, \(K_B\) is the complexification of \(V_B\)
in Proposition~\ref{filtered-weak-basis}; its dimension is
\(Bd_\Pi-s+s_{k-12B}\). For \(k-12B\le0\) the cusp space is zero
by holomorphy on a finite cover in weight zero and by the
negative-weight vanishing otherwise.
The map \(P\mapsto b_T(P)\) is injective, since its composite with
the tuple of fibre Betti classes is injective by
Theorem~\ref{fibre-detection}. Proposition~\ref{med-boundary}
identifies the desired subspace with the image under this injection
of \(\ker(Z_\kappa|_{K_B})\). Its dimension follows by rank and nullity.
\end{proof}

There may be further geometric boundaries among the same fibre tuples.
The distinction can be measured without constructing a higher Chow lift.
For \(B\ge k/12\), let \(\mathcal R_B\) consist of the \(P\in\mathcal P_B\)
for which \(b_T(P)\) is a boundary. Write
\(W_B=Z_\kappa(\mathcal P_B)\) and \(D_B=Z_\kappa(K_B)\).
The two tests on polar polynomials give a canonical exact sequence
\[
 0\longrightarrow K_B\cap\mathcal R_B
 \longrightarrow \mathcal R_B
 \xrightarrow{\ \beta\ } S^\vee
 \xrightarrow{\ \overline Z_B\ } W_B/D_B
 \longrightarrow 0 .
\]
Indeed, \(\beta\) maps \(\mathcal P_B\) onto \(S^\vee\). Two lifts of a
functional differ by an element of \(K_B\), so their cycle classes
agree modulo \(D_B\). If a lift \(P\) maps into \(D_B\), choose
\(Q\in K_B\) with \(Z_\kappa(Q)=Z_\kappa(P)\).
Then \(P-Q\in\mathcal R_B\) by localization and
\(\beta(P-Q)=\beta(P)\). Surjectivity at the last term follows from
the definition of \(W_B\).

In particular, once the weakly holomorphic relations vanish in Chow,
\(D_B=0\) and \(W_B=W_N\) by the low-coefficient bound.
The additional boundaries modulo the weakly holomorphic ones are then
canonically \(\ker(S^\vee\to W_N)\), of dimension \(s-\dim W_N\),
independent of \(B\). Thus increasing the pole bound adds
weakly holomorphic relations while the quotient measuring additional
geometric relations has already stabilized.

For a basis of the weakly holomorphic boundary space choose one higher Chow lift of
each tuple. The restrictions of the lifts to the generic fibre are
linearly independent modulo the correspondences of
Proposition~\ref{weight-residues}. Indeed, the latter have zero
projected Betti residues, whereas the chosen boundary tuples have
independent residues. In the level-one range of
Theorem~\ref{pointed-vanishing}, the full dimension \(B d_\Pi-s\)
also follows directly from the computation on pointed curves.
For \(d_\Pi>0\) these dimensions grow without bound with \(B\).
\subsection{Vertical cycles and admissible extensions}\label{sec:ks-integral}
Put \(a=\kappa-1\ge1\) and \(m=2a\). We use the Scholl projector
\(e\) on the \(m\)-fold elliptic power. At cyclic level we also take
invariants under the group preserving the cyclic subgroup, on the
union of the determinant components of a full fine cover.

\begin{lemma}[The good special fibre]\label{ks-good-fibre}
Let \(p\nmid h\). On the geometric special fibre at \(p\) of the
Kuga--Sato variety of cyclic level \(h\), the signed group
\(e\CH^a_{\mathbb Q}\) is zero. Here the cyclic summand is realized
on a fine cover of level prime to \(p\). The same vanishing holds
at the ordinary generic point of the modular curve after extending
its constant field to \(\overline{\mathbb F}_p\). The elliptic power
here is over the function field of that curve itself.
\end{lemma}
\begin{proof}
\emph{A relation which can be specialized.}
Use the geometric cyclic \(p\)-isogeny correspondence \(T_p\) of
Proposition~\ref{marked-cm-spectrum}. In characteristic zero the
proof of Proposition~\ref{md-t-spectrum} gives \(P_p(T_p)\in J\),
where the roots of \(P_p\) are classical eigenvalues in weight
\(m+2\), trivial character and level supported on \(h\).
For \(h=1\) use the full-level cuspidal Hecke algebra of
Section~\ref{med-u-five}. The inclusion in \(J\) gives a
presentation on isogeny curves with zero signed relative
realization, up to terms over finitely many pairs of base points.

A power of \(P_p(T_p)\) is supported on point pairs as a Chow
correspondence. Compose the presentations over a chain of normalized
isogeny curves. On the smooth modular locus the two maps to each
intermediate base curve are finite \'etale in characteristic zero.
Their fibre product is Tor-independent, so the refined diagonal for
the composite restricts at its generic point to composition of
the fibre correspondences. At that point, after
a finite extension if necessary, the isogenies identify all the
signed fibre motives with one motive \(M_e\). It is the alternating
summand of \(h^1(E)^{\otimes m}\), evenly finite dimensional of rank
\(m+1\), since its realization is \(\operatorname{Sym}^m H^1(E)\)
and \(m\) is even. The transported factors lie in the ideal of
numerically trivial endomorphisms of \(M_e\). This ideal has a
nilpotence bound depending only on the rank
\cite[Theorem~3.3 and Corollaries~3.4--3.5]{JannsenFinite}.
Thus one length \(L\) kills every generic composite. Finite trace
descends the vanishing, and localization on each base curve leaves
only classes above finitely many points. Paths passing through a
point pair had this support already, since both base maps are finite.

Two such point-pair correspondences compose to zero. Distinct
intermediate base points give an empty intersection. For a common
intermediate point \(x\), both classes are pushed forward from the
Cartier divisor defined by \(x\) in the middle factor of the triple
product. Its normal line is the constant tangent line of the base
at \(x\). The self-intersection formula therefore gives zero.
Thus \(P_p(T_p)^{2L}=0\) on the open square. On the compactified
isogeny curves both endpoints are cuspidal together. Restriction
of every composite Hecke graph to the open square therefore
removes precisely the paths through a cusp, by open base change.
Lemma~\ref{scholl-open-chow} now gives the same equality on the
compact square after projection at both ends. This equality involves only
the cyclic Hecke correspondence and algebraic scalars. If proved
over \(\mathbb C\), it descends to \(\overline{\mathbb Q}\) by
spreading a rational equivalence and specializing its parameters.

To change to a fine level prime to \(p\), first work over the
noncuspidal locus in characteristic zero and pass to a common
refinement. A finite étale covering there induces a map \(P\) of
elliptic powers with \(PP^t=\delta\,1\), where \(\delta=\deg P\).
Consequently \(T\mapsto P^tTP/\delta\) preserves composition,
normalized trace is its left inverse, and the identity maps to the
deck projector. These formulas compose in towers. They give the
polynomial relation on the open square at the prime-to-\(p\) level,
where a cyclic \(p\)-isogeny transports each source frame uniquely.
Averaging its graph at both ends gives \(P^tT_pP/\delta\).
Lemma~\ref{scholl-open-chow}, applied at both ends, gives the relation
on the compact square at this level. On its smooth proper model at
\(p\), specialization commutes with
composition of correspondences \cite[\S20.3]{Fulton}. We may
therefore apply the polynomial identity to cycles and restrict it
to the generic point of the special modular curve. Explicitly, a
cycle on the generic elliptic power extends by closure, followed by
the finite projectors. The base projections of the Hecke support
are finite and preserve the ordinary locus. Restricting the compact
action on such an extension therefore computes the ordinary
push--pull below, independently of the extension.

\emph{The two ordinary isogenies.}
We compute that action on cycles, with source pushforward and
target pullback. At ordinary points the two reduced components of
the degree-\(p\) correspondence, both of multiplicity one, are the
relative Frobenius and its dual Verschiebung
\cite[Theorem~13.4.7]{KatzMazurModuli}. On the first, the base map together with
the fibre Frobenius gives total Frobenius \(F\), contributing \(F^*\).
For the second, let \(u\) be the projection of the base-Frobenius
pullback of the elliptic power, let \(\varphi\) be its relative
fibre Frobenius, and let \(v\) be the fibrewise Verschiebung. Thus
\(u\varphi=F\), \(v\varphi=[p]\), and this branch contributes
\(u_*v^*\). The purely inseparable map \(u\) has degree \(p\).
On the signed, cyclic summand we obtain
\[
 F^* u_*v^*=p\,\varphi^*v^*=p[p]^*=p^{m+1}.
\]
Here each \(h^1\) contributes \(p\) under \([p]^*\)
\cite[Theorem~3.1]{DeningerMurre}. Multiplication by \(p\) preserves
the cyclic subgroup, so the change of a fine frame contributes no
diamond factor on this summand. Equivalently the opposite convention
for the Hecke graph gives the same correspondence by dual isogeny.
Consequently \(T_p=F^*+p^{m+1}(F^*)^{-1}\) over this ordinary locus.

A codimension-\(a\) cycle over this function field spreads on an
open subset over a finite constant field. Some power \((F^*)^b\)
therefore acts on it by \(p^{ab}\), by flat pullback under absolute
Frobenius on the smooth total space. On the span of its normalized
\(p^{-a}F^*\)-orbit the eigenvalues of \(F^*\) are \(p^a\zeta\), with
\(\zeta\) a root of unity. For the corresponding eigenvalue
\(\lambda\) of \(T_p\), the roots of \(X^2-\lambda X+p^{2a+1}\)
have absolute values \(p^a,p^{a+1}\). For a cuspidal root of
\(P_p\) both absolute values are \(p^{a+1/2}\), by Deligne's
theorem. For an Eisenstein root they are \(1,p^{2a+1}\).
The three pairs are distinct since \(a\ge1\). Thus \(P_p(T_p)\)
is invertible on this span and its \(2L\)-th power is zero by the
compact relation. The cycle therefore vanishes.

\emph{Closing the fibres.}
Localization over the smooth modular curve has remaining terms
\(e\CH^{a-1}(E_x^{m})\) at closed points \(x\). Over a finite
residue field of size \(q\), the Frobenius of \(h^1(E_x)\) satisfies
the elliptic Frobenius polynomial as an endomorphism of the elliptic
curve. All eigenvalues allowed by this polynomial have absolute
value \(\sqrt q\). Its tensor power on the signed motive thus
satisfies a polynomial with all roots of absolute value \(q^a\).
A codimension-\((a-1)\) cycle over a finite extension is multiplied
by \(q^{b(a-1)}\) by a power of Frobenius. Coprimality of these
polynomials kills the cycle, including at supersingular points.
At a cusp the Scholl boundary calculation in
Lemma~\ref{scholl-open-chow} has twist \(a-m-1<0\). The calculation
uses the split polygon and its toric resolution
\cite[\S\S2--3]{SchollMotives} and applies in characteristic \(p\)
with rational coefficients. This kills the remaining boundary.
\end{proof}

We next fix the integral model. Choose a full fine level \(n\)
divisible by \(N\) and by two relatively prime integers at least
three. Fix a finite Galois extension \(F/\mathbb Q\) containing
\(\mathbb Q(\zeta_n,\sqrt\Delta)\). Over \(\mathcal O_F\) take the
union of all Weil-pairing components of the full-level model,
base changing each from its cyclotomic base. Resolve the modular
base surface equivariantly over its supersingular points until the
reduced fibre there has smooth components with normal crossings.
Over its noncuspidal locus take the elliptic power. At the cusps
use the canonical Scholl resolution. Denote the resulting regular
proper arithmetic variety by \(\mathfrak Y\).
We use the componentwise construction of the fine models for CM heights in
\cite[\S\S2.1--2.3, 3.1]{Zhang} and
\cite[\S\S2.2.3--2.2.4]{QiuCM}, allowing further point blowups
of the base over supersingular points. Over the cyclotomic base the
ordinary locus is smooth, as is the formal base at a cusp. Thus
extension of the coefficient field leaves those charts smooth
\cite[Chs.~8, 12--13]{KatzMazurModuli}.

Call an arithmetic class admissible here if its curvature is zero
and the underlying cycle restricts to a cohomologically trivial
class on each geometric component of each finite fibre. Finite
cycle classes are taken in \(\ell\)-adic cohomology with \(\ell\)
different from the residue characteristic. We use rational
coefficients, with complexification for Fourier series. Let \(b_N\)
average the deck transformations that preserve the cyclic marking,
and take the character \(\chi_\Delta\) on \({\rm Gal}(F/\mathbb Q)\).
The Scholl, deck and Galois projectors commute.

\begin{proposition}[Admissible arithmetic representatives]
\label{ks-admissible-restriction}
On the \(e b_N\)-summand of character \(\chi_\Delta\), restriction
from the admissible arithmetic Chow group of \(\mathfrak Y\) to
\(\CH^\kappa(\mathfrak Y_F)_{\mathbb Q}\) is an isomorphism onto
the corresponding summand. In particular the admissible
arithmetic lifts of the Heegner coefficients preserve their
rational relations and form a cusp form of weight
\(\kappa+\tfrac12\) on this model.
\end{proposition}
\begin{proof}
\emph{Archimedean classes and strict components.}
Scholl's signed cohomology of the compact complex fibre is
concentrated in degree \(2a+1\). An \(e\)-projected generic
codimension-\(\kappa\) cycle thus admits a Green current of
curvature zero. For a fixed class on the total algebraic model the
remaining ambiguity is the image in arithmetic Chow of a form of
bidegree \((a,a)\) killed by \(\partial\bar\partial\), modulo
the images of \(\partial\) and \(\bar\partial\)
\cite[\S3.3]{GS}. On a compact K\"ahler
manifold this group of forms is isomorphic to \(H^{a,a}\) by the
\(\partial\bar\partial\)-lemma. Its \(e\)-part is zero
\cite[Theorem~2.3.1 and proof of Proposition~3.1.1]{Zhang}.
Hence a zero-curvature \(e\)-class with zero generic restriction
has a presentation by vertical algebraic cycles equipped with
zero Green currents, by ordinary localization.

The strict components of a geometric finite fibre are the elliptic
powers with Scholl resolution on the original reduced modular
components. At their cusps the model is the split Tate chart.
Elsewhere the elliptic power is smooth over the base, and the
proper transform of an original smooth base curve is isomorphic
to it. On each such component \(D\), Zhang's version of Scholl's
calculation gives
\[
        eH^{2a}(D)=eH^{2a+2}(D)=0
       \qquad\text{\cite[Proposition~2.1.2 and Theorem~2.3.1]{Zhang}} .
\]
In particular every projected extension is cohomologically
admissible on the strict components.

\emph{Removing the exceptional components.}
Over a supersingular point \(s\), an exceptional base curve \(C_i\)
gives the component \(C_i\times E_s^m\) by smooth base change.
After a finite residue extension the mixed graphs of a basis of
the quaternion endomorphism algebra of \(E_s\) span the four
mixed cohomology classes of \(E_s^2\). Their pairing is the
nondegenerate trace pairing. Its inverse splits four copies of
the Lefschetz motive off \(h^1(E_s)^{\otimes2}\). The complementary
projector is homologically trivial and vanishes by finite
dimensionality of elliptic motives \cite[Proposition~7.5]{KimuraFinite}.
Tensoring \(a\) times and applying the sign projector shows that
the signed fibre motive is a sum of copies of \(\mathbb L^a\).
Write \(W_s=e\CH^a(E_s^m)_{\mathbb Q}\). On the product with \(C_i\)
the signed Chow groups in codimensions \(a\) and \(a+1\) are
\(W_s\otimes\CH^0(C_i)_{\mathbb Q}\) and
\(W_s\otimes\CH^1(C_i)_{\mathbb Q}\), respectively. In the latter
group the degree on \(C_i\) suffices since its degree-zero Picard
group over the algebraic closure of a finite field is torsion.

Start with any projected closure of a generic class. The degrees
of its restrictions to the exceptional components form a column
\((r_j)_j\) in \(W_s\). Adding the pushforwards from \(C_i\times E_s^m\)
of \([C_i]\boxtimes w_i\) changes that column by
\((\sum_i (C_i\cdot C_j)w_i)_j\). Indeed, the normal line and
the intersections of distinct components are pulled back from the
base surface. The exceptional intersection matrix is negative
definite. To see this, write the whole connected fibre as
\(\sum m_i C_i\). For a divisor \(\sum d_i C_i\) supported on
the exceptional curves, the square is
\[
 -\sum_{i<j}m_i m_j(C_i\cdot C_j)
                (d_i/m_i-d_j/m_j)^2 .
\]
Equality forces a constant ratio on the whole fibre and the ratio
on a strict component is zero. There is therefore a unique
exceptional correction to the given closure. It descends from the
splitting residue field by trace and commutes with the deck and
Galois actions. Together with the zero-curvature Green current it
gives an admissible representative, proving surjectivity.

For a class with zero generic restriction, average a vertical
presentation over \(e\), \(b_N\) and the Galois character before
pushing forward the component classes. The contribution from a
strict component has zero cohomology by \(eH^{2a}(D)=0\), hence
its intersection contributes zero to the cohomology on every
exceptional component. The matrix just computed forces all
exceptional coefficients of an admissible presentation to vanish.
It remains to consider the averaged classes on the strict components.

\emph{Forgetting the ordinary Igusa generator.}
Fix a place above \(p\). If \(p\nmid n\), Lemma~\ref{ks-good-fibre}
applies directly after taking the cyclic invariants. Suppose now
\(p^r\Vert n\), and put \(R=\mathbb Z/p^r\mathbb Z\).
On the ordinary locus of a fixed determinant component the
Drinfeld \(p\)-marking has a kernel line \(K\subset R^2\) on
geometric points. For fixed \(K\), the remaining datum, up to a
purely inseparable change of base, is a generator of the
\'etale quotient on the ordinary elliptic curve
\cite[Theorems~12.6.1 and 13.10.3]{KatzMazurModuli}.
The scalar matrices \(uI\), \(u\in R^\times\), preserve \(K\)
and every cyclic marking and multiply this generator by \(u\).
An individual scalar also changes the Weil pairing by \(u^2\).

Choose a local inertia element \(\tau\) acting by \(u\) on
\(\mu_{p^r}\). Such an element exists because the cyclotomic
extension is totally ramified at \(p\). Its square \(\sigma=\tau^2\)
acts by \(u^2\), fixes the prime-to-\(p\) roots of unity, and
satisfies \(\chi_\Delta(\sigma)=1\). Use \(uI\) at the \(p\)-part
of the full level and the identity at its other primes. Composing
this deck map with the inverse Galois transport under \(\sigma\)
now preserves the chosen determinant branch. On the ordinary
chart, Galois transport between the two cyclotomic components
reduces to the identity on the underlying elliptic curve and
its Drinfeld sections: their residue field maps are identical.
Thus the composite acts on the surviving generator by \(u\)
and fixes the elliptic coordinates. This argument uses only the
ordinary charts, before any modification over the supersingular
points. It applies also at \(2\).

The component class in our averaged presentation is invariant
under each of these composites. Put \(k=\overline{\mathbb F}_p\)
and let \(L/K_0\) be the function fields of this ordinary component
and of the tame fine curve over \(k\). Write \(U_p\) for the effective
unit group on the generator. The extension \(L/L^{U_p}\) is Galois.
The finite map for \(L^{U_p}/K_0\) has one geometric generic point
above a given ordinary elliptic curve with its tame marking, since
the units act transitively on the generators. Thus this last
extension is purely inseparable. The elliptic power over \(L\)
is pulled back from \(K_0\). For the first extension, pullback
after trace is the sum of the unit actions on its Chow group.
For the purely inseparable extension normalized trace inverts
pullback. The invariant class therefore descends to the elliptic
power over \(K_0\).
Prime-to-\(p\) changes of frame preserving the tame cyclic marking
can be taken to have identity \(p\)-part. On a fixed Weil-pairing
component take their determinant-one subgroup. Averaging over
this group descends further to cyclic level \(h=N/p^{v_p(N)}\).
Indeed the tame deck group acts transitively on the determinant
components: its diagonal matrices realize every unit determinant.
Its invariants on their direct sum are precisely the invariants
under the determinant-one stabilizer on one component. Hence the
descended class lies in the summand of
Lemma~\ref{ks-good-fibre} over the function field of the tame
cyclic curve. The central automorphism \(-1\)
acts trivially on the even power. Lemma~\ref{ks-good-fibre} at
this tame level kills the descended class. On restoring the
closed points, the same elliptic Frobenius argument in that lemma
kills the signed groups \(e\CH^{a-1}(E_x^m)\). At each cusp the
Tate parameter on the cyclotomic chart is \(t=q^{1/n}\), and
the section represented by \(t\) meets component \(1\) of the
\(n\)-gon \cite[Theorems~8.11.10 and 10.9.1]{KatzMazurModuli}.
Even when \(p\)-power roots of unity reduce to one, this section
permutes all polygon components. On the torus inversion is
\(z\mapsto z^{-1}\), also in characteristic two. In Scholl's torus
calculation the isogeny with kernel \(\mu_2\) is then radicial of
degree two. Its pullback on rational Chow groups is still invertible
by finite flat push--pull and universal injectivity.
The split strata calculation of \cite[\S\S2.3--2.4, 3.1]{SchollMotives} therefore
gives the same negative-twist boundary with rational coefficients.
This proves that the strict component terms vanish and completes
the proof of injectivity.

Finally, Proposition~\ref{coef-galois-character} places all
generic coefficients in the indicated character space. Applying
the inverse restriction to Theorem~\ref{med-target} gives the
arithmetic modular form on this single model over \(F\).

One can also start with the integral CM cycles and admissible Green
currents of \cite[\S2.2.4]{QiuCM}. For each finite orbit choose a
finite Galois extension \(F'/F\) defining its representatives.
On a common refinement dominating the base change of the chosen
base surface, let \(\pi:\mathfrak Y'\to\mathfrak Y\) be the resulting
proper morphism. Such a refinement is obtained by resolving the
graph on the base surfaces over the supersingular points. At the
cusps the split toroidal construction commutes with base change.
Take normalized arithmetic trace \(\pi_* /[F':F]\), followed by
the deck and character projectors. Curvature zero is preserved.
For a strict or exceptional target component the refined pullback
of its defining divisor is the sum of the source components
above it with their multiplicities. Proper base change for
refined intersection therefore carries the zero restriction
classes on those components to a zero restriction class on the
target. By injectivity this trace is the unique lift of the
generic coefficient. The field \(F'\) may vary with the orbit.
\end{proof}

\begin{introremark}[The class and its intersection number]
The numerical pairing can already be formed after the existence
part of Proposition~\ref{ks-admissible-restriction}. Indeed, two
admissible lifts of one generic class differ in the signed summand
by a sum of vertical classes \(i_{D*}v_D\) with zero Green current.
For any admissible lift \(\widehat Y\), the contribution of such a
term to the arithmetic degree is proportional to
\(\deg_D(v_D\cdot i_D^*Y)\). This degree vanishes because the
cycle class of \(i_D^*Y\) is zero. Thus the number depends only
on the generic classes. Injectivity in the proposition gives the
stronger conclusion that the arithmetic Fourier coefficients
themselves satisfy the rational relations of the generic fibre.

For codimension two the injectivity argument also has a shorter
form. On each smooth strict component \(D\), a vertical term
comes from \(e\CH^1(D)_{\mathbb Q}\). The degree-two cohomology
of \(D\) has zero \(e\)-part. The divisor class map embeds the
rational N\'eron--Severi group in \(\ell\)-adic cohomology after
extension of scalars, and \(\operatorname{Pic}^0(D)\) over
\(\overline{\mathbb F}_p\) is torsion. Thus the strict terms
vanish already at this step. The exceptional intersection
matrix above then gives injectivity.
\end{introremark}

\section*{Acknowledgements}
The author thanks T.~Yang for suggestions and for help checking
details of the original arithmetic inner product manuscript.

\section*{AI Use Declaration}
This manuscript was reviewed and revised using GPT Codex.

\end{document}